\documentclass[12pt]{lecture}
\usepackage{amssymb,amsmath}
\usepackage{graphicx,rotating}
\usepackage{multirow}
\usepackage{colordvi}

\newtheorem{thm}{Theorem}[chapter]
\newtheorem{lemma}[thm]{Lemma}
\newtheorem{prop}[thm]{Proposition}
\newtheorem{cor}[thm]{Corollary}
\newtheorem{conj}[thm]{Conjecture}
\newtheorem{example}[thm]{Example}
\newtheorem{remark}[thm]{Remark}
\newtheorem{definition}[thm]{Definition}
\newenvironment{defn}{\begin{definition}\rm}{\end{definition}}
\newenvironment{rmk}{\begin{remark}\rm}{\end{remark}}
\newenvironment{ex}{\begin{example}\rm}{\end{example}}

\numberwithin{equation}{chapter}
\numberwithin{figure}{chapter}

\newcounter{FNC}[page]
\def\fauxfootnote#1{{\addtocounter{FNC}{2}$^\fnsymbol{FNC}$%
     \let\thefootnote\relax\footnotetext{$^\fnsymbol{FNC}$#1}}}
\newcommand{\DeCo}[1]{\Blue{#1}}
\newcommand{\col}{\mbox{\rm col}}
\newcommand{\conv}{\mbox{\rm conv}}
\newcommand{\gcf}{\mbox{\rm gcf}}
\newcommand{\Gr}{\mbox{\rm Gr}}
\newcommand{\sGr}{\mbox{\scriptsize\rm Gr}}
\newcommand{\Hom}{\mbox{\rm Hom}}
\newcommand{\ini}{\mbox{\rm in}}
\newcommand{\mdeg}{\mbox{\rm mdeg}}
\newcommand{\ord}{\mbox{\rm ord}}
\newcommand{\sign}{\mbox{\,\rm sign\,}}
\newcommand{\sgn}{\mbox{\rm sign}}
\newcommand{\Span}{\mbox{\rm span}}
\newcommand{\St}{\mbox{\rm St}}
\newcommand{\ubc}{\mbox{\rm ubc}}
\newcommand{\vol}{\mbox{\rm vol}}
\newcommand{\var}{\mbox{\rm var}}
\newcommand{\Wr}{\mbox{\rm Wr}\,}
\newcommand{\Edot}{E_\bullet}
\newcommand{\Fdot}{F_\bullet}
\newcommand{\ba}{{\bf a}}
\newcommand{\bg}{{\bf g}}
\newcommand{\e}{{\bf e}}
\newcommand{\bs}{{\bf s}}
\newcommand{\calA}{{\mathcal A}}
\newcommand{\calB}{{\mathcal B}}
\newcommand{\calH}{{\mathcal H}}
\newcommand{\calR}{{\mathcal R}}
\newcommand{\calX}{{\mathcal X}}
\newcommand{\A}{{\mathbb A}}
\newcommand{\C}{{\mathbb C}}
\newcommand{\G}{{\mathbb G}}
\newcommand{\N}{{\mathbb N}}
\renewcommand{\P}{{\mathbb P}}
\newcommand{\Q}{{\mathbb Q}}
\newcommand{\R}{{\mathbb R}}
\newcommand{\T}{{\mathbb T}}
\newcommand{\Z}{{\mathbb Z}}
\newcommand{\slm}{{\mathfrak{sl}_m}}
\newcommand{\bI}{{\includegraphics{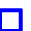}}}
\newcommand{\I}{{\includegraphics{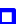}}}
\newcommand{\QED}{\hfill\raisebox{-2pt}{\includegraphics[height=12pt]{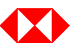}}\medskip}

\begin{document}
\title{\mbox{\ }\vspace{-100pt}\\Real Solutions to Equations From Geometry\vspace{35pt}}

\author{Frank Sottile\\
        Department of Mathematics\\
         Texas A\&M University\\
         College Station\\
         TX \ 77843\\
         USA \\
  {\tt sottile@math.tamu.edu}\\{\tt http://www.math.tamu.edu/\~{}sottile}
 \\\begin{picture}(10,10)\put(-230,-260){\includegraphics[height=240pt]{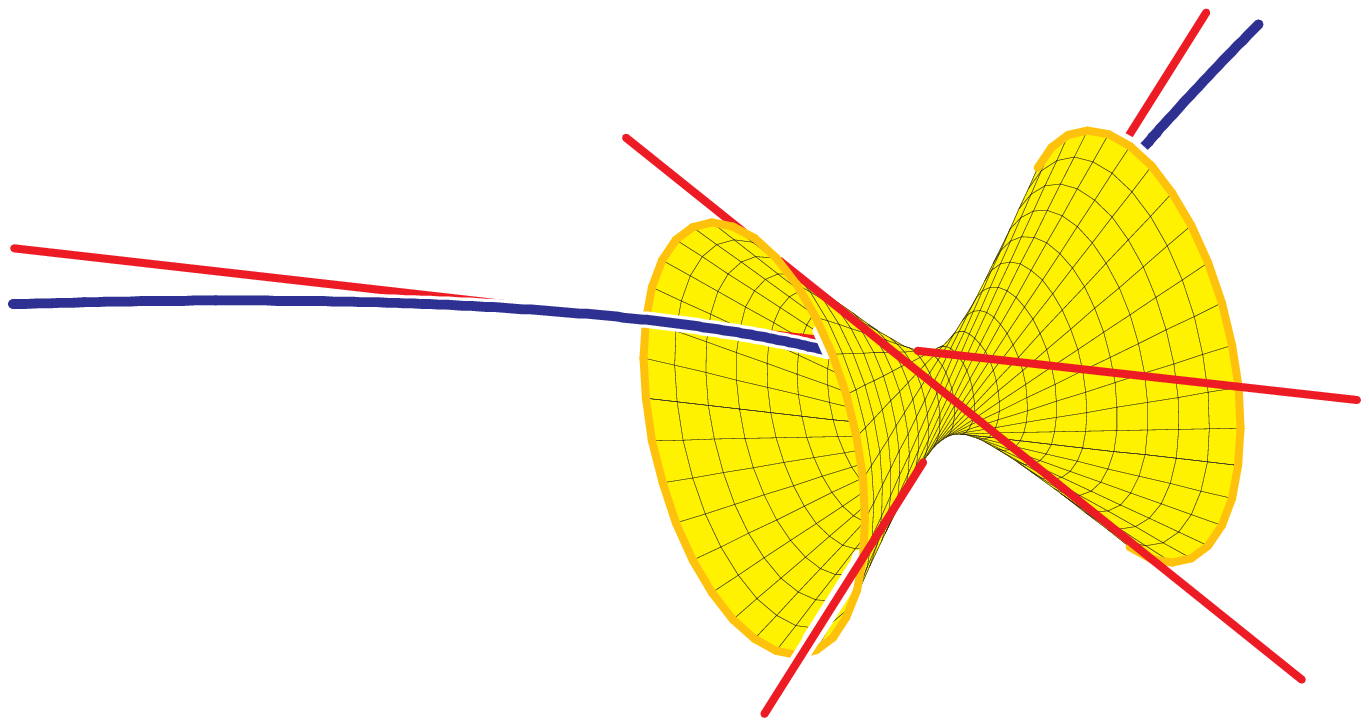}}\end{picture} 
}
%
%
%

\maketitle
\thispagestyle{empty} \mbox{\ }\newpage\setcounter{page}{1}
\tableofcontents
\begin{abstract}
   Understanding, finding, or even deciding on the existence of real solutions to a
   system of equations is a very difficult problem with many applications. While
   it is hopeless to expect much in general, we know a surprising amount about
   these questions for systems which possess additional structure. Particularly
   fruitful---both for information on real solutions and for applicability---are
   systems whose additional structure comes from geometry. Such equations from
   geometry for which we have information about their real solutions are the 
   subject of these notes. 

    We will focus on equations from toric varieties and homogeneous spaces,
  particularly Grassmannians. Not only is much known in these cases, but they
  encompass some of the most common applications. The results we discuss may be
  grouped into three themes: 
\begin{enumerate}
  \item[(I)] Upper bounds on the number of real solutions.
  \item[(II)] Geometric problems that can have all solutions be real.
  \item[(III)] Lower bounds on the number of real solutions.
\end{enumerate}
Upper bounds as in (I) bound the complexity of the set of real solutions---they
are one of the sources for the theory of o-minimal structures which are an
important topic in real algebraic geometry. 
The existence (II) of geometric problems that can have all solutions be real
was initially surprising, but this phenomena now appears ubiquitous.
Lower bounds as in (III) give an existence
proof for real solutions. Their most spectacular manifestation is the
nontriviality of the Welschinger invariant, which was computed via
tropical geometry.
One of the most surprising manifestations of this phenomenon is when the upper bound
equals the lower bound, which is the subject of the Shapiro conjecture.
\medskip

I thank the Institut Henri Poincar\'e, where a preliminary version of these notes were
produced during a course I taught there in November 2005. 
These notes were revised and expanded during course at Texas A\&M University in 2007 and a
lecture series at the Centre Interfacultaire Bernoulli at EPFL in 2008, and completed in 2010
with material from a lecture I gave in January 20-09 the Theorem of Mukhin, Tarasov, and
Varchenko, and lectures at the GAeL meeting in Leiden in June 2009.
I also thank Prof.~Dr.~Peter Gritzmann of the Technische Universit\"at M\"unchen, whose
hospitality enabled the completion of the first version of these notes.
During this period, my research was supported by NSF grants DMS-0701059 and CAREER grant
DMS-0538734. 
The point of view in these notes was developed through the encouragement of Bernd Sturmfels,
Askold Khovanskii, and Marie-Fran\c{c}oise Roy, and through my interactions with the many whose
work is mentioned were, including my collaborators from whom I have learned a great deal.
\vspace{30pt}

Frank Sottile

03.01.10, College Station, Texas.
\end{abstract}
\thispagestyle{empty} \mbox{\ }\newpage\setcounter{page}{1}

\chapter{Overview}\label{ch:1}

In mathematics and its applications, we are often faced with a system of polynomial
equations, and we need to study or find the solutions.
Such systems that arise naturally typically possess some geometric or combinatorial
structure that may be exploited to study their solutions.
Such structured systems are studied in enumerative algebraic geometry, which has given us
the deep and powerful tools of intersection theory~\cite{Fu} to count and analyze their
{\sl complex} solutions.
A companion to this theoretical work are algorithms, both symbolic (based on Gr\"obner
bases~\cite{CLO,GBCP}) and numerical (based on homotopy continuation~\cite{SW05}) for
solving and analyzing systems of polynomial equations.

Despite these successes, this line of research largely sidesteps the often primary goal of
formulating problems as solutions to systems of equations---namely to determine or study
their real solutions.
This deficiency is particularly acute in applications, from control~\cite{Byrnes},
Kinematics~\cite{BoRo}, statistics~\cite{PRW}, and computational biology~\cite{ASCB}, for
it is typically the real solutions that applications seek.
One reason that traditional algebraic geometry ignores the real solutions is that there
are few elegant theorems or general results available to study real solutions.
Nevertheless, the demonstrated importance of understanding the real solutions to systems
of equations demands our attention.

In the 19th century and earlier, many elegant and powerful methods were developed to study
the real roots of univariate polynomials (Sturm sequences, Budan-Fourier Theorem,
Routh-Hurwitz criterion), which are now standard tools in some applications of
mathematics. 
In contrast, it has only been in the pat few decades that serious attention has been paid
toward understanding the real solutions to systems of polynomial equations.

This work has concentrated on systems possessing some, particularly geometric, structure.
The reason for this is two-fold:
systems from nature typically possess some special structure that should be exploited in
their study, and it is unlikely that any results of substance hold for general or
unstructured systems.
In this period, a story has emerged of bounds (both upper and lower) on the number of real
solutions to certain classes of systems, as well as the discovery and study of systems
that have only real solutions.
This Overview will sketch this emerging landscape and the subsequent chapters will treat
these developments in more detail.

%
\section*{Introduction}
%

Our goal will be to say something
meaningful about the real solutions to a system of equations.
For example, consider a system 
 \begin{equation}\label{E1:system}
  f_1(x_1,\dotsc,x_n)\ =\ 
  f_2(x_1,\dotsc,x_n)\ =\ \dotsb\ =\ 
  f_N(x_1,\dotsc,x_n)\ =\ 0\,,
 \end{equation}
of $N$ real polynomials in $n$ variables.
Let \DeCo{$r$} be its number of real solutions and \DeCo{$d$} its number of complex
solutions\fauxfootnote{We shall always 
   assume that our systems are \DeCo{{\sl generic}} in the sense
  that all of their solutions occur with multiplicity 1, and the
  number \DeCo{$d$} of complex solutions is stable under certain
  \DeCo{allowed} perturbations of the coefficients.}.
Since every real number is complex, and since nonreal solutions come
in conjugate pairs, we have the following trivial inequality
\[
  d\ \geq\ r\ \geq\ d\!\mod 2\ \;\in\ \;\{0,1\}\,.
\]
We can say nothing more unless the equations have some structure,
and a particularly fruitful class of structures are those which come from
geometry.
The main point of this book is that we can identify
structures in equations that will allow us to do better
than this trivial inequality.

Our discussion will have three themes:
\begin{enumerate}
 \item[(I)] Sometimes, the upper bound $d$ is not sharp and there is a
       smaller bound for $r$.
 \item[(II)] For many problems from enumerative geometry, the
       upper bound is sharp.
 \item[(III)] The lower bound for $r$ may be significantly
       larger than $d\mod  2$.
\end{enumerate}
  A lot of time will be spent on the Shapiro Conjecture (Theorem of Mukhin, Tarasov, and
  Varchenko~\cite{MTV_Sh}) and its generalizations, which  
  is a situation where the upper bound of $d$ is also the lower
  bound---all solutions to our system are real.

We will not describe how to actually find the solutions to a system~\eqref{E1:system} and 
there will be no discussion of algorithms nor any 
complexity analysis.
The book of Basu, Pollack, and Roy~\cite{BPR03} is an excellent place to learn about 
algorithms for computing real algebraic varieties and finding real solutions.
%
%
We remark that some of the techniques we employ to study real solutions
underlie numerical algorithms to compute the solutions.\medskip

One class of systems that we will study are systems of sparse
polynomials.
Integer vectors $a=(a_1,\dotsc,a_n)\in\Z^n$ are exponents for (Laurent) monomials
 \[
   \Z^n\ni a\ \leftrightarrow\ \DeCo{x^a}\ :=\ x_1^{a_1}x_2^{a_2}\dotsb
   x_n^{a_n}\ \in\ \C[x_1,\dotsc,x_n,x_1^{-1},\dotsc,x_n^{-1}]\,.
 \]
Sometimes, we will just call elements of $\Z^n$ \DeCo{{\sl monomials}}.
Let $\DeCo{\calA}\subset\Z^n$ be a finite set of monomials.
A linear combination 
\[
   f\ =\ \sum_{a\in\calA} c_a x^a\qquad c_a\in\R
\]
of monomials from $\calA$ is a \DeCo{{\sl sparse polynomial}} with 
\DeCo{{\sl support} $\calA$}.
Sparse polynomials naturally define functions on the complex torus $(\C^\times)^n$.
A system~\eqref{E1:system} of $N=n$ polynomials in $n$ variables, where
each polynomial has support $\calA$, will be called a 
\DeCo{{\sl system}} (of polynomials) with \DeCo{{\sl support} $\calA$}.
These are often called \DeCo{{\sl unmixed systems}} in contrast to 
\DeCo{{\sl mixed systems}} where each polynomial may have different support.
While sparse systems occur naturally---multilinear or multihomogeneous polynomials are
an example---they also occur in problem formulations for the simple reason 
that we humans seek simple formulations of problems, and this may mean polynomials with
few terms.

A fundamental result about unmixed systems is the Kushnirenko bound
on their number of complex solutions.
The \DeCo{{\sl Newton polytope}} of a polynomial with support $\calA$
is the convex hull \DeCo{$\Delta_\calA$} of the set $\calA$ of monomials.
Write \DeCo{$\vol(\Delta)$} for the Euclidean volume of a polytope $\Delta$.

\begin{thm}[Kushnirenko~\cite{BKK}]\label{T1:Koushnirenko}
 A system of $n$ polynomials in $n$ variables with common support $\calA$ 
 has at most $n!\vol(\Delta_\calA)$ isolated solutions in 
 $(\C^\times)^n$, and exactly this number when the polynomials are generic 
 polynomials with support $\calA$.
\end{thm}

Bernstein generalized this to mixed systems.
The Minkowski sum $P+Q$ of two polytopes in $\R^n$ is their pointwise sum as sets of
vectors in $\R^n$.
Let $P_1,\dotsc,P_n\subset\R^n$ be polytopes.
The volume
\[
   \vol( t_1 P_1\ +\  t_2P_2\ +\ \dotsb\ +\ t_nP_n)
\]
is a homogeneous polynomial of degree $n$ in the variables
$t_1,\dotsc,t_n$~\cite[Exercise 15.2.6]{Grunbaum}. 
The \DeCo{{\sl mixed volume}  MV$(P_1,\dotsc,P_n)$} of $P_1,\dotsc,P_n$ is the coefficient
of the monomial $t_1t_2\dotsb t_n$ in this polynomial.

\begin{thm}[Bernstein~\cite{Be75}]\label{T1:Bernstein}
 A system of $n$ polynomials in $n$ variables where the polynomials have supports
 $\calA_1,\dotsc,\calA_n$ has at most $\mbox{\rm
 MV}(\Delta_{\calA_1},\dotsc,\Delta_{\calA_n})$  isolated solutions in  
 $(\C^\times)^n$, and exactly this number when the polynomials are generic for their given
 support.
\end{thm}

Since $\mbox{\rm MV}(P_1,\dotsc,P_n)=n!\vol(P)$ when $P_1=\dotsb=P_n=P$, this generalizes
Kushnirenko's Theorem. 
We will prove Kushnirenko's Theorem in Chapter~\ref{Ch:sparse},
but will not present a proof of Bernstein's Theorem.
Instead, we suggest two excellent sources by Sturmfels.
Both are similar, but the first is self-contained and superbly written.

\begin{itemize}

\item{\cite{St98}} \ \ 
    \emph{Polynomial equations and convex polytopes}, Amer. Math.
     Monthly \textbf{105} (1998), no.~10, 907--922.

\item{\cite{SPE}} \ \ Chapter 3 in 
    \emph{Solving systems of polynomial equations}, CBMS, vol.~97,
    American Mathematical Society, Providence, RI, 2002.

\end{itemize}

The bound of Theorem~\ref{T1:Koushnirenko} and its generalization 
 Theorem~\ref{T1:Bernstein} is often called the 
\DeCo{{\sl BKK bound}} for \DeCo{B}ernstein, \DeCo{K}hovanskii, 
and \DeCo{K}ushnirenko~\cite{BKK}.

\section{Upper bounds}

While the number of complex roots of a univariate polynomial is
typically equal to its degree, the number of real roots depends upon
the length of the expression for the polynomial.
Indeed, by Descartes's rule of signs~\cite{D1637} (see Section~\ref{S2:Descartes}),
a univariate polynomial with $m{+}1$ terms has at most $m$ positive roots,
and thus at most $2m$ nonzero real roots.
For example, the polynomial $x^d-a$ with $a\neq 0$ has 0, 1, or 2 real roots,
but always has $d$ complex roots.
Khovanskii generalized this type of a bound to multivariate polynomials with
his fundamental \DeCo{{\sl fewnomial bound}}.

\begin{thm}[Khovanskii~\cite{Kh80}]\label{T1:Khovanski}
 A system of\/ $n$ polynomials in $n$ variables having a total of \/
 $l{+}n{+}1$ distinct monomials has at most 
\[
   2^{\binom{l+n}{2}}(n+1)^{l+n}
\]
 nondegenerate positive real solutions.
\end{thm}

There are two reasons for this restriction to positive solutions.
Most fundamentally is that Khovanskii's proof requires this restriction.
This restriction also excludes the following type of trivial zeroes:
Under the substitution $x_i\mapsto x_i^2$, each positive solution becomes $2^n$ real
solutions, one in each of the $2^n$ orthants.
More subtle substitutions lead to similar extra trivial zeroes which differ 
from the positive solutions only by some sign patterns.

This is the first of many results which verified the principle of Bernstein and
Kushnirenko that the topological complexity of a set defined by real polynomials should
depend on the number of terms  in the polynomials and not on their degree.
Khovanskii's work was also a motivation for the notion of o-minimal structures~\cite{vdD,PiSt}.
The main point of Khovanskii's theorem is the existence of such a bound and not the actual
bound itself.
For each $l,n\geq 1$, we define the \DeCo{{\sl Khovanskii number} $X(l,n)$} to be the
maximum number of nondegenerate positive solutions to a system of $n$ polynomials in $n$
variables with $l+n+1$ monomials.
Khovanskii's Theorem gives a bound for $X(l,n)$, but that bound is enormous.
For example, when $l=n=2$, the bound is 5184.
Because of this, it was expected to be far from sharp.
Despite this expectation, the first nontrivial improvement was only given in 2003.

\begin{thm}[Li, Rojas, and Wang~\cite{LRW03}]\label{T1:LRW}
 Two trinomials in two variables have at most $5$ nondegenerate positive real
 solutions. 
\end{thm}

This bound sharp.
Haas~\cite{Ha02} had shown that the system of two
trinomials in $x$ and $y$ 
 \begin{equation}\label{E1:Haas}
  10x^{106} + 11 y^{53} - 11y\ =\ 
  10y^{106} + 11 x^{53} - 11x\ =\ 0\,,
 \end{equation}
has 5 positive solutions.

Since we may multiply one of the trinomials in~\eqref{E1:Haas} by an arbitrary monomial
without changing the solutions, we can assume that the two trinomials~\eqref{E1:Haas} share a
common monomial, and so there are at most $3+3-1=5=2+2+1$ monomials between the two
trinomials, and so two trinomials give a fewnomial system with $l=n=2$.
While 5 is less than 5184, Theorem~\ref{T1:LRW} does
not quite show that $X(2,2)=5$ as two trinomials do not constitute a general
fewnomial system with $l=n=2$.
Nevertheless, Theorem~\ref{T1:LRW} gave strong evidence that Khovanskii's bound may be
improved. 
Such an improved bound was given in~\cite{BS07}.

\begin{thm}\label{T1:New_Bounds}
 ${\displaystyle X(l,n)\ <\ \tfrac{e^2+3}{4}2^{\binom{l}{2}}n^l}$.
\end{thm}

For small values of $l$, it is not hard to improve this.
For example, when $l=0$, the support $\calA$ of the system is a simplex, and there will be
at most 1 positive real solution, so $X(0,n)=1$.
Theorem~\ref{T1:New_Bounds} was inspired by the sharp bound of
Theorem~\ref{Th1:BBS} when $l=1$~\cite{BBS}. 
A set $\calA$ of exponents is \DeCo{{\sl primitive}} if $\calA$ affinely spans
the full integer lattice $\Z^n$.

\begin{thm}\label{Th1:BBS}
 If $l=1$ and the set $\calA$ of exponents is primitive, then 
 there can be at most $2n{+}1$ nondegenerate nonzero real solutions, and this
 is sharp in that for any $n$ there exist systems with $n{+}2$ monomials and $2n{+}1$
 nondegenerate real solutions whose exponent  vectors affinely span $\Z^n$. 
\end{thm}

Observe that this bound is for all real solutions, not just positive solutions.
We will discuss this in Section~\ref{S4:circuit}.
Further analysis gives the sharp bound for $X(1,n)$.

\begin{thm}[Bihan~\cite{Bihan}]
  $X(1,n)=n+1$.
\end{thm}

In contrast to these results establishing absolute upper bounds for the
number of real solutions which improve the trivial bound of the number
$d$ of complex roots, there are a surprising number of problems that
come from geometry for which all solutions can be real.
For example, Sturmfels~\cite{St94} proved the following.

\begin{thm}\label{T1:Sturm}
  Suppose that a lattice polytope $\Delta\subset\Z^n$ admits a regular 
  triangulation with each simplex having minimal volume $\frac{1}{n!}$.
  Then there is a system of sparse polynomials with support
  $\Delta\cap\Z^n$ having all solutions real.
\end{thm}

For many problems from enumerative geometry, it is similarly possible 
that all solutions can be real. 
This will be discussed in Chapter~\ref{S:ERAG}.
The state of affairs in 2001 was presented in~\cite{So03}.

%
\section{The Wronski map and the Shapiro Conjecture}
%

The \DeCo{{\sl Wronskian}} of univariate polynomials 
$f_1(t), f_2(t), \dotsc,f_m(t)$ is the determinant
\[
   \DeCo{\Wr(f_1,\dotsc,f_m)}\ :=\ \det \bigl( 
    ({\textstyle \frac{\partial}{\partial t}})^{j-1} f_i(t)
     \bigr)_{i,j=1,\dotsc,m}
  \ .
\]
When the polynomials $f_i$ have degree $m{+}p{-}1$ and are linearly independent, the
Wronskian has degree at most $mp$.
For example, if $m=2$, then $W(f,g)=f'g-fg'$, which has degree $2p$ as the
coefficients of $t^{2p+1}$ in this expression cancel.
Up to a scalar, the Wronskian depends only upon the linear span of the
polynomials $f_1,f_2,\dotsc,f_m$.
Removing these ambiguities gives the \DeCo{{\sl Wronski map}},
 \begin{equation}\label{E1:Wronski}
  \DeCo{\Wr}\ :=\ \Gr(m,\C_{m{+}p{-}1}[t])\ \longrightarrow\ \P^{mp}\,,
 \end{equation}
where \DeCo{$\Gr(m,\C_{m{+}p{-}1}[t])$} is the \DeCo{{\sl Grassmannian}} of $m$-dimensional
subspaces of the linear space $\C_{m{+}p{-}1}[t]$ of polynomials of degree $m{+}p{-}1$ in 
the variable $t$, and $\P^{mp}$ is the projective space of
polynomials of degree at most $mp$, which has dimension equal to the
the dimension of the Grassmannian.

Work of Schubert in 1886~\cite{Sch1886b}, combined with a result of
Eisenbud and Harris in 1983~\cite{EH83} shows that the Wronski map is surjective
and the general polynomial $\Phi\in\P^{mp}$ has
 \begin{equation}\label{E1:WronDeg}
   \DeCo{\#_{m,p}}\ :=\ 
   \frac{1!2!\dotsb(m{-}1)!\cdot(mp)!}%
        {m!(m{+}1)!\dotsb (m{+}p{-}1)!}
 \end{equation}
preimages under the Wronski map.
These results concern the complex Grassmannian and complex projective space.

Boris Shapiro and Michael Shapiro made a 
conjecture in 1993/4 about the Wronski map from 
the real Grassmannian to real projective space.
This was proven when $\min(k,d{+}1{-}k)=2$ by Eremenko and Gabrielov~\cite{EG02}, and finally
settled by Mukhin, Tarasov, and Varchenko~\cite{MTV_Sh}. 
They have subsequently found another proof~\cite{MTV_R}.

\begin{thm}\label{T1:Shap_Conj}
 If the polynomial $\Phi\in\P^{mp}$ has only real zeroes, then every point in 
 $\Wr^{-1}(\Phi)$ is real.
 Moreover, if $\Phi$ has $mp$ simple real zeroes 
 then there are $\#_{m,p}$ real points in $\Wr^{-1}(\Phi)$.
\end{thm}

We will sketch the proof of Mukhin, Tarasov, and Varchenko in Chapter~\ref{ch:MTV}.
This \DeCo{{\sl Shapiro Conjecture}} has appealing geometric interpretations,  enjoys links
to several areas of mathematics, and has many theoretically satisfying generalizations
which we will discuss in Chapters~\ref{Ch:ScGr},~\ref{Ch:EG}, and~\ref{Ch:Frontier}.
We now mention two of its interpretations.

%
\subsection{The problem of four lines}
A geometric interpretation of the Wronski map and the Shapiro Conjecture
when $m=p=2$ is a variant of the classical problem of the
lines in space which meet four given lines.
Points in $\Gr(2,\C_3[t])$ correspond to lines in $\C^3$ as follows.
The \DeCo{{\sl moment curve} $\gamma$} in $\C^3$ is the curve with
parameterization 
\[
   \gamma(t)\ :=\ (t, t^2, t^3)\,.
\]
A cubic polynomial $f$ is the composition of $\gamma$ and an affine-linear map
$\C^3\to\C$, and so a two-dimensional space of cubic polynomials is to a
two-dimensional space of affine-linear maps whose common kernel is the corresponding line
in $\C^3$.
(This description is not exact, as some points in $\Gr(2,\C_3[t])$ correspond to lines at
infinity.)

Given a polynomial $\Phi(t)$ of degree 4 with distinct real roots, 
points in the fiber $\Wr^{-1}(\Phi)$ correspond to the lines in space
which meet the four lines tangent to the moment curve $\gamma$ at its
points coming from the roots of $\Phi$.
There will be two such lines, and the Shapiro conjecture asserts that
both will be real.

\begin{figure}[htb]
 \[
  \begin{picture}(372,185.2)(0,1)
   \put(0,0){\includegraphics[height=190pt]{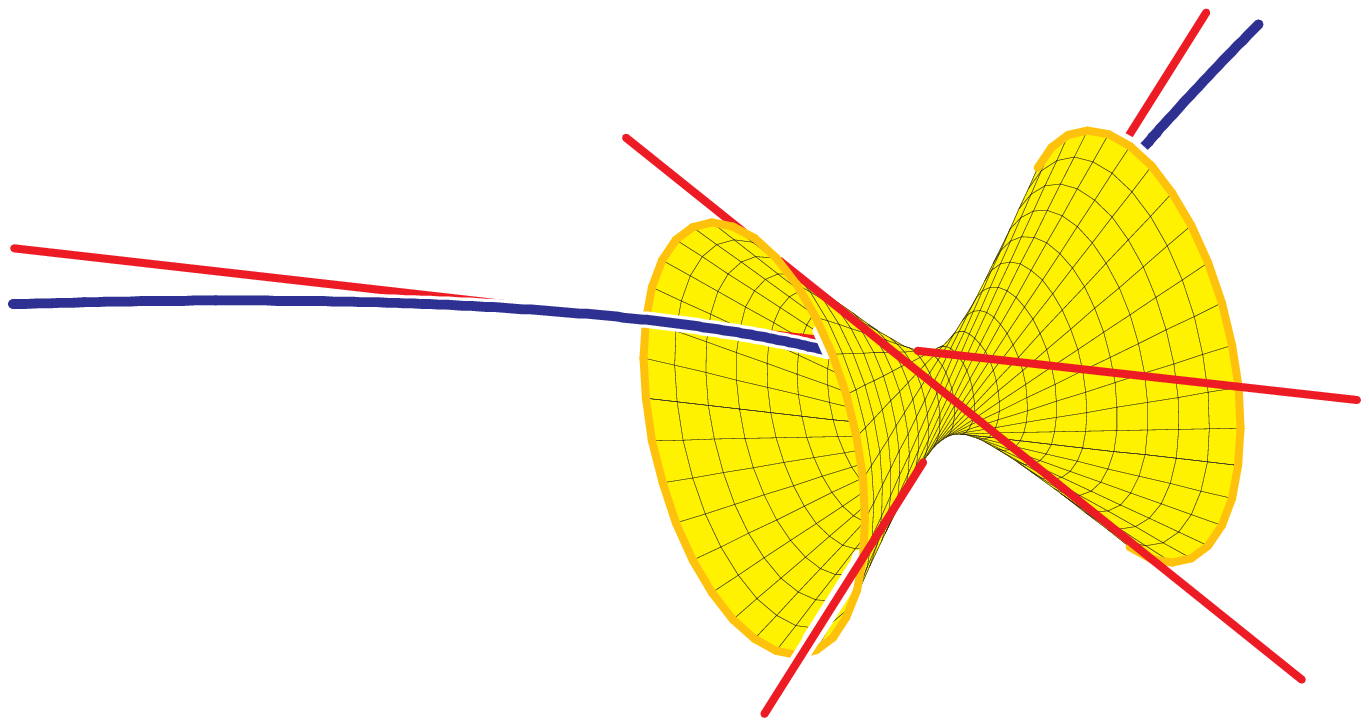}}
   \put(170,0){\Red{$\ell(-1)$}} \put(354,0){\Red{$\ell(0)$}}
   \put(0,136){\Red{$\ell(1)$}} \put(0,94){\Blue{$\gamma$}}
   \put(280,178){\Red{$\ell(-1)$}} \put(337,174){\Blue{$\gamma$}}
   \put(320,135){\Brown{$Q$}}
  \end{picture}
 \]
 \caption{Quadric containing three lines tangent to $\gamma$\label{F1:hyperboloid}.} 
 \end{figure}

It is not hard to see this directly.
Any fractional linear change of parameterization of the moment curve
is realized by a projective linear transformation of 3-dimensional space which
stabilizes the image of the moment curve.
Thus we may assume that the polynomial $\Phi(t)$ is equal to $(t^3-t)(t-s)$,
which has roots $-1,0,1$, and $s$, where $s\in (0,1)$.
Applying an affine transformation to 3-dimensional space, the moment curve becomes
the curve with parameterization
\[
   \gamma\ :\ t\ \longmapsto\ 
     (6t^2-1, \ \tfrac{7}{2}t^3+\tfrac{3}{2}t,\ \tfrac{3}{2}t-\tfrac{1}{2}t^3)\,.
%
%
\]
Then the lines tangent to $\gamma$ at the roots $-1,0,1$ of $\Phi$ have
parameterizations
\[
  (-5-s,\;5+s,\;-1)\,,\ (-1,\;s,\;s)\,,\ (5+s\;,5+s\;,1)\qquad s\in\R\,.
\]
These lie on a hyperboloid \Brown{$Q$} of one sheet, which is defined by
 \begin{equation}\label{Eq1:Hyperboloid}
   1 - x_1^2 + x_2^2 - x_3^2\ =\ 0\,.
 \end{equation}
We display this geometric configuration in Figure~\ref{F1:hyperboloid}.
There, $\ell(i)$ is the line tangent to $\gamma$ at the point $\gamma(i)$.
The quadric \Brown{$Q$} has two rulings.
One ruling contains our three tangent lines and 
the other ruling (which is drawn on \Brown{$Q$}) consists of the lines which
meet our three tangent lines. 

Now consider the fourth line \ForestGreen{$\ell(s)$} which is tangent to
\Brown{$\gamma$} at the point \Brown{$\gamma(s)$}.
This has the parameterization
\[
  \ForestGreen{\ell(s)}\ =\ 
   \bigl( 6s^2-1\,,\ 
       \tfrac{7}{2}s^3+\tfrac{3}{2}s\,,\ 
        \tfrac{3}{2}s-\tfrac{1}{2}s^3 \bigr)\ +\ 
   t \bigl(12s\,,\ 
       \tfrac{21}{2}s^2+\tfrac{3}{2}\,,\ 
        \tfrac{3}{2} -\tfrac{3}{2}s^2\bigr)\,.
\]
We compute the intersection of the fourth line with $Q$.
Substituting its parameterization into~\eqref{Eq1:Hyperboloid} and dividing by $-12$ gives
the equation  
\[
   (s^3-s)(s^3-s + t(6s^2-2)+9st^2)\ =\ 0\,.
\]
The first (nonconstant) factor $s^3-s$ vanishes when \ForestGreen{$\ell(s)$} is equal to
one of \Red{$\ell(-1)$}, \Red{$\ell(-0)$}, or \Red{$\ell(-1)$}--for these values of $s$ 
every point of \ForestGreen{$\ell(s)$} lies on the quadric \Brown{$Q$}.
The second factor has solutions
\[
  t\ =\ -\; \frac{3s^2-1\pm\sqrt{3s^2+1}}{9s}\ .
\]
Since $3s^2+1>0$ for all $s$, both solutions will be real.

\begin{figure}[htb]
 \[
  \begin{picture}(285, 159)(-34.5,-11.5)
   \put(0,0){\includegraphics[height=130pt]{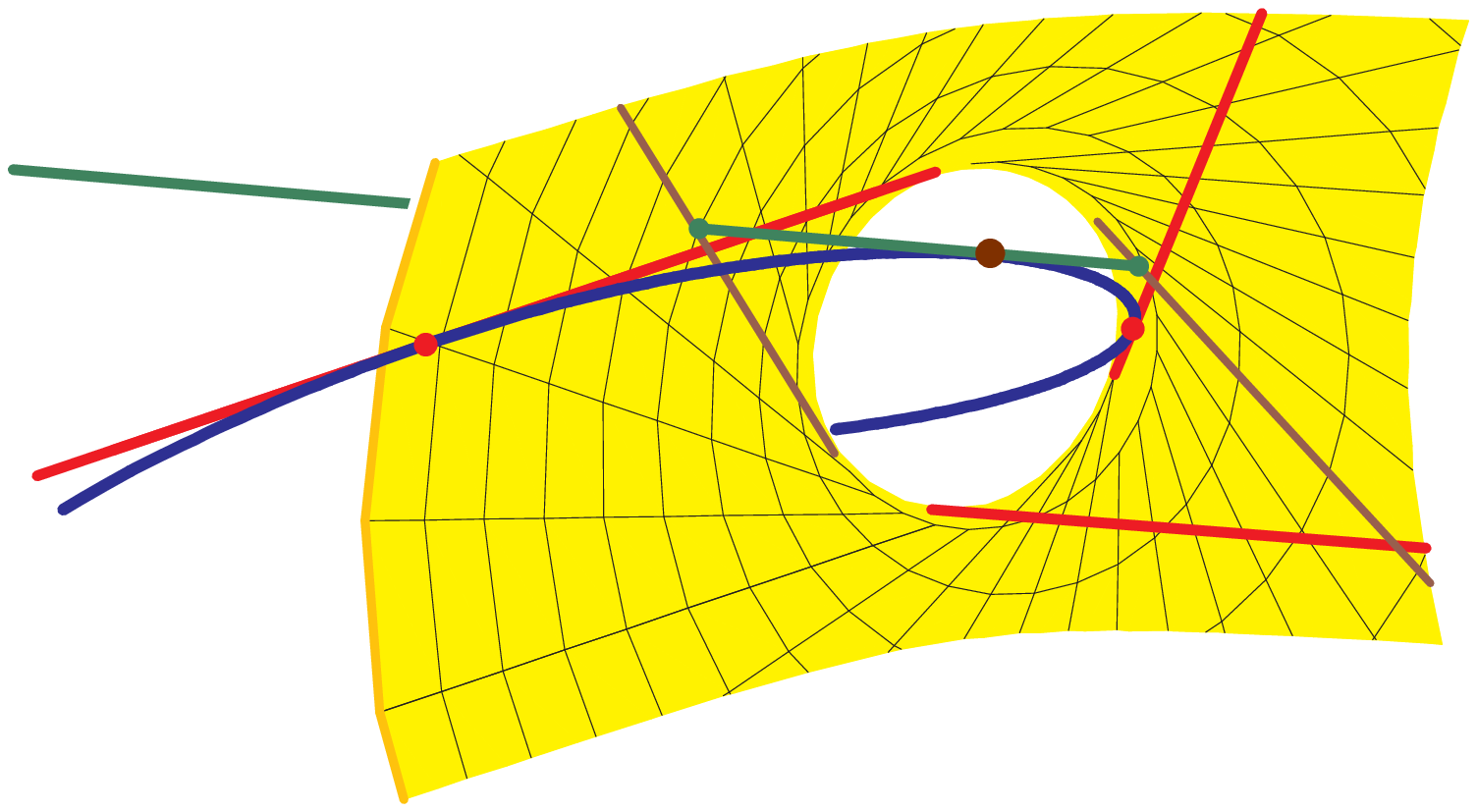} }
   \put(-12,63.5){\Red{$\ell(1)$}} \put(237,34.6){\Red{$\ell(-1)$}} 
   \put(193,138){\Red{$\ell(0)$}}
   \put(8.4,34.5){\Blue{$\gamma$}} \put(2.4,111){\ForestGreen{$\ell(s)$}}
   \put(117.6,131){\Brown{$Q$}}     
  \thicklines
    \put(141.6,-13){\Brown{$\gamma(s)$}}
    \put(156,2.4){\Brown{\vector(0,1){82}}} 
  \end{picture}
\]
 \caption{The fourth tangent line meets \Brown{$Q$} in two real
   points\label{F1:line4}.} 
 \end{figure}

We may also see this geometrically.
Consider the fourth line \ForestGreen{$\ell(s)$} for $0<s<1$.
In Figure~\ref{F1:line4}, we look down the throat of the hyperboloid
at the interesting part of this configuration.
This picture demonstrates that \ForestGreen{$\ell(s)$} must meet \Brown{$Q$} in
two real points. 
Through each point, there is a real line in the second ruling which
meets all four tangent lines, and this proves Shapiro's conjecture
for $m=p=2$.

%
\subsection{Rational functions with real critical points}
When $m=2$, the Shapiro conjecture may be interpreted in terms of
rational functions.
A rational function $\rho(t)=f(t)/g(t)$ is a quotient of two univariate
polynomials, $f$ and $g$.
This defines a map $\rho\colon \P^1\to\P^1$ whose critical points are those $t$ for which
$\rho'(t)=0$. 
Since $\rho'(t)=(f'g-g'f)/g^2$, we see that the critical points are
the roots of the Wronskian of $f$ and $g$.
Composing the rational function $\rho\colon\P^1\to\P^1$ with an automorphism of the target
$\P^1$ gives an equivalent rational function, and the equivalence class of $\rho$ is
determined by the linear span of its numerator and denominator.
Thus Shapiro's conjecture asserts that a rational function having only real critical  
points is equivalent to a real rational function.

Eremenko and Gabrielov~\cite{EG02} proved exactly this statement in 2002, and
thereby established the Shapiro Conjecture in the case $m=2$.

\begin{thm}
 A rational function with only real critical points is equivalent to a
 real rational function.
\end{thm}

In Chapter~\ref{Ch:EG} we will present an elementary proof of this result
that Eremenko and Gabrielov found in 2005~\cite{EG05}.

%
\section{Lower bounds}\label{S1:lower}
%

We begin with some of perhaps the most exciting recent
development in real algebraic geometry.
It begins with the fundamental observation of Euclid that two points determine a line.
Many people who have studied geometry know that five points on the plane determine a
conic.  In general, if you have $m$ random points in the plane and you
want to pass a rational curve of degree $d$ through all of them,
there may be no solution to this interpolation problem (if $m$ is
too big), or an infinite number of solutions (if $m$ is too small),
or a finite number of solutions (if $m$ is just right).  It turns
out that ``$m$ just right'' means $m=3d{-}1$ ($m=2$ for lines and
$m=5$ for conics). 

A harder question is, if $m=3d{-}1$, {\it how many\/} rational curves of degree
$d$ interpolate the points?  Let's call this number $N_d$, so that 
$N_1=1$ and $N_2=1$ because the line and conic of the previous
paragraph are unique.  It has long been known that $N_3=12$, and
in 1873 Zeuthen~\cite{Ze1873} showed that $N_4=620$.
That was where matters stood until 1989, when Ran~\cite{R89} gave a recursion for these
numbers. 
In the 1990's, Kontsevich and Manin~\cite{KM} used associativity in quantum
cohomology of $\mathbb{P}^2$ to give the elegant recursion
 \begin{equation}\label{E1:Konts}
  N_d\ =\   \sum_{a+b=d} N_a N_b \left( a^2b^2\binom{3d-4}{3a-2} - 
           a^3b\binom{3d-4}{3a-1}\right)\ ,
 \end{equation}
which begins with the Euclidean declaration that two points determine a line ($N_1=1$).
These numbers grow quite fast, for example $N_5=87304$.

The number of real rational curves which interpolate 
a given $3d-1$ points in the real plane $\R\P^2$ will depend rather
subtly on the configuration of the points.
To say anything about the real rational curves  would seem impossible.
However this is exactly what Welschinger~\cite{W} did.
He found an invariant which does not depend upon the
choice of points.

A rational curve in the plane is necessarily singular---typically it has 
$\binom{d-1}{2}$ nodes.
Real curves have three types of nodes.
Only two types are visible in $\R\P^2$, and we are familiar with them from
rational cubics.
The curve on the left below has a node with two real branches, and the curve
on the right has a \DeCo{{\sl solitary point}} `{\small$\bullet$}', where two complex conjugate
branches meet. 
 \[
   \includegraphics[height=70pt]{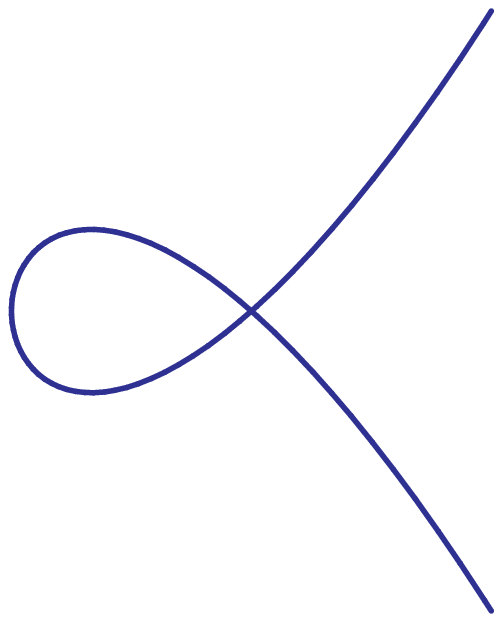}\qquad\qquad\qquad\qquad\qquad
   \includegraphics[height=70pt]{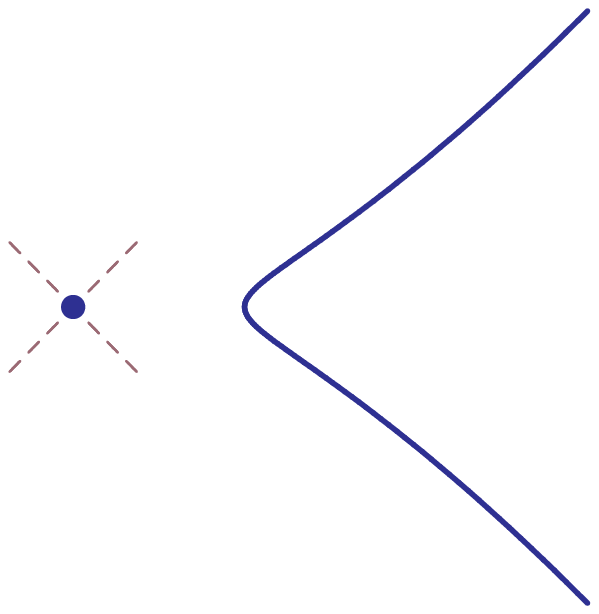}
 \]
The third type of node is a pair of complex conjugate nodes, which are 
not visible in $\R\P^2$.

\begin{thm}[Welschinger~\cite{W}]
  The sum,  
 \begin{equation}\label{Eq1:Wel}
   \sum (-1)^{\#\{\mbox{\scriptsize solitary points, }\bullet, \mbox{\scriptsize in }C\}}\,,
 \end{equation}
  over all real rational curves $C$ of degree $d$ interpolating $3d{-}1$
  general points in  $\R\P^2$ does not depend upon the choice of the points.
\end{thm}

Set $\DeCo{W_d}$ to be the sum~\eqref{Eq1:Wel}.
The absolute value of this \DeCo{{\sl Welschinger invariant}} is then a lower
bound for the number of real rational curves of degree $d$ interpolating $3d-1$
points in $\R\P^2$.
Since $N_1=N_2=1$, we have $W_1=W_2=1$.
Prior to Welschinger's discovery, Kharlamov~\cite[Proposition 4.7.3]{DeKh00} (see also
Section~\ref{sec:ratcubics}) showed that $W_3=8$.
The question remained whether any other Welschinger invariants were
nontrivial. 
This was settled in the affirmative by Itenberg, Kharlamov, and
Shustin~\cite{IKS03,IKS04}, who used Mikhalkin's 
Tropical Correspondence Theorem~\cite{Mi05} to show
\begin{enumerate}
 \item If $d>0$, then $W_d\geq \frac{d!}{3}$.  (Hence $W_d$ is positive.) \smallskip
 \item ${\displaystyle \lim_{d\to\infty} 
            \frac{\log N_d}{\log W_d}=1}$. \ 
       (In fact, $\log N_d \sim 3d\log d\sim \log N_d$.)
\end{enumerate}
In particular, there are always quite a few real rational curves
of degree $d$ interpolating $3d{-}1$ points in $\R\P^2$.
Since then, Itenberg, Kharlamov, and Shustin~\cite{IKS09} gave a recursive formula for
the Welschinger invariant which is based upon Gathmann and Markwig's~\cite{GM} 
tropicalization of the Caporaso-Harris~\cite{CH98} formula.
This shows that $W_4=240$ and $W_5=18264$.
Solomon~\cite{Sol} has also found an intersection-theoretic interpretation for these invariants.  

These ideas have also found an application.
Gahleitner, J\"uttler, and Schicho~\cite{GJH} proposed a method to compute an
approximate parametrization of a plane curve using rational cubics.
Later, Fiedler-Le Touz\'e~\cite{FLT} used the result of Kharlamov
(that $W_3=8$), and an analysis of pencils of plane cubics to prove that this
method works. 
\smallskip

While the story of this interpolation problem is fairly well-known, it
was not the first instance of lower bounds in enumerative real
algebraic geometry. 
In their investigation of the Shapiro conjecture, Eremenko and
Gabrielov found a similar invariant $\sigma_{m,p}$ which gives a lower bound on the
number of real points in the inverse image $\Wr^{-1}(\Phi)$ under the Wronski
map of a real polynomial $\Phi\in\R\P^{mp)}$.
Assume that $m\leq p$.
If $m+p$ is odd, set
\[
  \DeCo{\sigma_{m,p}}\ :=\ 
   \frac{1!2!\dotsb(m{-}1)!(p{-}1)!(p{-}2)!\dotsb(p{-}m{+}1)!(\frac{mp}{2})!}
   {(p{-}m{+}2)!(p{-}m{+}4)!\dotsb(p{+}m{-}2)!\left(\frac{p-m+1}{2}\right)!%
     \left(\frac{p-m+3}{2}\right)!\dotsb\left(\frac{p+m-1}{2}\right)!}\ .
\]
If $m+p$ is even, then set $\sigma_{m,p}=0$.
If $m>p$, then set $\sigma_{m,p}:=\sigma_{p,m}$.

\begin{thm}[Eremenko-Gabrielov~\cite{EG01}]\label{T1:E-G01}
  If\/ $\Phi(t)\in\R\P^{mp}$ is a generic real polynomial of degree
  $mp$ (a regular value of the Wronski map), then there are at least
  $\sigma_{m,p}$ real $m$-dimensional subspaces of polynomials of degree $m{+}p{-}1$
  with Wronskian $\Phi$. 
\end{thm}

\begin{rmk}
Recall that the number of complex points in $\Wr^{-1}(\Phi)$ is
$\#_{m,p}$~\eqref{E1:WronDeg}.
It is instructive to compare these numbers.
We show them for $m{+}p=11$ and $m=2,\dotsc,5$.
\[
   \begin{tabular}{|c|r|r|r|r|}\hline
      $m$ & 2 & 3 & 4 & 5\\\hline
      $\sigma_{m,p}$& 14 & 110 & 286 & 286 \\\hline
      $\#_{m,p}$ &4862 &23371634&13672405890&396499770810 \\\hline
    \end{tabular}
\]
We also have $\sigma_{7,6}\approx 3.4\times 10^4$ and 
$\#_{7,6}\approx 9.5\times 10^{18}$.
Despite this disparity in their magnitudes, the asymptotic ratio of 
$\log(\sigma_{m,p})/\log(\#_{m,p})$ appears to be close to $1/2$.
We display this ratio in the table below, for different values of $m$ and $p$.
 \[
   \begin{tabular}{|c|r|c|c|c|c|c|c|}\hline
    \multicolumn{2}{|c|}{\multirow{2}{80pt}{
       ${\displaystyle \frac{\log(\sigma_{m,p})}{\log(\#_{m,p})}}$}}&
    \multicolumn{6}{|c|}{$m$\rule{0pt}{13 pt}}\\\cline{3-8}
    \multicolumn{2}{|c|}{\ }
    &2& $\frac{m{+}p{-}1}{10}$&$2\frac{m{+}p{-}1}{10}$&$3\frac{m{+}p{-}1}{10}$
         &$4\frac{m{+}p{-}1}{10}$&$5\frac{m{+}p{-}1}{10}$\raisebox{-6pt}{\rule{0pt}{19pt}}\\\hline
  \multirow{6}{15pt}{\hspace{10pt}\raisebox{-28pt}{\begin{rotate}{90}$m{+}p{-}1$\end{rotate}}\hspace{-10pt}}
    &     100& 0.47388& 0.45419& 0.43414& 0.41585& 0.39920& 0.38840\\\cline{2-8}
    &    1000& 0.49627& 0.47677& 0.46358& 0.45185& 0.44144& 0.43510\\\cline{2-8}
    &   10000& 0.49951& 0.48468& 0.47510& 0.46660& 0.45909& 0.45459\\\cline{2-8}
    &  100000& 0.49994& 0.48860& 0.48111& 0.47445& 0.46860& 0.46511\\\cline{2-8}
    & 1000000& 0.49999& 0.49092& 0.48479& 0.47932& 0.47453& 0.47168\\\cline{2-8}
    &10000000& 0.50000& 0.49246& 0.48726& 0.48263& 0.47857& 0.47616\\
\hline
  \end{tabular}
 \]
 Thus, the lower bound on the number of real points in a fiber of the Wronski map 
 appears asymptotic to the square root of the number of complex solutions.

 It is interesting to compare this to the the result of Shub and Smale~\cite{ShSm93} that
 the expected number of real solutions to a system of $n$ Gaussian random polynomials in
 $n$ variables of degrees $d_1,\dotsc,d_n$ is $\sqrt{d_1\dotsb d_n}$, which is the square
 root of the number of complex solutions to such a system of polynomials. 
\end{rmk}

The idea behind the proof of Theorem~\ref{T1:E-G01} is to compute the topological degree
of the real Wronski map, which is the restriction of the Wronski map to real subspaces of
polynomials, 
\[
  \Wr\!_\R\ :=\ \Wr\!|_{\sGr(m,\R_{m{+}p{-}1}[t])}\ \colon\ 
  \Gr(m,\R_{m{+}p{-}1}[t])\ \longrightarrow\ \R\P^{mp}\,.
\]
This maps the Grassmannian of real subspaces to the space of real
Wronski polynomials.
Recall that the topological degree of a
map $f\colon X\to Y$ between two oriented manifolds $X$ and $Y$ of the same dimension 
is the number $d$ such that $f_*[X]=d[Y]$, where $[X]$ and $[Y]$ are the fundamental
homology cycles of $X$ and $Y$, respectively, and $f_*$ is the functorial map in homology.
When $f$ is differentiable, this may be computed as follows.
Let $y\in Y$ be a regular value of $f$ so that at any point $x$ in the fiber $f^{-1}(y)$
above $y$ the derivative $d_xf\colon T_xX \to T_yY$ is an isomorphism.
Since $X$ and $Y$ are oriented, the isomorphism $d_xf$ either preserves the orientation or it reverses
it.
Let $P$ be the number of points $x\in f^{-1}(y)$ at which $d_xf$ preserves the
orientation and $R$ be the number of points where the orientation is reversed.
Then the degree of $f$ is the difference $P-R$.

There is a slight problem in computing the degree of $\Wr\!_\R$, as neither the real
Grassmannian nor the real projective space is orientable when $m{+}p$ is odd, and thus the 
topological degree of $\Wr\!_\R$ is not defined when $m{+}p$ is odd.
Eremenko and Gabrielov get around this by computing the degree of the restriction of the
Wronski map to open cells of $\Gr_\R$ and $\R\P^{mp}$, which is a proper map.
They also show that it is the degree of a lift of the 
Wronski map to oriented double covers of both spaces.
The degree bears a resemblance to the Welschinger invariant
as it has the form $|\sum \pm 1|$, the sum over all real points in 
$\Wr_\R^{-1}(\Phi)$, for $\Phi$ a regular value of the Wronski map.
This resemblance is no accident.
Solomon~\cite{Sol} showed how to orient a moduli space of rational curves with marked
points so that the Welschinger invariant is indeed the degree of a map.
\medskip

While both of these examples of geometric problems possessing a lower bound
on their numbers of real solutions are quite interesting, they are
rather special.
The existence of lower bounds for more general problems or for more general
systems of polynomials would be quite important in applications, as these lower 
bounds guarantee the existence of real solutions.

With Soprunova, we~\cite{SS} set out to develop a theory of lower bounds
for sparse polynomial systems, using the approach 
of Eremenko and Gabrielov via topological degree.
This is a first step toward practical applications of these ideas.
Chapter~\ref{Ch:lower} will elaborate this theory.
Here is an outline:
\begin{enumerate}
 \item[({\it i})] Realize the solutions to a system of polynomials as the fibers
            of a map from a toric variety.
 \item[({\it ii})] Give a condition which implies that the degree of this map (or a
            lift to double covers) exists.
 \item[({\it iii})] Develop a method to compute the degree in some (admittedly
            special) cases.
 \item[({\it iv})] Give a nice family of examples to which this theory applies.
 \item[({\it v})] Use the sagbi degeneration of a Grassmannian to a toric variety~\cite[Ch.~11]{GBCP} and
             the systems of ({\it iv}) to recover the result of Eremenko and Gabrielov.
\end{enumerate}

\begin{ex}\label{Ex1:P2P2}
We close this Chapter with one example from this theory.
Let $\Blue{w},\Blue{x},\Brown{y},\Brown{z}$ be indeterminates, and consider a
sparse polynomial of the form 
\begin{eqnarray}
& c_4\, \Blue{wx}\Brown{yz}\hspace{0.2em}&\nonumber\\
&+\,c_3(\Blue{wx}\Brown{z} + \Blue{x}\Brown{yz})\hspace{0.9em}&\nonumber\\
&+\,c_2(\Blue{wx} + \Blue{x}\Brown{z} + \Brown{yz})\hspace{1em}&\label{E1:P2P2}\\
&+\,c_1(\Blue{x}+\Brown{z})\hspace{0.74em}&\nonumber\\
&\hspace{.25em}\,+\,c_0\,,&\nonumber
\end{eqnarray}
where the coefficients $c_0,\dotsc,c_4$ are real numbers.\medskip

\begin{thm}\label{T1:P2P2}
 A system of four equations involving polynomials of the
  form~\eqref{E1:P2P2} has six solutions, at least \Blue{two} of which are real.
\end{thm}

We make some remarks to illustrate the ingredients of this theory.
First, the monomials in the sparse system~\eqref{E1:P2P2} are the
integer points in the \Blue{{\sl order polytope}} of the poset $P$, 
\[
  \begin{picture}(90,50)(-35,0)
    \put(-40,20){$P\ :=$}
    \put(14,3){\includegraphics{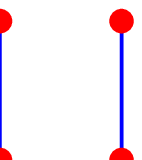}}
    \put(0,40){$\Blue{x}$}   \put(56,40){$\Brown{z}$}
    \put(0, 0){$\Blue{w}$}   \put(56, 0){$\Brown{y}$}
    \put(65,20){.}
  \end{picture}
\]
That is, each monomial corresponds to an order ideal of $P$ (a subset which is
closed upwards).
The number of complex roots is the number of linear extensions of the
poset $P$.
There are six, as each is a permutation of the word $\Blue{wx}\Brown{yz}$ where
$\Blue{w}$ precedes $\Blue{x}$ and $\Brown{y}$ precedes $\Brown{z}$.

One result ({\it ii}) gives conditions on the Newton polytope which imply that
any polynomial system with that Newton polytope has a
lower bound, and order polytopes satisfy these conditions.
Another result ({\it iv}) computes that lower bound for certain families of 
polynomials with support an order polytope.
Polynomials in these families have the form~\eqref{E1:P2P2} in that
monomials with the same total degree have the same coefficient.
For such polynomials, the lower bound is the absolute value of the sum of the signs of the
permutations underlying the linear extensions.
We list these for $P$.
\begin{center}
 \begin{tabular}{|r||c|c|c|c|c|c||c|}\hline
  permutation&$\Blue{wx}\Brown{yz}$
             &$\Blue{w}\Brown{y}\Blue{x}\Brown{z}$
             &$\Brown{y}\Blue{wx}\Brown{z}$
             &$\Blue{w}\Brown{yz}\Blue{x}$
             &$\Brown{y}\Blue{w}\Brown{z}\Blue{x}$
             &$\Brown{yz}\Blue{wx}$&sum\\\hline
  sign & $+$ & $-$  & $+$  & $+$ & $-$  & $+$ & 2 \\\hline
 \end{tabular}
\end{center}
This shows that the lower bound in Theorem~\ref{T1:P2P2} is 2.

We record the frequency of the different root counts in each of 
10,000,000 instances of this polynomial system, where the coefficients were chosen
uniformly from $[-200,200]$.
\[
   \begin{tabular}{|r|c|c|c|c|}\hline
      number of real roots&0&2&4&6\\\hline
      frequency &0&9519429&0&480571\\\hline
   \end{tabular}
\]
This computation took 19,854 CPU-seconds, or 5:30:54 on a AMD Athlon 64-bit dual core
processor 4600+ (2.4 GHz, 512KB cache), using a single core.
\QED 
\end{ex}
%
%
%
%
The apparent gap in the numbers of real solutions (4 does not seem a possible number of
real solutions) is proven for the system of Example~\ref{Ex1:P2P2}.
This is the first instance we have seen of this phenomena of gaps in the numbers of real
solutions.
More are found in~\cite{SS} and~\cite{RSSS}.
%
%
\chapter{Real solutions to univariate polynomials}\label{ch:Sturm}

%
%

Before we study the real solutions to systems of multivariate polynomials, we will 
review some of what is known for univariate polynomials.
The strength and precision of results concerning real solutions to univariate polynomials
forms the gold standard in this subject of real roots to systems of 
polynomials. 
We will discuss two results about univariate polynomials: 
Descartes' rule of signs and Sturm's Theorem.
Descartes' rule of signs, or rather its generalization in the Budan-Fourier Theorem, gives a
bound for the number of roots in an interval, counted with multiplicity.
Sturm's theorem is topological---it simply counts the number of roots of a univariate polynomial
in an interval without multiplicity.
From Sturm's Theorem we obtain a symbolic algorithm to count the number of real solutions
to a system of multivariate polynomials in many cases.
We underscore the topological nature of Sturm's Theorem by presenting a new and very
elementary proof due to Burda and Khovanskii~\cite{BKh}.
These and other fundamental results about real roots of univariate polynomials were
established in the 19th century.
In contrast, the main results about real solutions to multivariate polynomials were
only established in recent decades.

%
\section{Descartes' rule of signs}\label{S2:Descartes}
%

Descartes' rule of signs~\cite{D1637} is fundamental for real algebraic
geometry. 
Suppose that $f$ is a univariate polynomial and write its terms in increasing order of
their exponents,
 \begin{equation}\label{Eq2:Descartes}
   f\ =\  c_0 x^{a_0} +  c_1 x^{a_1} + \dotsb +  c_m x^{a_m}\,,
 \end{equation}
where $ c_i\neq 0$ and $a_0<a_1<\dotsb<a_m$.

\begin{thm}[Descartes' rule of signs]\label{T2:Descartes}
 The number, $r$, of positive roots of $f$, counted with multiplicity, is at most 
 the variation in sign of the coefficients of $f$,
\[
   \#\{i\mid 1\leq i\leq m\mbox{ and }  c_{i-1} c_i<0\}
    \ \leq\ r\,,
\]
 and the difference between the variation and $r$ is even.
\end{thm}

We will prove a generalization, the Budan-Fourier Theorem, which provides a
similar estimate for any interval in $\R$.
We first formalize this notion of variation in sign that appears in Descartes' rule.

The \DeCo{{\sl variation} $\var(c)$} in a finite sequence $ c$ of real numbers is
the number of times that consecutive elements of the sequence have opposite signs, after
we remove any 0s in the sequence.
For example, the first sequence below has four variations, while the
second has three.
 \[
   8\DeCo{,}-4,-2,-1\DeCo{,}2,3\DeCo{,}-5\DeCo{,}7,11,12
   \qquad\qquad
    -1,0\DeCo{,}1,0,1\DeCo{,}-1\DeCo{,}1,1,0,1\,.
 \]

Suppose that we have a sequence $F=(f_0,f_1,\dotsc,f_k)$ of polynomials and a real number
$a\in \R$. 
Then \DeCo{$\var(F,a)$} is the variation in the sequence 
$f_0(a),f_1(a),\dotsc,f_k(a)$.
This notion also makes sense when $a=\pm\infty$:
We set $\var(F,\infty)$ to be the variation in the sequence of leading coefficients of the
$f_i(t)$, which are the signs of $f_i(a)$ for $a\gg0$, and set $\var(F,-\infty)$ to be the
variation in the leading coefficients of $f_i(-t)$.
\medskip

Given a univariate polynomial $f(t)$ of degree $k$, let $\delta f$ be the sequence of its
derivatives,
\[
  \DeCo{\delta f}\ :=\ (f(t), f'(t), f''(t),\dotsc, f^{(k)}(t))\,.
\]
For $a,b\in\R\cup\{\pm\infty\}$, let $r(f,a,b)$ be the number of roots of $f$ in the
interval $(a,b]$, counted with multiplicity.
We prove a version of Descartes' rule due to Budan~\cite{Budan} and
Fourier~\cite{Fourier}.

\begin{thm}[Budan-Fourier]
 Let $f\in\R[t]$ be a univariate polynomial and $a<b$ two numbers in
 $\R\cup\{\pm\infty\}$. 
 Then
\[
    \var(\delta f,a)\ -\ \var(\delta f,b)\ \geq\ r(f,a,b)\,,
\]
 and the difference is even.
\end{thm}

We may deduce Descartes' rule of signs from the Budan-Fourier Theorem once we 
observe that for the polynomial $f(t)$~\eqref{Eq2:Descartes}, 
$\var(\delta f,0)=\var(c_0, c_1,\dotsc, c_m)$,
while $\var(\delta f,\infty)=0$, as the leading coefficients of $\delta f$ all have the
same sign.\medskip

\begin{ex}\label{Ex2:BF}
 The the sextic $f=5t^6-4t^5-27t^4+55t^2-6$ whose graph is displayed below
\[
   \begin{picture}(300,175)
    \put(0,0){\includegraphics{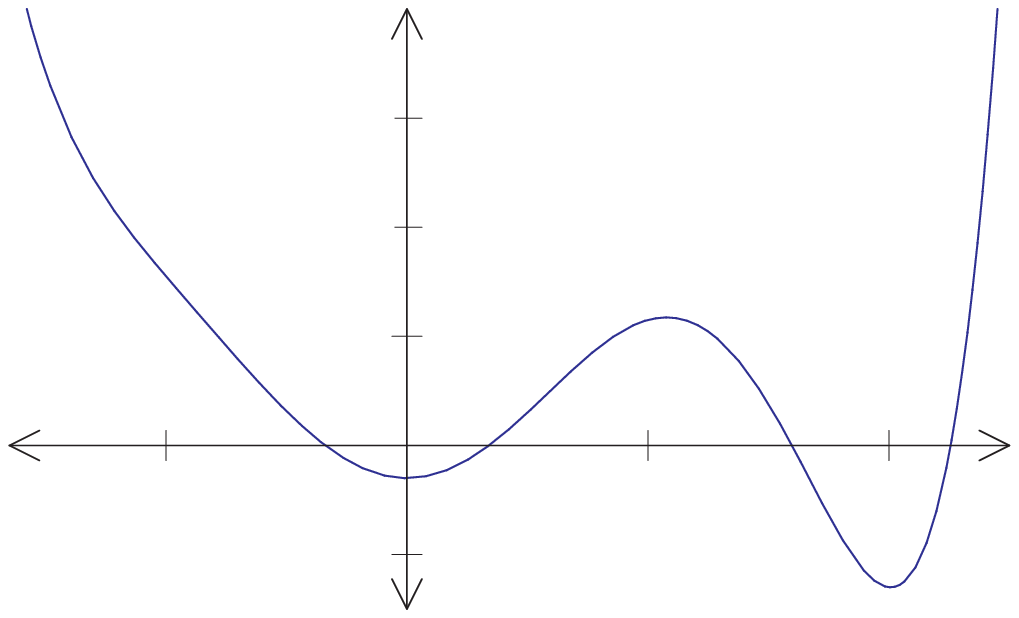}}

    \put(33,30){$-1$}    \put(182,30){$1$}    \put(252,30){$2$}
    \put(88,13){$-20$}   \put(97,76){$20$}   
    \put(97,108){$40$}   \put(97,139){$60$}
    \put(273,163){$f$}   \put(290,53){$t$}
   \end{picture}
\]
 has four real zeroes at approximately $-0.339311$, $0.340401$, $1.59753$, $2.25615$.
 If we evaluate the derivatives of $f$ at 0 we obtain
 \[
    \delta f(0)\ =\ -6,\, 0\DeCo{,}\, 110,\, 0\DeCo{,}\, -648,\, -480\DeCo{,}\, 3600\,,
 \]
 which has 3 variations in sign.
 If we evaluate the derivatives of $f$ at $2$, we obtain
 \[
     \delta f(2)\ =\ -26,\, -4\DeCo{,} \,574,\, 2544,\, 5592, \,6720,\, 3600\,,
 \]
 which has one sign variation.
 Thus, by the Budan-Fourier Theorem, $f$ has either 2 or 0 roots in the interval $(0,2)$,
 counted with multiplicity.
 This agrees with our observation that $f$ has 2 roots in the interval $[0,2]$. 
 \QED
\end{ex}

\noindent{\it Proof of Budan-Fourier Theorem.}
 Observe that $\var(\delta f,t)$ can only change when $t$ passes a root $c$ of 
 some polynomial in the sequence $\delta f$ of derivatives of $f$.
 Suppose that $c$ is a root of some derivative of $f$ and 
 let $\epsilon>0$ be a positive number such that
 no derivative $f^{(i)}$ has a root in the interval $[c-\epsilon,c+\epsilon]$, except
 possibly at $c$.  
 Let $\DeCo{m}$ be the order of vanishing of $f$ at $c$.
 We will prove that 
 \begin{equation}\label{E2:claims}
  \begin{tabular}{cl}
   (1) &$\var(\delta f,c)=\var(\delta f,c+\epsilon)$, \ and\\
   (2) &$\var(\delta f,c-\epsilon)\geq \var(\delta f,c)+m$, \  \rule{0pt}{15pt}
     and the difference is even.
  \end{tabular}\hspace{80pt}
 \end{equation}

 We deduce the Budan-Fourier theorem from these conditions.
 As $t$ ranges from $a$ to $b$, both $r(f,a,t)$ and $\var(\delta f,t)$ only change when
 $t$ passes a root $c$ of $f$.
 (These could change at a root of a derivative of $f$, but in fact do not.)  
 At such a point, $r(f,a,t)$ jumps by the multiplicity $m$ of that root of $f$, while 
 $\var(\delta f,t)$ drops by $m$, plus a nonnegative even integer.
 Thus the sum $r(f,a,t)+\var(\delta f,t)$ can only change at roots $c$ of $f$, where it
 drops by an even integer.
 Since this sum equals $\var(\delta f,a)$ when $t=a$, the Budan-Fourier Theorem
 follows.\smallskip 

 Let us now prove our claim about the behavior of $\var(\delta f,t)$ in a neighborhood of 
 a root $c$ of some derivative $f^{(i)}$.
 We argue by induction on the degree of $f$.
 When $f$ has degree 1, then we are in one of the following two cases, depending upon the
 sign of $f'$
\[
   \begin{picture}(122,92)(-22,0)   
    \put(1,0){\includegraphics{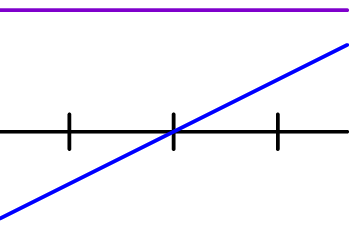}}
    \put(-22,17){$f(t)$}   \put(-24,78){$f'(t)$}  
    \put(7.5,54){$c-\epsilon$}  \put(47,30){$c$} \put(68,30){$c+\epsilon$}
   \end{picture}
    \qquad\quad
   \begin{picture}(122,92)(-22,0)      
    \put(1,0){\includegraphics{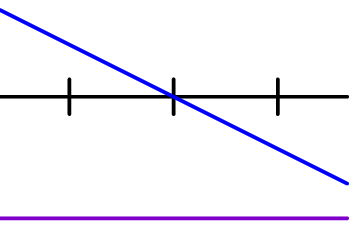}}
    \put(-22,67){$f(t)$}   \put(-24,7){$f'(t)$}  
    \put(7.5,30){$c-\epsilon$}  \put(47,30){$c$} \put(68,54){$c+\epsilon$}
   \end{picture}
\]
 In both cases, $\var(\delta f,c-\epsilon)=1$, but 
 $\var(\delta f, c)=\var(\delta f,c+\epsilon)=0$, which proves the claim when $f$ is linear.

 Now suppose that the degree of $f$ is greater than $1$ and let
 $m$ be the order of vanishing of $f$ at $c$.
 We first treat the case when $f(c)=0$, and hence $m>0$ so that $f'$ vanishes at $c$ to
 order $m{-}1$.  
 We apply our induction hypothesis to $f'$ and obtain that
\[
  \var(\delta f',c)\ =\ \var(\delta f', c+\epsilon),
   \qquad\mbox{and}\qquad
  \var(\delta f',c-\epsilon)\ \geq\ \var(\delta f', c)\ +\ (m-1)\,,
\]
 and the difference is even.
 By Lagrange's Mean Value Theorem applied to the intervals $[c-\epsilon,c]$ and 
 $[c,c+\epsilon]$, $f$ and $f'$ must have opposite signs at $c-\epsilon$,
 but the same signs at $c+\epsilon$, and so 
 \begin{eqnarray*}
   \var(\delta f,c)&=&\var(\delta f',c)\ =\ 
   \var(\delta f',c+\epsilon)\ =\ \var(\delta f, c+\epsilon)\,,\\
   \var(\delta f,c-\epsilon)&=&\var(\delta f',c-\epsilon)+1\ \geq\ 
   \var(\delta f',c)+(m-1)+1\ =\ \var(\delta f, c)+m\,,
 \end{eqnarray*} 
 and the difference is even.
 This proves the claim when $f(c)=0$.

 Now suppose that $f(c)\neq 0$ so that $m=0$.
 Let \DeCo{$n$} be the order of vanishing of $f'$ at $c$.
 We apply our induction hypothesis to $f'$ to obtain that 
\[
  \var(\delta f',c)\ =\ \var(\delta f', c+\epsilon),
   \qquad\mbox{and}\qquad
  \var(\delta f',c-\epsilon)\ \geq\ \var(\delta f', c)\ +\ n\,,
\]
 and the difference is even.
 We have $f(c)\neq 0$, but
 $f'(c)=\dotsb=f^{(n)}(c)=0$, and $f^{(n+1)}(c)\neq 0$.
 Multiplying $f$ by $-1$ if necessary, we may assume that $f^{(n+1)}(c)>0$.
 There are four cases: $n$ even or odd, and $f(c)$ positive or negative.
 We consider each case separately.

 Suppose that $n$ is even.
 Then both $f'(c-\epsilon)$ and $f'(c+\epsilon)$ are positive and so for each
 $t\in\{c-\epsilon,c,c+\epsilon\}$ the first nonzero term in the sequence 
 \begin{equation}\label{E2:val_seq}
  f'(t),\ f''(t),\ \dotsc,\ f^{(k)}(t)
 \end{equation}
 is positive.
 When $f(c)$ is positive, this implies that  $\var(\delta f,t)=\var(\delta f',t)$ 
 and when $f(c)$ is negative, that  $\var(\delta f,t)=\var(\delta f',t)+1$.
 This proves the claim as it implies that $\var(\delta f,c)=\var(\delta f,c+\epsilon)$ and
 also that 
\[
   \var(\delta f,c-\epsilon)\ -\ \var(\delta f,c)\ =\ 
   \var(\delta f',c-\epsilon)\ -\ \var(\delta f', c)\,,
\]
 but this last difference exceeds $n$ by an even number, and so is even as $n$ is even.

 Now suppose that $n$ is odd.
 Then $f'(c-\epsilon)<0<f'(c+\epsilon)$ and so 
 the first nonzero term in the sequence~\eqref{E2:val_seq} has sign 
 $-,+,+$ as $t=c-\epsilon,c,c+\epsilon$.
 If $f(c)$ is positive, then 
 $\var(\delta f,c-\epsilon)=\var(\delta f',c-\epsilon)+1$ and the other two variations
 are unchanged, but if $f(c)$ is negative, then the variation at $t=c-\epsilon$ is
 unchanged, but it increases by 1 at $t=c$ and $t=c+\epsilon$.
 This again implies the claim,
 as $\var(\delta f, c)=\var(\delta f,c+\epsilon)$, but
\[
   \var(\delta f,c-\epsilon)-\var(\delta f,c)\ =\ 
   \var(\delta f',c-\epsilon)-\var(\delta f',c)\ \pm\ 1\,.
\]
 Since the difference $ \var(\delta f',c-\epsilon)-\var(\delta f',c)$ is equal to
 the order $n$ of the vanishing of $f'$ at $c$ plus a nonnegative even number, if we add
 or subtract 1, the difference is a nonnegative even number.
 This completes the proof of the Budan-Fourier Theorem.
\QED

%
\section{Sturm's Theorem}
%

Let $f,g$ be univariate polynomials.
Their \DeCo{{\sl Sylvester sequence}} is the sequence of polynomials
\[
   f_0\;:=\;f, \ \ f_1\;:=\;g, \ \ f_2,\ \dotsc,\ f_k\,,
\]
where $f_k$ is a greatest common divisor of $f$ and $g$, and 
\[
   \DeCo{-f_{i+1}}\ :=\ \mbox{remainder}(f_{i-1},f_i)\,,
\]
the usual remainder from the Euclidean algorithm.
Note the sign.
We remark that we have polynomials $q_1,q_2,\dotsc,q_{k-1}$ such that
 \begin{equation}\label{E2:quotient_polys}
   f_{i-1}(t)\ =\ q_i(t)f_i(t)\ -\ f_{i+1}(t)\,,
 \end{equation}
and the degree of $f_{i+1}$ is less than the degree of $f_i$.
The \DeCo{{\sl Sturm sequence}} of a univariate polynomial $f$ is the
Sylvester sequence of $f, f'$.

\begin{thm}[Sturm's Theorem]\label{T2:Sturm}
 Let $f$ be a univariate polynomial and $a,b\in\R\cup\{\pm\infty\}$ with 
 $a<b$ and $f(a),f(b)\neq 0$.
 Then the number of zeroes of $f$ in the interval $(a,b)$ is the difference
\[
   \var(F,a)\ -\ \var(F,b)\,,
\]
 where $F$ is the Sturm sequence of $f$.
\end{thm}

\begin{ex}\label{Ex2:Sturm}
 The sextic $f$ of Example~\ref{Ex2:BF} has Sturm sequence
 \begin{eqnarray*}
   f&=& 5t^6-4t^5-27t^4+55t^2-6\\
   f_1:=f'(t)&=& 30t^5-20t^4-108t^3+110t\\
   f_2&=&\tfrac{84}{9}t^4 +\tfrac{12}{5}t^3 -\tfrac{110}{3}t^2-\tfrac{22}{9}t+6\\
   f_3&=&\tfrac{559584}{36125}t^3+\tfrac{143748}{1445}t^2
                -\tfrac{605394}{7225}t-\tfrac{126792}{7225}\\
   f_4&=&\tfrac{229905821875}{724847808}t^2+\tfrac{1540527685625}{4349086848}t
           +\tfrac{7904908625}{120807968}\\
   f_5&=&-\tfrac{280364022223059296}{58526435357253125}t
          +\tfrac{174201756039315072}{292632176786265625}\\
    f_6&=&-\tfrac{17007035533771824564661037625}{162663080627869030112013128}\,.
 \end{eqnarray*}
 Evaluating the Sturm sequence at $t=0$ gives
\[
   -6,\ 0\DeCo{,}\ 6\DeCo{,}\ -\tfrac{126792}{7225}\DeCo{,}\ 
    \tfrac{174201756039315072}{292632176786265625}\DeCo{,}\ 
   -\tfrac{17007035533771824564661037625}{162663080627869030112013128}\ ,
\]
 which has $4$ variations in sign, while evaluating the Sturm sequence at $t=2$ gives
\[
  -26,\ -4\DeCo{,}\ \tfrac{1114}{45},\ \tfrac{3210228}{36125}\DeCo{,}\ 
  -\tfrac{1076053821625}{2174543424},\  
   -\tfrac{2629438466191277888}{292632176786265625},\ 
   -\tfrac{17007035533771824564661037625}{162663080627869030112013128}\ , 
\]
 which has $2$ variations in sign.
 Thus by Sturm's Theorem, we see that $f$ has $2$ roots in the interval $[0,2]$,
 which we have already seen by other methods.
 \QED
\end{ex}

An application of Sturm's Theorem is to isolate real solutions to a univariate
polynomial $f$ by finding intervals of a desired width that contain a unique root of 
$f$.
When $(a,b)=(-\infty,\infty)$, Sturm's Theorem gives the total number of real roots of a
univariate polynomial.
In this way, it leads to an algorithm to investigate the number of real
roots of generic systems of polynomials.
We briefly describe this algorithm here.
This algorithm was used in an essential way to get information on real solutions which 
helped to formulate many results discussed in later chapters.

Suppose that we have a system of real multivariate polynomials
 \begin{equation}\label{Eq2:systemN}
   f_1(x_1,\dotsc,x_n)\ =\ 
   f_2(x_1,\dotsc,x_n)\ =\ \dotsb\ =\ 
   f_N(x_1,\dotsc,x_n)\ =\ 0\,,
 \end{equation}
whose number of real roots we wish to determine.
Let $I\subset\R[x_1,\dotsc,x_n]$ be the ideal generated by the polynomials
$f_1,f_2,\dotsc,f_N$. 
If~\eqref{Eq2:systemN} has finitely many complex zeroes, then the dimension of 
the quotient ring (the \DeCo{{\sl degree}} of $I$) is finite, and for each 
variable $x_i$, there is a univariate polynomial $g(x_i)\in I$ of minimal degree, called
an \DeCo{{\sl eliminant}} for $I$.

\begin{prop}
 The roots of $g(x_i)=0$ form the set of $i$th coordinates of solutions
 to~\eqref{Eq2:systemN}. 
\end{prop}

The main tool here is a consequence of the Shape Lemma~\cite{BMMT}.

\begin{thm}[Shape Lemma]
 Suppose that $I$ has an eliminant $g(x_i)$ whose degree is equal to the degree of $I$.
 Then the number of real solutions to~\eqref{Eq2:systemN} is equal to the number of 
 real roots of $g$.
\end{thm}

Suppose that the coefficients of the polynomials $f_i$ in the
system~\eqref{Eq2:systemN} lie in a computable subfield of $\R$, for example, $\Q$
(e.g.~if the coefficients are integers). 
Then the degree of $I$ may be computed using Gr\"obner bases, and we may also use
Gr\"obner bases to compute an eliminant $g(x_i)$.
Since Buchberger's algorithm does not enlarge the field of the coefficients, 
$g(x_i)\in\Q[x_i]$ has rational coefficients, and so we may use Sturm sequences to compute the
number of its real roots.\medskip

\noindent{\bf Algorithm}

\noindent{\sc Given:} $I=\langle f_1,\dotsc,f_N\rangle\subset\Q[x_1,\dotsc,x_n]$

\begin{enumerate}

 \item Use Gr\"obner bases to compute the degree $d$ of $I$.

 \item Use Gr\"obner bases to compute an eliminant $g(x_i)\in I\cap \Q[x_i]$ for $I$.

 \item If $\deg(g)=d$, then use Sturm sequences to compute the number $r$ of real roots of
        $g(x_i)$, and output ``The ideal $I$ has $r$ real solutions.''

 \item Otherwise output ``The ideal $I$ does not satisfy the hypotheses of the Shape
           Lemma.'' 
\end{enumerate}

If this algorithm halts with a failure (step 4), it may be called again to compute an eliminant
for a different variable.
Another strategy is to apply a random linear transformation before eliminating.
An even more sophisticated form of elimination is Roullier's rational univariate
representation~\cite{Rou99}.

\subsection{Traditional Proof of Sturm's Theorem}

Let $f(t)$ be a real univariate polynomial with Sturm sequence $F$.
We prove Sturm's Theorem by looking at the variation
$\var(F,t)$ as $t$ increases from $a$ to $b$.
This variation can only change when $t$ passes a number $c$ where some member $f_i$ of the 
Sturm sequence has a root, for then the sign of $f_i$ could change.
We will show that if $i>0$, then this has no effect on the variation of the sequence, but
when $c$ is a root of $f=f_0$, then the variation decreases by exactly $1$ as $t$ passes
$c$.
Since multiplying a sequence by a nonzero number does not change its variation, we 
will at times make an assumption on the sign of some value $f_j(c)$ to reduce the number
of cases to examine.

Observe first that by~\eqref{E2:quotient_polys}, if $f_i(c)=f_{i+1}(c)=0$, then 
$f_{i-1}$ also vanishes at $c$, as do the other polynomials $f_j$.
In particular $f(c)=f'(c)=0$, so $f$ has a multiple root at $c$.
Suppose first that this does not happen, either that $f(c)\neq 0$ or that $c$ is a simple
root of $f$.

Suppose that $f_i(c)=0$ for some $i>0$.
The vanishing of $f_i$ at $c$, together with~\eqref{E2:quotient_polys} implies that
$f_{i-1}(c)$ and $f_i(c)$ have opposite signs.
Then, whatever the sign of $f_i(t)$ for $t$ near $c$, there is exactly one variation in
sign coming from the subsequence $f_{i-1}(t),f_i(t),f_{i+1}(t)$, and so the vanishing of
$f_i$ at $c$ has no effect on the variation as $t$ passes $c$.
Note that this argument works equally well for any Sylvester sequence.

Now we consider the effect on the variation when $c$ is a simple root of $f$.
In this case $f'(c)\neq 0$, so we may assume that $f'(c)>0$.
But then $f(t)$ is negative for $t$ to the left of $c$ and positive for $t$ to the right
of $c$.
In particular, the variation $\var(F,t)$ decreases by exactly 1 when $t$ passes a simple
root of $f$ and does not change when $f$ does not vanish.

We are left with the case when $c$ is a multiple root of $f$.
Suppose that its multiplicity is $m+1$.
Then $(t-c)^m$ divides every polynomial in the Sturm sequence of $f$.
Consider the sequence of quotients,
 \[
   G\ =\ (g_0,\dotsc,g_k)\ :=\ 
    \left(f/(t-c)^m,\ f'/(t-c)^m,\ f_2/(t-c)^m,\ \dotsb,\ f_k/(t-c)^m\right)\,.
 \]
Note that $\var(G,t)=\var(F,t)$ when $t\neq c$, as multiplying a sequence by a nonzero
number does not change its variation.
Observe also that $G$ is a Sylvester sequence.
Since $g_1(c)\neq 0$, not all polynomials $g_i$ vanish at $c$.
But we showed in this case that there is no contribution to a change in the variation by 
any polynomial $g_i$ with $i>0$.

It remains to examine the contribution of $g_0$ to the variation as $t$ passes $c$.
If we write $f(t)=(t-c)^{m+1}h(t)$ with $h(c)\neq 0$, then 
\[
  f'(t)\ =\ (m+1)(t-c)^mh(t)\ +\ (t-c)^{m+1}h'(t)\,.
\]
In particular,
\[
  g_0(t)\ =\ (t-c)h(t)
  \qquad\mbox{and}\qquad
  g_1(t)\ =\ (m+1)h(t)\ +\ (t-c)h'(t)\,.
\]
If we assume that $h(c)>0$, then $g_1(c)>0$ and $g_0(t)$ changes from negative to positive
as $t$ passes $c$.
Once again we see that the variation $\var(F,t)$ decreases by 1 when $t$ passes a root of
$f$.
This completes the proof of Sturm's Theorem.
\QED

%
\section{A topological proof of Sturm's Theorem}
%
We present a second, very elementary, proof of Sturm's Theorem due to Burda and
Khovanskii~\cite{BKh} whose virtue is in its tight connection to topology.
We first recall the definition of topological degree of a continuous function 
$\varphi\colon \R\P^1\to\R\P^1$.
Since $\R\P^1$ is isomorphic to the quotient $\R/\Z$, we may 
pull $\varphi$ back to the interval $[0,1]$ to obtain a map
$[0,1]\to\R\P^1$.
This map lifts to the universal cover of $\R\P^1$ to obtain a map
$\psi\colon[0,1]\to \R$.  
Then the \DeCo{{\sl mapping degree}},  \DeCo{$\mdeg(\varphi)$}, of $\varphi$ is simply
$\psi(1)-\psi(0)$, which is an integer. 
We call this mapping degree to distinguish it from 
the usual algebraic degree of a polynomial or rational function.

The key ingredient in this proof is a formula to compute the mapping degree of a
rational function $\varphi\colon \R\P^1\to\R\P^1$.
Any rational function $\varphi=f/g$ where $f,g\in\R[t]$ are polynomials has a continued
fraction expansion of the form
 \begin{equation}\label{E2:Cfrac}
    \varphi\ =\ q_0\ +\ \cfrac{1}{q_1+\cfrac{1}{q_2+
         \cfrac{1}{\begin{array}{cr}\raisebox{-4pt}{$\ddots$}&\\
            &\,+\,\dfrac{1}{q_k}\end{array}\hspace{-6pt}}}}
 \end{equation}
where $q_0,\dotsc,q_k$ are polynomials.
Indeed, this continued fraction is constructed recursively.
If we divide $f$ by $g$ with remainder $h$, so that $f=q_0g+h$ with the degree of $h$ less
than the degree of $g$, then 
\[
    \varphi\ =\ q_0 + \frac{h}{g}\ =\  q_0 + \cfrac{\quad 1\quad}{\cfrac{g}{h}}\,.
\]
We may again divide $g$ by $h$ with remainder, $g=q_1h+k$ and obtain
\[
    \varphi\ =\ q_0 + \cfrac{1}{q_1+\cfrac{\quad 1\quad}{\cfrac{h}{k}}}\,.
\]
As the degrees of the numerator and denominator drop with each step, this 
process terminates with an expansion~\eqref{E2:Cfrac} of $\varphi$.

For example, if $f=4t^4-18t^2-6t$ and $g=4t^3+8t^2-1$, then
\[
   \frac{f}{g}\ =\ t-2\ +\ \cfrac{1}{-2t+1\ +\ \cfrac{1}{-2t-3\ +\ \cfrac{1}{t+1}}}
\]
This continued fraction expansion is just the Euclidean algorithm in disguise. 

Suppose that $q=c_0+c_1t+\dotsb+c_dt^d$ is a real polynomial of degree $d$.
Define
\[
   \DeCo{[q]}\ :=\ \sign(c_d)\cdot (d\mod 2)\ \in\ \{\pm1,0\}\,.
\]

\begin{thm}\label{T2:top-degree}
  Suppose that $\varphi$ is a rational function with continued fraction
  expansion~\eqref{E2:Cfrac}.
  Then the mapping degree of $\varphi$ is
\[
   [q_1] - [q_2] + \dotsb + (-1)^{k-1} [q_k]\,.
\]
\end{thm}

We may use this to count the roots of a real polynomial $f$ by 
the following lemma.

\begin{lemma}\label{L2:root-deg}
  The number of roots of a polynomial $f$, counted without multiplicity is 
  the mapping degree of the rational function $f/f'$.
\end{lemma}

We deduce Sturm's Theorem from Lemma~\ref{L2:root-deg}.
Let \DeCo{$f_0,f_1,f_2\dotsc,f_k$} be the Sturm sequence for $f$.
Then $f_0=f$, $f_1=f'$, and for $i>1$, $-f_{i+1}:=\mbox{remainder}(f_{i-1},f_i)$.
That is, $\deg(f_i)<\deg(f_{i-1})$ and there are univariate polynomials 
$g_1,g_2,\dotsc,g_k$ with
\[
   f_{i-1}\ =\ g_if_i-f_{i+1}\qquad\mbox{for}\quad i=1,2,\dotsc,k{-}1\,.
\]
We relate these polynomials to those obtained from the Euclidean algorithm applied to
$f,f'$ and thus to the continued fraction expansion of $f/f'$.
It is clear that the $f_i$ differ only by a sign from the remainders in the Euclidean
algorithm.
Set $r_0:=f$ and $r_1=f'$, and for  $i>1$, $r_i:=\mbox{remainder}(r_{i-2},r_{i-1})$.
Then $\deg(r_i)<\deg(r_{i-1})$, and there are univariate polynomials 
$q_1,q_2,\dotsc,q_k$ with
\[
   r_{i-i}\ =\ q_ir_i+r_{i+1}\qquad\mbox{for}\quad i=1,\dotsc,k{-}1\,.
\]
We leave the proof of the following lemma as an exercise for the reader.

\begin{lemma}
  We have $g_i=(-1)^{i-1}q_i$ and $f_i=(-1)^{\lfloor\frac{i}{2}\rfloor}r_i$, for
  $i=1,2,\dotsc,k$. 
\end{lemma}

Write $F$ for the Sturm sequence $(f_0,f_1,f_2\dotsc,f_k)$ for $f$.
Then $\var(F,\infty)$ is the variation in the leading coefficients
$(f_0^{\rm top},f_1^{\rm top},\dotsc,f_k^{\rm top})$ of the polynomials in $F$.
Similarly, $\var(F,-\infty)$ is the variation in the sequence
\[
  ((-1)^{\deg(f_0)}f_0^{\rm top},(-1)^{\deg(f_1)}f_1^{\rm top},\dotsc,
   (-1)^{\deg(f_k)}f_k^{\rm top})\,.
\]
Note that the variation in a sequence $(c_0,c_1,\dotsc,c_k)$ is just the sum of the
variations in each subsequence $(c_{i-1},c_i)$ for $i=1,\dotsc,k$.
Thus
 \begin{multline}\label{Eq:variation_sum}
  \qquad \var(F,-\infty)-\var(F,\infty)\\ =\ 
   \sum_{i=1}^k \left(\var( (-1)^{\deg(f_{i-1})}f_{i-1}^{\rm top},(-1)^{\deg(f_i)}f_i^{\rm top})
           \ -\  \var(f_{i-1}^{\rm top},f_i^{\rm top})\right)\,.\qquad
 \end{multline}
Since $f_{i-1}=g_if_i-f_{i+1}$ and $\deg(f_{i+1})<\deg(f_i)<\deg(f_{i-1})$, 
we have
\[
   f_{i-1}^{\rm top}\ =\ g_i^{\rm top}f_i^{\rm top}
       \qquad\mbox{and}\qquad
   \deg(f_{i-1})\ =\ \deg(g_i)+\deg(f_i)\,.
\]
Thus we have
 \begin{eqnarray*}
   \var(f_{i-1}^{\rm top},f_i^{\rm top})&=&\var(g_i^{\rm top}, 1)\,,\qquad\mbox{and}
   \\
   \var( (-1)^{\deg(f_{i-1})}f_{i-1}^{\rm top},(-1)^{\deg(f_i)}f_i^{\rm top})&=&
   \var( (-1)^{\deg(g_i)}g_i^{\rm top}, 1)\,.
 \end{eqnarray*}
Thus the summands in~\eqref{Eq:variation_sum} are
 \begin{eqnarray*}
   \var( (-1)^{\deg(g_i)}g_i^{\rm top}, 1)\ -\ \var(g_i^{\rm top}, 1)
   &=& \sign(g_i^{\rm top}) (\deg(g_i)\mod 2)\\
   &=& [g_i]\ =\ (-1)^{i-1}[q_i]\,,
 \end{eqnarray*}
This proves that
 \begin{eqnarray*}
   \var(F,-\infty)-\var(F,\infty)&=& [g_1]+[g_2]+\dotsb+[g_k]\\
      &=& [q_1]-[q_2]+ \dotsb + (-1)^{k-1}[q_k]\,.
 \end{eqnarray*}
But this proves Sturm's Theorem, as this is the number of roots of $f$, by
Lemma~\ref{L2:root-deg}. \QED\smallskip

The key to the proof of Lemma~\ref{L2:root-deg} is an alternative formula for the mapping
degree of a continuous function $\varphi\colon\R\P^1\to\R\P^1$.
Suppose that $p\in\R\P^1$ is a point with finitely many inverse images $\varphi^{-1}(p)$. 
To each inverse image we associate an index that records the behavior of 
$\varphi(t)$ as $t$ increases past the inverse image.
The index is $+1$ if $\varphi(t)$ increases past $p$, it is $-1$ if 
$\varphi(t)$ decreases past $p$, and it is $0$ if $\varphi$ stays on the same side of $p$.
(Here, increase/decrease are taken with respect to the orientation of $\R\P^1$.)
For example, here is a graph of a function $\varphi$ in relation to the value $p$ with the
indices of inverse images indicated.
\[
  \begin{picture}(270,118)(-30,0)
   \put(-20,0){\includegraphics{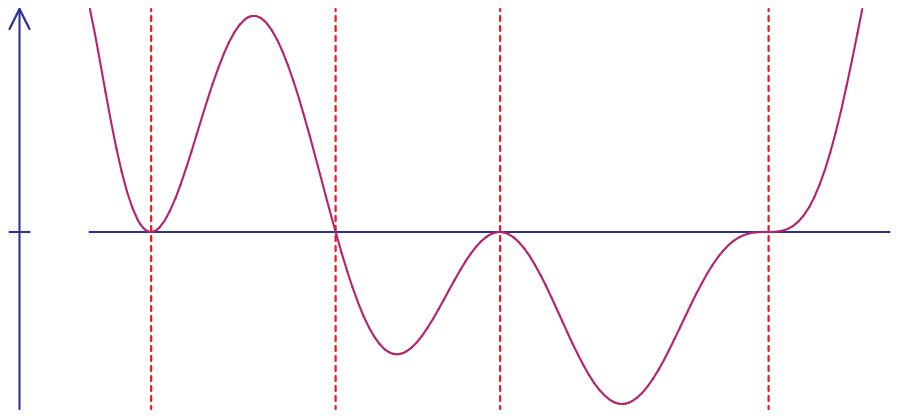}}

   \put(230,110){$\varphi$}  \put(-30,49){$p$}
   \put(22,40){$0$}  \put(74,53){$-1$}
   \put(122,53){$0$} \put(182,53){$+1$}
  \end{picture}
\]
With this definition, the mapping degree of $\varphi$ is the sum of the indices of the
points in a fiber $\varphi^{-1}(p)$, whenever the fiber is finite.
That is, 
\[
   \mdeg(\varphi)\ =\ \sum_{a\in\varphi^{-1}(p)} \mbox{index of $a$}\,.
\]

\noindent{\it Proof of Lemma~$\ref{L2:root-deg}$.}
  The zeroes of the rational function $\varphi:=f/f'$ coincide with the zeroes of $f$.
  Suppose $f(a)=0$ so that $a$ lies in $\varphi^{-1}(0)$.
  The lemma will follow once we show that $a$ has index $+1$.
  Then we may write $f(t)=(t-a)^d h(t)$, where $h$ is a polynomial with $h(a)\neq 0$.
  We see that $f'(t)=d(t-a)^{d-1}h(t) + (t-a)^d h'(t)$, and so
\[
   \varphi(t)\ =\ \frac{f(t)}{f'(t)}\ =\ 
    \frac{(t-a)h(t)}{dh(t)+ (t-a)h'(t)}\ \approx\ \frac{t-a}{d}\,,
\]
 the last approximation being valid for  $t$ near $a$ as $h(t)\neq 0$.
 Since $d$ is positive, we see that the index of the point $a$ in the fiber
 $\varphi^{-1}(0)$ is $+1$. 
\QED

\noindent{\it Proof of Theorem~$\ref{T2:top-degree}$.}
 Suppose first that $\varphi$ and $\psi$ are rational functions with no common poles.
 Then
\[
   \mdeg(\varphi+\psi)\ =\ \mdeg(\varphi) + \mdeg(\psi)\,.
\]
 To see this, note that $(\varphi +\psi)^{-1}(\infty)$ is just the union of the sets 
 $\varphi^{-1}(\infty)$ and $\psi^{-1}(\infty)$, and the index of a pole of 
 $\varphi$ equals the index of the same pole of  $(\varphi +\psi)$.

 Next, observe that $\mdeg(\varphi)=-\mdeg(1/\varphi)$.
 For this, consider the behavior of $\varphi$ and $1/\varphi$ near the level set $1$.
 If $\varphi>1$ than $1/\varphi<1$ and vice-versa.
 The two functions have index $0$ at the same points, and opposite index
 at the remaining points in the fiber $\varphi^{-1}(1)=(1/\varphi)^{-1}(1)$.

 Now consider the mapping degree of $\varphi=f/g$ as we construct its continued fraction 
 expansion.
 At the first step $f=f_0g+h$, so that $\varphi=f_0 + h/g$.
 Since $f_0$ is a polynomial, its only pole is at $\infty$, but
 as the degree of $h$ is less than the degree of $g$, $h/g$ does not have a pole at
 $\infty$.
 Thus the mapping degree of $\varphi$ is 
\[
   \mdeg(f_0 + h/g)\ =\ 
    \mdeg(f_0) + \mdeg\bigl(\tfrac{h}{g}\bigr)\ =\ 
    \mdeg(f_0) - \mdeg\bigl(\tfrac{g}{h}\bigr)\,.
\]
 The theorem follows by induction, as $\mdeg(f_0)=[f_0]$.
\QED

We close this chapter with an application of this method.
Suppose that we are given two polynomials $f$ and $g$, and we wish to count the
zeroes $a$ of $f$ where $g(a)>0$.
If $g=(x-b)(x-c)$ with $b<c$, then this will count the zeroes of $f$ in the interval
$[b,c]$, which we may do with either of the main results of this chapter.
The question is much more general, and it is not {\it a priori} clear how to use 
the methods in the first two sections of this chapter to solve this problem. 

A first step toward solving this problem is to 
compute the mapping degree of the rational function
\[
  \varphi\ :=\ \frac{f}{gf'}\,.
\]
We consider the indices of its zeroes.
First, the zeroes of $\varphi$ are those zeroes of $f$ that are not zeroes of $g$,
together with a zero at infinity if $\deg(g)>1$.
If $f(a)=0$ but $g(a)\neq 0$, then $f=(t-a)^dh(t)$ with $h(a)\neq 0$.
For $t$  near $a$,
\[
   \varphi(t)\ \approx\ \frac{t-a}{d\cdot g(a)}\,,
\]
and so the preimage $a\in\varphi^{-1}(0)$ has index $\sign(g(a))$.
If $\deg(g)=e>1$ and $\deg(f)=d$ then the asymptotic expansion of $\varphi$ for $t$ near
infinity is
\[
   \varphi(t)\ \approx\ \frac{1}{d g_e t^{e-1}}\,,
\]
where $g_e$ is the leading coefficient of $g$.
Thus the index of $\infty\in\varphi^{-1}(0)$ is 
$\sign(g_e)(e{-}1\mod 2)\ =\ [g']$.
We summarize this discussion.

\begin{lemma}\label{L:sign_counting}
 If $\deg(g)>1$, then 
\[
   \sum_{\{a\mid f(a)=0\}} \sign(g(a))\ =\ \mdeg(\varphi)\ -\ [g']\,,
\]
 and if $\deg(g)=1$, the correction term $-[g']$ is omitted.
\end{lemma}
Since $\mdeg(\varphi)=-\mdeg(1/\varphi)$, we have the alternative expression
for this sum.

\begin{lemma}
 Let $q_1,q_2,\dotsc,q_k$ be the successive quotients in the Euclidean algorithm 
 applied to the division of $f'g$ by $f$.
 Then
\[
   \sum_{\{a\mid f(a)=0\}} \sign(g(a))\ =\ [q_2]-[q_3]+\dotsb+(-1)^k[q_k]\,.
\]
\end{lemma}

\noindent{\it Proof.}
  We have
\[
   \mdeg\frac{f}{f'g}\ =\ -\mdeg\frac{f'g}{f}\ =\ 
    -[q_1]+[q_2]- \dotsb +(-1)^k[q_k]\,,
\]
 by Theorem~\ref{T2:top-degree}.
 Note that we have $f'g=q_1f+r_1$.
 If we suppose that $\deg(f)=d$ and $\deg(g)=e$, then $\deg(q_1)=e-1$.
 Also, the leading term of $q$ is $dg_e$, where $g_e$ is the leading term of $g$,
 which shows that $[q_1]=[g']$.
 Thus the lemma follows from Lemma~\ref{L:sign_counting}, when $\deg(g)\geq 2$.

 But it also follows when $\deg(g)<2$ as $[q_1]=0$ in that case.
\QED

Now we may solve our problem.
For simplicity, suppose that $\deg g>1$.
Note that 
\[
   \tfrac{1}{2} \bigl( \sign(g^2(a))+ \sgn(g(a))\ =\ 
   \left\{ \begin{array}{rcl} 1&&\mbox{if }g(a)>0\\0&&\mbox{otherwise}\end{array}\right.\ .
\]
And thus
\[
  \#\{a\mid f(a)=0,\ g(a)>0\}\ =\ 
  \frac{1}{2} \mdeg\left(\frac{f}{g^2f'}\right)\ +\ 
  \frac{1}{2} \mdeg\left(\frac{f}{gf'}\right)\,,
\]
which solves the problem.

%
%
%
%
\chapter{Sparse Polynomial Systems} \label{Ch:sparse}
%
%
%
Consider a system of $n$ polynomials in $n$ variables
 \begin{equation}\label{E3:PolySystem}
  f_1(x_1,\dotsc,x_n)\ =\ f_2(x_1,\dotsc,x_n)\ =\ 
   \dotsb\ =\ f_n(x_1,\dotsc,x_n)\ =\ 0\,,
 \end{equation}
where the polynomial $f_i$ has total degree $d_i$.
By B\'ezout's Theorem~\cite{Bez}, this system has
at most $d_1 d_2\dotsb d_n$ isolated complex solutions, and exactly that number if the
polynomials are generic among all polynomials with the given degrees.

Polynomials in nature (e.g.~from applications) are
not necessarily generic---they often have some additional structure which we
would like our count of solutions to reflect.

\begin{ex}\label{E3:sparse}
 Consider the system of two polynomials in the variables $x$ and $y$
\[
   f\ :=\ x^2y+2xy^2+xy-1\ =\ 0
     \quad\mbox{and}\quad 
   g\ :=\ x^2y-xy^2-xy+2\ =\ 0\,.
\]
Since these equations have the algebraic consequences,
 \begin{eqnarray*}
  f\cdot(y-x+1) + g\cdot(x+2y+1)& = &  3 y + 3 x + 1\,,\qquad\mbox{ and}\\
  f\cdot(3y^3-3xy^2+5y^2-2xy+2y-3)\qquad\qquad&&\\
  +\ g\cdot(6y^3+3xy^2+7y^2+2xy+2y+3)& =& 9y^3+9y^2+2y+9\,,
 \end{eqnarray*}
we see that the system has three solutions.
(The  linear consequence allows us to recover $f$ and $g$ from the cubic
in $y$.)

Both polynomials $f$ and $g$ have degree three, but they only have three common solutions,
which is fewer than the nine predicted by B\'ezout's Theorem.
The key idea behind this deficit of $6=9-3$ is illustrated by plotting the exponent
vectors of monomials which occur in the polynomials $f$ and $g$.
 \begin{equation}\label{E3:triangle}
 \begin{minipage}{3cm}
  \begin{eqnarray*}
   x^2y&\leftrightarrow&(2,1)\\
   xy^2&\leftrightarrow&(1,2)\\
   xy&\leftrightarrow&(1,1)\\
   1&\leftrightarrow&(0,0)\\\ \\
  \end{eqnarray*}
 \end{minipage}\qquad\qquad
 \raisebox{-49pt}{
  \begin{picture}(115,115)(0,-3)
   \put( 0, 0){\includegraphics{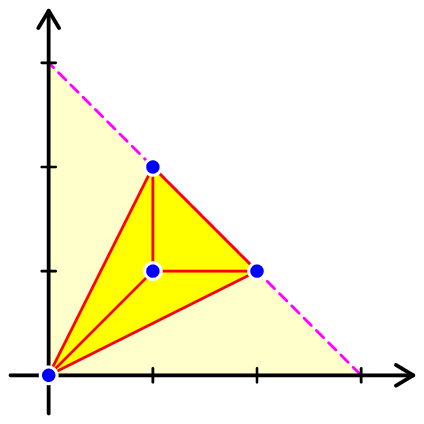}}
   \put(38,-3){$1$}   \put(0,37){$1$}
   \put(68,-3){$2$}   \put(0,67){$2$}
   \put(98,-3){$3$}   \put(0,97){$3$}
  \end{picture}}
 \end{equation}
The \DeCo{{\sl Newton polytope}} of $f$ (and of $g$) is the convex hull of these
exponent vectors.
This triangle has area $\frac{3}{2}$ and is covered
by \DeCo{three} lattice triangles. 
We will see why this number of lattice triangles equals the number of
solutions to the original system. \QED
\end{ex}

%
\section{Kushnirenko's  Theorem}
%
The polynomial system in Example~\ref{E3:sparse} is a sparse system
whose support is the set of integer points in the triangle of~\eqref{E3:triangle}. 
Let $\DeCo{\calA}\subset\Z^n$ be a finite set of exponent vectors which affinely spans
$\R^n$.
A \DeCo{{\sl sparse polynomial}} $f$ with \DeCo{{\sl support}} $\calA$
is a linear combination 
\[
   f\ =\ \sum_{a\in\calA} c_a x^a\qquad c_a\in\R
\]
of monomials from $\calA$.
Let $\Delta_\calA\subset\R^n$ be the convex hull of the vectors in $\calA$.
Let $\DeCo{\T}:=\C^\times:=\C-\{0\}$ be the algebraic torus.
We recall Kushnirenko's Theorem from Chapter~1.\medskip

\noindent{\bf Kushnirenko's Theorem.}
{\it
 A system~\eqref{E3:PolySystem} of $n$ polynomials in $n$ variables 
 with common support $\calA$ has at most $n!\vol(\Delta_\calA)$ isolated
 solutions in $\T^n$, and exactly this number if the polynomials are generic
 given their support $\calA$. 
}\medskip

We will give two proofs of this result, one using algebraic geometry that is due to
Khovanskii and another which is algorithmic. 
Each proof introduces some important geometry related to sparse systems of polynomials.

It is worth remarking that while sparse polynomials occur naturally---multilinear or
multihomogeneous polynomials are an example---they also occur in problem formulations due
to human psychology.
Most of us are unable to write down or reason with polynomials having thousands of terms,
and we instead seek problem formulations with fewer terms.

%
\subsection{The geometry of sparse polynomial systems}\label{S3:sparse}
%
Consider the map
\[
   \DeCo{\varphi_{\calA}}\ \colon\ \T^n\ni x\ 
   \longmapsto\ [x^a\mid a\in \calA]\ \in\ \P^{\calA}\,,
\]
where \DeCo{$\P^{\calA}$} is the projective space with homogeneous coordinates 
$[z_a\mid a\in\calA]$ indexed by $\calA$.
This map factors
\[
  \T^n\ \longrightarrow\ \T^\calA\ \joinrel\relbar\twoheadrightarrow\ 
  \T^\calA/\delta(\T)\ \subset\ \P^{\calA}\,,
\]
where $\DeCo{\T^\calA}=(\C^\times)^{|\calA|}$ is the torus with coordinates indexed
by $\calA$ and \DeCo{$\delta(\T)\subset\T^\calA$} is the diagonal torus.
The quotient \DeCo{$\T^\calA/\delta(\T)$} is the dense torus in the projective space
$\P^\calA$. 
It consists of those points $[z_a\mid a\in\calA]\in\P^\calA$ with no coordinate zero.
Notice that $\varphi_\calA$ is a homomorphism into this dense torus.
It is often convenient to identify $\calA$ with the $n\times|\calA|$ matrix whose columns
are the exponent vectors in $\calA$.

\begin{ex}\label{E3:Hex1}
  Suppose that $\calA$ consists of the seven exponent vectors
 $(0,0)$, $(1,0)$, $(0,1)$, $(1,1)$, $(-1,0)$, $(0,-1)$, and 
 $(-1,-1)$.  
 Here is the corresponding matrix,
\[
   \left(\begin{array}{rrrrrrr}
         0&1&1&0&-1&-1& 0\\
         0&0&1&1& 0&-1&-1\end{array}\right)\,.
\]
 The convex hull $\Delta_\calA$ of these points is the hexagon,
\[
  \begin{picture}(122,122)(0,0)
   \put( 0, 0){\includegraphics[height=122.5pt]{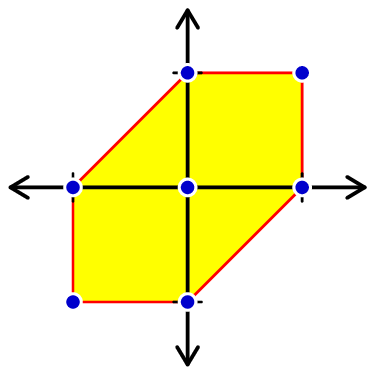}}
   \put(98.5,42){$1$}     \put(46.8,96){$1$}
   \put(68.4,18.6){$-1$}  \put(8,72){$-1$}
  \end{picture}
\]
 and the map $\varphi_\calA$ is
\[
  \varphi_\calA\ \colon (x,y)\in \T^2\ 
      \longmapsto\ 
    [1,x,y,xy,x^{-1},y^{-1},x^{-1}y^{-1}]\in\P^\calA\ =\ \P^6\,. 
   \makebox[.1in][l]{\hspace{1.8cm}\raisebox{-2pt}{\includegraphics[height=12pt]{figures/HSBC.eps}}}
\]
\end{ex}

Consider a linear form $\Lambda$ on $\P^\calA$, 
\[
  \Lambda\ =\ \sum_{a\in\calA} c_a z_a\,.
\]
Its pullback $\varphi^*_{\calA}(\Lambda)$ 
along $\varphi_{\calA}$ is a polynomial with support $\calA$,
\[
   \varphi^*_{\calA}(\Lambda)\ =\ 
    \sum_{a\in\calA} c_a x^a\,.
\]

This provides a bijective correspondence between linear forms
on $\P^{\calA}$ and sparse polynomials with support $\calA$.
Under $\varphi_{\calA}$, the zero set of a sparse polynomial is mapped 
to a hyperplane section of $\varphi_\calA( \T^n)$
(the hyperplane is where the corresponding linear form vanishes).
Since $n$ general linear forms on $\P^\calA$ cut out a linear subspace $L$ of
codimension $n$, a general system~\eqref{E3:PolySystem} of sparse polynomials with
support $\calA$ is the pullback along $\varphi_{\calA}$ of a codimension $n$
linear section of $\varphi_{\calA}(\T^n)$. 
That is, it equals $\varphi_{\calA}^{-1}(L)$ or 
$\varphi_{\calA}^{-1}(\varphi(\T^n)\cap L)$.

The closure of the image of $\varphi_\calA$ 
is the \DeCo{{\sl toric variety}}\fauxfootnote{There are competing
  notions of toric variety.  We follow the convention from symbolic
  computation~\cite{GBCP}, rather than from algebraic geometry~\cite{Fu93}.
  In particular, we do not assume that $X_\calA$ is normal.}
\DeCo{$X_\calA$} parameterized by the monomials $\calA$. 
Since $\varphi_\calA$ is a homomorphism, the number $d_\calA$ of
solutions to a general sparse system with support $\calA$ is the
product
 \begin{equation}\label{E3:Da_Formula}
   d_\calA\ =\ |\ker(\varphi_\calA)|\cdot \deg(X_\calA)\ .
 \end{equation}
Indeed, if $L$ is a general linear subspace of $\P^\calA$ of
codimension $n$, then Bertini's theorem implies that 
$\varphi_\calA(\T^n)\cap L=X_\calA\cap L$ and this intersection is
transverse\fauxfootnote{This transversality will also follow from the arguments given in our  
  second proof of Kushnirenko's Theorem.}.
When the intersection is not transverse, $d_\calA$ will be the sum of the
multiplicities of the solutions.
The number of points in such a linear section is the degree \DeCo{$\deg(X_\calA)$} of
$X_\calA$, and each point pulls back under 
$\varphi_\calA$ to $|\ker(\varphi_\calA)|$ solutions of the sparse
system corresponding to the linear section.

%
\subsection{Algebraic-geometric proof of Kushnirenko's Theorem}
%
We prove Kushnirenko's Theorem by showing that 
\[
  n!\cdot \vol(\Delta_\calA)\ =\ 
   |\ker(\varphi_\calA)|\cdot \deg(X_\calA)\ =\ d_\calA \ .
\]
This proof is due to Khovanskii~\cite{Kh95}
The same idea of proof is used by Khovanskii and Kaveh~\cite{KK}, where they extend the theory 
of Newton polyhedra for functions on $\T^n$ to convex bodies associated to functions on an
arbitrary affine variety.

We first determine the kernel of the map $\varphi_\calA$, which is the composition
 \begin{eqnarray*}
  \T^n & \longrightarrow& \T^\calA\ \longrightarrow\ \T^\calA/\delta(\T)\ 
         \subset\ \P^\calA\\
   x&\longmapsto& (x^a\mid a\in\calA)\ .
 \end{eqnarray*}
To facilitate this computation, we assume that $0\in\calA$.
This is no loss of generality, for if $0\not\in\calA$, then we simply translate
$\calA$ so that one of its exponent vectors is the origin.
This has the effect of multiplying each point in $\varphi_\calA(\T^n)$ by a
scalar, and so it does not change $X_\calA$.
It also multiplies each polynomial in~\eqref{E3:PolySystem} by a common
monomial, which affects neither the solutions in $\T^n$ nor their number
$d_\calA$. 
By the relation~\eqref{E3:Da_Formula}, this translation does not change the
cardinality of the kernel of $\varphi_\calA$. 

Suppose that $0\in\calA$.
Then the $z_0$-coordinate of $\varphi_\calA$ is constant ($x^0=1$) and so the map 
which sends $x\in\T^n$ to $(x^a\mid a\in\calA)$ maps
$\T^n$ into $1\times\T^N\subset\T^\calA$, where $|\calA|=N+1$.
The composition of the two maps 
\[
  1\times \T^N\ \longrightarrow\ \T^\calA\ 
   \longrightarrow\ \T^\calA/\delta(\T)\ \left(\subset\ \P^\calA\right)
\]
is an isomorphism.
Thus it is sufficient to compute the kernel of the map
 \begin{eqnarray*}
  \psi_\calA\ \colon\ \T^n & \longrightarrow& \T^\calA\,,\\
   x&\longmapsto& (x^a\mid a\in\calA)\,,
 \end{eqnarray*}
when $0\in\calA$.

Let $\Z\calA\subset\Z^n$ be the sublattice spanned by the exponent vectors in
$\calA$.
It has full rank, by our assumption that $\calA$ affinely spans $\R^n$.
Then the quotient $\Z^n/\Z\calA$ is a finite abelian group whose order is the 
\DeCo{{\sl lattice index}} $[\Z^n : \Z\calA]$.
Its group $\Hom(\Z^n/\Z\calA,\T)$ of characters (homomorphisms to $\T$) is 
also called its \DeCo{{\sl Pontryagin dual}}.
We have the sequence of abelian groups,
\[
   0\ \longrightarrow\ \Z\calA\ \longrightarrow\ \Z^n\ 
    \relbar\joinrel\twoheadrightarrow\ \Z^n/\Z\calA\ \longrightarrow\ 0\,.
\]
This gives rise to the sequence of Pontryagin duals,
\[
  0\ \longrightarrow\ \Hom(\Z^n/\Z\calA, \T)\ 
      \longrightarrow\ \Hom(\Z^n, \T)
      \longrightarrow\ \Hom(\Z\calA, \T)\,,
\]
which is exact in that the image of $\Hom(\Z^n/\Z\calA, \T)$ in $\Hom(\Z^n, \T)$ is the
kernel of the last map.
If we identify $\Hom(\Z^n, \T)$ with $\T^n$, then this sequence of groups, or direct
calculation, shows that 
\[
  \Hom(\Z^n/\Z\calA, \T)\ =\ 
   \{ x\in \T^n \mid x^a=1 \ \ \forall a\in\calA\}\ =\ 
   \ker(\varphi_\calA)\,.
\]
Thus $\ker(\varphi_\calA)$ is the Pontryagin dual to the quotient
$\Z^n/\Z\calA$.
Since $\Z^n/\Z\calA$ is a finite abelian group, we see that 
 \begin{equation}\label{E3:kernel}
  |\ker(\varphi_\calA)|\ =\  [\Z^n : \Z\calA]\ =\ |\Z^n/\Z\calA|\,.\medskip
 \end{equation}

The \DeCo{{\sl Hilbert polynomial} $h_X(d)$} a projective variety $X$ is the
polynomial which is eventually equal to the dimension of the
$d$th graded piece of the homogeneous coordinate ring $\C[X]$ of $X$,
\[
   h_X(d)\ =\ \dim_\C \C_d[X]\,,  \mbox{ for all $d$ sufficiently large.}
\]
The Hilbert polynomial contains many numerical invariants of $X$.
For example, the degree of the Hilbert polynomial is the dimension $n$ of $X$
and its leading coefficient is~$\deg(X)/n!$.
For a discussion of Hilbert polynomials, see Section 9.3 of~\cite{CLO}.

We determine the Hilbert polynomial of the toric variety $X_\calA$.
For this, it is helpful to consider the homogeneous version of the
parametrization map $\varphi_\calA$.
We lift $\calA\subset\Z^n$ to a homogenized set of exponent
vectors $\calA^+\subset 1\times\Z^n$ by prepending a component of $1$ to each
vector in  $\calA$. 
That is, 
\[
   \DeCo{\calA^+}\ :=\ \{(1,a)\mid a\in\calA\}\,.
\]
The matrix $A^+$ is obtained from the matrix $A$ by
adding a new first row of 1s.
For the hexagon of Example~\ref{E3:Hex1}, this is 
 \[
   \calA^+\ =\ \left(\begin{array}{ccccrrr}1&1&1&1&1&1&1\\
                     0&1&1&0&-1&-1&0\\
                    0&0&1&1&0&-1&-1\end{array}\right)\ .
 \]
Here is the lifted hexagon, where the first coordinate is vertical.
\[
  \includegraphics[height=3cm]{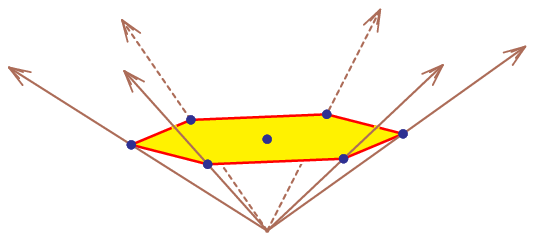}
\]

The map $\varphi_{\calA^+}$ on $\T^n$ has the same image
in $\P^\calA$ as does $\varphi_\calA$.  
The advantage is that the image of the map $\varphi_{\calA^+}$ in $\C^\calA$ is
stable under multiplication by scalars---this is built into it via the new first
coordinate of $\calA^+$. 
If we let $t$ be the first coordinate of $\T^{1+n}$, then 
the pullback of the coordinate ring of $\P^\calA$ to the ring of Laurent
polynomials (the coordinate ring of $\T^{1+n}$) is 
\[
   \DeCo{S_\calA}\ :=\ \C[t x^a\mid a\in\calA]\ \simeq\ \C[\N\calA^+]\,.
\]
This is also the homogeneous coordinate ring of the toric variety 
$X_\calA=\overline{\varphi_\calA(\T^{1+n})}$ (the closure taken in $\P^\calA$).

The grading on $S_\calA$ is given by the
exponent of the variable $t$. 
It follows that the $d$th graded piece of $S_\calA$ has a basis of monomials
\[
    \{ t^d\cdot x^a\mid a\in d\calA\}\,,
\]
where \DeCo{$d\calA$} is the set of $d$-fold sums of vectors in $\calA$.
This is just the set of integer points in $d\Delta_\calA$ which lie in the 
nonnegative integer span $\N\calA$ of $\calA$.
If we let $H_\calA(d)$ be the dimension of the $d$th graded piece of the
homogeneous coordinate ring of $X_\calA$ (also called the 
\DeCo{{\sl Hilbert function}} of $X_\calA$), then these
arguments show that 
\[
  H_\calA(d)\ =\ |d\Delta_\calA \cap \N\calA|.
\]
We will estimate this Hilbert function, which will enable us to determine the
leading coefficient of the Hilbert polynomial, as the Hilbert function and Hilbert
polynomial agree for $d$ sufficiently large.

Let $\Lambda\simeq\Z^n$ be a lattice in $\R^n$ and $\Delta$ be a polytope with
vertices in $\Lambda$.
Ehrhart~\cite{Eh62} showed that the counting function
\[
    P_\Delta\ \colon\ \N\ni d\ \longmapsto \ |d\Delta\cap\Lambda|
\]
for the points of $\Lambda$ contained in positive integer multiples of the polytope
$\Delta$ is a polynomial in $d$.
This polynomial is called the \DeCo{{\sl Ehrhart polynomial}} of the
polytope $\Delta$, and its degree is the dimension of the affine span of
$\Delta$.
When $\Delta$ has dimension $n$, its leading coefficient is the volume of
$\Delta$, normalized so that a fundamental domain of the lattice $\Lambda$ has
volume 1.
That is, it is the Euclidean volume divided by the lattice index $[\Z^n : \Z\calA]$.
When $\Lambda=\Z^n$, this is the ordinary Euclidean volume of $\Delta$.

Now suppose that $\Delta=\Delta_\calA$, the convex hull of $\calA$.
Since $d\calA\subset d\Delta_\calA\cap \Z\calA$, if $\Lambda=\Z\calA$, we have 
 \begin{equation}\label{Eq3:Hilbert_compare}
   P_{\Delta_\calA}(d)\ \geq\ H_{\calA}(d)\,.
 \end{equation}

We give a lower bound for $H_{\calA}(d)$.
Let $\calB$ be the points $b$ in $\Z\calA$ which may be written as
\[
   b\ =\ \sum_{a\in\calA} c_a a\,,
\]
where $c_a$ is a rational number in $[0,1)$.
Fix an expression for each $b\in\calB$ as an integer linear combination of
elements of $\calA$, and let $-N$ with $N\geq 0$ be an integer lower bound for the
coefficients in these expressions for the finitely many elements of $\calB$.

For $d\geq N|\calA|$ we claim that translation by the vector
$N\sum_{a\in\calA}a$ defines a map
\[
   \Z\calA\cap (d- N|\calA|)\Delta_\calA
    \ \longrightarrow\ \N\calA \cap d\Delta_\calA\,.
\]
Indeed, a point $v\in \Z\calA\cap (d-N|\calA|)\Delta_\calA$ is a nonnegative
rational combination of the vectors in $\calA$.
Taking fractional parts gives $v=b+c$, where $b\in\calB$ and $c\in \N\calA$.
Adding $N\sum_{a\in\calA}a$ to the fixed integral expression of $b$ gives a
positive integral expression, which proves the claim.
This shows that 
\[ 
    H_{\calA}(d)\ \geq\ 
     P_{\Delta_\calA}(d-N|\calA|)\,.
\]
If we combine this estimate with~\eqref{Eq3:Hilbert_compare}, and use the result that the 
Hilbert function equals the Hilbert polynomial for $d$
large enough, then we have shown that the Hilbert polynomial $h_{\calA}$
of $X_\calA$ has the  same degree and leading coefficient as the Ehrhart
polynomial $P_{\Delta_\calA}$.

Thus the Hilbert polynomial has degree $n$ and its leading coefficient is
$\vol(\Delta_\calA)/[\Z^n\colon \Z\calA]$, which is the normalized volume of the polytope 
$\Delta_\calA$ with respect to the lattice $\Z\calA$.
We conclude that the degree of $X_\calA$ is
\[
   n!\; \frac{\vol(\Delta_\calA)}{[\Z^n\colon \Z\calA]}\,.
\]

Recall~\eqref{E3:kernel} that the kernel of $\varphi_\calA$ has order
$[\Z^n\colon \Z\calA]$.
Then the formula~\eqref{E3:Da_Formula} for the number $d_\calA$ of solutions to a
sparse system~\eqref{E3:PolySystem} with support $\calA$ becomes 
\[
   d_\calA\ =\ |\ker(\varphi_\calA)|\cdot \deg(X_\calA)\ =\ 
   [\Z^n\colon \Z\calA]\cdot n!\;\frac{\vol(\Delta_\calA)}{[\Z^n\colon \Z\calA]}
   \ =\ n!\, \vol(\Delta_\calA)\,,
\]
which proves Kushnirenko's Theorem.\QED

%
\section{Algorithmic proof of Kushnirenko's Theorem}
%

We present a second proof of Kushnirenko's Theorem
whose advantage is that it introduces more geometry connected to toric varieties
(which will be useful in subsequent chapters), in particular, the proof uses 
\DeCo{{\sl toric degenerations}} (sometimes called Gr\"obner degenerations).
These toric degenerations underlie the method of Viro~\cite{Viro,St94b}, which is
important for many constructions in real algebraic geometry.
We first analyze the simplest case of Kushnirenko's Theorem---when $|\calA|=n+1$.

%
\subsection{Kushnirenko's Theorem for a simplex}\label{S3:simplex}
%
Suppose that $|\calA|=n+1$ so that $\Delta_\calA$ is a simplex with vertices
$\calA$. 
Then $X_\calA = \P^n=\P^\calA$.
As we saw, the solutions to any sparse system~\eqref{E3:PolySystem} with support $\calA$ 
have the form $\varphi_\calA^{-1}(L)$, where $L\subset\P^n$ is the codimension-$n$ plane
cut out by the linear forms which define the polynomials of the system.
In this case $L$ is simply a point $\beta\in\P^n$ (which lies in the dense torus as the
equations are general), and the solutions have the form $\varphi_\calA^{-1}(\beta)$.
Since $\varphi_\calA$ is a homomorphism to the dense torus of $\P^\calA$, these solutions
form a single coset of $\ker(\varphi_\calA)$.

We may determine these solutions explicitly.
Assume that $0\in\calA$.
Since $|\calA|=n+1$, we may write the sparse system~\eqref{E3:PolySystem} as
 \begin{equation}\label{E3:matrix_system}
   C\cdot (x^{a_1}, x^{a_2}, \dotsc, x^{a_n})^T\ =\ b\,,
 \end{equation}
where $\calA-\{0\}=\{a_1,a_2,\dotsc,a_n\}$, $C$ is the $n$ by $n$ matrix of
coefficients, and $b\in\C^n$.
If our system is generic, then $C$ is invertible and we may perform row
operations on $C$ and hence on the system~\eqref{E3:matrix_system} to obtain an
equivalent binomial system
\[
  x^{a_i}\ =\ \beta_i\qquad\mbox{for }i=1,\dotsc,n\,.
  \eqno{\eqref{E3:matrix_system}'}
\]
where $\beta_1,\dotsc,\beta_n\in\T$,
and so the system has the form $\varphi_\calA^{-1}(\beta)$.
(In fact, the requirements that $C$ be invertible and that the resulting
 constants $\beta_i\in\T$ are the conditions for genericity of this system.)

Let $A$ be the $n$ by $n$ matrix whose columns are the exponent vectors in
$\calA-\{0\}$. 
We will see how  the \DeCo{{integer linear algebra}} of the matrix
$A$ is used to solve the system~\eqref{E3:matrix_system}$'$.

\begin{ex}\label{Ex3:simplex}
 Consider the system of equations with support
 $\calA=\{(16,14),  (22,18), (0,0)\}$
 \begin{equation}\label{Eq3:simplexEq}
  \begin{array}{rcl}
    23 x^{16}y^{14}\ -\ x^{22}y^{18}& =& -27\,,\qquad\mbox{and}\\
    35 x^{16}y^{14}\ -\ x^{22}y^{18}& =&  -9\,.\rule{0pt}{14pt}
  \end{array}
 \end{equation}
 Subtracting the two equations gives the binomial, 
\[
   12 x^{16}y^{14}\ =\ 36
   \qquad\mbox{or}\qquad
    x^{16}y^{14}\ =\ 3\,,
\]
 and solving for $x^{22}y^{18}$ gives $x^{22}y^{18}=96$.
 Thus we have the equivalent binomial system 
 \begin{equation}\label{E3:binomial_system}
  x^{16}y^{14}\ =\ 3
   \qquad\mbox{and}\qquad
  x^{22}y^{18}\ =\ 96\,.
 \end{equation}
 Under the invertible substitution  (the inverse is given by $u=x^8y^7$ and $v=(xy)^{-1}$)
 \begin{equation}\label{E3:Coordinate_Change}
   x\ =\ uv^7\qquad\mbox{and}\qquad y\ =\ u^{-1}v^{-8}\,,
 \end{equation}
 our equations become triangular
 \begin{equation}\label{E3:triangular}
  \begin{array}{rlll}
   (uv^7)^{16}(u^{-1}v^{-8})^{14}& =\ 
    u^{16}v^{112}u^{-14}v^{-112}  &  =\ 
    u^2 & =\  3\,,\quad\mbox{ and}\\
   (uv^7)^{22}(u^{-1}v^{-8})^{18}& =\ 
    u^{22}v^{154}u^{-18}v^{-144}  &  =\ 
    u^4v^{10} & =\  96\,.
  \end{array}
 \end{equation}
While the solution is now immediate via back substitution, we make one further 
simplifying substitution.
Write $f$ for the first equation and $g$ for the second.
Replacing $g$ by $gf^{-2}$ yields a diagonal system
which is now completely trivial to solve
\[
  \begin{array}{rlll}
    f& \colon\ 
    u^2 & =\  3\,,\\
    gf^{-2}& \colon\ 
    v^{10} & =\  32/3\,.
  \end{array}
\]
The solutions are $u=\pm\sqrt{3}$ and $v=\zeta\sqrt{2}$, where $\zeta$ runs over
all 10th roots of $\frac{1}{3}$.  
Substituting these into~\eqref{E3:Coordinate_Change}, gives the 20 solutions to our
original system of equations~\eqref{Eq3:simplexEq}. \smallskip\QED
\end{ex}

Underlying these simplifications is the relation between the integer linear
algebra of $n$ by $n$ matrices and (multiplicative) coordinate changes
in $\T^n$. 
Since $\T^n=\Hom(\Z^n,\T)$, its automorphism group is ${\it GL}(n,\Z)$, the
group of invertible $n$ by $n$ integer matrices, and this is the source of that
relation. 

\begin{ex}
The monomials in~\eqref{E3:binomial_system} correspond to the columns of the
matrix
\[
   A\ =\ \left(\begin{matrix}16&22\\14&18\end{matrix}\right)\ ,
\]
and the coordinate change~\eqref{E3:Coordinate_Change} corresponds to 
\DeCo{left} multiplication (hence row operations) by the matrix
\[
   \left(\begin{matrix}1&-1\\7&-8\end{matrix}\right)\ .
\]
Indeed,
\[
  \left(\begin{matrix}1&-1\\7&-8\end{matrix}\right)
  \left(\begin{matrix}16&22\\14&18\end{matrix}\right)\  =\ 
  \left(\begin{matrix}2&4\\0&10\end{matrix}\right)\ ,
\]
which corresponds to the exponent vectors in the triangular
system~\eqref{E3:triangular}. 
This upper triangular matrix is the 
\DeCo{{\sl Hermite normal form}} of the matrix $A$---the row reduced
echelon form over $\Z$.
This notion makes sense for matrices
whose entries lie in any principal ideal domain.

Multiplicative reductions using the \DeCo{equations} correspond to multiplicative
coordinate changes in the target torus $\T^2$ and are represented by column
operations, or \DeCo{right} multiplication by integer matrices.
Indeed
\[
  \left(\begin{matrix}2&4\\0&10\end{matrix}\right)
  \left(\begin{matrix}1&-2\\0&1\end{matrix}\right)\ =\ 
  \left(\begin{matrix}2&0\\0&10\end{matrix}\right)\ ,
\]
which is the \DeCo{{\sl Smith normal form}} of the integer matrix $A$.\QED
\end{ex}

The Smith normal form of an $n$ by $n$ matrix $A$ is the diagonal matrix
with entries $d_1, d_2, \dotsc, d_n$ where $d_i$ is the greatest common divisor
of all $i$ by $i$ subdeterminants (minors) of $A$.
These are called the \DeCo{{\sl invariant factors}} of $A$, and we have
$d_1|d_2|\dotsb|d_n$.
Then
\[
   \Z^n/\Z\calA\ \simeq\ 
    \Z/d_1\Z\, \times \,\Z/d_2\Z\, \times\, \dotsb\,\times\,\Z/d_n\Z \,, 
\]
the canonical way to write a finite abelian group.
As with the Hermite normal form, the Smith normal form makes sense for 
presentations of modules over any principal ideal domain.

%
\subsection{Regular triangulations and toric degenerations}
%

Let $\omega\colon\calA\to \N$.
This is a weight for $\T$ acting diagonally on the space $\C^\calA$,
if $z=(z_a\mid a\in\calA)\in\C^\calA$, then $t.z=(t^{\omega(a)}z_a\mid a\in\calA)$.
This induces a similar action on $\P^\calA$, and a dual action on the 
homogeneous coordinate ring $\C[z_a\mid a\in \calA]$ of $\P^\calA$, namely 
\[
   t.z_a\ =\ t^{-\omega(a)}z_a\quad\mbox{ for a}\in \calA\,,
\]
where $z_a$ is a variable (coordinate function on $\C^\calA$).
Let $\calX_\calA\subset \C\times\P^\calA$ be the closure of the family over
$\T$ of deformations of the toric variety $X_\calA$ under this action, 
\[
   \DeCo{\calX_\calA}\ :=\ 
   \overline{\{(t,\,t.z)\subset \C\times \P^\calA\mid t\in\T\,,\ z\in X_\calA\}}\ .
\]
This is a \DeCo{flat family} over $\C$~\cite[Ch.~15]{E95}.
(This technical fact implies that all fibers have the
same Hilbert polynomial.)
The fiber of $\calX_\calA$ over a point $t\in\T$ is the translated toric variety
$t.X_\calA$, and the fiber over $0\in\C$ is called the \DeCo{{\sl
scheme-theoretic limit}} of the family $t.X_\calA$, which is written
 \[
   \lim_{t\to 0} t. X_\calA\,.
 \]
We use geometric combinatorics to study this limit.
The passage from $X_\calA$ to such a scheme-theoretic limit of an action of $\T$
is a \DeCo{{\sl toric degeneration}}.

We may also use the weight $\omega\colon\calA\to \N$ to lift the vector configuration
$\calA$ into $\R^{1+n}$.
Consider the convex hull $P_\omega$ in $\R^{1+n}$ of the lifted vectors
 \begin{equation}\label{Eq3:Lifted}
  \DeCo{P_\omega}\ :=\ \conv\{(w(a),a)\, \mid\, a\in\calA\}\,.
 \end{equation}
The \DeCo{{\sl lower facets}} of this polytope are those facets of $P_\omega$ whose
outward-pointing normal vector has a negative first coordinate.
Projecting these lower facets back to $\R^n$ gives the facets in the
\DeCo{{\sl regular polyhedral subdivision $\Delta_\omega$}} of 
the convex hull of $\calA$ induced by the lifting function $\omega$. 
The vertices in this subdivision are \DeCo{some} of the vectors in $\calA$. 
While much of the following makes sense for general polyhedral subdivisions, we
shall henceforth assume that $\Delta_\omega$ is a \DeCo{{\sl triangulation}} 
in that each lower facet is a simplex.
(This is a mild genericity assumption on $\omega$.)
We display the lower facets and the resulting triangulation $\Delta_\omega$ for a weight
function $\omega$ on where $\calA=\{(0,0),(1,0),(0,1),(2,0),(1,1),(0,2)\}$.
\[
   \includegraphics[height=140pt]{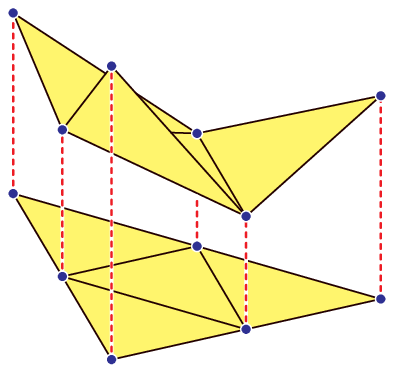}
\]

Exponent vectors $\alpha$ of monomials in the coordinate ring
$\C[z_a\mid a\in\calA]$ of $\P^\calA$ are elements of $\N^\calA$, and 
$\omega$ induces a linear form on $\N^\calA$,
\[
  \DeCo{\omega\cdot\alpha}\ :=\  \omega\cdot (\alpha_a\mid a\in\calA)\ =\ 
   \sum_{a\in\calA}\omega(a)\alpha_a\ .
\]

The \DeCo{{\sl initial form} $\ini_\omega(g)$} of a 
homogeneous form $g(z)\in\C[z_a\mid a\in\calA]$
is the sum of the terms $c_\alpha z^\alpha$ of $g$ 
for which $\omega\cdot\alpha$ is maximal among all terms of $g$. 
Let $\omega(g)$ be this maximal value.
Multiplying $t.g$ by $t^{\omega(g)}$  we see that
 \[
   t^{\omega(g)}(t. g(z))\ =\ \ini_{\omega}g(z)\ +\ h\,,
 \]
where $h$ is divisible by $t$.
Thus
 \[
   \lim_{t\to 0}\,t^{\omega(g)}(t. g(z))\ =\ \ini_{\omega}g(z)\,.
 \]

Write $I_\calA$ for the ideal defining $X_\calA$.
(This {\DeCo{\it toric ideal}} has a linear basis of
binomials $z^\alpha-z^\beta$ 
such that $A^+\alpha=A^+\beta$, where $A$ is the matrix whose columns are the
exponent vectors in $\calA^+$~\cite[Lemma 4.1]{GBCP}.) 
These binomials have the following geometric interpretation.
The product $A^+\alpha$ is a positive linear combination of the vectors in $\calA^+$,
so binomials $z^\alpha-z^\beta$ in $I_\calA$ record vectors in $\N\calA$ that have 
two distinct representations as positive linear combinations of the vectors in $\calA^+$. 
If we divide by the the initial coordinate of the vector $A^+\alpha$, we obtain a vector
of the form $(1,x)$, where $x$ is a rational point lying in the convex hull of $\calA$.
When $A^+\alpha=A^+\beta$, this pont $x$ has two distinct rational representations as a
convex combination of elements of $\calA$.
Thus $x$ lies in the convex hulls of two different subsets of $\calA$.
Conversely, any such point gives rise to a binomial $z^\alpha-z^\beta\in I_\calA$. 

When $\omega\cdot\alpha>\omega\cdot\beta$, then $z^\alpha$ is the initial term
$\ini_\omega(z^\alpha-z^\beta)$.
The \DeCo{{\sl initial ideal}} of $I_\calA$ is 
\[
  \DeCo{\ini_\omega(I_\calA)}\ =\ 
  \{ \ini_\omega(g)\mid g\in I_\calA\}\,.
\]
Since the ideal $I(t. X_\calA)$ of $t. X_\calA$ is 
\[
   I(t. X_\calA)\ =\ \left\{ t. g(y) \mid g\in I(X_\calA)\right\}
    \ =\ t.I_\calA\,,
\]
we see that 
\[
  \lim_{t\to 0} t.I_\calA\ =\ \ini_\omega(I_\calA)\,.
\]
Thus this initial ideal is the ideal of the scheme-theoretic limit of the 
family $t.X_\calA$.

%
\subsection{Kushnirenko's Theorem via toric degenerations}
%

Since the family ${\mathcal X}_\calA\to \C$ is flat, every fiber has
the same degree, and so Kushnirenko's Theorem follows if we can prove that the
degree of the limit scheme $\ini_\omega(X_\calA)$ is 
$n!\vol(\Delta_\calA)$ divided by the degree $[\Z^n\colon \Z\calA]$ of
$\varphi_\calA$. 
It is proved in Chapter 8 of~\cite{GBCP} that
 \begin{equation}\label{E3:rad_decomp}
  \sqrt{\ini_\omega(I_\calA)\,}\ =\ 
    \bigcap_{\tau} \langle z_a\ \mid\ a\not\in\tau \rangle\ ,
 \end{equation}
the intersection over all facet $n$-simplices $\tau$ of the regular
triangulation of $\Delta_\calA$.
This result is not so hard.
If the segment $\overline{ab}$ for $a,b\in\calA$ is not a face of the triangulation
$\Delta_\omega$, then it crosses a minimal face conv$(\tau)$ of the triangulation.
This implies that there is a binomial $z_a^Nz_b^M-z^\gamma$, where 
$M,N$ are positive integers and the monomial $z^\gamma$ involves the variables in $\tau$.
By the construction of the triangulation $\Delta_\omega$, the
corresponding lifted segment $\overline{(\omega(a),a),(\omega(b),b)}$
of $P_\omega$ lies above the lift of the face conv$(\tau)$,
and thus the initial term of this binomial is $z_a^Nz_b^M$,
and so $z_az_b$ lies in the radical $\sqrt{\ini_\omega(I_\calA)}$ of the initial ideal.

It follows that the limit scheme $\ini_\omega(X_\calA)$ is supported on the union
of coordinate $n$-planes $\P^\tau$, one for each facet $n$-simplex $\tau$ of the
regular triangulation.
(Here $\P^\tau$ is the coordinate plane which is spanned by the coordinates indexed by
$\tau$.) 
The degree of this initial scheme is then
\[
   \deg(\ini_\omega(X_\calA))\ =\ 
   \sum_\tau \mbox{\rm mult}_{\P^\tau}(\ini_\omega(X_\calA))\ ,
\]
the sum over these facets $\tau$ of the algebraic multiplicity of 
the limit scheme $\ini_\omega(X_\calA)$ along the coordinate $n$-plane
$\P^\tau$.

In~\cite[Chapter 8]{GBCP}, and under the (mild) assumption that $\calA$ is
primitive ($\Z\calA=\Z^n$), Sturmfels shows that this multiplicity is
$n!\vol(\Delta_\tau)$. 
Since these facets cover $\Delta_\calA$, the degree of the limit scheme
$\ini_\omega(X_\calA)$ is $n!\vol(\Delta_\calA)$.
As the family $\calX_\calA$ is flat, this degree is the degree of
$X_\calA$, and we may deduce Kushnirenko's Theorem from this.

Before we continue, we deduce a corollary from this.
A triangulation of a polytope in $\R^n$ is \DeCo{{\sl unimodular}} is every facet has
minimal volume $1/n!$.

\begin{cor}\label{C3:reg_triangulation}
 Suppose that $\Delta_\omega$ is a regular unimodular triangulation.
 Then the limit scheme of the corresponding flat toric degeneration is a union of
 coordinate $n$-planes, one for every facet $\tau$ of $\Delta_\omega$, 
\[
  \lim_{t\to 0} t.x_\calA\ =\ \bigcup_{\tau\in\Delta_\omega}\P^\tau\,.
\]
\end{cor}

We give a different proof which is adapted from~\cite{SS} and does not assume
primitivity.
Its limitation is that it only works over $\C$ as it uses metric properties of $\C$ and not an
arbitrary algebraically closed field.\smallskip

\begin{figure}[htb]
  \[
   \begin{picture}(310,160)(-2,0)
    \put(0,0){\includegraphics[width=4in]{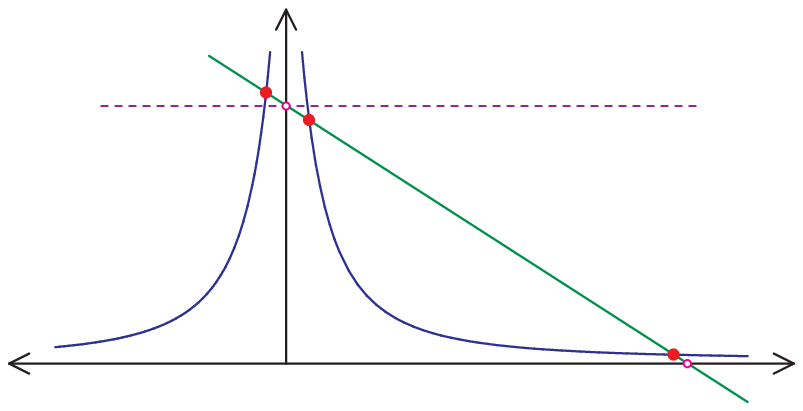}}
    \put(98,150){$\P^\tau$}
    \put(280,25){$\P^\sigma$}
    \put(168,75){\ForestGreen{$L$}}    \put(33,37){\Blue{$t.X_\calA$}}
    \put(215,118){\Violet{$\pi_\tau^{-1}(p_\tau)$}}
    \put(57,93){$p_\tau$}  \put(70,98.3){\vector(3,1){30}}
    \put(256,52){$p_\sigma$} \put(260,48){\vector(-1,-3){10}}
   \end{picture}
  \]
  \caption{Points in $L\cap t.X_\calA$ near $\P^\tau$ for small
  $t$.}\label{F:odd_even} 
 \end{figure}

The main idea is to fix a general linear subspace $L$ of codimension $n$
in $\P^\calA$  and consider the family of linear sections 
$L\cap t.X_\calA$ of fibers of the flat family $\mathcal{X}_\calA$, for $t$ near 0. 
This is illustrated in Figure~\ref{F:odd_even}.
%
The subspace $L$ will meet each facet $n$-plane $\P^\tau$ in a single point
$p_\tau$, and the points of $L\cap t.X_\calA$ for $t$ small will be clustered
near the different points $p_\tau$.

We could determine the number of points clustered near one of the $p_\tau$, which 
is the algebraic multiplicity, $\mbox{\rm mult}_{\P^\tau}(\ini_\omega(X_\calA))$.
It is in fact easier to determine the number of points in $\T^n$ of the form
$\varphi^{-1}_\calA(t^{-1}.x)$, for $x$ a point clustered near one of the $p_\tau$ and
when $t$ is small.
This is more direct, and it bypasses computing this algebraic multiplicity.
This is also where we avoid the assumption of primitivity, but must work over $\C$.

Note that in a neighborhood of $\P^\tau$, the linear space $L$ is isotopic to 
$\pi_\tau^{-1}(p_\tau)$, which is a fiber of the coordinate projection 
$\pi_\tau\colon \P^\calA\ - \to \P^\tau$.
This is a rational map not defined on the linear span of coordinates
$\{z_a\mid a\not\in\tau\}$.
It follows that the number of points in $\T^n$ coming from points in the linear
section $L\cap t.X_\calA$  near $p_\tau$ is 
equal to the number of points in $\T^n$ coming from points in the linear
section $\pi_\tau^{-1}(p_\tau)\cap t.X_\calA$  near $p_\tau$.
This is simply the degree of the map which is the composition of the parametrization 
$\varphi_\calA$ of $X_\calA$, the map $z\mapsto t.z$ on $\P^\calA$, and this
projection $\pi_\tau$,  
\[
  \T^n\ \xrightarrow{\ \varphi_\calA\ }\ X_\calA\ 
  \xrightarrow{\ t\ }\ t.X_\calA\ 
  \xrightarrow{\ \pi_\tau\ }\ \P^\tau\,.
\]
Since multiplication by $t$ is homotopic to the identity and it commutes with
the projection $\pi_\tau$, we may assume now
that $t=1$, and so this composition is nothing other than the 
parametrization $\varphi_\tau$ of $\P^\tau$ by the monomials
corresponding to integer points of $\tau$.
The degree of this map is the order of the kernel of $\varphi_\tau$, which
is $n!\cdot\vol(\Delta_\tau)$.
Summing this quantity over all facets $\tau$ of the triangulation shows that
there are
\[
   \sum_\tau n!\cdot\vol(\Delta_\tau)\ =\ n!\vol(\Delta_\calA)
\]
points in $\T^n$ which are pullbacks under $\varphi_\calA$ of the linear section
$L\cap{X}_\calA$.
This completes our algorithmic proof of Kushnirenko's Theorem.
\QED\medskip

The reason that it is algorithmic is that it (more-or-less) counts the solutions
to the system $L\cap t.X_\calA$, for $t$ small, while also giving enough
information on their location and structure to determine them numerically so
that they may become the input to the polyhedral homotopy method~\cite{HS95,VVC94}
for computing the solutions at $t=1$.

Note that this proof also shows that the intersection $X_\calA\cap L$ is transverse when
$L$ is general, and thus gives a proof of Bertini's Theorem in this context.
Since, for $t$ small, the intersection near $p_\tau$ may be deformed to the intersection
of $X_\calA$ with the horizontal subspaces $\pi_\tau^{-1}(p_\tau)$, and this is deformed
to the system $\varphi_\tau^{-1}(p_\tau)$, which consists of $n~\vol(\Delta_\tau)$
distinct points, the general such intersection is transverse.

%
\subsection{A brief aside about real solutions}
%

Observe that if the triangulation $\Delta_\omega$ is 
\DeCo{{\sl unimodular}}, in that each facet has minimal volume
$1/n!$, then near each point $p_\tau$ there will be exactly \Blue{one} point
of $L\cap t.X_\calA$ and one corresponding solution in $\T^n$.
If both $L$ and $t$ are real, then each $p_\tau$ and each nearby point in 
$t.X_\calA$ will be real.
Since
\[
   t^{-1}.(L\cap t.X_\calA)\ =\ \bigl(t^{-1}.L\bigr) \cap X_\calA\,,
\]
and the points in the left hand side are all real, so are the points in the right
hand side.
This right hand side corresponds to a system of real polynomials with support
$\calA$.
This proves a theorem of Sturmfels~\cite{St94}, and gives his argument in a
nutshell.\medskip 

\noindent{\bf Theorem~\ref{T1:Sturm}}\ 
{\em 
  If a lattice polytope $\Delta\subset\Z^n$ admits a regular unimodular
  triangulation, then there exist real polynomial systems with support
  $\Delta\cap\Z^n$ having all solutions real.
  }\medskip

A more careful analysis, which begins by examining real solutions when $|\calA|=n+1$,
leads to the more refined result for not necessarily
unimodular triangulations that appears in Sturmfels's paper.

%
%
%
%
\chapter{Upper Bounds}\label{Ch:Khov}
%

%
%

Recall Descartes rule of signs~(\cite{D1637} or Section~\ref{S2:Descartes}), which gives a
bound for the number of positive solutions to a univariate polynomial.\medskip

\noindent{\bf Theorem~\ref{T2:Descartes} (Descartes's rule of signs)} 
{\it 
 The number, $r$, of positive roots of a univariate polynomial
 \begin{equation}\label{E4:Univariate}
   f(x)\ =\ c_0 x^{a_0} + c_1 x^{a_1} + \dotsb + c_m x^{a_m}\,,
 \end{equation}
 counted with multiplicity, is bounded above by the number of variations in sign of the
 coefficients of $f$, 
 \[
   \#\{i\mid 1\leq i\leq m\mbox{ and }  c_{i-1} c_i<0\}
    \ \leq\ r\,,
\]
 and the difference between the variation and $r$ is even.\medskip
}

In~\eqref{E4:Univariate} we assume that $a_0<a_1<\dotsb<a_m$ and no
coefficient $c_i$ vanishes.  

Thus a univariate polynomial with $m+1$ monomials has at most 
$m$ positive roots.
This bound is sharp, as the polynomial
 \begin{equation}\label{E4:Descartes}
   (x-1)(x-2)\dotsb (x-m)
 \end{equation}
has $m+1$ distinct terms and $m$ positive roots.
Replacing $x$ by $x^2$ gives a polynomial with $m{+}1$ terms and $2m$ nonzero real roots.

This chapter and the next will discuss 
extensions of this Descartes bound to systems of multivariate polynomials.

%
\section{Khovanskii's fewnomial bound}
%

Descartes's rule of signs suggests that the number of real roots to a system of
polynomials depends not on its degree, but rather on the complexity of its
description.
D.~Bernstein and A.~Kushnirenko formulated the principle that the topological complexity
of a set in $\R^n$ defined by real polynomials is controlled by the complexity of the
description of the polynomials, rather than by their degree or Newton polytopes.
This is exactly what Khovanskii found in 1980 with his celebrated 
\DeCo{{\sl fewnomial bound}}.

\begin{thm}[Khovanskii~\cite{Kh80}]\label{T4:Fewnomial}
  A system of $n$ real polynomials in $n$ variables involving $l{+}n{+}1$ distinct
  monomials will have at most 
 \begin{equation}\label{E4:Fewnomial_Bound}
    2^{\binom{l+n}{2}}\cdot (n+1)^{l+n}
 \end{equation}
 nondegenerate positive solutions.
\end{thm}

We remark that nondegenerate solutions are isolated, and there are finitely
many of them.
This bound, like other bounds in this part of the subject, 
considers solutions in the \DeCo{{\sl positive orthant} $\R_>^n$}.
A consequence of Khovanskii's bound is that for each $l$ and $n$, there is a number 
\DeCo{$X(l,n)$} which is equal to the maximum number of positive solutions to a system of
$n$ polynomials in $n$ variables having $l{+}n{+}1$ distinct monomials.
A central question in this area is to determine the \DeCo{{\sl Khovanskii number}} 
$X(l,n)$ exactly, or give good bounds.
Khovanskii's Theorem shows that $X(l,n)$ is bounded above by the
quantity~\eqref{E4:Fewnomial_Bound}. 

A complete proof of Theorem~\ref{T4:Fewnomial} may be found in Khovanskii's
book~\cite{Kh91}, where much else is also developed.
Chapter 1 of that book contains an accessible sketch.
Benedetti and Risler~\cite[\S 4.1]{Benedetti_Risler} have a
careful and self-contained exposition of Khovanskii's fewnomial
bound.
We give a sketch of the main ideas in the exposition of
Benedetti and Risler, to which we refer for further details
(this is also faithful to Khovanskii's sketch).
We remark that Sturmfels has also sketched~(\cite[pp.~39--40]{SPE} and
in~\cite{St98}) a version of the proof.
This omits some contributions to the root count and is therefore regrettably incorrect.

Khovanskii looks for solutions in the positive orthant $\R_>^n$, proving a far
more general result involving solutions in $\R^n$ of polynomial functions in logarithms of 
the coordinates and monomials.
For this, he took logarithms of the coordinates.
Set
 \begin{equation}\label{E4:substitution}
   \DeCo{z_i}\ :=\ \log(x_i)
  \qquad\mbox{and}\qquad
   \DeCo{y_j}\ :=\ e^{z\cdot a_j}\ =\ x^{a_j}\,,
 \end{equation}
where $i=1,\dotsc,n$, $j=1,\dotsc,k$, and $a_j\in\R^n$ can be real
exponents. 
Consider a system of functions of the form
 \begin{equation}\label{Eq4:Khovanski_poly}
   F_i(z_1,\dotsc,z_n,\, y_1,\dotsc,y_k)\ =\ 0
   \qquad i=1,\dotsc,n\,,
 \end{equation}
where each $y_j=y_j(z)$ is an exponential function $e^{z\cdot a_j}$
and the $F_i$ are polynomials in $n{+}k$ indeterminates.

\begin{thm}[Khovanskii's Theorem]\label{T4:Kh_real_thm}
 The number of nondegenerate real solutions to the
 system~$\eqref{Eq4:Khovanski_poly}$ is at most
 \begin{equation}\label{Eq4:Real_fewnomial bound}
    \Bigl(\prod_{i=1}^n \deg F_i\Bigr) \cdot
    \Bigl( 1 + \sum_{i=1}^n \deg F_i \Bigr)^k\cdot
     2^{\binom{k}{2}} \,.
 \end{equation}
\end{thm}

\noindent{\it Proof of Theorem~$\ref{T4:Fewnomial}$.}
 Given a system of $n$ real polynomials in $n$ variables involving
 $l{+}n{+}1$ distinct monomials, we may assume that one of the monomials is 1.
 Under the substitution~\eqref{E4:substitution}, this becomes a
 system of the form~\eqref{Eq4:Khovanski_poly}, where each $F_i$ is
 a degree 1 polynomial in $k=l+n$ variables.
 Then $\deg F_i=1$ and the bound~\eqref{Eq4:Real_fewnomial bound} reduces
 to~\eqref{E4:Fewnomial_Bound}. 
\QED\medskip

\noindent{\it Sketch of proof of Theorem~$\ref{T4:Kh_real_thm}$.}
 We proceed by induction on $k$, skipping some technicalities involving Sard's Theorem.
 When $k=0$, there are no exponentials, and the system is just a
 system of $n$ polynomials in $n$ variables, whose number of
 nondegenerate isolated solutions is bounded above by the B\'ezout
 number,
\[
   \prod_{i=1}^n \deg F_i\,,
\]
 which is the bound~\eqref{Eq4:Real_fewnomial bound} when $k=0$.

 Suppose that we have the bound~\eqref{Eq4:Real_fewnomial bound} for
 systems of the form~\eqref{Eq4:Khovanski_poly} with $k$ exponentials, and
 consider a system with $k{+}1$ exponentials and one added variable \DeCo{$t$}. 
 \begin{eqnarray}
  \label{Eq4:subsystem}
   G_i(z,t)\ :=\ F_i(z_1,\dotsc,z_n,\; y_1,\dotsc,y_k,\,t\cdot y_{k+1})&=& 0\
   \qquad i=1,\dotsc,n\\
   \nonumber
     t&=& 1
 \end{eqnarray}
 The subsystem~\eqref{Eq4:subsystem} defines an analytic curve
 $C$ in $\R^{n+1}$, which we assume is smooth and transverse to
 the hyperplane at $t=1$.

 Write $z_{n+1}$ for $t$ and consider the vector field $\xi$ in
 $\R^{n+1}$ whose $r$th component is
 \begin{equation}\label{Eq4:zeta_r}
  \xi_r\ :=\ (-1)^{n+1-r}\det\left(\frac{\partial G_i}{\partial
  z_j}\right)^{i=1,\dotsc,n}_{j=1,\dotsc,\hat{r},\dotsc,n+1}\ .
 \end{equation}
 This vector field is tangent to the curve $C$, and we write
 $\xi_t=\xi_{n+1}$ for its component in the $t$-direction.
 An important ingredient in our proof of Theorem~\ref{T4:Kh_real_thm}
 is a special case of the Khovanskii-Rolle Theorem~\cite[pp.~42--51]{Kh91}.

\begin{thm}[Khovanskii-Rolle Theorem]\label{T4:Kh-R}
 The number of points of\/ $C$ where $t=1$ is bounded 
 above by
\[
   N\ +\ q\,,
\]
 where $N$ is the number of points of\/ $C$ where $\xi_t=0$ and 
 $q$ is the number of unbounded components of\/ $C$.
\end{thm}
   
\noindent{\it Proof.}
 Note that $\xi_t$ varies continuously along $C$.
 Suppose that $a$ and $b$ are consecutive points along an arc of
 $C$ where $t=1$.
 Since $C$ is transverse to the hyperplane $t=1$, we have
 $\xi_t(a)\cdot\xi_t(b)<0$, and so there is a point $c$ of $C$
 on the arc between $a$ and $b$ with $\xi_t(c)=0$.

 The hyperplane $t=1$ cuts a compact connected component of $C$
 into the same number of arcs as points where $t=1$.
 Since the endpoints of each arc lie on the hyperplane $t=1$,
 there is at least one point $c$ on each arc with $\xi_t(c)=0$.
 Similarly, the hyperplane $t=1$ cuts a noncompact component into
 arcs, and each arc with two endpoints in the  hyperplane $t=1$ gives a
 point $c$ with $\xi_t(c)=0$.  
 However, there will be one more point with $t=1$ on this component
 than such arcs.
\QED\medskip

We illustrate the argument in the proof below.
\[
  \begin{picture}(220,160)(0,-7)
   \put(0,0){\includegraphics[height=150pt]{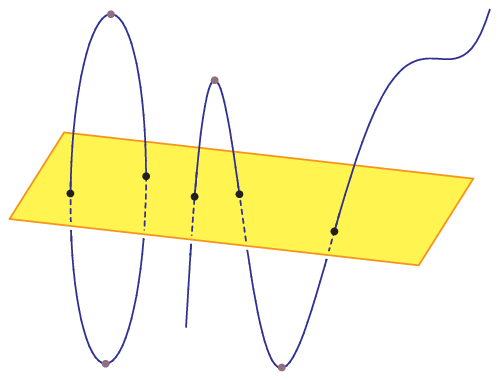}}
   \put(197,76){$t=1$}
   \put(72,144){$\xi_t=0$}
   \put(68,147){\vector(-1,0){20}} \put(85.8,140){\vector(0,-1){17}}
   \put(62,-8){$\xi_t=0$}
   \put(59,-2.5){\vector(-4,1){15}}   \put(95,-3.5){\vector(4,1){15}}
   \put(198,130){$C$}
  \end{picture}
\]

 The key to the induction in the proof of Khovanskii's formula is to replace the last
 exponential by a new variable. 
 This substitution is omitted in Sturmfels's argument, which also does not use
 the Khovanskii-Rolle Theorem~\ref{T4:Kh-R}.
 Since we have 
 \begin{eqnarray*}
   \frac{\partial G_i}{\partial z_r}
   &=&
  \frac{\partial F_i}{\partial z_r}
     (z_1,\dotsc,z_n,\; y_1,\dotsc,y_k,t y_{k+1})\\
   &&\ +\ \sum_{j=1}^k \frac{\partial F_i}{\partial y_j}
     (z_1,\dotsc,z_n,\;y_1,\dotsc,y_k,t\, y_{k+1})\cdot
      a_{j,r}\, y_j\\
   &&\ +\ \frac{\partial F_i}{\partial y_{k+1}}
     (z_1,\dotsc,z_n,\;y_1,\dotsc,y_k,t\, y_{k+1})\cdot
      a_{k+1,r}\  t\, y_{k+1}\,,
 \end{eqnarray*}
if we set $\DeCo{u}:= t y_{k+1}$ and define \DeCo{$\phi(z,u)$} to be the expression for
$\xi_t=\xi_{n+1}$~\eqref{Eq4:zeta_r} written in terms of $z$ and $u$, then 
the total degree (in $z_1,\dotsc,z_n,y_1,\dotsc,y_k,\DeCo{u}$) of 
$\phi_t(z,u)$ is at most $\sum_{i=1}^n \deg F_i$.

Thus number $N$ of Theorem~\ref{T4:Kh-R} is the number of
solutions to the system
 \begin{equation}\label{Eq4:Sys_A}
 \begin{array}{rcl}
   F_i(z_1,\dotsc,z_n,\; y_1,\dotsc,y_k,\,\DeCo{u})&=& 0\
   \qquad i=1,\dotsc,n\\
    \phi_t(z,\DeCo{u})&=& 0\,.\rule{0pt}{15pt}
 \end{array}
 \end{equation}
This has the form~\eqref{Eq4:Khovanski_poly} with $k$ exponentials.
Given any solution to the system~\eqref{Eq4:Sys_A}, we use the substitution $u=t
y_{k+1}$ to solve for $t$ and get a corresponding point $c$ on 
the curve $C$ with $\xi_t(c)=0$.
We apply our induction hypothesis to the system~\eqref{Eq4:Sys_A} (which has $k$ exponentials
and $n{+}1$ equations in $n{+}1$ variables) to obtain
\[
    N\ \leq \ 
    \prod_{i=1}^n \deg F_i \cdot
    \Bigl( \sum_{i=1}^n \deg F_i \Bigr)\cdot
    \Bigl( 1 + 2\sum_{i=1}^n \deg F_i \Bigr)^k\cdot
     2^{\binom{k}{2}} \,.
\]

We similarly estimate the number $q$ of noncompact components of
$C$.
We claim that this is bounded above by the maximum number of points of
intersection of $C$ with a hyperplane.
Indeed, since each noncompact component has two infinite branches,
there are  
$2q$ points (counted with multiplicity) on the sphere $S^{n}$
corresponding to directions of accumulation points of the branches of $C$ 
at infinity.
Any hyperplane through the origin not meeting these points will have
at least $q$ of these points in one of the hemispheres into which it
divides the sphere.
If we translate this hyperplane sufficiently far toward infinity, it
will meet the branches giving these accumulation points, and thus
will meet $C$ in at least $q$ points.

Thus $q$ is bounded by the number of solutions to a system of the form
 \begin{equation}\label{Eq4:Sys_B}
  \begin{array}{rcl}
    F_i(z_1,\dotsc,z_n,\; y_1,\dotsc,y_k,\,\DeCo{u})&=& 0\
    \qquad i=1,\dotsc,n\\
    l_0+ l_1z_1 + l_2z_2 +\dotsb+l_nz_n + l_u \DeCo{u} &=& 0 \ ,\rule{0pt}{15pt}
  \end{array}
 \end{equation}
where $l_0,\dotsc,l_n,\,l_u$ are some real numbers.
This again involves only $k$ exponentials, and the last equation has degree 1, so we have
\[ 
   q\ \leq\ 
    \prod_{i=1}^n \deg F_i \cdot
    1\cdot
    \Bigl(1 + \sum_{i=1}^n \deg F_i \ +1\Bigr)^k\cdot
     2^{\binom{k}{2}} \,.
\]
Combining these estimates gives
\[
   N+q\ \leq\   
    \prod_{i=1}^n \deg F_i \cdot    2^{\binom{k}{2}} \cdot 
   \Bigl[
      \Bigl( \sum_{i=1}^n \deg F_i \Bigr)\cdot
     \Bigl( 1 + 2\sum_{i=1}^n \deg F_i \Bigr)^k
    \ +\ 
    \Bigl( 2 + \sum_{i=1}^n \deg F_i \Bigr)^k
    \Bigr]  \,.
\]

We can obtain a simpler (but larger) estimate by bounding $N+q$  by the number of
solutions to the single system of equations,
 \begin{eqnarray*}
   F_i(z_1,\dotsc,z_n,\; y_1,\dotsc,y_k,\,\DeCo{u})&=& 0\
   \qquad i=1,\dotsc,n\\
   \DeCo{F_{n+1}}\ :=\  (l_0+ l_1z_1 + l_2z_2 +\dotsb+l_nz_n + l_u \DeCo{u})
        \cdot\phi_t(z,\DeCo{u}) &=& 0 \ .
 \end{eqnarray*}
By our induction hypothesis, we have 
\[
   N+q\ \leq\ 
    \prod_{i=1}^{n+1}\deg F_i \cdot
    \Bigl( 1 + \sum_{i=1}^{n+1} \deg F_i \Bigr)^k\cdot
     2^{\binom{k}{2}} \,.
\]
But we saw that $F_{n+1}$ has degree at most $\sum_{i=1}^n \deg F_i\ +1$, and so 
the number, $M$, of solutions to the system with $k+1$
exponentials is bounded by 
 \begin{eqnarray*}
  M \ \leq \ N+q&\leq& \prod_{i=1}^n \deg F_i \cdot
    \Bigl(1+ \sum_{i=1}^n \deg F_i \Bigr)\cdot
    \Bigl(1 + \sum_{i=1}^n \deg F_i +1+\sum_{i=1}^n \deg F_i\Bigr)^k \cdot 2^{\binom{k}{2}} \\
   &=& \prod_{i=1}^n \deg F_i \cdot
        \Bigl(1 + \sum_{i=1}^n \deg F_i \Bigr)^{k+1} \cdot 2^{\binom{k+1}{2}} \,,
 \end{eqnarray*}
which completes the proof of Theorem~\ref{T4:Kh_real_thm}.
\QED\medskip

We see that the result of Theorem~\ref{T4:Kh_real_thm} is much more
general than the statement of Theorem~\ref{T4:Fewnomial}.
Also, the bound is not sharp.
While no one believed that Khovanskii's bound~\eqref{E4:Fewnomial_Bound} was anywhere near
the actual upper bound $X(l,n)$, it was been extremely hard to improve it.
We discuss the first steps in this direction.

%
\section{Kushnirenko's conjecture}
%

One of the first proposals of a more reasonable bound than Khovanskii's for the number of
positive solutions to a system of polynomials was due to 
Kushnirenko,  and for many years experts believed that this may
indeed be the truth.

\begin{conj}[Kushnirenko]\label{C4:Kushnirenko}
 A system  $f_1=f_2=\dotsb=f_n=0$ of real polynomials where each $f_i$ has
 $m_i{+}1$ terms has at most $m_1m_2\dotsb m_n$ nondegenerate positive
 solutions. 
\end{conj}

This generalizes the bound given by Descartes's rule of signs.
It easy to use the example~\eqref{E4:Descartes} for the sharpness of Descartes's
rule to construct systems of the form
\[
  f_1(x_1)\ =\ f_2(x_2)\ =\ \dotsb\ =\ f_n(x_n)\ =\ 0\,,
\]
which achieve the bound of Conjecture~\ref{C4:Kushnirenko}.

Soon after Kushnirenko made this conjecture, K.~Sevostyanov found a counterexample which was
unfortunately lost.
Nevertheless, this conjecture passed into folklore until Haas~\cite{Ha02}
found a simple example of two trinomials ($3=2+1)$ in variables $x$ and $y$ with  
$5\,(>4=2\cdot 2)$ isolated nondegenerate positive solutions.
 \begin{equation}\label{E4:Haas}
   10 x^{106}\ +\ 11 y^{53}\ -\ 11y\ =\ 
   10 y^{106}\ +\ 11 x^{53}\ -\ 11x\ =\ 0\,.
 \end{equation}

There have been other attempts to find better bounds than the
Khovanskii bound. 
Sturmfels~\cite{St94} used the toric degenerations introduced in Chapter~\ref{Ch:sparse}
to show how to construct systems with many real roots (the root count
depends upon a mixture the geometry of the Newton polytopes and some
combinatorics of signs associated to lattice points)\fauxfootnote{Actually, he
  used the toric proof of Bernstein's Theorem.}.
This inspired Itenberg and Roy~\cite{IR96} to propose a multivariate 
version of Descartes's rule of signs, which was later found to be too
optimistic~\cite{LW98}. 
An interesting part of this story is told in the cheeky paper of Lagarias and
Richardson~\cite{LR97}. 

More recently, Li, Rojas, and Wang looked closely at Haas's
counterexample to Kushnirenko's conjecture, seeking to obtain 
realistic bounds for the number of positive solutions which depended only on the number of
monomials in the different polynomials.
For example, they showed that Haas's counterexample was the best
possible.

\begin{thm}[Li, Rojas, and Wang~\cite{LRW03}]
  A system consisting of two trinomials in two variables has at 
  most $5$ nondegenerate positive solutions.
\end{thm} 

Dickenstein, Rojas, Rusek,  and Shih~\cite{DRRS} used exact
formulas for $\calA$-dis\-criminants~\cite{DFS} to study systems of two trinomials in two 
variables which achieve this bound of five positive solutions.
They gave the following example, which indicates how difficult it is to find systems with 
many real solutions.

\begin{ex}
 Consider the family of systems of bivariate sextics,
 \begin{equation}\label{Eq:DRRS_system}
    x^6\ +\ a y^3\ -\ y\ =\ y^6\ +\ bx^3\ -\ x\ =\ 0\,,
 \end{equation}
 where $a,b$ are real numbers.
 When $a=b=78/55$, this has five positive real solutions
\[
   (0.8136, 0.6065)\,,\ 
   (0.7888, 0.6727)\,,\ 
   (0.7404, 0.7404)\,,\ 
   (0.6727, 0.7888)\,,\ 
   (0.6065, 0.8136)\,.
\]
 We now investigate the set of parameters $(a,b)$ for which this achieves the trinomial
 bound of five positive solutions.
 This turns out to be a single connected component in the complement of the discriminant
 for this family of systems.
 This discriminant is a polynomial of degree 90 in $a,b$ with 58 terms whose leading and
 trailing terms are 
 \begin{multline*}
  1816274895843482705708030487016037960921088a^{45}b^{45} \ +\ 
   \dotsb\ \mbox{56 terms}\ \dotsb\ +\ \\
   1102507499354148695951786433413508348166942596435546875.
 \end{multline*}
 We display this discriminant in the square $[0,4]\times[0,4]$, as well as three successive
 magnifications, each by a factor of 11.
 The shaded region in the last picture is the set of pairs $(a,b)$ for
 which~\eqref{Eq:DRRS_system} achieves the trinomial bound of five positive real
 solutions.
\[
  \includegraphics[height=110pt]{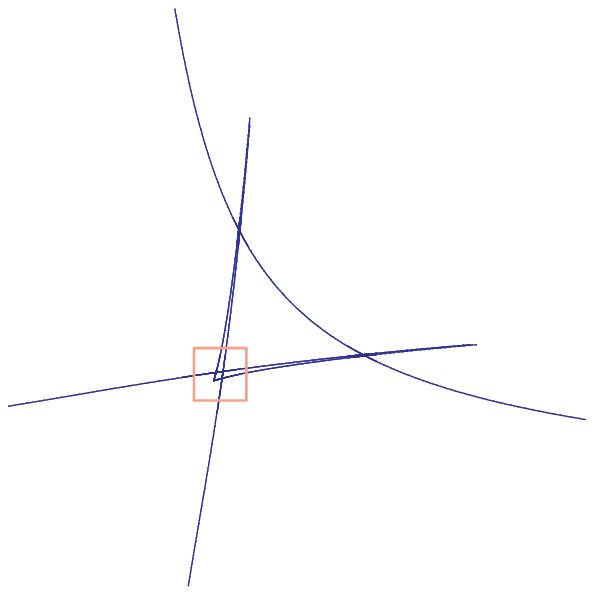}
  \includegraphics[height=110pt]{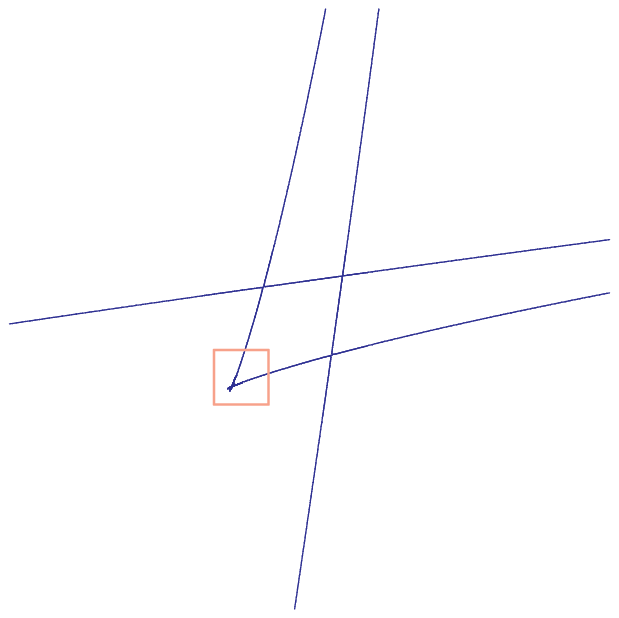}
  \includegraphics[height=110pt]{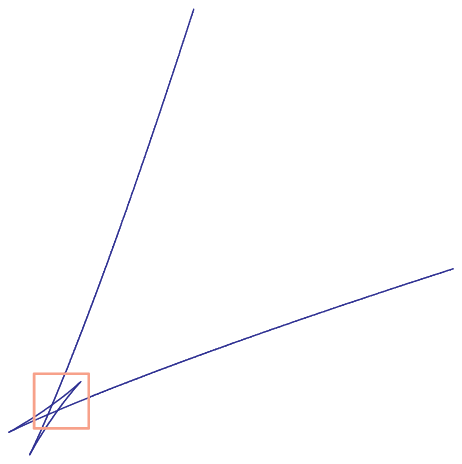}
  \includegraphics[height=110pt]{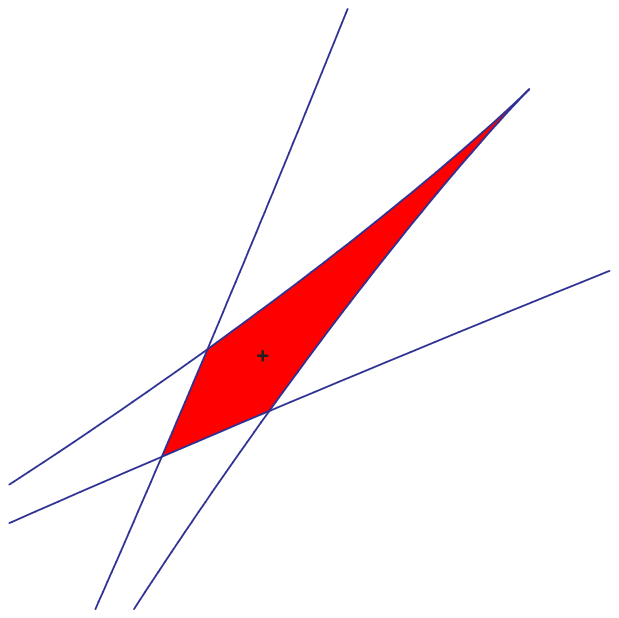}
\]
\end{ex}

To compare the trinomial bound in~\cite{LRW03} to the fewnomial
bound~\eqref{E4:Fewnomial_Bound}, note that we may multiply one of the polynomials by a
monomial so that the two trinomials share a monomial. 
Then there are at most $5=2+2+1$ distinct monomials occurring in the two trinomials.
The fewnomial bound for $l=n=2$ is 
\[
   X(2,2)\ \leq \  2^{\binom{4}{2}}\cdot (2+1)^4\ =\ 5184\,.
\]
We remark that a trinomial system is not quite a general
fewnomial system with $l=n=2$.
Still, the bound of five real solutions lent credence to the belief $X(2,2)$ is closer to
five than to 5184 and that Khovanskii's fewnomial bound~\eqref{E4:Fewnomial_Bound} could be
improved. 

In addition to providing the counterexample to Kushnirenko's conjecture, Sevostyanov
also established the first result of fewnomial-type.
He showed the existence of an absolute bound $\sigma(d,m)$ for the number of real
solutions  to a system
\[
   f(x,y)\ =\ g(x,y)\ =\ 0\,,
\]
where $f$ is a polynomial of degree $d$ and $g$ has $m$ terms.
The proof of this result, like his counterexample, has unfortunately been lost.
This result however, was the inspiration for Khovanskii to develop his theory of
fewnomials.

Recently, Avenda\~{n}o~\cite{Av}, established a precise version of a special case of
Sevostyanov's theorem.

\begin{thm}
 Suppose that $f(x,y)$ is linear and $g(x,y)$ has $m$ terms.
 Then the system
\[
   f(x,y)\ =\ g(x,y) =\ 0\,,
\]
 has at most $6m-4$ real solutions.
\end{thm}

%
\section{Systems supported on a circuit}\label{S4:circuit}
%
Restricting the analysis of Section~\ref{S3:simplex} to real solutions shows that 
$X(0,n)=1$.
Recently, it was shown that $X(1,n)=1{+}n$.
We discuss this here.

A collection $\calA$ of $n{+}2$ vectors in $\Z^n$ which affinely spans $\R^n$ is
called a \DeCo{{\sl circuit}}.
The circuit is \DeCo{{\sl primitive}} if its $\Z$-affine span is all of $\Z^n$.
When $0\in\calA$, this means that $\Z^n=\Z\calA$.

\begin{thm}[Bertrand, Bihan, and Sottile~\cite{BBS}]\label{T4:BBS}
  A polynomial system supported on a primitive circuit has at most $2n+1$
  nondegenerate nonzero real solutions.
\end{thm}

\begin{thm}[Bihan~\cite{Bihan}]\label{T4:Bihan}
  A polynomial system supported on a circuit has at most $n{+}1$
  nondegenerate positive solutions, and there exist systems supported on a circuit having
  $n{+}1$ positive solutions.
\end{thm}

This can be used to construct fewnomial systems with relatively many positive
solutions.

\begin{cor}[Bihan, Rojas, Sottile~\cite{BRS}]\label{C4:BRS}
  There exist systems of $n$ polynomials in $n$ variables having $l{+}n{+}1$ monomials and 
  at least 
  $\lceil\frac{n}{l}\rceil^l$ positive solutions.
\end{cor}

This gives a lower bound for $X(l,n)$ of $\lceil\frac{n}{l}\rceil^l$, and is the best
construction when $l$ is fixed and $n$ is large.
It remains an open problem to give constructions with more solutions, or 
constructions with many solutions when $l$ is not fixed.

The construction establishing Corollary~\ref{C4:BRS} is quite simple.
Suppose that $n=kl$ is a multiple of $l$, and let 
\[
   f_1(x_1,\dotsc,x_k)\ =\ f_2(x_1,\dotsc,x_k)\ =\ \dotsb\ =\ 
    f_k(x_1,\dotsc,x_k)\ =\ 0\,,
\]
be a system with $k{+}2$ monomials and $k{+}1$ positive solutions.
Such systems exist, by Theorem~\ref{T4:Bihan}.
Write $F(x)=0$ for this system and assume that one of its monomials is a constant.
For each $i=1,\dotsc,l$, let $y^{(i)}=(y_1^{(i)},\dotsc,y^{(i)}_k)$ be a set of $k$ variables.
Then the system
\[
   F(y^{(1)})\ =\ F(y^{(2)})\ =\ \dotsb\ =\ F(y^{(l)})\ =\ 0\,,
\]
has $(k+1)^l$ solutions, $kl$ variables, and $1+l+kl$ monomials.
When $n=kl+r$ with $r<;$, adding extra variables $y_i$ and equations 
$y_i=1$ for $i=1,\dotsc,r$ gives a system with $(\lfloor\frac{n}{l}\rfloor +1)^l$ 
positive solutions and $1+l+n$ monomials.\medskip

When $l=1$ the fewnomial bound~\eqref{E4:Fewnomial_Bound} becomes
\[
   2^{\binom{1+n}{2}} \cdot (n+1)^{1+n}\,,
\]
which is considerably larger than Bihan's bound of $n{+}1$.
Replacing $l{+}n$ by $n$ in the fewnomial bound, it becomes equal to Bihan's bound when
$l=1$.
When $l=0$, this  same substitution in~\eqref{E4:Fewnomial_Bound} yields $1$, which is the
sharp bound when $l=0$.
In Chapter~\ref{Ch:Gale}, we will outline generalizations of Theorems~\ref{T4:BBS}
and~\ref{T4:Bihan} to arbitrary $l$, giving the bound 
 \begin{equation}\label{Eq4:newfewnomial}
   X(l,n)\ <\ \frac{e^2+3}{4}2^{\binom{l}{2}}n^l
 \end{equation}
for positive solutions and, when $\calA$ is primitive, a bound for all real solutions,
\[
     \frac{e^{\Red{4}}+3}{4}2^{\binom{l}{2}}n^l\,.
\]
This is only slightly larger---the difference is in the exponents $2$ and $4$ of $e$.
These are proven in~\cite{BaBiSo,BS07}.

By Corollary~\ref{C4:BRS} and the bound~\eqref{Eq4:newfewnomial},
\[
   l^{-l}n^l\ <\ \left\lceil\frac{n}{l}\right\rceil^l\ 
   <\ X(n,l)\ <\ \frac{e^2+3}{4}2^{\binom{l}{2}}n^l\,.
\]
This reveals the correct asymptotic information for $X(n,l)$, when $l$ is fixed,
$X(n,l)=\Theta(n^l)$.

Theorems~\ref{T4:BBS} and~\ref{T4:Bihan} are related, and we will outline 
their proofs, following the papers in which they appear, where more details may be  
found. 
To begin, let 
 \begin{equation}\label{E4:system}
  f_1(x_1,x_2,\dotsc,x_n)\ =\ 
  f_2(x_1,x_2,\dotsc,x_n)\ =\ \dotsb\ =\ 
  f_n(x_1,x_2,\dotsc,x_n)\ =\ 0
 \end{equation}
be a system with support a circuit $\calA$.
Suppose that $0\in\calA$ and list the elements of the circuit 
$\calA=\{0,a_0,a_1,\dotsc,a_n\}$.  
After a multiplicative change of coordinates (if necessary), we may assume that
$a_0=\ell \e_n$, where $\e_n$ is the $n$th standard basis vector.
Since the system~\eqref{E4:system} is generic, row operations on the equations
put it into diagonal form
 \begin{equation}\label{E4:Diagonal_circuit}
   x^{a_i}\ =\ w_i\ +\ v_i x_n^\ell
   \qquad\mbox{for }i=1,\dotsc,n\,.
 \end{equation}

When $\calA$ was a simplex we used integer linear algebra to reduce the
equations to a very simple system in Section~\ref{S3:simplex}.
We use (different)  integer linear algebra to simplify this system
supported on a circuit.

%
\subsection{Some arithmetic for circuits}
%
Suppose that $\{0, a_0, a_1,\dotsc,a_n\}\subset\Z^n$ is a primitive circuit. 
We assume here that it is nondegenerate---there is no affine dependency
involving a subset.  
(The bounds in the degenerate case are lower, replacing $n$ by the size of this
smaller circuit.)
After possibly making a coordinate change, we may assume that 
$a_0=\ell\cdot \e_n$, where $\e_n$ is the $n$th standard basis vector.

For each $i$, we may write $a_i=b_i+k_i\cdot \e_n$, where $b_i\in\Z^{n-1}$.
Removing common factors from a nontrivial integer linear relation among 
the $n$ vectors $\{b_1,\dotsc,b_n\}\subset\Z^{n-1}$ gives us the
\DeCo{{\sl primitive relation}} among them (which is well-defined up to
multiplication by $-1$),
\[
   \sum_{i=1}^p \lambda_i b_i\ =\ 
   \sum_{i=p+1}^n \lambda_i b_i\,.
\]
Here, each $\lambda_i>0$, and we assume that the vectors are ordered
so that the relation has this form.
We further assume that
\[
   N\ :=\ \sum_{i=p+1}^n \lambda_i k_i\ -\ 
          \sum_{i=1}^p \lambda_i k_i\ >\ 0\,.
\]
Then we have
 \[
   N \e_n\ +\ \sum_{i=1}^p \lambda_i a_i\ -\ 
             \sum_{i=p+1}^n \lambda_i a_i\ =\ 0\,,
 \]
and so
 \begin{equation}\label{E4:primitiverelation}
   x_n^N \cdot \prod_{i=1}^p (x^{a_i})^{\lambda_i} \ -\ 
             \prod_{i=p+1}^n (x^{a_i})^{\lambda_i}\ =\ 0\,.
 \end{equation}
%

%
\subsection{Elimination for circuits}
%

 Using~\eqref{E4:Diagonal_circuit} to 
 substitute for $x^{a_i}$ in~\eqref{E4:primitiverelation} 
 gives the univariate consequence
 of~\eqref{E4:Diagonal_circuit}
 \begin{equation}\label{Eq4:eliminant}
  \DeCo{f}(x_n)\ :=\ x_n^N \prod_{i=1}^p (w_i+v_i x_n^\ell)^{\lambda_i}
        \ -\ \prod_{i=p+1}^n (w_i+v_i x_n^\ell)^{\lambda_i}\ .
 \end{equation}
 Some further arithmetic of circuits (which may be found in~\cite{BBS}) shows 
 that $f(x_n)$ has degree equal to $n!\vol(\Delta_\calA)$.
 This is in fact the eliminant of the system.

\begin{lemma}\label{L4:eliminant}
 The association of a solution $x$ of~\eqref{E4:Diagonal_circuit} to its $n$th
 coordinate $x_n$ gives a bijection between the solutions
 of~\eqref{E4:Diagonal_circuit} and the roots of $f~\eqref{Eq4:eliminant}$ which restricts to a
 bijection between their real solutions/roots. 
\end{lemma}

 While $f$ is the eliminant of the system, we do not have a Gr\"obner basis or
 even a triangular system to witness this fact, and the proof proceeds by
 explicitly constructing a solution~\eqref{E4:Diagonal_circuit} from a root 
 $x_n$ of $f$.

 The upper bound is found by writing $f=F-G$ as in~\eqref{Eq4:eliminant} and then
 perturbing $f$, 
\[
   f_t(y)\ =\ t\cdot F(y)\ -\ G(y)\,.
\]
 We simply estimate the number of changes in the the real roots of $f_t$ as
 $t$ passes from $-\infty$ to $1$, which can occur only at the singular roots of
 $f_t$. 
 While similar to the proof of Khovanskii's theorem, this is not inductive, but
 relies on the form of the Wronskian $F'G-G'F$ whose roots are the singular
 roots of $f_t$.
 This may also be seen as an application of Rolle's Theorem.
 We note that this estimation also uses Viro's construction
 for $t$ near $0$ and $\infty$.

 These estimates prove the bounds in Theorems~\ref{T4:BBS} and~\ref{T4:Bihan}.
 Sharpness comes from construction.
 In~\cite{BBS} Viro's method for univariate polynomials is used to 
 construct polynomials $f$ are constructed with  $2n{+}1$ real solutions, for special primitive
 circuits. 
 Bihan~\cite{Bihan} constructs a system with $n{+}1$ positive solutions using
 Grothendieck's {\it dessins d'enfants}.

%
\subsection{A family of systems with a sharp bound}\label{S4:exact}
%

We give a family of systems that illustrate the result of Theorem~\ref{T4:BBS}
(actually of an extension of it) and which may be treated by hand.
These systems come from a family of polytopes $\Delta\subset\Z^n$ for which we
prove a nontrivial upper bound on the number of real solutions to polynomial
systems with primitive support $\calA:=\Delta\cap\Z^n$. 
That is, the integer points $\calA$ in $\Delta$ affinely span $\Z^n$,
so that general systems supported on $\Delta$ have $n!\vol(\Delta)$ complex solutions, 
but there are fewer than $n!\vol(\Delta)$ real solutions to polynomial systems with
support $\calA$.
This is intended to not only give a glimpse of the more general results
in~\cite{BBS}, but also possible extensions which are not treated in~\cite{BS07}.

Let $l > k>0$ and $n \ge 3$ be integers and 
$\DeCo{\epsilon} = (\epsilon_1,\dotsc,\epsilon_{n-1})\in\{0,1\}^{n-1}$ have at least one
nonzero coordinate. 
The polytope $\DeCo{\Delta_{k,l}^\epsilon} \subset \R^n$ is the convex hull of the points 
 \[
   (0,\dots,0),\ (1,0,\dots,0),\dots, (0,\dots,0,1,0),\ (0,\dots,0,k),\
   (\epsilon_1,\dots, \epsilon_{n-1},l)\,.
 \]
 The configuration $\DeCo{\calA_{k,l}^\epsilon}=\Delta_{k,l}^\epsilon\cap\Z^n$ 
 also includes the points along the last axis
\[
  (0,\dotsc,0,1),\ (0,\dotsc,0,2),\ \dotsc,\ (0,\dotsc,0,k{-}1)\,.
\]
These points include the standard basis and the origin, so $\calA_{k,l}^\epsilon$
is primitive in that $\Z\calA=\Z^n$.

Set $|\epsilon|:=\sum_i\epsilon_i$.  Then the volume of
$\Delta_{k,l}^\epsilon$ is $(l+k|\epsilon|)/n!$. Indeed,
the configuration $\calA_{k,l}^\epsilon$ can be triangulated into
two simplices $\Delta_{k,l}^\epsilon \setminus \{(\epsilon_1,\dots, \epsilon_{n-1},l)\}$
and $\Delta_{k,l}^\epsilon \setminus \{0\}$ with volumes $k/n!$ and $(l-k+k|\epsilon|)/n!$, 
respectively. One way to see this is to apply the affine transformation
\[ 
  (x_1,\dotsc,x_n)\ \longmapsto\ 
  (x_1,\dotsc,x_{n-1},x_n-k+k\sum_{i=1}^{n-1}x_i)\,.
\]

\begin{thm}\label{T4:aBound}
 The number, $r$, of real solutions to a generic system of $n$ real
 polynomials with support $\calA_{k,l}^\epsilon$ lies in the interval
\[
    0 \ \leq\ r\ \leq \ k + k|\epsilon|+2\,,
\]
 and every number in this interval with the same parity as $l+k|\epsilon|$
 occurs. 
\end{thm}

This upper bound does not depend on $l$ and, since
$k < l$, it is smaller than or equal to the number $l+k|\epsilon|$
of complex solutions.
We use elimination to prove this result.

\begin{ex}\label{E:km-polytope}
Suppose that $n=k=3$, $l=5$, and $\epsilon = (1,1)$.  
\[
  \begin{array}{rclcl}
  \makebox[210pt][l]{Then the system}&&&&
   \multirow{9}*{
    \begin{picture}(80,105)(0,15)
        \put( 4,0){\includegraphics[height=4.5cm]{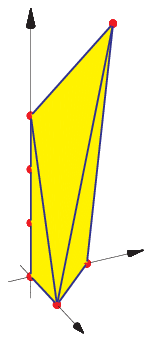}}
        \put(-1,24){$1$}   \put(16, 3){$x$}   \put(39,20){$y$}
        \put(-1,42){$z$}   \put(-1,63){$z^2$} \put(-1,85){$z^3$} 
        \put(50,120){$xyz^5$} 
        \put(55,60){$\Delta_{3,5}^{(1,1)}$}
    \end{picture}}
  \\
    x + y + xyz^5 +1+z+z^2+z^3&=&0&\quad&\rule{0pt}{16pt}\\
    x + 2y + 3xyz^5 +5+7z+11z^2+13z^3&=&0&\quad&\\
   2x + 2y + xyz^5 +4+8z+16z^2+32z^3&=&0&\quad&\\
     \makebox[210pt][l]{is equivalent to }&&\rule{0pt}{16pt}\\
    x    -(5 +11z +23z^2 +41z^3)&=&0\rule{0pt}{16pt}\\
     y    +( 8 +18z +38z^2 +72z^3)&=&0\\
    xyz^5 -(2 +6z  +14 z^2 +30 z^3)&=&0\\
    \mbox{\ }
   \end{array}
\]
And thus its number of real roots equals the number of real roots of
\[
  z^5(5 +11z +23z^2 +41z^3)(8 +18z +38z^2 +72z^3) -
  (2 +6z  +14 z^2 +30 z^3)\,,
\]
%
%
which, as we invite the reader to check, is 3.  \QED
\end{ex}

\noindent{\it Proof of Theorem~$\ref{T4:aBound}$.}
 A generic real polynomial system with support $\calA_{k,l}^\epsilon$ has the
 form 
 \[
    \sum_{j=1}^{n-1} c_{ij}x_j \ + c_{in}x^\epsilon x_n^l\ \ 
    +\ f_i(x_n)\ =\ 0\ \quad 
    {\rm for\  }i=1,\dots,n\,,
 \]
 where each polynomial $f_i$ has degree $k$ and 
 $x^\epsilon$ is the monomial 
 $x_1^{\epsilon_1}\dotsb x_{n-1}^{\epsilon_{n-1}}$.

 Since all solutions to our system are simple, we may perturb the
 coefficient matrix $(c_{ij})_{i,j=1}^n$ if necessary and then use
 Gaussian elimination to obtain an equivalent system
 \begin{equation}\label{system}
  x_1-g_1(x_n)\ =\ \dotsb \ = \ x_{n-1}-g_{n-1}(x_n) \ =\ 
  x^\epsilon x_n^l -g_{n}(x_n)\ =\ 0\,,
 \end{equation}
 where each polynomial $g_i$ has degree $k$.
 Using the first $n-1$ polynomials to eliminate the
 variables $x_1, \dots , x_{n-1}$ gives the univariate polynomial
 \begin{equation}\label{E:eliminantex}
   x_n^l \cdot g_1(x_n)^{\epsilon_1} \dotsb 
    g_{n-1}(x_n)^{\epsilon_{n-1}}\ -\ g_n(x_n)\,,
 \end{equation}
 which has degree $l+k|\epsilon|=v(\Delta_{k,l}^\epsilon)$.  Any zero
 of this polynomial leads to a solution of the original
 system~\eqref{system} by back substitution.  This implies that the
 number of real roots of the polynomial~\eqref{E:eliminantex} is equal
 to the number of real solutions to our original
 system~\eqref{system}.

The eliminant~\eqref{E:eliminantex} has no terms of degree $m$ for $k<m<l$, and
so it has at most $k+k|\epsilon|+2$
nonzero real roots, by Descartes's rule of signs.
This proves the upper bound. 

We omit the construction which shows that this bound is sharp.
\QED

%
%
%
\chapter{Fewnomial upper bounds from Gale dual polynomial systems}\label{Ch:Gale}
%

\Magenta{{\sf This needs a proper Introduction}}

\begin{figure}[htb]
 \[
   \begin{picture}(280,180)
    \put(0,0){\includegraphics[height=180pt]{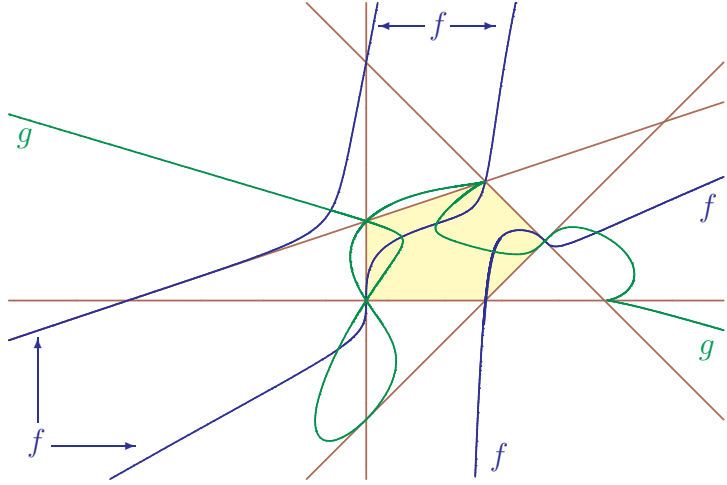}}
    \put(4,128){$\ForestGreen{g}$}
    \put(262,48){$\ForestGreen{g}$}
    \put( 8,11){$\Blue{f}$}  \put(183,6){$\Blue{f}$}
    \put(160,169){$\Blue{f}$}  \put(262,100){$\Blue{f}$}
    \put(17,13){\Blue{\vector(1,0){32}}}
    \put(12,22){\Blue{\vector(0,1){32}}}

    \put(158,172){\Blue{\vector(-1,0){17}}}
    \put(168,172){\Blue{\vector(1,0){17}}}

   \end{picture} 
 \]
\caption{Curves and lines}\label{F:Gale_pic1}
\end{figure}

Suppose that we have the system of polynomials,
 \begin{eqnarray}
    v^2w^3 - 11uvw^3   - 33uv^2w + 4v^2w + 15u^2v + 7 &=& 0\,,\nonumber\\
    v^2w^3\hspace{53pt}+\hspace{5.6pt}5uv^2w 
         - 4v^2w -\hspace{5.6pt}3u^2v + 1 &=& 0\,,\label{Eq5:system1}\\
    v^2w^3 - 11uvw^3   - 31uv^2w + 2v^2w + 13u^2v + 8 &=& 0\,.\nonumber
 \end{eqnarray}
If we solve them for the monomials $v^2w^3$, $v^2w$, and $uvw^3$, 
we obtain
 \begin{eqnarray}
    v^2w^3 &=&               1    - u^2v - uv^2w\,,\nonumber\\   \label{Eq5:system2}
    v^2w   &=&       \tfrac{1}{2} - u^2v + uv^2w\,,\\
   uv  w^3 &=&\tfrac{10}{11}(1    + u^2v - 3uv^2w)\,.\nonumber
 \end{eqnarray}
Since 
 \begin{eqnarray*}
  \left(uv^2w\right)^3 \cdot \left(v^2w^3\right) &=& u^3v^8w^6\ =\ 
   \left(u^2v\right)\cdot \left(v^2w\right)^3\cdot \left(uvw^3\right)\,
\qquad\mbox{and}\\
  \left(u^2v\right)^2 \cdot \left(v^2w^3\right)^3
   &=& u^4v^8w^9\ =\ 
  \left(uv^2w\right)^2\cdot \left(v^2w\right)\cdot \left(uvw^3\right)^2\,,
 \end{eqnarray*}
we may substitute the expressions on the right hand sides of~\eqref{Eq5:system2} 
for the monomials $v^2w^3$, $v^2w$, and $uvw^3$ in these expressions to obtain the system 
 \begin{eqnarray*}
  \left(uv^2w\right)^3 \cdot \left(1 - u^2v - uv^2w\right)
   &=&  \left(u^2v\right)\cdot \left(\tfrac{1}{2} - u^2v + uv^2w\right)^3
       \cdot \left(\tfrac{10}{11}(1+u^2v - 3uv^2w)\right)\,
    \qquad\mbox{and}\\
  \left(u^2v\right)^2 \cdot \left(1 - u^2v - uv^2w\right)^3
   &=& 
  \left(uv^2w\right)^2\cdot \left(\tfrac{1}{2} - u^2v + uv^2w\right)
        \cdot \left(\tfrac{10}{11}(1+u^2v - 3uv^2w)\right)^2\,.
 \end{eqnarray*}
Writing \DeCo{$x$} for $u^2v$ and \DeCo{$y$} for $uv^2w$ and solving for $0$, these become
 \begin{equation}\label{Eq5:Gale1}
  \begin{array}{rcrcl}
   \Blue{f}&:=&
       y^3(1-x-y)\ \ -\ \       x(\tfrac{1}{2}-x+y)^3\left(\tfrac{10}{11}(1+x-3y)\right) &=&
        0\,,\quad\mbox{and}\\ 
   \ForestGreen{g}&:=&
        x^2(1-x-y)^3\ \ -\ \ y^2(\tfrac{1}{2}-x+y)\left(\tfrac{10}{11}(1+x-3y)\right)^2 &=&
                  0\,.\rule{0pt}{15pt} 
  \end{array}
 \end{equation}
Figure~\ref{F:Gale_pic1} shows the curves these define and the 
lines given by the linear factors in~\eqref{Eq5:Gale1}.

It is clear that the solutions to~\eqref{Eq5:Gale1} in the complement of the lines are
consequences of solutions to~\eqref{Eq5:system1}.
More, however, it true.
The two systems define isomorphic schemes as complex or as real varieties, with the
positive solutions to~\eqref{Eq5:system1} corresponding to the solutions
of~\eqref{Eq5:Gale1} lying in the central pentagon.
Gale duality, which generalizes this isomorphism, is a first step towards the new
fewnomial bounds of~\cite{BaBiSo,BS07}. 

We remark that these new bounds are derived using the general method that Khovanskii
developed in~\cite{Kh80}.
However, they take advantage of special geometry (encoded in Gale duality) available to
systems of polynomials in a way that the proof of Khovanski's bound
(Theorem~\ref{T4:Fewnomial}) did not.
Their main value is that they are sharp, in the asymptotic sense described after
Corollary~\ref{C4:BRS}.

%
\section{Gale duality for polynomial systems}
%

Gale duality is an alternative way to view a sparse system of polynomials.
It was developed in~\cite{BS_Gale} in more generality than we will treat here.
Let us work over the complex numbers.
Let $\calA=\{0,a_1,\dotsc,a_{l+n}\}\subset\Z^n$ be integer vectors which span $\R^n$ and
suppose that we have a system
 \begin{equation}\label{Eq5:sp_system}
  f_1(x_n,\dotsc,x_n)\ =\ f_2(x_1,\dotsc,x_n)\ =\ 
   \dotsb \ =\ f_n(x_1,\dotsc,x_n)\ =\ 0
 \end{equation}
of polynomials with support $\calA$.
As in Section~\ref{S3:sparse}, the solutions may be interpreted geometrically as
$\varphi_\calA^{-1}(L)$, where $\varphi_\calA$ is the map
\[
  \begin{array}{rcl}
   \varphi_\calA\ \colon\ \T^n&\longrightarrow& \T^{l+n}\ \subset\
     \C^{l+n}\\
   x&\longmapsto& (x^{a_1},x^{a_2},\dotsc,x^{a_{l+n}})\,,\rule{0pt}{14pt}
  \end{array}
\]
and $L\subset\C^{l+n}$ is a codimension $n$ plane defined by degree 1 polynomials
$\Lambda_1,\dotsc,\Lambda_n$
corresponding to the polynomials $f_i$.
While we work here in $\C^{n+l}$, but used $\P^\calA=\P^{n+l}$ in Section~\ref{S3:sparse}, 
there is no essential difference.

Suppose that $\calA$ is primitive in that $\Z\calA=\Z^n$, so that the homomorphism
$\varphi_\calA$ is injective. 
Then the subscheme of $\T^n$ defined by~\eqref{Eq5:sp_system} is isomorphic to the
subscheme $\DeCo{X}:=\varphi(\T^n)\cap L$ of $\T^{l+n}$ or $\C^{l+n}$.
If we change our perspective and view $X$ as the basic object, then the parametrization
$\varphi_\calA$ of $\varphi_\calA(\T^n)$ realizes $X$ (or rather $\varphi_\calA^{-1}(X)$) as the
subscheme of $\T^n$ defined by~\eqref{Eq5:sp_system}.

The main idea behind Gale duality for polynomial systems is to instead parameterize $L$
with a map $\psi_p\colon \C^l\to L$ and then consider the subscheme $\psi_p^{-1}(X)$ of
$\C^l$, which is isomorphic to $X$. 
This is in fact what we did in 
transforming~\eqref{Eq5:system1} into~\eqref{Eq5:Gale1}. 
We will show that $\psi_p^{-1}(X)$ is defined in $\C^l$ by a system of master functions,
which we define in the next section.

%
\subsection{Master functions}
%

Let \DeCo{$p_1(y),\dotsc,p_{l+n}(y)$} be pairwise non-proportional degree 1 polynomials on
$\C^l$. 
Their product $\prod_i p_i(y)=0$ defines a hyperplane arrangement \DeCo{$\calH$} in
$\C^l$.
Let $\DeCo{\beta}\in\Z^{l+n}$ be an integer vector, called a \DeCo{{\sl weight}} for the
arrangement $\calH$.
We use this to define a rational function $p^\beta$,
\[
  \DeCo{p^\beta}\ =\ p(y)^\beta\ :=\ p_1(y)^{b_1} p_2(y)^{b_2}\dotsb
  p_{l+n}(y)^{b_{l+n}}\,,
\]
where $\beta=(b_1,\dotsc,b_{l+n})$.
This rational function $p(y)^\beta$ is a \DeCo{{\sl master function}} for the arrangement
$\calH$. 
As the components of $\beta$ may be negative, its natural domain of definition is the
complement $\DeCo{M_\calH}:=\C^l\setminus\calH$ of the arrangement.

A \DeCo{{\sl system of master functions}} in $M_\calH$ with weights
$\calB=\{\beta_1,\dotsc,\beta_l\}$ is the system of equations in $M_\calH$,
 \begin{equation}\label{Eq5:master}
   p(y)^{\beta_1}\ =\ p(y)^{\beta_2}\ =\ \dotsb\ =\ p(y)^{\beta_l}\ =\ 1\,.
 \end{equation}
More generally, we could instead consider equations of the form
$p(y)^\beta=\alpha$, where $\alpha\in\T$ is an arbitrary nonzero complex number.
We may however absorb such
constants into the polynomials $p_i(y)$, as there are $l{+}n$ such polynomials but only
$l$ constants in a system of master functions.
This may be viewed as the source of the factor $\frac{10}{11}$ in~\eqref{Eq5:Gale1}.
We further assume that the system~\eqref{Eq5:master} defines a zero-dimensional scheme in
$M_\calH$.
This implies that the weights $\calB$ are linearly independent, and that the
suppressed constants multiplying the $p_i(y)$ are sufficiently general.

As with sparse systems, a system of master functions may be realized geometrically through
an appropriate map.
The polynomials $p_1(y),\dotsc,p_{l+n}(y)$ define an affine map
\[
  \begin{array}{rcl}
   \DeCo{\psi_p}\ \colon\ \C^l&\longrightarrow& \C^{l+n}\\
   y&\longmapsto& (p_1(y), p_2(y),\dotsc,p_{l+n}(y))\,.\rule{0pt}{14pt}
  \end{array}
\]
This map is injective if and only if the polynomials $\{1,p_1(y),\dotsc,p_{l+n}(y)\}$
span the space of degree 1 polynomials on $\C^l$, in which case the hyperplane arrangement
$\calH$ is called essential.
The hyperplane arrangement $\calH$ is the pullback along $\psi_p$ of the coordinate
hyperplanes $z_i=0$ in $\C^{l+n}$, and its complement $M_\calH$ is the pullback of the
torus $\T^{l+n}$ which is the complement of the coordinate hyperplanes in $\C^{l+n}$.

The weights $\calB$ are \DeCo{{\sl saturated}} if they are linearly independent and span a
saturated subgroup of 
$\Z^{l+m}$, that is, if $\Z\calB=\Q\calB\cap\Z^{l+n}$.
Linear independence of $\calB$ is equivalent to the subgroup 
\DeCo{$\G$} of the torus $\T^{l+n}$ defined by the equations 
 \begin{equation}\label{Eq5:Torus_eqs}
   z^{\beta_1}\ =\ z^{\beta_2}\ =\ \dotsb\ =\ z^{\beta_l}\ =\ 1
 \end{equation}
having dimension $n$ and saturation is equivalent to $\G$ being connected.
(Here, $z_1,\dotsc,z_{l+n}$ are the coordinates for $\C^{l+n}$.)
In this way, we see that the equations~\eqref{Eq5:master} describe the pullback
$\psi_p^{-1}(\G)$ of this subgroup $\G$.
We summarize this discussion.

\begin{prop}\label{P5:Master}
  A system of master functions~$\eqref{Eq5:master}$ in $M_\calH$ is the pullback along
  $\psi_p$ of the intersection of the linear space $\psi_p(\C^l)$ with a subgroup $\G$ of\/
  $\T^{l+n}$ of dimension $n$, and any such pullback 
  defines a system of master functions in $M_\calH$.
  When $\psi_p$ is injective, it gives a scheme-theoretic isomorphism between the
  solutions of the system of master functions and the intersection
  $\G\cap \psi_p(\C^l)$.
\end{prop}

%
\subsection{Gale duality}
%
Proposition~\ref{P5:Master} is the new ingredient needed for the notion of Gale duality.
Suppose that $\G\subset\T^{l+n}$ is a connected subgroup of dimension $n$ and that
$L\subset\C^{l+n}$ is a linear subspace of dimension $l$ not parallel to any coordinate
hyperplane. 
Then the intersection $\G\cap L$ has dimension $0$.

\begin{defn}\label{D5:Gale_duality}
 Suppose that we are given
 \begin{enumerate}
  \item
        An isomorphism $\varphi_\calA\colon\T^n\to\G$ for $\calA=\{0,a_1,\dotsc,a_{l+n}\}$
        and equations~\eqref{Eq5:Torus_eqs} defining $\G$ as a subgroup of
        $\T^{l+n}$.
        Necessarily $\calA$ is primitive and $\calB=\{\beta_1,\dotsc,\beta_l\}$ is saturated.

  \item
        A linear isomorphism $\psi_p\colon\C^l\to L$ and degree 1 polynomials
        $\Lambda_1,\dotsc,\Lambda_n$ on $\C^{l+n}$ defining $L$.
 \end{enumerate}

 Let $\calH\subset\C^l$ be the pullback of the coordinate hyperplanes of $\C^{l+n}$ along 
 $\psi_p$.
 We say that the system of sparse polynomials on $\T^n$
 \begin{equation}\label{Eq5:GD_def_poly}
   \varphi^*(\Lambda_1)\ =\ \varphi^*(\Lambda_2)\ =\ \dotsb\ =\ 
      \varphi^*(\Lambda_n)\ =\ 0
 \end{equation}
 with support $\calA$ is \DeCo{{\sl Gale dual}} to the system of master functions
 \begin{equation}\label{Eq5:GD_def_master}
   p(y)^{\beta_1}\ =\ p(y)^{\beta_2}\ =\ \dotsb\ =\ p(y)^{\beta_l}\ =\ 1
 \end{equation}
 with weights $\calB$ on the hyperplane complement $M_\calH$ and vice-versa.\QED
\end{defn}

The following is immediate.

\begin{thm}\label{T5:Gale_dual}
  A pair of Gale dual systems~\eqref{Eq5:GD_def_poly} and~\eqref{Eq5:GD_def_master} define
  isomorphic schemes. 
\end{thm}

This notion of Gale duality involves two different linear algebraic dualities
in the sense of linear functions annihilating vector spaces.
In the first,  the degree 1 polynomials $p_i(y)$
defining the map $\psi_p$ are annihilated by the degree 
1 polynomials $\Lambda_i$ which define the system of sparse
polynomials~\eqref{Eq5:GD_def_poly}.
The second is integral, as the weights $\calB$ form a basis for the free
abelian group of integer linear relations among the nonzero elements of $\calA$.
Writing the elements of $\calB$ as the rows of a matrix, the $l{+}n$ columns form the 
\DeCo{{\sl Gale dual}} or Gale transform~\cite[\S 5.4]{Grunbaum} of the vector
configuration $\calA$---this is the source of our terminology.

If we restrict the domain of $\varphi_\calA$ to the real numbers or to the positive
reals, then we obtain the two forms of Gale duality which are relevant to us.
Set $\DeCo{\T_\R}:=\R\setminus\{0\}$, the real torus and \DeCo{$\R_>$} to be the positive
real numbers.
Suppose that $\calA$ is not necessarily primitive, but that the lattice index
$[\Z^n\colon\Z\calA]$ is odd.
Then $\varphi_\calA\colon \T_\R^n\to\T_\R^{l+n}$ is injective.
Similarly, if $\Z\calB$ has odd index in its saturation $\Q\calB\cap\Z^{l+n}$,
which is the group of integer linear relations holding on $\calA$, 
then the equations~\eqref{Eq5:Torus_eqs} define a not necessarily connected subgroup
$\G\subset\T^{l+n}$ whose real points $\G_\R$ lie in its connected component containing
the identity. 
When the linear polynomials $\Lambda_i$ of~\eqref{Eq5:GD_def_poly} and $p_i$
of~\eqref{Eq5:GD_def_master} are real and annihilate each other, then these two
systems---which do not necessarily define isomorphic schemes in $\T^n$ and
$M_\calH$---define isomorphic real analytic sets in $\T^n_\R$ for~\eqref{Eq5:GD_def_poly}
and in the complement $\DeCo{M^\R_\calH}:=\R^l\setminus \calH_\R$ of the real hyperplanes
defined by the $p_i$ for~\eqref{Eq5:GD_def_master}.

In the version valid for the positive real numbers, we may suppose that the exponents
$\calA$ are real vectors, for if $r\in\R_>$ and $a\in\R$, then
$\DeCo{r^a}:=\exp(a\log(x))$ is well-defined. 
In this case, the weights $\calB$ should be a basis for the vector space of linear
relations holding on $\calA$, and the degree 1 polynomials $\Lambda_i$ and $p_i$ are again
real and dual to each other. 
The equations~\ref{Eq5:Torus_eqs} for $z\in\R^{l+n}_>$ define a connected subgroup of
$\R^{l+n}_>$ which equals $\varphi_\calA(\R^n_>)$.
In this generality, the polynomial system~\eqref{Eq5:GD_def_poly} makes sense only for 
$x\in\R^n_>$ and the system of master functions~\eqref{Eq5:GD_def_master}  only makes
sense for $y$ in the \DeCo{{\sl positive chamber}} $\Delta_p$ of the hyperplane complement 
$M^\R_\calH$,
\[
  \DeCo{\Delta_p}\ :=\ 
   \{y\in\R^l\mid p_i(y)>0\quad i=1,\dotsc,l{+}n\}\,,
\]
and the two systems define isomorphic real analytic sets in $\R^n_>$
for~\eqref{Eq5:GD_def_poly} and in $\Delta_p$ for~\eqref{Eq5:GD_def_master}.

%
\subsection{Algebra of Gale duality}
%

The description of Gale duality in Definition~\ref{D5:Gale_duality} lends itself
immediately to an algorithm for converting a system of sparse polynomials into an
equivalent system of master functions.
We describe this over $\C$, but it works equally well over $\R$ or over $\R_>$. 
Suppose that $\calA\subset\Z^n$ is a primitive collection of integer vectors and suppose
that
 \begin{equation}\label{Eq5:PCI}
   f_1(x_1,\dotsc,x_n)\ =\ f_2(x_1,\dotsc,x_n)\ =\ 
   \dotsb\ =\ f_n(x_1,\dotsc,x_n)\ =\ 0
 \end{equation}
defines a zero dimensional subscheme of $\T^n$.
In particular, the polynomials $f_i$ are linearly independent.
We may solve these equations for some of the monomials to obtain
 \begin{equation}\label{Eq5:Diag_system}
  \begin{array}{rclcl}
  x^{a_1} &=& g_1(x)& =:& \DeCo{p_1}(x^{a_{n+1}},\dotsc,x^{a_{l+m+n}})\\\rule{0pt}{14pt}
          &\vdots&&&\ \DeCo{\vdots}\\\rule{0pt}{14pt}
  x^{a_n} &=& g_n(x)& =:& \DeCo{p_n}(x^{a_{n+1}},\dotsc,x^{a_{l+m+n}})
  \end{array}
 \end{equation}
Here, for each $i=1,\dotsc,n$, $g_i(x)$ 
is a polynomial with support $\{0,\, a_{n+1},\dotsc,a_{l+n}\}$ which is a 
degree 1 polynomial function $p_i(x^{a_{n+1}},\dotsc,x^{a_{l+n}})$ in the given $l$
arguments. 
For $i=n{+}1,\dotsc,l{+}n$, set 
$\DeCo{p_i}(x^{a_{n+1}},\dotsc,x^{a_{l+n}}):= x^{a_i}$.

An integer linear relation among the exponent vectors in $\calA$,
\[
  b_1 a_1\;+\; b_2a_2\;+\;\dotsb\;+\;
   b_{l+n}a_{l+n}\ =\ 0\,,
\]
is equivalent to the monomial identity
\[
   (x^{a_1})^{b_1}\cdot (x^{a_2})^{b_2}\dotsb(x^{a_{l+n}})^{b_{l+n}}
    \ =\ 1\,,
\]
which gives the consequence of the system~\eqref{Eq5:Diag_system}
\[
   \bigl(p_1(x^{a_{n+1}},\dotsc,x^{a_{l+n}})\bigr)^{b_1}\ \dotsb\ 
   \bigl(p_{l+n}(x^{a_{n+1}},\dotsc,x^{a_{l+n}})\bigr)^{b_{l+n}}
    \ =\ 1\,.
\]

Define \DeCo{$y_1,\dotsc,y_{l}$} to be new variables which are coordinates for 
$\C^{l}$.
The degree 1 polynomials $p_i(y_1,\dotsc,y_{l})$ define a hyperplane arrangement
$\calH$ in $\C^{l}$.
Let  $\calB:=\{\beta_1,\dotsc,\beta_l\}\subset\Z^{l+n}$ be a basis for the 
$\Z$-module of integer linear relations among the nonzero vectors in $\calA$.
These weights $\calB$ define a system of master functions 
 \begin{equation}\label{Eq5:MFCI}
   p(y)^{\beta_1}\ =\ p(y)^{\beta_2}\ =\ \dotsb\ =\ 
  p(y)^{\beta_l}\ =\ 1
 \end{equation}
in the complement $M_\calH$ of the arrangement.

\begin{thm}\label{Thm5:Alg-GD}
  The system of polynomials~\eqref{Eq5:PCI} in $\T^n$
  and the system of master functions~\eqref{Eq5:MFCI} in  
  $M_\calH$ define isomorphic schemes.
\end{thm}

\noindent{\it Proof.}
 Condition (1) in Definition~\ref{D5:Gale_duality} holds as $\calA$ and $\calB$ are both
 primitive and annihilate each other.
 The linear forms $\Lambda_i$ that pull back along $\varphi_\calA$ to 
 define the system~\eqref{Eq5:Diag_system} are
\[
   \Lambda_i(z)\ =\ z_i\ -\ p_i(z_{n{+}1},\dotsc,z_{l+n})\,,
\]
 which shows that condition (2) holds, and so the statement follows from
 Theorem~\ref{T5:Gale_dual}.
\QED

The example at the beginning of this chapter illustrated Gale duality, but 
the equations~\eqref{Eq5:Gale1} are not of the form $p^\beta=1$.
They are, however, easily transformed into such equations, and we obtain
 \begin{equation}\label{Eq5:Gale-as-master}
  \frac{y^3(1-x-y)}{x(\tfrac{1}{2}-x+y)^3\left(\tfrac{10}{11}(1+x-3y)\right)}\ =\ 
  \frac{x^2(1-x-y)^3}{y^2(\tfrac{1}{2}-x+y)\left(\tfrac{10}{11}(1+x-3y)\right)^2}\ =\ 1\,.
 \end{equation}
Systems of the form~\eqref{Eq5:Gale1} may be obtained from systems of master functions by 
multiplying $p^\beta=1$ by the terms of $p^\beta$ with negative exponents to clear the
denominators, to obtain $p^{\beta^+}=p^{\beta^-}$ and thus $p^{\beta^+}-p^{\beta^-}=0$,
where $\DeCo{\beta_{\pm}}$ is the componentwise maximum of the vectors
$(0,\pm\beta)$. 

%
\section{New fewnomial bounds}
%

The transformation of Gale duality is the key step in establishing the new fewnomial bounds.

\begin{thm}\label{T5:NewBounds}
  A system~$\eqref{Eq5:sp_system}$ of $n$ polynomials in $n$ variables having a total of
  $l{+}n{+}1$ monomials with exponents $\calA\subset\R^n$ has at most 
\[
   \frac{e^2+3}{4}\, 2^{\binom{l}{2}}n^l
\]
 positive nondegenerate solutions.

 If $\calA\subset\Z^n$ and $\Z\calA$ has odd index in $\Z^n$,
 then the system has at most
\[
   \frac{e^{\Red{4}}+3}{4}\, 2^{\binom{l}{2}}n^l
\]
 nondegenerate real solutions.
\end{thm}

The first bound is proven in~\cite{BS07} and the second in~\cite{BaBiSo}. 
By Gale duality, Theorem~\ref{T5:NewBounds} is equivalent to the nest Theorem.

\begin{thm}\label{T5:Gale_bound}
  Let $p_1(y),\dotsc,p_{l+n}(y)$ be degree $1$ polynomials on $\R^l$ that, with $1$, span
  the space of degree $1$ polynomials.
  For any linearly independent vectors $\calB=\{\beta_1,\dotsc,\beta_l\}\subset\R^{l+n}$,
  the number of solutions to
\[
   1\ =\ p(y)^{\beta_j}\qquad\mbox{for}\quad j=1,\dotsc,l
\]
  in the positive chamber $\Delta_p$ is less than
\[
   \frac{e^2+3}{4}\,2^{\binom{l}{2}}n^l\,.
\]
  If $\calB\subset\Z^{l+n}$ and has odd index in its saturation, then the number of solutions
  in $M^\R_\calH$  is less than
\[
   \frac{e^{\Red{4}}+3}{4}\,2^{\binom{l}{2}}n^l\,.
\]
\end{thm}

\Magenta{{\sf Should replace the hypotheses on the degree 1 polynomials by the notion of
    essential hyperplane arrangement, and have that arise in the discussion of Gale duality.}} 

We outline the proof of Theorem~\ref{T5:Gale_bound} in the next three sections.

%
\subsection{One idea}
%
The basic idea behind the proof of Theorem~\ref{T5:Gale_bound} is to use the
Khovanskii-Rolle Theorem, but in a slightly different 
form than given in Theorem~\ref{T4:Kh-R}.
Using it in this way to establish bounds for real solutions to equations was first done
in~\cite{GNS}.  
Given functions $g_1,\dotsc,g_m$ defined on a domain $D$, let 
\DeCo{$V(g_1,\dotsc,g_m)$} be their set of common zeroes.
If $C$ is a curve in $D$, let \DeCo{$\ubc(C)$} be its number of unbounded
components in $D$.

\begin{thm}[Khovanskii-Rolle]\label{T5:Kh-Ro}
  Let $g_1,\dotsc,g_l$ be smooth functions defined on a domain $D\subset\R^l$
  with finitely many common zeroes and suppose that $\DeCo{C}:=V(g_1,\dotsc,g_{l-1})$ is a
  smooth curve. 
  Set $\DeCo{J}$ to be the Jacobian determinant, $\det(\partial g_i/\partial y_j)$, 
  of $g_1,\dotsc,g_l$.
  Then we have 
 \begin{equation}\label{Eq5:KhRo}
   |V (g_1,\dotsc,g_l)|\ \leq\ \ubc(C) + |V(g_1,\dotsc,g_{l-1},\,J)| \,.    
 \end{equation}
\end{thm}

This form of the Khovanskii-Rolle Theorem follows from the from the usual Rolle Theorem.
Suppose that $g_l(a)=g_l(b)=0$, for points $a,b$ on the same component of $C$.
Let $s(t)$ be the arclength along this component of $C$, measured from a point $t_0\in C$,
and consider the map, 
 \begin{eqnarray*}
   C&\longrightarrow&\R^2\\
   t&\longmapsto& (s(t),g_l(t))\,.
 \end{eqnarray*}
This is the graph of a differentiable function $g(s)$ which vanishes when $s=s(a)$ and
$s=s(b)$, so there is a point $s(b)$ between $s(a)$ and $s(b)$ where its derivative aso
vanishes, by the usual Rolle Theorem.
But then $c$ lies between $a$ and $b$ on that component of $C$, and the vanishing of
$g'(s(c))$ is equivalent to the Jacobian determinant $J$ vanishing at $c$.

Thus along any arc of $C$ connecting two zeroes of $g_l$, the Jacobian
vanishes at least once.  
\[
  \begin{picture}(262,130)(-20,0)
   \put(0,0){\includegraphics[height=120pt]{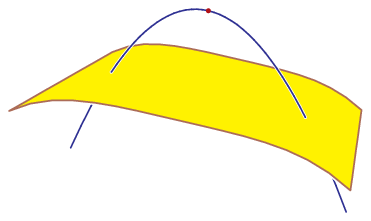}}
    \put(75,3){$C=V_\Delta(g_1,\dotsc,g_{l-1})$}
    \put(-20,65){$g_l=0$}
    \put(148,120){$V_\Delta(g_1,\dotsc,g_{l-1},J)$}
    \put(146,122){\vector(-3,-1){30}}
  \end{picture}
\] 

As in the proof of Theorem~\ref{T4:Kh-R}, the estimate of~\eqref{Eq5:KhRo} follows from
the observation concerning consecutive zeroes of $g_l$ along $C$.

%
\subsection{A generalization and two reductions}\label{S5:reduction}
%

We first make an adjustment to the system of master 
functions in Theorem~\ref{T5:Gale_bound}, replacing  each master function
$p(y)^\beta=p_1(y)^{b_1}\dotsb p_{l+n}(y)^{b_{l+n}}$ by 
 \[
   \DeCo{|p(y)|^\beta}\ := |p_1(y)|^{b_1}\dotsb |p_{l+n}(y)|^{b_{l+n}}\,.
 \]
For example, if we take absolute values in the system of master
functions~\eqref{Eq5:Gale-as-master}, we obtain
 \begin{equation}\label{Eq5:absolute_system}
   \frac{|y|^3|1-x-y|}{|x||\tfrac{1}{2}-x+y|^3|\tfrac{10}{11}(1+x-3y)|}\ =\ 
   \frac{|x|^2|1-x-y|^3}{|y|^2|\tfrac{1}{2}-x+y||\tfrac{10}{11}(1+x-3y)|^2}\ =\ 1\,.
 \end{equation}
This new system with absolute values will still have the same number of solutions in the
positive chamber $\Delta_p$ as the original system, since $|p_i(y)|=p_i(y)$ for
$i=1,\dotsc,l{+}n$ and $y\in\Delta_p$. 
Its solutions in the hyperplane complement $M^\R_\calH$ will include the solutions to the
system of master functions from Theorem~\ref{T5:Gale_bound}, but there may be more solutions. 

\begin{figure}[htb]
 \[
   \begin{picture}(250,200)
    \put(0,0){\includegraphics[height=200pt]{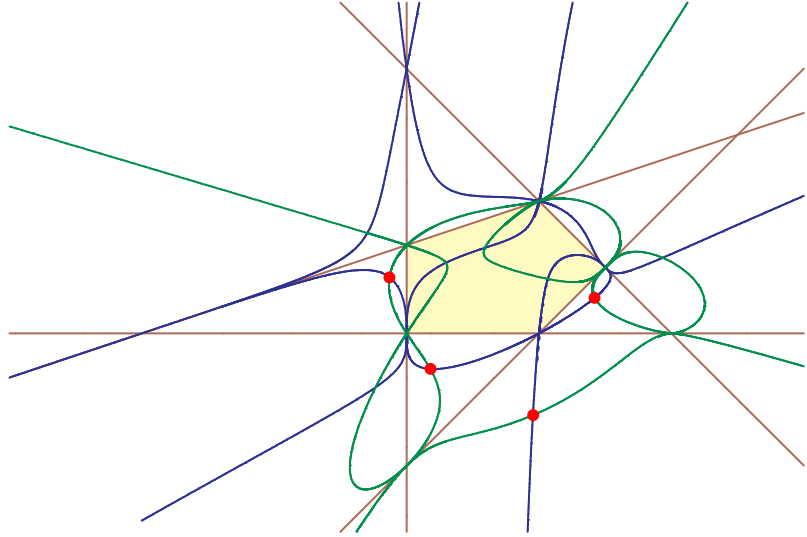}}
   \end{picture}
 \]
\caption{Curves from absolute values}\label{F:Gale_pic2}
\end{figure}

We illustrate this for the system~\eqref{Eq5:absolute_system} in Figure~\ref{F:Gale_pic2},
which we may compare to Figure~\ref{F:Gale_pic1} as the system of master
functions~\eqref{Eq5:Gale-as-master} is equivalent to the system~\eqref{Eq5:Gale1} 
in the complement of the lines.
In particular, among the solutions to the system~\eqref{Eq5:absolute_system} are the 
three solutions to~\eqref{Eq5:Gale1} in the positive chamber (which is
shaded) as well as the three solutions to~\eqref{Eq5:Gale1} outside the positive chamber.
The system~\eqref{Eq5:absolute_system} has four additional solutions outside the positive
chamber, which are marked in Figure~\ref{F:Gale_pic2}.

We give a proof of Theorem~\ref{T5:Gale_bound} for systems of the form
$|p(y)|^{\beta_j}=1$ for $j=1,\dotsc,l$.
This will imply the bound for systems of master functions.
Taking absolute values allows nonintegral (real number) exponents in 
$|p(y)|^{\beta_j}=1$, and so we need not require that exponents
are integral.\medskip

We promised two reductions.
\begin{enumerate}

\item  The degree 1 polynomials $p_i(y)$ are in general position
       in that the hyperplanes in the arrangement $\calH$ are in linear general position.
       That is, any $j$ of them meet in an affine linear subspace of codimension $j$, if $j\leq l$,
       and their intersection is empty if $j>l$.
       We may do this, as we are bounding nondegenerate solutions, which cannot be
       destroyed if the $p_i(y)$ are perturbed to put the hyperplanes into this general
       position. 

\item  Let $B$ be the matrix whose rows are $\beta_1,\dotsc,\beta_l$.  
       We may assume that every minor of $B$ is non-zero.
       This may be done by perturbing the real-number exponents in the functions
       $|p(y)|^{\beta_j}$. 
       This will not reduce the number of nondegenerate solutions.

       Perturbing exponents is not as drastic of a measure as it first seems.
       Note that in the hyperplane complement, $|p(y)|^\beta=1$ defines the
       same set as $\log(|p(y)|^\beta)=0$.
       If $\beta=(b_1,\dotsc,b_{l+n})$, then this is simply
 \begin{equation}\label{Eq5:logs}
      b_1\log|p_1(y)|\ +\  b_2\log|p_2(y)|\ +\ \dotsb \ +\  b_{l+n}\log|p_{l+n}(y)|\ =\
       0\,.
 \end{equation}
       Expressing the equations in this form shows that we may perturb the exponents. 
\end{enumerate}

We first look at these reductions in the context of the system of master
functions~\eqref{Eq5:Gale-as-master}. 
The hyperplane arrangement $\calH$ is an arrangement of lines in which no three meet and 
no two are parallel, and thus they are in general position.
The matrix of exponents is
\[
    B\ =\  \left(\begin{array}{rrrrr}
       2 & -2 & 3 & -1 & -2 \\ 1 & -3 & -1 & 3 & 1
    \end{array}\right)\,.
\]
No entry and no minor of $B$ vanishes.

\begin{figure}[htb]
\[
  \begin{picture}(380,235)
   \put(0,0){\includegraphics[height=235pt]{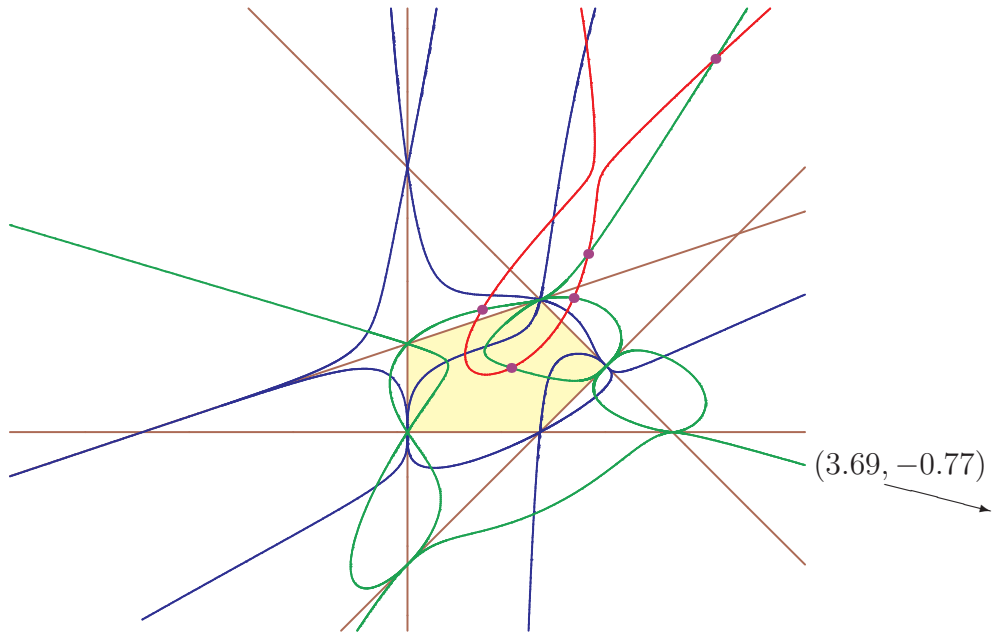}}
   \put(305,60){$(3.69,-0.77)$} 
   \put(332,56){\vector(4,-1){40}}
  \end{picture}
\]
\caption{Gale system and Jacobian $J_2$.}\label{F6:KhRo_proof}
\end{figure}

Let us now see how the Khovanskii-Rolle Theorem applies to the 
system~\eqref{Eq5:absolute_system} of Figure~\ref{F:Gale_pic2}, restricted to the positive
chamber.
First, take logarithms and rearrange to obtain 
 \begin{eqnarray*}
   \ForestGreen{2\log |x| - 2 \log |y| + 3\log |1{-}x{-}y|  -\log |\tfrac{1}{2}{-}x{+}y| -
   2\log|\tfrac{10}{11}(1{+}x{-}3y)|}
      &=& 0\\
   \Blue{\log |x| - 3\log |y| + 1 \log |1{-}x{-}y|  + 3\log |\tfrac{1}{2}{-}x{+}y|
   +\log |\tfrac{10}{11}(1{+}x{-}3y)|} &=&0
 \end{eqnarray*}
Call these functions \ForestGreen{$f_1$} and \Blue{$f_2$}, respectively.
Their Jacobian is the rational function
\[
    \frac{2x^3-16x^2y+12xy^2+6y^3 
      -\tfrac{31}{2}x^2+26xy-\tfrac{53}{2}y^2+\tfrac{9}{2}x+\tfrac{15}{2}y-2}%
    {xy(1-x-y)(\tfrac{1}{2}{-}x{+}y)(1+x-3y)}
\]
whose denominator is the product of the linear factors defining the lines in
Figure~\ref{F:Gale_pic2}. 
Clearing the denominator and multiplying by 2, we obtain a cubic polynomial
\[
   \DeCo{J_2}\ :=\  \Red{4x^3-32x^2y+24y^2x+12y^3-31x^2+52xy-53y^2+9x+15y-4}\,.
\]
Its zero set meets the curve \ForestGreen{$C_1$} (which is defined by $\ForestGreen{f_1}=0$) in
6 points, five of which we display in Figure~\ref{F6:KhRo_proof}---the sixth is at
$(3.69,-0.77)$. 
By the Khovanskii-Rolle Theorem, the number of solutions to $f_1=f_2=0$ is at most this
intersection number, $|V(f_1,J_2)|$, plus the number of unbounded components of $C_1$.

We see that \ForestGreen{$C_1$} has 14 unbounded components, which gives the inequality
\[
   10\ =\ |V(f_1,f_2)|\ \leq\  \ubc(C_1) + |V(f_1,J_2)| \ =\ 14 +6\ =\ 18\,.
\]

%
\subsection{Three lemmata}
%

We follow the suggestion in the second reduction in Section~\ref{S5:reduction} and replace
the master functions by the logarithms of their absolute values.
For each $j=1,\dotsc,l$, define
\[
   \DeCo{g_j(y)}\ :=\ \log|p(y)|^{\beta_j}\ =\ 
    \sum_{i=1}^{l+n} \beta_{i,j} \log|p_i(y)|\ ,
\]
where $\beta_{i,j}$ is the $i$th component of $\beta_j$.
Observe that both $f_j=0$ and $|p(y)|^{\beta_j}=1$ have the same solutions in the hyperplane
complement $M^\R_\calH$.
For each $j=1,\dotsc,l$, define $\DeCo{\mu_j}:= V(f_1,\dotsc,f_j)$.

Recall that the connected components of the complement $M^\R_\calH$ are called
chambers. 
A \DeCo{{\sl flat}} of the arrangement $\calH$ is an affine subspace which is an
intersection of some hyperplanes in $\calH$.
By our assumption that the hyperplane arrangement $\calH$ is in general position, a flat
of $\calH$ has codimension $j$ exactly when it is the intersection of $j$ hyperplanes in
$\calH$. 

\begin{lemma}\label{L5:one}
  For each $j=1,\dotsc,l{-}1$, $\mu_j$ is a smooth submanifold of $M^\R_\calH$
  of codimension $j$.
  The closure $\overline{\mu_j}$ of $\mu_j$ in $\R^l$ meets the arrangement $\calH$ in a union
  of codimension $j{+}1$ flats.
  In the neighborhood of point on a codimension $j{+}1$ flat meeting $\overline{\mu_j}$, 
  $\mu_j$ has at most one branch in each chamber incident on that point.
\end{lemma}

Define functions $J_l,J_{l-1},\dotsc,J_1$ by recursion,
\[
   J_j\ :=\ \mbox{Jacobian of } g_1,\dotsc,g_j,\ J_{j+1},\dotsc,J_l\,.
\]
The following is proven in~\cite[\S~3.1]{BS07} using the Cauchy-Binet Theorem.

\begin{lemma}\label{L5:Jacobians}
 ${\displaystyle J_j\cdot\Bigl(\prod_{i=1}^{l+n} p_i(y) \Bigr)^{2^{l-j}}}$ is a polynomial
 of degree $2^{l-j}\cdot n$.
\end{lemma}

For each $j=1,\dotsc,l$, define $\DeCo{C_j}:=\mu_{j-1}\cap V(J_{j+1},\dotsc,J_l)$.
By our assumptions that the polynomials $g_i$ and the exponents $\beta_j$ are general,
this will be a smooth curve in the hyperplane complement $M^\R_\calH$.

We now iterate the Khovanskii-Rolle Theorem~\ref{T5:Kh-Ro} to estimate the number of
solutions to a system of master functions as in Theorem~\ref{T5:Gale_bound}.
 \begin{eqnarray}
   |V(g_1,g_2,\dotsc,g_l)| &\leq& \ubc(C_l)\,+\,   \nonumber
    |V(g_1,g_2,\dotsc,g_{l-1},\,J_l)|\\
  &\leq& \ubc(C_l)\,+\,\ubc(C_{l-1})\, +\,         \label{Eq:Kp_form}
     |V(g_1,\dotsc,g_{l-2},\,J_{l-1},J_l)|\\
  &\leq& \ubc(C_l)\,+\,\dotsb\,+\,\ubc(C_1)\, +\,
    |V(J_1,J_2,\dotsc,J_l)|\,.  \nonumber
 \end{eqnarray}
Here, $V(\dotsb)$ is the common zeroes in the hyperplane
complement $M^\R_\calH$.
Let $\ubc_\Delta(C)$ count the number of unbounded components of the curve $C$ in
$\Delta_p$ and $V_\Delta(\dotsb)$ be the common zeroes in $\Delta_p$.
Then the analog of~\eqref{Eq:Kp_form} holds in $\Delta_p$.

\begin{lemma}\label{L5:estimates}
 With these definitions, we have the estimates
 \begin{enumerate}
  \item[$(1)$] $|V_\Delta(J_1,\dotsc,J_l)|\ \leq\ |V(J_1,\dotsc,J_l)|\ \leq\ 2^{\binom{l}{2}}n^l$.
  \item[$(2)$] $\ubc_\Delta(C_j)\ \leq\ \tfrac{1}{2} \binom{1+l+n}{j}\cdot 2^{\binom{l-j}{2}}
        n^{l-j}$.
  \item[$(3)$] $\ubc(C_j)\ \leq\ \tfrac{1}{2} \binom{1+l+n}{j}\cdot 2^{\binom{l-j}{2}}
        n^{l-j}\cdot \Red{2^j}$.
 \end{enumerate}
\end{lemma}

The first statement follows from Lemma~\ref{L5:Jacobians} and B\'ezout's Theorem.
For the second, nore that $C_j=\mu_{j-1}\cap V(J_{j+1},\dotsc,J_l)$.
Since each unbounded component has two ends, we estimate the number of unbounded
components by $\frac{1}{2}$ of the number of points in the boundary of $\Delta$ that are
limits of points of $C_j$.
By Lemma~\ref{L5:one}, these will be points in the codimension skeleta where 
$J_{j+1},\dotsc,J_l$ vanish.
The bound in the second statement of~\ref{L5:estimates} is simply $\frac{1}{2}$ multiplied
by the product of $\binom{1+l+n}{j}$ and $2^{\binom{l-j}{2}}n^{l-j}$.
That is by the number of codimension $J$ flats in $\calH$ (some of which meet the
boundary of $\Delta$) multiplied by the B\'ezout number of the
system $J_{j{+}1}\ =\ \dotsb\ =J_{l}= 0$.

Note that the bound in (2)~holds for any chamber of $M^\R_\calH$.
We get the bound in (3)~by noting that in the neighborhood of any
codimension $j$ stratum of $\calH$, the complement has at most $2^j$ chambers,
and so each point of $C_j$ at the boundary of $M^\R_\calH$ can contribute at most
$2^j$ such ends.

The complement $M^\R_\calH$ of the hyperplane arrangement consists of 
many chambers.
The first bound of Theorem~\ref{T5:Gale_bound} is in fact a bound for the number of
solutions in any chamber, while the second bound is for the number of solutions in all
chambers.
This is smaller than what one may naively expect.
The number of chambers in a generic arrangement of $1{+}n{+}l$ hypersurfaces in $\R^l$ is
\[
   \tbinom{2+n+l}{l}+\tbinom{2+n+l}{l-2}+\dotsb+\left\{ \begin{array}{rcl}
       \tbinom{2+n+l}{2}+1&\ &\mbox{if $1+n+l$ is even}\\
       1+n+l &\ &\mbox{if $1+n+l$ is odd} \end{array}\right.\ .
\]
Thus, we would naively expect that the ratio between  $\ubc_\Delta(C_j)$ and 
$\ubc(C_j)$ to be this number, rather than the far smaller $2^j$.
This is the source for the mild difference between the two estimates in
Theorem~\ref{T5:NewBounds}.

We may combine the estimates of Lemma~\ref{L5:estimates} with~\eqref{Eq:Kp_form} to estimate
$|V_\Delta(J_1,\dotsc,J_l)|$ and $|V(J_1,\dotsc,J_l)|$,
 \begin{eqnarray*}
  |V_\Delta(J_1,\dotsc,J_l)|&\leq& 
    \frac{1}{2} \sum_{j=1}^l \tbinom{1+l+n}{j}\cdot 2^{\binom{l-j}{2}}n^{l-j}\ +\ 
              2^{\binom{l}{2}}n^l\ ,\quad\mbox{and}\\
  |V(J_1,\dotsc,J_l)|&\leq& 
    \frac{1}{2} \sum_{j=1}^l \tbinom{1+l+n}{j}\cdot 2^{\binom{l-j}{2}}n^{l-j}\cdot 2^j\ +\ 
              2^{\binom{l}{2}}n^l\,.
 \end{eqnarray*}
It is not hard to show the estimate~\cite[Eq.(3.4)]{BS07}
 \[
   \tbinom{1+l+n}{j}\cdot 2^{\binom{l-j}{2}}n^{l-j}\ \leq\ 
   \frac{2^{j-1}}{j!} 2^{\binom{l}{2}}n^l\,,
 \]
so that these estimates become
 \[
  \begin{array}{rclcl}
  |V_\Delta(J_1,\dotsc,J_l)|&\leq& 
    \Bigl(\frac{1}{2} \sum_{j=1}^l \frac{2^{j-1}}{j!} \ \ +\ 1\Bigr) 2^{\binom{l}{2}}n^l
   &\leq& \frac{e^2+3}{4}  2^{\binom{l}{2}}n^l\ ,\quad\mbox{and}   \\
  |V(J_1,\dotsc,J_l)|&\leq& 
    \Bigl(\frac{1}{2} \sum_{j=1}^l \frac{2^{2j-1}}{j!} \ \ +\ 1\Bigr) 2^{\binom{l}{2}}n^l
   &\leq& \frac{e^4+3}{4}  2^{\binom{l}{2}}n^l\,. 
  \end{array}
 \]
This implies Theorem~\ref{T5:Gale_bound} and thus the new fewnomial bounds,
Theorem~\ref{T5:NewBounds}. \QED


%
%
%
%
%
%
%
\chapter{Lower Bounds for Sparse Polynomial Systems} \label{Ch:lower}
%

%
%
In Chapter~\ref{ch:1}, we mentioned how work of Welschinger~\cite{W},
Mikhalkin~\cite{Mi05}, and of Kharlamov, Itenberg, and
Shustin~\cite{IKS03,IKS04} combined to show that there is a 
nontrivial lower bound $W_d$ for the number of real rational curves
of degree $d$ interpolating $3d{-}1$ points in 
$\R\P^2$.
If $N_d$ is the Kontsevich number~\eqref{E1:Konts} of complex rational curves
of degree $d$ interpolating $3d-1$ points in $\C\P^2$, then we have
 \begin{enumerate}
  \item ${\displaystyle W_d\geq  \frac{d!}{3}}$, \quad and
  \item ${\displaystyle \lim_{d\to\infty}\frac{\log(N_d)}{\log(W_d)}\ =\ 1}$. \ 
        (In fact, $\log(N_d)\sim 3d\log(d)\sim \log(W_d)$.)
 \end{enumerate}
An exposition of the beginning of this story written for a general mathematical audience
is given in~\cite{So04}.
Similar results were found by Solomon~\cite{Sol} for rational curves on real Calabi-Yau threefolds.
For example, there are at least $15$ real lines on a smooth quintic hyperurface in $\P^4$.

Eremenko and Gabrielov~\cite{EG01} have a similar result 
for the number of real solutions to the inverse Wronski problem.
They gave numbers $\sigma_{m,p}>0$ for $0<m\leq p$ with $m{+}p$ odd,
and proved that if $\Phi$ is a real polynomial of degree $mp$ then there
are at least $\sigma_{m,p}$ different $m$-dimensional subspaces of the
form Span$\{f_1,f_2,\dotsc,f_m\}$ where each $f_i$ is a real polynomial of
degree $m{+}p{-}1$, and their Wronskian $\Wr(f_1,f_2,\dotsc,f_m)$ is a scalar multiple
of $\Phi$.  
We will discuss this in Section~\ref{S6:SAGBI}.

While these results are spectacular, they are but the beginning of what
we believe will be a bigger and more important story (at least for the
applications of mathematics).
These are examples of natural geometric problems
possessing a \DeCo{{\sl lower bound}} on their numbers of real solutions.
It would be a very important development for some applications 
if this phenomenon were widespread,  if there were
methods to detect when such lower bounds existed, and if there were also methods to
compute or estimate these lower bounds.
The point is that nontrivial lower bounds imply the existence of real
solutions to systems of equations, or to interesting geometric problems.
A beginning of the interaction between applications and this new theory of
lower bounds is found in work of Fiedler-Le Touz\'e~\cite{FLT} and discussed
in Section~\ref{S1:lower} of Chapter~\ref{ch:1}.

This chapter will report on the first steps toward a
theory of lower bounds for sparse polynomial systems as given in~\cite{SS}.
There are three papers where one may read more about this subject.

\begin{enumerate}

\item[{\cite{EG01}}]
A.~Eremenko and A.~Gabrielov, \emph{Degrees of real {W}ronski maps}, Discrete
  Comput. Geom. \textbf{28} (2002), no.~3, 331--347.

  \begin{itemize}
   \item Establishes a lower bound for the Wronski map, realized as
    a degree of its lift to oriented double covers.
  \end{itemize}

\item[{\cite{SS}}]
E.~Soprunova and F.~Sottile, \emph{Lower bounds for real solutions to sparse
  polynomial systems}, Advances in Math., \textbf{204} (2006), no.~1, 116--151.

  \begin{itemize}
   \item Begins the theory of lower bounds to sparse polynomial systems.
  \end{itemize}

\item[{\cite{JW07}}] M.~Joswig and N.~Witte, \emph{Products of foldable
    triangulations},  Advances in Math., {\bf 210} (2007), no.~2, 769--796.

  \begin{itemize}
   \item  Uses geometric combinatorics to give many more examples of sparse
    polynomial systems with a lower bound on their number of real solutions.
  \end{itemize}
\end{enumerate}

The last two papers study lower bounds for unmixed systems, such as those covered by
Kouchnirenko's Theorem.
It remains an important open problem to develop a theory for unmixed systems such as those
which appear in Bernstein's Theorem.\fauxfootnote{Put a statement like this in the chapter
  on Gale duality!} 
Example~\ref{Ex6:mixed} at the end of this chapter is a first step in this direction.

\section{Polynomial systems from posets}\label{S6:PolyPoset}
 
Let $P$ be a finite partially ordered set (\DeCo{{\sl poset}}), whose elements
we take to be our variables.
For convenience, assume that $P$ has $n$ elements.
The order ideals $I\subset P$ are its up sets---subsets that are closed upwards.
Specifically, a subset $I\subset P$ is an \DeCo{{\sl order ideal\,}} if whenever $x<y$ in
$P$ with $x\in I$, then we must have $y\in I$.
Our polynomials will have monomials indexed by these order ideals.
For $I\subset P$ set
\[
   x^I\ :=\ \prod_{x\in I} x\,.
\]
A \DeCo{{\sl linear extension}} $w$ of $P$ is a permutation 
$w\colon x_1,x_2,\dotsc,x_n$ of the elements of $P$
such that $x_i<x_j$ in $P$ implies that $i<j$.
Let $\lambda(P)$ denote the number of linear extensions of a poset $P$.
The \DeCo{{\sl sign-imbalance $\sigma(P)$}} of $P$ is
 \begin{equation}\label{Eq6:sign_imbalance}
   \sigma(P)\ :=\ \left|\sum \mbox{sgn}(w)\right|\ ,
 \end{equation}
the sum over all linear extensions $w$ of $P$.
Here, $\mbox{sgn}(w)$ is the sign of the permutation $w$.
(To be precise, this makes sense only if we first fix one linear extension $\pi$
 of $P$ and measure the sign of $w\pi^{-1}$.
 Taking the absolute value in~\eqref{Eq6:sign_imbalance} removes the effect of that choice.) 
Lastly, a \DeCo{{\sl Wronski polynomial}} for $P$ is a polynomial of the
form
 \begin{equation}\label{E6:Wr_simple}
  \sum_I c_{|I|} x^I\,,
 \end{equation}
the sum over all order ideals $I$ of $P$, where $c_0, c_1,\dotsc,c_n\in\R$.
Note that the coefficient of the monomial $x^I$ depends only upon its degree,
$|I|$. 
The reason for the term Wronski polynomial is that this comes from a linear
projection with a special form that is also shared by the Wronski map.
(This is discussed just before Section~\ref{S6:Schub_Vars} in Chapter~\ref{Ch:ScGr}.)

\begin{thm}\label{T6:poset}
 Suppose that $P$ is a finite poset in which every maximal chain has the same
 parity.
 Then a generic system of Wronski polynomials for $P$ has $\lambda(P)$ complex
 solutions, at least $\sigma(P)$ of which will be real.
\end{thm}

\begin{rmk}
 A variant of this is to first fix real numbers $\alpha_I$ for each order ideal $I$ and
 then consider polynomials of the form
 \begin{equation}\label{E6:Wronski_general}
  \sum_I c_{|I|} \alpha_I x^I\,.
 \end{equation}
 That is, the coefficients $c_0,c_1,\dotsc,c_n$ vary in $\R^n$, but the numbers
 $\alpha_I$ are fixed.
 The same statement holds about such systems.
 In fact, this variant is closest to the Wronski map in Schubert calculus.
 \QED
\end{rmk}

\begin{ex}
 Let $P$ be the incomparable union of two chains, each of length 2,
\[
  \begin{picture}(95,50)(-35,0)
    \put(-35,20){$P\ :=$}
    \put(14,3){\includegraphics{figures/6/P2P2.eps}}
    \put(0, 0){\DeCo{$w$}}   \put(61, 0){\Brown{$y$}}
    \put(0,40){\DeCo{$x$}}   \put(61,40){\Brown{$z$}}
    \put(75,20){.}
  \end{picture}
\]
 Here are the monomials corresponding to the order ideals of $P$.
\[
  \{\emptyset,\ \DeCo{x},\ \Brown{z},
     \ \DeCo{wx},\ \DeCo{x}\Brown{z},\ \Brown{yz},
     \ \DeCo{x}\Brown{yz},\ \DeCo{wx}\Brown{z},\ \DeCo{wx}\Brown{yz}\}\,.
\]
 A Wronski polynomial for $P$ has the form
 \begin{eqnarray}
  &c_4\, \DeCo{wx}\Brown{yz}\hspace{.2em}&\nonumber\\
  &+\ c_3(\DeCo{x}\Brown{yz}\ +\ \DeCo{wx}\Brown{z})\hspace{1em}&\nonumber\\
  &+\ c_2(\DeCo{wx}\ +\ \DeCo{x}\Brown{z}\ +\ \Brown{yz})\hspace{1.2em}&\label{E6:P2P2}\\
  &+\ c_1(\DeCo{x}\ +\ \Brown{z})\hspace{1.2em}&\nonumber\\
  &\hspace{.2em}+\ c_0\,,&\nonumber
 \end{eqnarray}
where the coefficients $c_0,\dotsc,c_4$ are real numbers.

There are six linear extensions of $P$, as each is a permutation of the word
$\DeCo{wx}\Brown{yz}$ where $\DeCo{w}$ precedes $\DeCo{x}$ 
and $\Brown{y}$ precedes $\Brown{z}$.
The sign-imbalance of $P$ is seen to be 2, as computed in the following table.
\begin{center}
 \begin{tabular}{|r||c|c|c|c|c|c||c|}\hline
  permutation&$\DeCo{wx}\Brown{yz}$
             &$\DeCo{w}\Brown{y}\DeCo{x}\Brown{z}$
             &$\Brown{y}\DeCo{wx}\Brown{z}$
             &$\DeCo{w}\Brown{yz}\DeCo{x}$
             &$\Brown{y}\DeCo{w}\Brown{z}\DeCo{x}$
             &$\Brown{yz}\DeCo{wx}$&$\sigma(P)$\\\hline
  sign & $+$ & $-$  & $+$  & $+$ & $-$  & $+$ & 2\\\hline
 \end{tabular}
\end{center}

Since $P$ has two maximal chains and each contains 2 elements, it satisfies the
hypotheses of Theorem~\ref{T6:poset}, and so we conclude that:
\begin{center} 
 A system of four equations involving polynomials of the\\
 form~\eqref{E6:P2P2} has six solutions, at least two of which are real.
 \makebox[.1in][l]{\hspace{2.1cm}\raisebox{-2pt}{\includegraphics[height=12pt]{figures/HSBC.eps}}}
\end{center}
\end{ex}

The rest of this chapter will discuss the steps in the proof of
Theorem~\ref{T6:poset}, the more general results contained
in~\cite{SS}, and some examples that illustrate this phenomenon of gaps.

\section{Orientability of real toric varieties}

Recall from Chapter~\ref{Ch:sparse} that a sparse polynomial system with support 
$\calA\subset\Z^n$ is
equivalent to a linear section $X_\calA\cap L$ of the toric variety $X_\calA$.
In fact, its solutions are 
 \[
  \varphi_\calA^{-1}(X_\calA\cap L)\,,
 \]
where 
$\DeCo{\varphi_\calA}\colon\T^n\ni x\mapsto[x^a\mid a\in\calA]\in\P^\calA$ is the
parametrization map and $X_\calA$ is the closure of its image.\medskip

\noindent\Brown{{\bf Issue 1.}} In order for there to be a bijective 
correspondence between real solutions of the original system and real points in
the linear section, we need that the map $\varphi_\calA$ be injective on 
$(\R^\times)^n$, which is the condition that the index $[\Z\calA\colon\Z^n]$ be
odd.
This should be regarded as a minor issue and a mild assumption.
We have already seen this in Chapters~\ref{Ch:sparse} and~\ref{Ch:Gale}.
\medskip

\begin{figure}[htb]
 \[
  \begin{picture}(230,150)(-20,-2)
   \put(0,0){\includegraphics{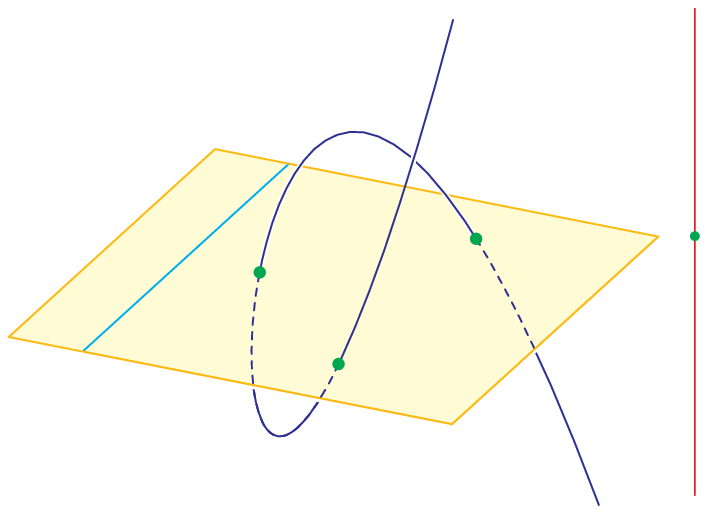}} 

   \put(10,130){$\pi$}\put(20,133){\vector(1,0){20}}

   \put(-18,72){$E$}\put(-6,74){\vector(2,-1){20}}
   
   \put(141,113){$L$}
   \put(145,109){\vector(-1,-4){5}}

   \put(207,79){$y$}\put(204,81){\vector(-1,0){20}}
   \put(207,119){$M$}\put(204,121){\vector(-1,0){20}}

   \put(155,8){$X_\calA$}

   \put(77,14){\vector(0,1){26}}
   \put(83,14){\vector(1,2){31}}
   \put(73,14){\vector(-1,3){17}}
   \put(67,2){$\pi^{-1}(y)$}

  \end{picture}
 \]
\caption{A linear projection $\pi$ with center $E$.\label{F6:cubic_Projection}}
 
\end{figure}

The key idea at the beginning of this theory is to realize the intersection
$X_\calA\cap L$ as the fiber of a map.
To that end, let $E\subset L$ be a hyperplane in $L$ that does not meet $X_\calA$ and
$M\simeq\P^n$ a linear space that is disjoint from $E$.
Then $E$ has codimension $n{+}1$ and the set of codimension $n$ planes containing $E$ is
naturally identified with $M\simeq\P^n$,
as each codimension $n$ plane containing $E$ meets $M$ in a unique point.
Define the \DeCo{{\sl linear projection}}
 \begin{equation}\label{Eq6:linProj}
   \DeCo{\pi} = \pi_E\ \colon\ \P^\calA-E\ \ \longrightarrow\ M\simeq \P^n
 \end{equation}
by sending a point $p\in \P^\calA-E$ to the intersection of $M$ with span of $E$ and $p$.
Figure~\ref{F6:cubic_Projection} illustrates this in $\P^3$, where $E$ and $M$ are lines.
We sometimes write $\pi\colon\P^\calA-\to\P^n$, using the broken arrow $-\to$ to indicate
that the map is not defined on all of $\P^\calA$.

Write $\pi$ for the restriction of the linear projection to the toric variety $X_\calA$.
If $\DeCo{y}:=L\cap M$ is the point where $L$ meets $M$, then 
\[
   X_\calA \cap L \ =\ \pi^{-1}(y)\,.
\]
This is also illustrated in Figure~\ref{F6:cubic_Projection}, where
$X_\calA$ is a cubic curve.
The reason that we reformulate our system of polynomials as the fiber of a map is so that
we may use topological methods to study its solutions.

That is, we work over $\R$, define $Y_\calA:=X_\calA\cap\R\P^\calA$, and consider
the real linear section $Y_\calA\cap L_\R$.
Equivalently, we restrict the projection $\pi$ further to a map
\[
   f\ \colon\ Y_\calA\ \to\ \R\P^n\,,
\]
and consider points in the fiber $f^{-1}(p)$, where $p\in\R\P^n$.
Since $Y_\calA$ and $\R\P^n$ have the same dimension, the map $f$ may have a topological
degree.
For this, $Y_\calA$ and $\R\P^n$ must be orientible and we must fix orientations of
$Y_\calA$ and of $\R\P^n$. 
Then, for every regular regular value $y\in\R\P^n$ of $f$, the differential map $df$
is a bijection on the tangent spaces of every point $x\in f^{-1}(y)$.
Define $\sgn(x):=1$ if the differential preserves the orientations at $x$ and $y$ and 
$\sgn(x):=-1$ if the orientations are reversed.
The topological degree of $f$ is the sum
\[
   \sum_{x\in f^{-1}(y)} \sgn(x)\,.
\]
This definition does not depend upon the choice of a regular value $y$, as $\R\P^n$ is
connected.\fauxfootnote{\Magenta{Give a reference for this.}}
The value of the notion of topological degree is the
following.\fauxfootnote{\Red{Reconcile this with the use of degree in Chapter 2.}}

\begin{thm}\label{T6:degree}
   The number of points in a fiber $f^{-1}(y)$, for $y\in\R\P^n$ a regular value of $f$,
   is at least the absolute value of the topological degree of $f$.
\end{thm}

Since both $\R\P^n$ and $Y_\calA$ are not always orientable, the topological degree of $f$
is not always defined, and so we do not get a bound as in Theorem~\ref{T6:degree}.
To remedy this, we consider also the lift of the linear projection $\pi$ to 
the double covers coming from spheres:
 \begin{equation}\label{E6:lifts}
  \raisebox{-20pt}{
   \begin{picture}(160,50)(-5,0)
   \put( -3, 0){$f$}
   \put( 11, 0){$\colon$}
   \put( 23, 0){$Y_\calA$}
   \put( 45, 0){$\subset$}
   \put( 57, 0){$\R\P^\calA$}
   \put( 84, 0){$\stackrel{\pi}{---\to}$}
   \put(135, 0){$\R\P^n$}

   \put( -3,34){$f^+$}
   \put( 11,34){$\colon$}
   \put( 23,34){$Y^+_\calA$}
   \put( 45,34){$\subset$}
   \put( 62,34){$S^\calA$}
   \put( 84,34){$\stackrel{\pi^+}{---\to}$}
   \put(140,34){$S^n$}

   \put( 28,28){\vector(0,-1){15}}
   \put( 68,28){\vector(0,-1){15}}
   \put(144,28){\vector(0,-1){15}}
    \end{picture}}
 \end{equation}
Then the topological degree of the map $f^+$ provides a lower bound for the number of
solutions, as in Theorem~\ref{T6:degree}.
We seek criteria which imply that the map $f$  
or its lift $f^+$ to the double covers has a well-defined degree.
That is, either $Y_\calA$ or $Y^+_\calA$ is orientable.

For this, we use Cox's construction of $X_\calA$ as a quotient
of a torus acting on affine space, as detailed in~\cite[Theorem 2.1]{Co95}.
(In the symplectic category, this is due to Delzant.)
The convex hull of $\calA$ is the polytope $\Delta:=\Delta_\calA$.
We need the dual description of $\Delta$ in terms of intersections of
half-spaces, or facet inequalities, which have the form
\[
   \Delta\ =\ \{z\in\R^n\ \mid\ B\cdot z\ \geq\ -b\}\,,
\]
where $B$ is an integer $N$ by $n$ matrix (N is the number of facets of
$\Delta$) whose rows are the inward-pointing normals to the facets of $\Delta$, and
$b\in\Z^N$ measures the (signed) lattice distance of each facet from the origin.
Then there is a natural parametrization 
$\psi_\Delta\colon\C^N-Z_\Delta\twoheadrightarrow X_\calA$,
realizing $X_\calA$ as a geometric invariant theory quotient of $\C^N$ by a
subtorus of $\T^N$.
(We do not give the details here, see~\cite{Co95} or~\cite{SS}.)
\medskip

\noindent\Brown{{\bf Issue 2.}} 
 While $\psi_\Delta$ restricts to a map 
 $\psi_\Delta\colon\R^N-Z_\Delta\to Y_\calA$ (and also to $Y^+_\Delta$), these
 maps are not necessarily surjective.
 This map $\psi_\Delta$ is surjective if and only if the column span 
 $\col(B)$ of the matrix $B$ has odd index in its saturation, 
 $\Z^N\cap (\Q\otimes_\Z\col(B))$.
\medskip

When the standard orientation of $\R^N$ drops to an orientation either of 
$Y_\calA$ or of $Y^+_\calA$ under the Cox quotient map $\psi_\Delta$, we say
that $Y_\calA$ is \DeCo{{\sl Cox-orientable}}.

\begin{thm}
 With the assumptions outlined in \Brown{Issue $1$} and  \Brown{Issue $2$}, if
 $\col(B)+\Z\cdot b$ contains a vector, all of whose coordinates are odd, then
 $Y_\calA$ is Cox-orientable.
 If $\col(B)$ has such a vector, then $Y_\calA$ receives the orientation,
 otherwise $Y^+_\calA$ receives the orientatation.
\end{thm}

\noindent{\bf Remarks.}
\begin{enumerate}

\item[1.]
   Given a projection map $\pi$~\eqref{Eq6:linProj} whose center $E$ is disjoint
   from $Y_\calA$, write $f$ for its restriction to $Y_\calA$.
   Suppose that $Y_\calA$ is Cox-orientable.
   Then lift $f^+\colon Y^+_\calA\to S^n$ of $f$ to $Y^+_\calA$ has a well-defined degree,
   which is a lower bound for the number of 
   real solutions to polynomial systems arising as fibers of the map $f$.
   Call this number the \DeCo{{\sl real degree}} of the map $f$.

\item[2.]  
   If $P$ is a poset and $\calA\subset\{0,1\}^P$ consists of the indicator functions of
   its order ideals, so that a Wronski polynomial for $P$~\eqref{E6:Wr_simple} has support
   $\calA$, then it is not hard to show that $\Z\calA=\Z^n$ and also $Y_\calA$ is
   Cox-orientable if all maximal chains of $P$ have the same parity.
   
   This parity condition is sufficient, but not necessary, for the Wronski
   polynomial system on $Y_\calA$ to have a lower bound.
   The poset
\[
  \raisebox{11pt}{P\ =\ }\ \ \includegraphics[height=40pt]{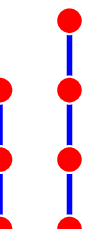}\
  \raisebox{11pt}{\ ,}
\]
   has two maximal chains of lengths 3 and 4, and so it is not necessarily Cox-oriented.
   Nevertheless, its Wronski polynomial systems have $\binom{7}{3}=35$ 
   solutions, at least $\binom{3}{1}=3$ of which are real. 
   See Section~\ref{S6:posets}\QED
\end{enumerate}

\section{Degree from foldable triangulations}

These results provide us with a challenge:
compute the real degree of a (or any) map $f$ arising as a linear projection
of a Cox-orientable toric variety $Y_\calA$.
We give a methods that uses the toric degenerations of Chapter~\ref{Ch:sparse} to provide
an answer to this question.
It is by no means the only answer.

By Corollary~\ref{C3:reg_triangulation}, if $\calA$ consists of the integer points in a
lattice polytope $\Delta$ which has a regular unimodular triangulation $\Delta_\omega$,
then there exists a flat deformation of $X_\calA$ within the projective space $\P^\calA$
into a limit scheme which is the union of coordinate $n$-planes
\[
  \lim_{t\to0} t.X_\calA\ =\ \bigcup_{\tau\in\Delta_\omega}\P^\tau\,.
\]
The idea here is to find conditions on the triangulation that allow us to compute
the real degree of some map $f$ from this deformation.
Some triangulations $\Delta_\omega$ have a naturally defined linear projection
with a nicely defined degree.

\begin{ex}
Consider the triangulation of the hexagon (the HSBC Bank symbol rotated
$45^\circ$ anti-clockwise) shown below with the
vertices labeled \DeCo{$a$}, \Magenta{$b$}, and \PineGreen{$c$}.
Mapping these vertices to the corresponding vertices of the simplex
defines a piecewise-linear `folding' map \DeCo{$\pi_\omega$} whose degree is
$4-2=2$. 
\[
  \begin{picture}(170,88)(0,-2)
   \put(10,10){\includegraphics[height=60pt]{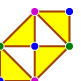}}
   \put( 0, 0){\DeCo{$a$}} \put(30,30){\DeCo{$a$}} \put(75,75){\DeCo{$a$}}
   \put( 0,37){\PineGreen{$c$}}  \put(75,37){\PineGreen{$c$}}
   \put(38,-3){\Magenta{$b$}}    \put(38,75){\Magenta{$b$}}
   \put(89, 37){$\xrightarrow{\ \pi_\omega\ }$}
   \put(130,25){\includegraphics[height=30pt]{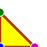}}
   \put(120,15){\DeCo{$a$}} 
   \put(120,55){\PineGreen{$c$}}
   \put(162,12){\Magenta{$b$}}  
   
   \put(295,0){\includegraphics[height=12pt]{figures/HSBC.eps}}
  \end{picture}
\]
\end{ex}

Suppose now that $\calA=\Delta\cap\Z^n$, where $\Delta$ is a lattice polytope.
A triangulation $\Delta_\omega$ is \DeCo{{\sl foldable}} if its facet simplices
may be properly 2-colored, which is equivalent to there being a coloring of the
vertices $\calA$ of the triangulation with $n+1$ labels, where each simplex
receives all $n{+}1$ labels~\cite{J02}. 
Both the 2-coloring and the vertex-labeling are unique up to permuting the
colors and labels. 
The difference in the number of simplices of different colors is the 
\DeCo{{\sl signature} $\sigma(\omega)$} of the balanced triangulation.
Up to a sign, it is the degree of the combinatorial folding map 
from $\Delta$ to an $n$-simplex given by the labeling of $\calA$.

The vertex labels define a linear projection 
$\pi_\omega\colon\P^\calA-\to \P^n$.
Let $\e_a$ be the standard basis vector of $\P^\calA$ corresponding to the monomial
$a\in\calA$ and $\e_j$ for $j=0,1,\dotsc,n$ the standard basis vector of $\P^n$.
We assume that the labels in the balanced triangulation take values
$0,1,\dotsc,n$.
Then $\pi_\omega(\e_a)=\e_j$, where the vertex $a$ of has label $j$.
The restriction of this linear projection to the coordinate spheres
$\pi_\omega\colon \cup_\tau S^\tau \to S^n$ is the geometric counterpart of this
combinatorial folding map.

From  the geometry of this map, we can show that for $t>0$ and small,
the restriction $f$ of $\pi_\omega$ to $t.Y^+_\calA$ has degree
$\sigma(\omega)$. 
This can be deduced from the real version of the arguments used in the
algorithmic proof of Koushnirenko's Theorem in Lecture~2.
In fact, Figure~\ref{F:odd_even} from Lecture~2 is reproduced from this argument
in~\cite{SS}. 
\Magenta{{\sf expand on this??}}

The restrictions of $\pi_\omega$ to $t.Y^+_\calA$ for $t$ near zero and $t=1$
are homotopic, \DeCo{{\bf if}} \,
$t.Y_\calA\cap\ker(\pi_\omega)=\emptyset$, for $t\in(0,1)$.
This condition is not hard to check on
specific examples.

\section{Reprise: polynomial systems from posets}\label{S6:posets}
 Suppose now that $P$ is a poset and $\calA\subset\{0,1\}^P$ is the set of indicator
 functions of order ideals $I$ of $P$, so that a Wronski polynomial for $P$ has support
 $\calA$. 
 Then $\Delta_P:=\conv(\calA)$ is the \DeCo{{\sl order polytope}} of $P$.
 This was studied by Stanley~\cite{St86b}.
 Its facet inequalities are
\[
  B\ \colon\ \left\{\begin{array}{rrrcl}
    y& \geq&  0&\quad&\mbox{if }y\; \in\; P\ \mbox{is mininmal}\\
   -y& \geq& -1&\quad&\mbox{if }y\; \in\; P\ \mbox{is maximal}\\
  x-y& \geq&  0&&\mbox{if } y\; \lessdot\; x\ 
        \mbox{is a cover in }P\end{array}\right.\ .
\]

 \Magenta{{\sl Give more details in these brief paragraphs}}
 It has a unimodular triangulation $\Delta_\omega$ whose simplices correspond to
 linear extensions of $P$, and it is easy to see that it is regular.
 Unimodular implies that $\Z\calA=\Z^n$.
 Recently, Stanley~\cite{St05} and others introduced the notion of 
 \DeCo{{\sl sign-(im)balance}} of posets, $\sigma(P)$, and it is 
 tautological that $\sigma(\omega)=\sigma(P)$.
 It is also easy to use the description of the facet inequalities to show that
 $\col(B)$ is saturated, and if every maximal chain in $P$ has the same parity,
 then $Y_\calA$ is Cox-oriented.
 This explains Theorem~\ref{T6:poset}. \QED

 \Magenta{{\sl Discuss the phenomenon of gaps.}}

\section{Sagbi degenerations}\label{S6:SAGBI}
 This method of computing a degree through a limiting process applies to the
original result of Eremenko and Gabrielov on the degree of the Wronski map in
Schubert calculus~\cite{EG01}.
The Grassmannian $\Gr(p,m{+}p)$ of $p$-planes in $(m{+}p)$-space admits a
\DeCo{{\sl sagbi degeneration}} to the toric 
variety associated to the poset $P_{p,m}$ which is the product of a chain with
$p$ elements and a chain with $m$ elements, and one may show that the degree of
the Wronski map for the Grassmannian $\Gr(p,m{+}p)$ is equal to the degree of the
Wronski map $f$ (the restriction of $\pi_\omega$) on $Y^+_\calA$, where $\calA$
is the set of monomials coming from order ideals of $P_{p,m}$.
(For this, we use a more general version of the Wronski projection $\pi_\omega$
which corresponds to the more general Wronski polynomial
of~\eqref{E6:Wronski_general}.) 

Using the sagbi degeneration, we can recover the results of Eremenko and
Gabrielov.
In fact, this connection between the Grassmannian and its toric degeneration,
as well as between Eremenko and Gabrielov's formula for the degree of the
Wronski map, the sign-imbalance of $P_{p,m}$, and the geometric folding map
was the genesis of this work~\cite{SS}.


\section{Open problems}
 There is much more to be done in this area.
 Here are some suggestions.

 \begin{enumerate}
   \item Give more comprehensive conditions which imply that $f$, or its lift
   to some (not necessarily the one given above) double cover of $Y_\calA$ is
   orientable. 

   \item Find other methods to give polynomial systems whose degree may be
   computed or estimated.

   \item Find more balanced triangulations (see~\cite{JW07}).

   \item Apply these ideas to specific problems from the applied sciences.

   \item Extend any of this from unmixed systems (all polynomials have the same
         Newton polytope) to more general mixed systems (those whose
         polynomials have different Newton polytopes).
         We end with an example in this direction which is due to Chris Hillar.
 \end{enumerate}

\begin{ex}\label{Ex6:mixed}
 Let $P$ and $Q$ be the two lattice polytopes given below
\[
   \begin{picture}(348,90)(-27,-8)
    \put(0,0){\includegraphics{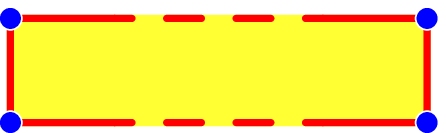}}
    \put(59,43){$P$}
    \put(-27,38){$(0,1)$}    \put(130,38){$(m,1)$}
    \put(-27,-8){$(0,0)$}    \put(130,-8){$(m,0)$}

    \put(230,0){\includegraphics{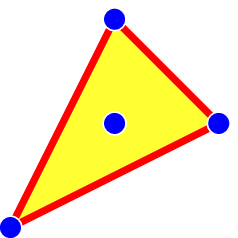}}
    \put(232,38){$Q$}
               \put(250, 73){$(1,2)$}
                               \put(300, 31){$(2,1)$}
    \put(203,-8){$(0,0)$}

   \end{picture}
\]
 A polynomial with support $P$ has the form
\[
   g\ :=\ A(x)\ +\ y B(X)\,,
\]
 where $A$ and $B$ are univariate polynomials in $x$ with degree $m$.
 Suppose that their coefficients are $a_0,\dotsc,a_m$ and $b_0,\dotsc,b_m$, and let
 $h$ be the polynomial with support $Q$,
\[
   h\ :=\ c\ +\ dxy\ +\ ex^2y\ +\ fxy^2\,.
\]
 By Bernstein's Theorem (Theorem~\ref{T1:Bernstein}), the mixed system $g(x,y)=h(x,y)=0$
 will have $2m+2$ solutions in $(\C^\times)^2$, as $2m+2$ is the mixed volume of $P$ and
 $Q$.
 (For polygons $P,Q$, the mixed volume is $\vol(P+Q)-\vol(P)-\vol(Q)$.)
 We can compute an eliminant for this system by substituting $-A(x)$ for $yB(x)$ in 
 $h\cdot B(x)^2$, to obtain
\[
    cB(x)^2\ -\ dxA(x)B(x)\ -\ ex^2A(x)B(x) + fxA(x)^2\,.
\]
 This has constant term $cb_0^2$ and leading term $-ea_mb_m$.
 If $ce>0$ and $a_mb_m>0$, then these have different signs, which implies that
 the mixed system has at least one positive root (and hence at least two real roots).
 This may be ensured by the condition that none of the coefficients vanish and
 the two linear equations, 
 $c+e=a_m+b_m=0$.
\QED
\end{ex}

%
%
\renewcommand{\QED}{\hfill\raisebox{-4.5pt}{\includegraphics[height=17pt]{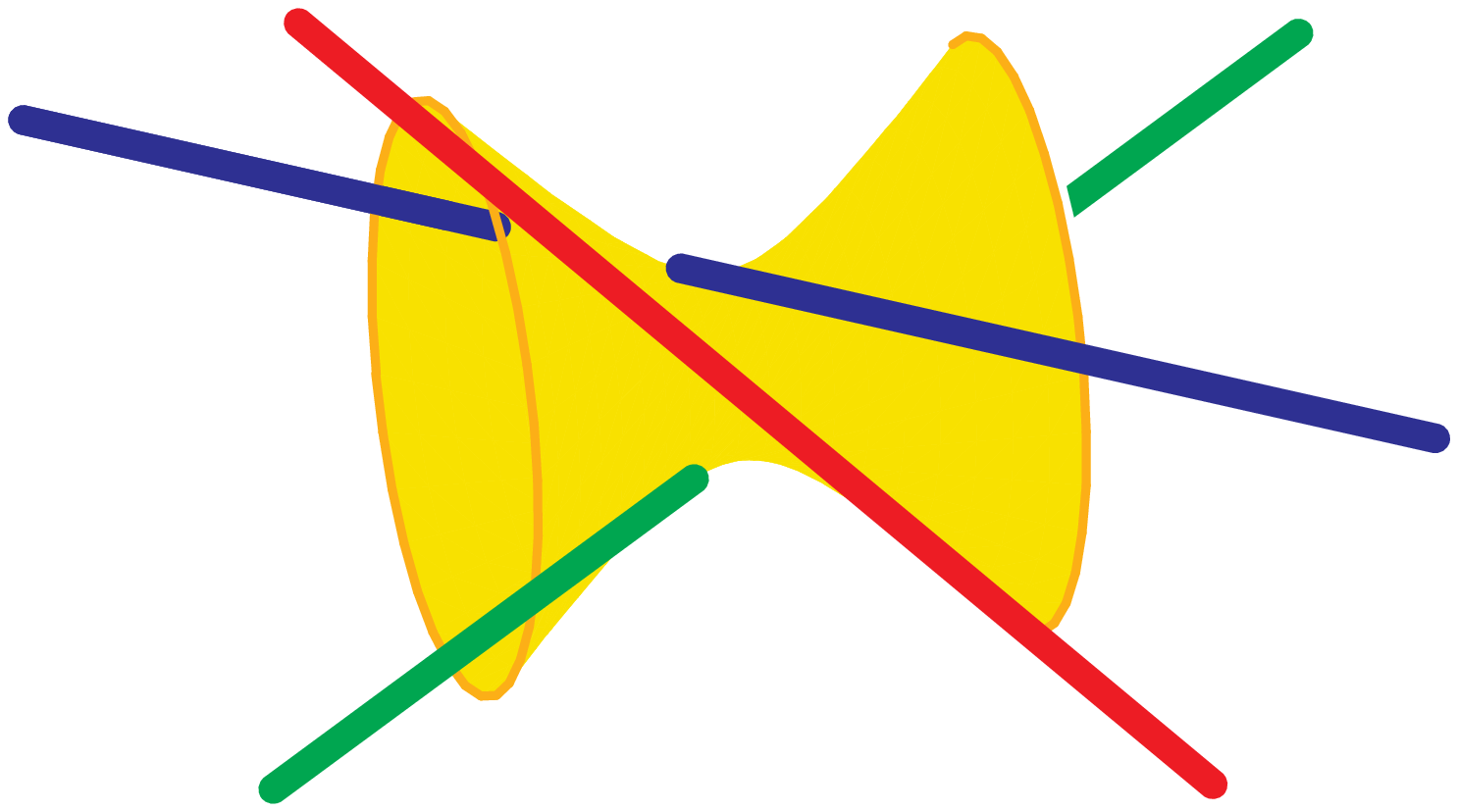}}\medskip}
%
%
%
%
%
\chapter{Enumerative Real Algebraic Geometry}\label{S:ERAG}

%
%

Enumerative geometry is the art of counting geometric figures satisfying 
conditions imposed by \DeCo{other}, fixed, geometric figures.  
For example, in 1848, Steiner~\cite{St1848} asked how many plane conics are tangent to
five given conics?
His answer, $6^5=7776$, turned out to be incorrect, and in 1864
Chasles~\cite{Ch1864} gave the correct answer of 3264.
These methods were later systematized and used to great effect by
Schubert~\cite{Sch1879}, who codified the field of enumerative geometry.

This classical work always concerned \DeCo{{\sl complex}} figures.
It was only in 1984 that the question of reality was
posed by 
Fulton~\cite[p.~55]{Fu96b}: ``The question of how many solutions of real
equations can be real is still very much open, particularly for enumerative
problems.'' 
He goes on to ask: ``For example, how many of the 3264 conics tangent to five
general conics can be real?''
He later determined that all can be real, but did not publish that
result.
Ronga, Tognoli, and Vust~\cite{RTV97} later gave a careful argument that 
\DeCo{all} 3264 can be real.

Since this work of Ronga, Tognoli, and Vust, there have been many geometric problems for
which it was shown that all solutions may be real.
This means that the upper bound
of $d$ (= number of complex solutions) is sharp for these problems.
In this chapter, we describe some of these problems, beginning with the problem of conics,
and concentrating on the Schubert calculus.
A survey of these questions (circa 2002) is given in~\cite{So03}.

\section{3264 real conics}

The basic idea of the arguments of Fulton and of Ronga, Tognoli, and Vust is to deform the
same special configuration. 
We will sketch the idea in Fulton's construction.

Suppose that $\ell_1,\dotsc,\ell_5$ are the lines supporting the edges of a
convex pentagon and $p_i\in\ell_i$, $i=1,\dotsc,5$ are points in the interior of
the corresponding edge. 
\[
  \begin{picture}(155,110)(-4,3)
   \put(-4,0){\includegraphics[height=110pt]{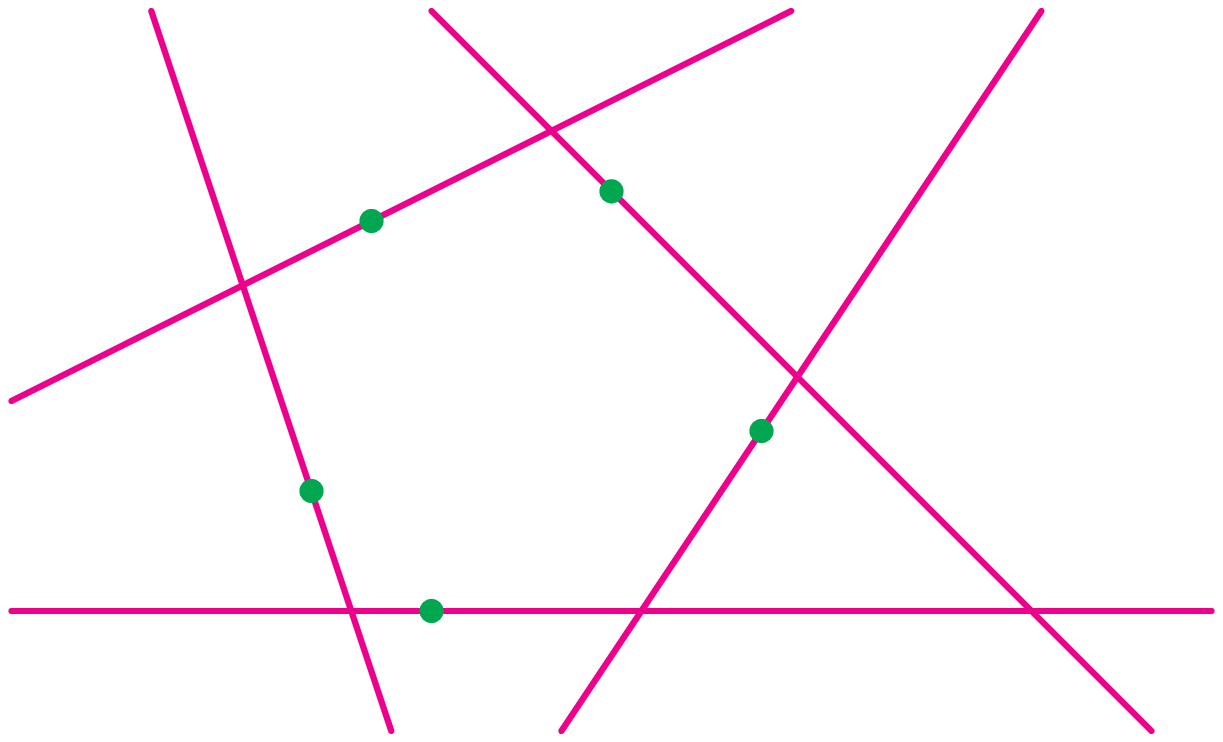}}
   \put(7,94){$\ell_5$}   \put(96,109){$\ell_4$}  \put(142,31){$\ell_3$}
   \put(1,8){$\ell_1$}    \put(88,0){$\ell_2$}
   \put(27,39){$p_5$}   \put(39,83){$p_4$}  \put(93,84){$p_3$}
   \put(59,9){$p_1$}    \put(95,45){$p_2$}
  \end{picture}
\]
The points in this example are $\{(0,0), (\frac{11}{4},\frac{3}{2}),
 (\frac{3}{2},\frac{7}{2}),(-\frac{1}{2},\frac{13}{4}), (-1,1)\}$, and the corresponding
 slopes of the lines are  $0,\frac{3}{2},-1,\frac{1}{2},-3$.

\begin{figure}[htb]
\[
   \begin{picture}(450,270)
    \put(0,0){\includegraphics[height=270pt]{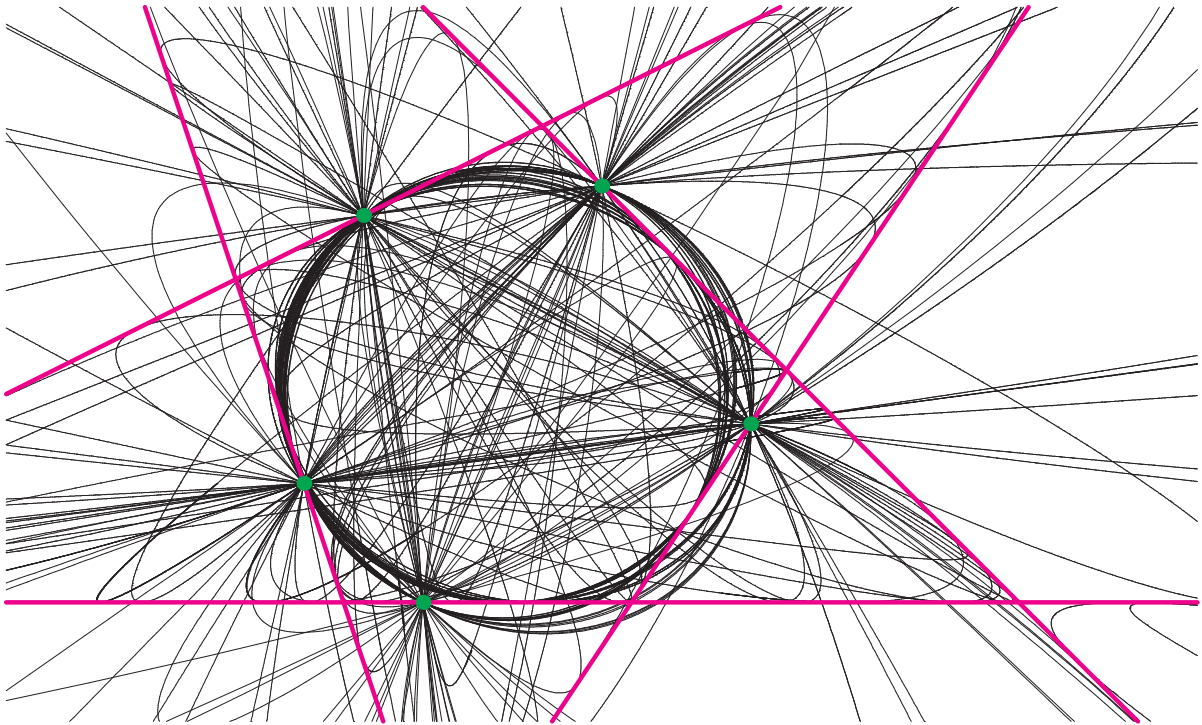}}
  \end{picture}
\]
\caption{102 conics}\label{F7:102_conics}
\end{figure}

For every subset $S$ of the lines, there are \DeCo{$2^{\min\{|S|, 5-|S|\}}$}
conics that are tangent to the lines in $S$ and that meet the $5{-}|S|$ points
not on the lines of $S$.
This is the number of complex conics, and it does not depend upon the configuration of
(generic) points and lines.
However, when the points and lines are chosen in convex position,
then all such conics will be real.
Altogether, this gives
\[
   2^0\tbinom{5}{0}\ +\ 2^1\tbinom{5}{1}\ +\ 2^2\tbinom{5}{2}\ +\ 
   2^2\tbinom{5}{3}\ +\ 2^1\tbinom{5}{4}\ +\ 2^0\tbinom{5}{5}\ =\ 102
\]
\DeCo{real} conics, that, for each $i=1,\dotsc,5$ either meet $p_i$ or are tangent to
$\ell_i$. 
We draw these in Figure~\ref{F7:102_conics}.
Since our pentagon was asymmetric, exactly 51 of these conics meet each point $p_i$ and 
none of the 51 conics tangent to $\ell_i$ are tangent at $p_i$. 

The idea now is to replace each pair $p_i\in\ell_i$ by a 
hyperbola \DeCo{$h_i$} that is close to the pair $p_i\in\ell_i$,
in that $h_i$ lies close to its asymptotes, which are two lines close to $\ell_i$ that meet at
$p_i$. 
If we do this for one pair $p_i\in\ell_i$, then, for every conic in our configuration,
there will be two nearby conics tangent to $h_i$. 
To see this, suppose that $i=1$.
Then the set \DeCo{$C$} of conics which satisfy one of the conditions ``meet $p_j$'' or 
``tangent to $\ell_j$'' for each $j=2,3,4,5$ will form an irreducible curve $C$.
For each conic in $C$ that meets $p_1$, there will be two nearby conics in $C$ tangent to
$h_1$ near $p_1$, and for each conic in $C$ tangent to $\ell_1$, there will be two nearby 
conics in $C$ tangent to each of the two nearby branches of $h_1$.
We illustrate this when $C$ is the curve of conics tangent to $\ell_2,\ell_3,\ell_4$, and
$\ell_5$, showing the  conics in $C$,
\[
  \begin{picture}(295,240)(-35,0)
   \put(0,0){\includegraphics[height=240pt]{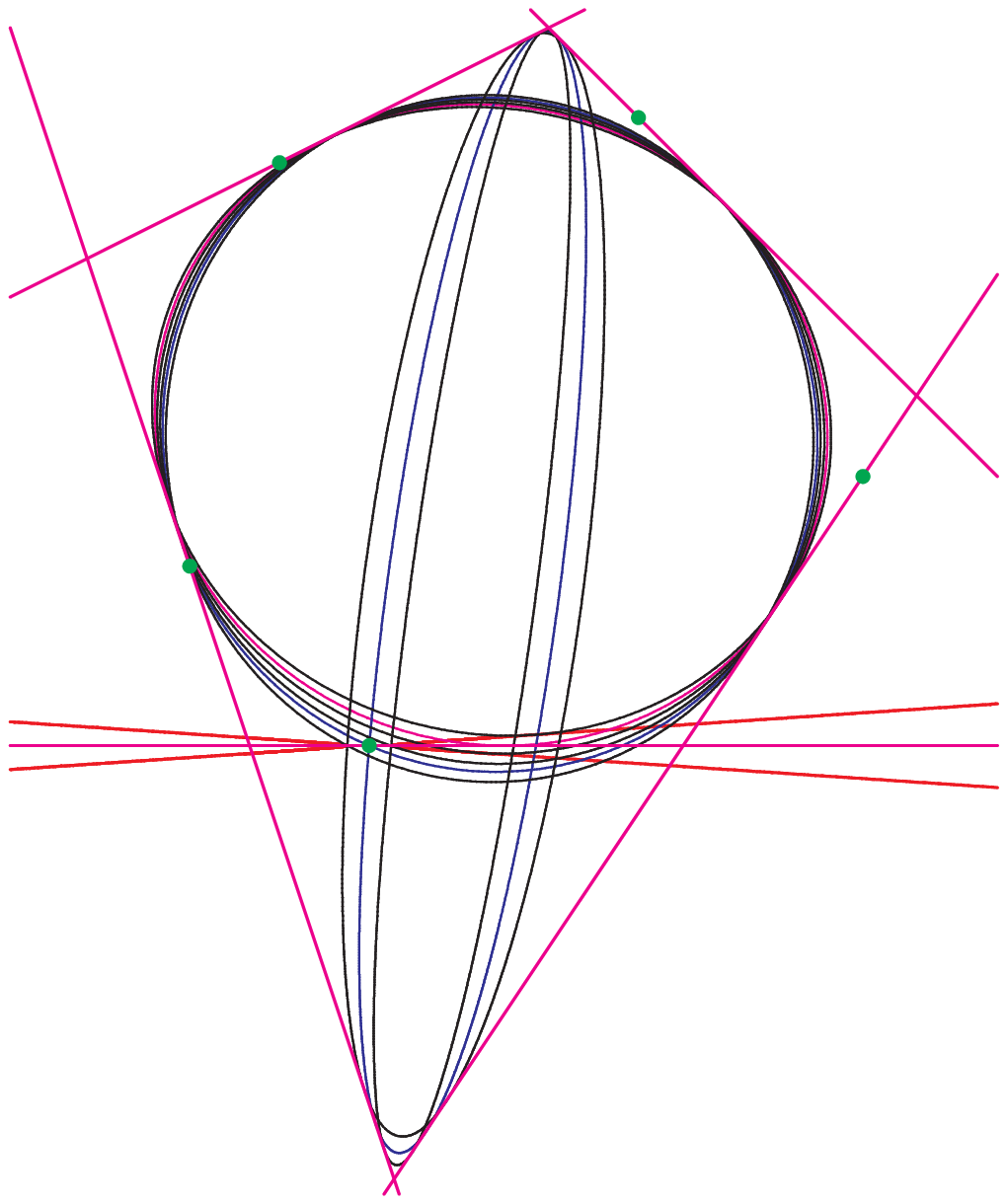}}
   \put(-35,87){$\ell_1$}\put(-23,91.2){\vector(1,0){20}}
   \put(230,87){$h_1$}
    \put(227,92.1){\vector(-4,1){25}} \put(227,90.1){\vector(-4,-1){25}}

   \put(8,60){$p_1$}  \put(19,64){\vector(2,1){45}}

   \put(45,216){$p_4$} \put(131,223){$p_3$}  
   \put(23,122){$p_5$} \put(177,140){$p_2$}
   
   \put(-13,176){$\ell_4$}  \put(195,153){$\ell_3$} 
   \put(52,40){$\ell_5$}    \put(110,40){$\ell_2$} 
  \end{picture}
\]
 and then a closeup near $l_1$.
\[
   \begin{picture}(408,140)(-28,0)
    \put(-28,30){$\ell_1$}\put(-17,33.5){\vector(1,0){15}}

   \put(0,0){\includegraphics[height=80pt]{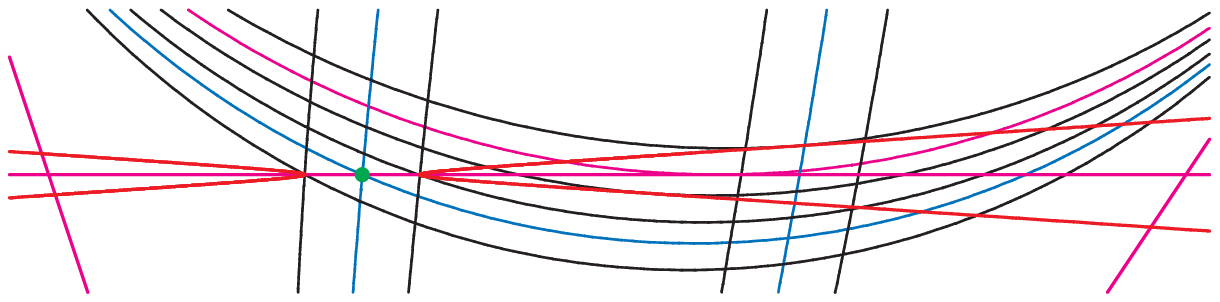}}
   \put(9,125){Conics in $C$}\put(11,114){meeting $p_1$}
    \put(32,109){\line(0,-1){31.5}}    \put(40,109){\line(2,-1){64.7}}

   \put(91,125){Nearby}\put(93,114){conics}
    \put(93,109){\line(-1,-1){48.5}}  \put(103,109){\line(-1,-1){45.4}}  
    \put(112,109){\line(-3,-4){24.1}} \put(120,109){\line(0,-1){48}}

   \put(224,125){Conic in $C$}\put(220,114){tangent to $\ell_1$}
    \put(255,109){\line(0,-1){69.7}}  

   \put(304,125){Nearby}\put(306,114){conics}
    \put(318,109){\line(0,-1){42}}  \put(328,109){\line(0,-1){45.5}}  

   \put(383,29){$h_1$}
    \put(381,35.4){\vector(-3,1){37.5}} \put(381,31.6){\vector(-3,-1){37.5}}
   \end{picture}
\]

\begin{figure}[htb]
\[
    \includegraphics[height=270pt]{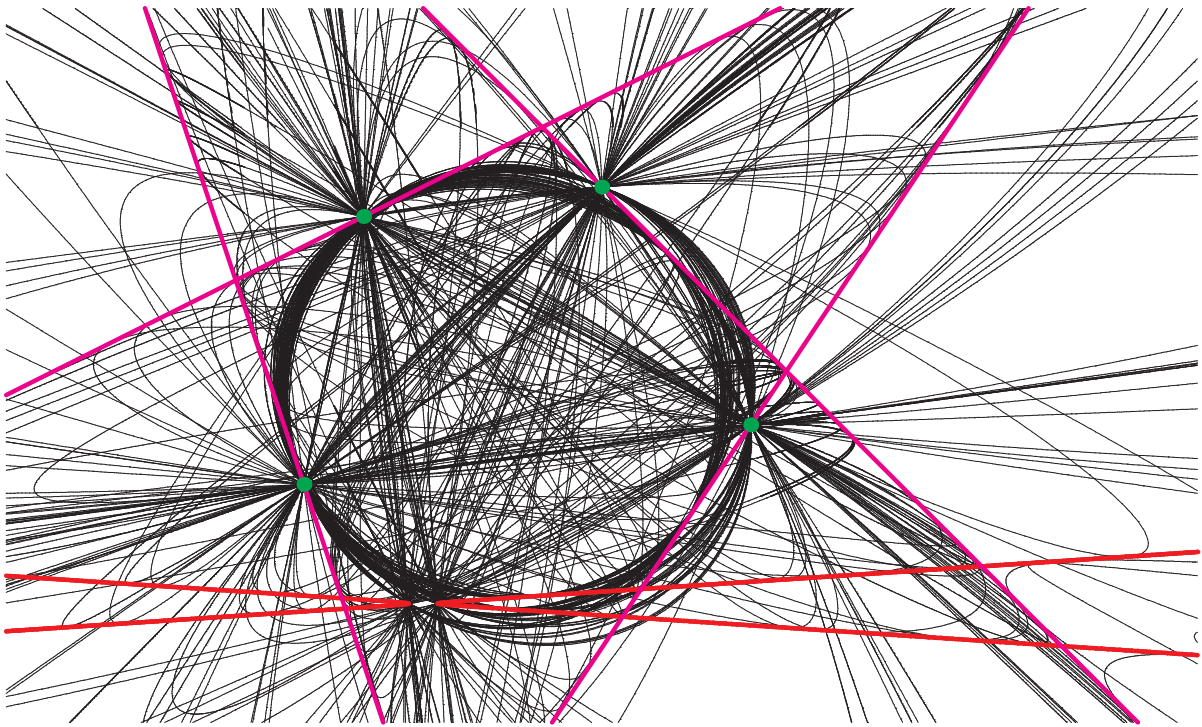}
\]
\caption{204 conics}\label{F7:204_conics}
\end{figure}

For the configuration of 102 conics of figure~\ref{F7:102_conics}, the
hyperbola
\[
   \DeCo{h_1}\ \colon\qquad
     (y\;-\;\tfrac{1}{15}x) (y\;+\;\tfrac{1}{15}x)  \ +\ \tfrac{1}{15000}\ =\ 0\,,
\]
is sufficiently close to its asymptotes, which meet at $p_1$ and are sufficiently close to
$\ell_1$, and so each of our 102 conics 
that meet $p_1$ or are tangent to $\ell_1$ becomes two conics tangent to $h_1$.
We  first show the  configuration of 102 conics
of Figure~\ref{F7:102_conics} in the neighborhood of $\ell_1$,
\[
  \begin{picture}(434,110)
    \put(0,0){\includegraphics[height=110pt]{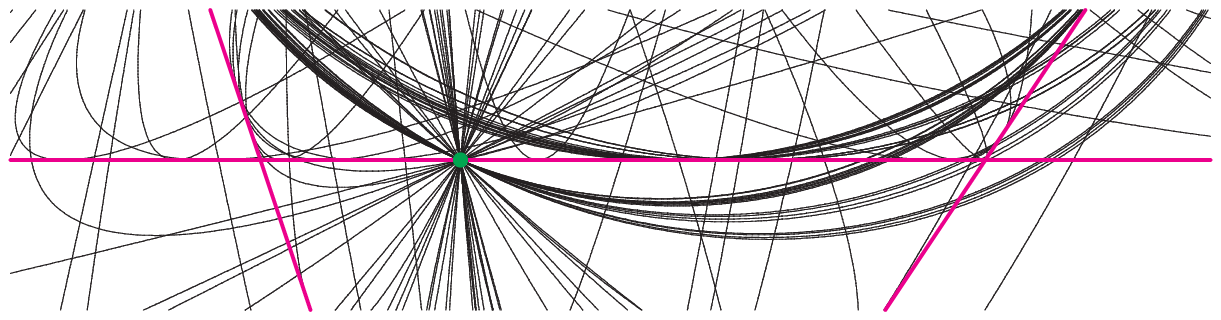}}
    \put(415,10){$\ell_1$}\put(419,21){\vector(0,1){32}}
   \end{picture}
\]
and then the resulting 204 conics in the same region.
\[
    \includegraphics[height=110pt]{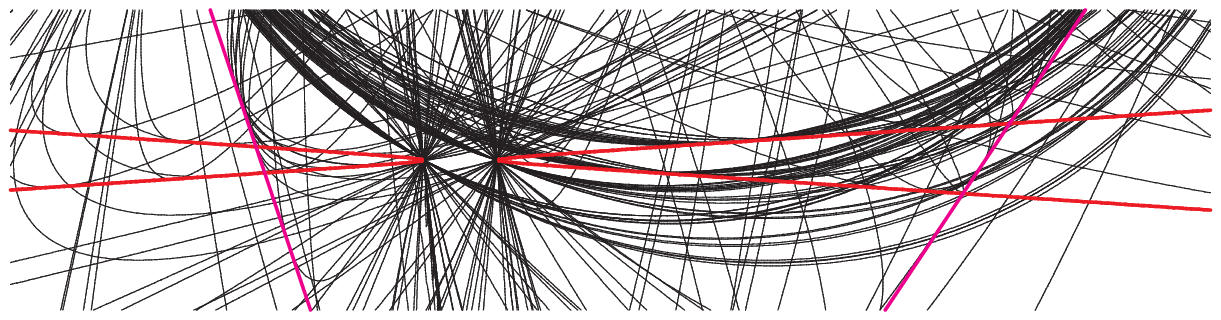}
\]
The key to the existence of this construction was that no tangent direction to a 
conic through $p_1$ met $h_1$, which is possible as no conic was tangent to
$\ell_1$ at $p_1$.
Figure~\ref{F7:204_conics} shows the resulting 204 conics that are tangent to 
$h_1$ and, for each $i=2,3,4,5$ either contain $p_i$ or are tangent to $\ell_i$.

If we now replace $p_2\in\ell_2$ by a similar nearby hyperbola, then the 204 conics
become 408.
Replacing $p_3\in\ell_3$ by a nearby hyperbola, will give 816 conics.
Continuing with $p_4\in\ell_4$ gives 1632, and finally replacing $p_5\in\ell_5$ with 
a hyperbola gives five hyperbolae, $h_1,\dotsc,h_5$ for which there are 
$2^5\cdot 102=3264$ real conics tangent to each $h_i$.
In this way, the classical problem of 3264 conics can have all of its solutions be real.
Observe that this discussion also gives a derivation of the number 3264 without reference
to intersection theory~\cite{Fu}.

%
%
%
\section{Some geometric problems}

We discuss some other geometric problems that can have all their solutions be real.

%
%
%
\subsection{The Stewart-Gough platform}
\begin{figure}[htb]
$$
  \setlength{\unitlength}{1.3pt}
  \begin{picture}(330,220)(-7,0)
   \put(  9.9, 0){\includegraphics[width=319pt]{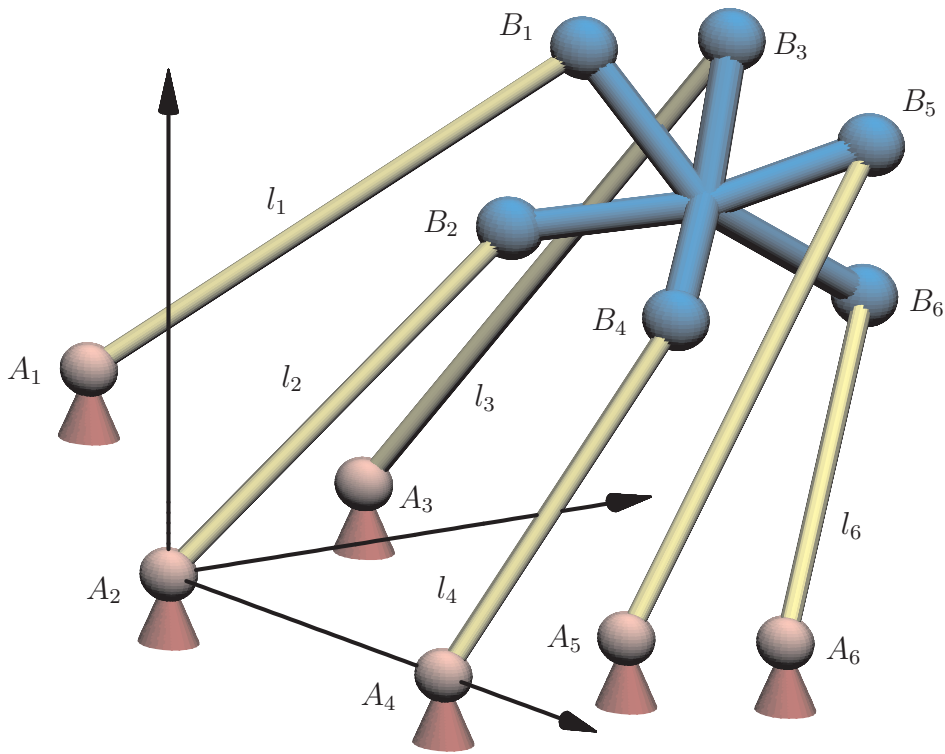}}
   \put( -5,109){$A_1$}   \put( 18,45){$A_2$}
   \put(108.9,71.5){$A_3$}   \put(98,14){$A_4$}
   \put(152,30.8){$A_5$}   \put(233.2,27.5){$A_6$}

   \put(138,209){$B_1$} \put(116,151.8){$B_2$}
   \put(217.8,203.5){$B_3$}   \put(165,123.2){$B_4$}
   \put(255.2,187){$B_5$}   \put(257.4,128.7){$B_6$}

   \put( 70.4,158.4){$l_1$}   \put( 74.8,106.7){$l_2$}
   \put(130.9,101.2){$l_3$}   \put(120, 45.1){$l_4$}
   \put(1204.6, 71.5){$l_5$}   \put(237.6, 63.8){$l_6$}

  \end{picture}
$$
\caption{A Stewart platform\label{fig:stewart}.}
\end{figure}

The position of a rigid body in $\R^3$ has 6 degrees of freedom (3 rotations and
3 translations). 
This is exploited in robotics, giving rise to the \DeCo{{\sl Stewart-Gough
platform}}~\cite{Go57,St65}:
Suppose that we have 6 fixed points $A_1,A_2,\ldots, A_6$ in space
and 6 points $B_1,B_2,\ldots, B_6$ on a rigid body $B$ (the framework of
Figure~\ref{fig:stewart}).
The body is controlled by varying each distance $l_i$ between the
fixed point $A_i$ and the point $B_i$ on $B$.
This may be accomplished by attaching rigid actuators between 
spherical joints located at the points $A_i$ and $B_i$, or by suspending the
platform from a ceiling with cables.
(Apparently, this configuration is often used in factories.)

Given a position of the body $B$, the distances $l_1,l_2,\ldots, l_6$ are
uniquely determined.
A fundamental problem is the inverse question (in robotics, this is called the
forward problem): 
Given a platform (positions of the $A_i$ fixed and the relative positions of
the $B_i$ specified) and a sextuple of distances $l_1,l_2,\ldots,l_6$, what is the
position of the platform?

It had long been understood that several positions were possible for a given
sextuple of lengths.
This led to the following enumerative problem.

\begin{center}
 For a given (or general) Stewart platform, how many (complex) positions are
 there\\  for a
 generic choice of the distances $l_1,l_2,\ldots,l_6$? \ 
 How many of these can be real?
\end{center}

In the early 1990's, several approaches (including a nice interaction between
theory and computer experimentation~\cite{La93,RV95})
showed that there are 40 complex positions of a general Stewart platform.
The obviously practical question of how many positions could be real 
remained open until 1998, when Dietmaier introduced a novel method involving
numerical homotopy to find a platform and value of the distances $l_1,l_2,\ldots,l_6$ with   
all 40 positions real. 

\begin{thm}[Dietmaier~\cite{Di98}]
  \mbox{{\DeCo{All}} $40$ positions can be real.}
\end{thm}

%
%
%
\subsection{Real rational cubics through 8 points in 
            ${\mathbb P}^2_\R$}\label{sec:ratcubics}
Recall from Section~\ref{S1:lower} that there are 12 singular (rational) cubic curves
containing 8 general points in the plane.
Kharlamov studied this over the real numbers.

\begin{thm}[{\cite[Proposition 4.7.3]{DeKh00}}]\label{thm:cubics}
 Given $8$ general points in ${\mathbb P}^2_\R$, there are at least $8$
 real rational cubics containing them, and there are choices of the $8$ points
 for which \DeCo{all}\/ $12$ rational cubics are real.
\end{thm}

The bound is a nice exercise in Euler characteristic.
A homogeneous cubic has 10 coefficients, so the set of plane cubics is naturally
identified with $9$-dimensional projective space.
Let $p_1,\dotsc,p_8$ be general points in $\R\P^2$.
As the condition for a cubic to contain a point $p_i$ is linear in the coefficients of the
cubic, there is a pencil ($\P^1$) of cubics through these 8 points.
Let $P, Q$ be two distinct cubics in this pencil, which is then parametrized by $sP+tQ$
for $[s,t]\in\P^1$. 
By B\'ezout's Theorem, the cubics $P$ and $Q$, and hence every cubic in the pencil, 
vanish at a ninth point, $p_9$. 

It is not hard to see that there is a unique cubic in the
pencil that vanishes at any point $p\in\R\P^2-\{p_1,\dotsc,p_9\}$.
A little harder, but still true, is that there is a unique cubic in the pencil with any
given tangent direction at some point $p_i$.
In this way, we have maps
\[
   \begin{picture}(113,50)(-7,0)
    \put(3,38){$\DeCo{Z}\ :=\ \text{Bl}_{\{p_1,\dotsc,p_9\}}\R\P^2$}
    \put(-7,20){$C$}  \put(7,32){\vector(0,-1){20}}
    \put(70,32){\vector(0,-1){20}}  \put(74,20){$\pi$}
    \put(0,0){$\R\P^1$}   \put(64,0){$\R\P^2$}
   \end{picture}
\] 
where $\text{Bl}_{\{p_1,\dotsc,p_9\}}\R\P^2$ is the blow-up of $\R\P^2$ in the 9 points,
which is obtained by removing each point $p_i$ and replacing it with the tangent directions
$\R\P^1\simeq S^1$ to $\R\P^2$ at $p_i$.
The map $\pi$ is the blow-down, and the map $C$ associates a point of $Z$ to the
unique curve in the pencil which contains that point.

Because $Z$ is a blow-up, we may compute its Euler characteristic to obtain
 \begin{eqnarray*}
   \chi(Z)&=& \chi(\R\P^2)\ -\ 9\cdot\chi(\text{pt})\ +\ 9\cdot\chi(S^1)\\
          &=& 1\ -\ 9\ +\ 0\\
          &=& -8\,.
 \end{eqnarray*}
The key to Theorem~\ref{thm:cubics} is to compute the Euler characteristic of $Z$ a second way
using the map $C\colon Z\to S^1$.  
The fibers of this map are the cubic curves in the pencil.
Smooth real cubics either have one or two topological components, 
\[
    \includegraphics[height=70pt]{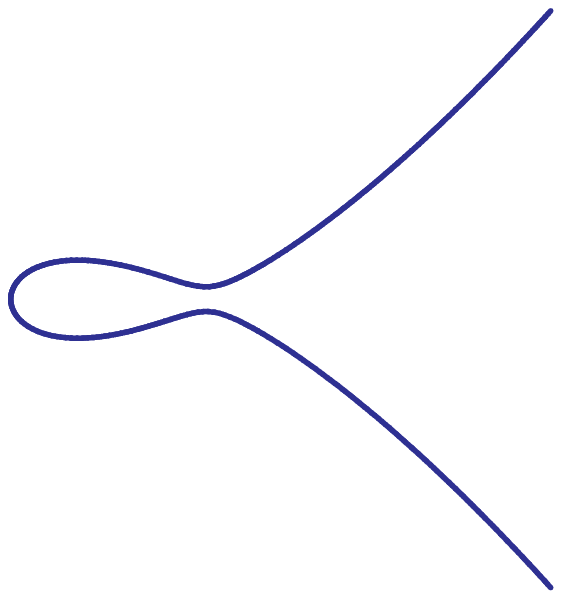}
    \qquad\qquad\qquad
    \includegraphics[height=70pt]{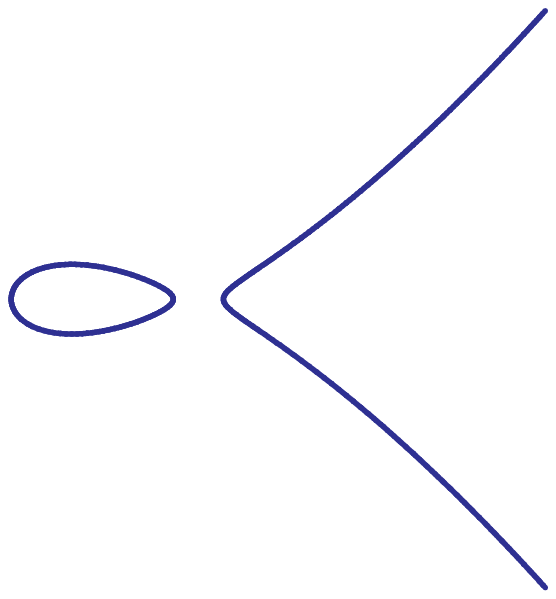}
\]
and hence are homeomorphic to one or two copies of $S^1$.
In either case, their Euler characteristic is zero, and so the
Euler characteristic of $Z$, $-8$, is the sum of the Euler characteriatics of singular fibers
of the map $C$. 

Because the points $p_1,\dotsc,p_8$ were general, the only possible singular fibers are
nodal cubics, and there are two types of real nodal cubics.
 \[
   \includegraphics[height=70pt]{figures/1/renodes.eps}\qquad\qquad\qquad
   \includegraphics[height=70pt]{figures/1/conodes.eps}
 \]
The first is the topological join of two circles and has Euler characteristic $-1$,
while the second is the disjoint union of a circle with a point and therefore has Euler
characteristic $1$.
Thus, if $n$ is the number of nodal cubics in the pencil and $s$ is the number with a
solitary point, then we have $n+s\leq 12$ (as there are 12 complex rational cubics), and 
$s-n=-8$.

There are three solutions to this system, $(n,s)\in\{(8,0),(9,1),(10,2)\}$.
Thus there are at least $8$ real cubics through the $8$ points, and this is the derivation
of the Welschinger number $W_3=8$.
Moreover, if there are 2 cubics in the pencil with solitary points, then all 12 rational
cubics will be real.
Such a pencil is generated by the two cubics given below, whose curves we also display.
\[
   \begin{picture}(300,100)
    \put(0,0){\includegraphics[height=100pt]{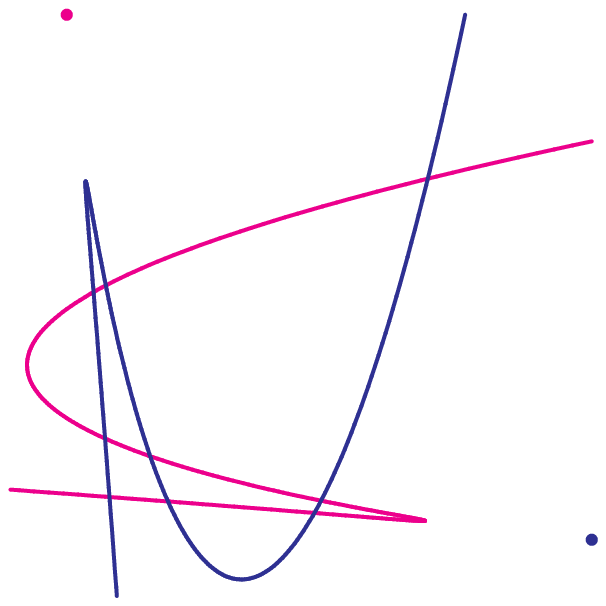}}
    \put(123,74){\Magenta{$(y+9x-28)^2\ =\ 4x^2(x-1)$}}
    \put(118,76){\line(-1,0){17}}    \put(118,79){\line(-6,1){104}}
    \put(92,40){\Blue{$(x+9y-28)^2\ =\ 4y^2(y-1)$}}
    \put(87,43){\line(-1,0){20}}  \put(97,33){\line(0,-1){20}}
   \end{picture}
\]
We conclude that there will be 12 real rational cubics interpolating any subset of 8 of
the 9 points where these two cubics meet.

The question of how many of the $N_d$~\eqref{E1:Konts} rational curves of
degree $d$ which interpolate $3d{-}1$ points in $\R\P^2$ can be real remains
open.
(This was posed in~\cite{So97c}.)

\subsection{Common tangent lines to $2n{-}2$ spheres in 
 $\R^n$}\label{sec:tangent} 
\begin{figure}[htb]
$$
  \makebox{\includegraphics[height=3.6in]{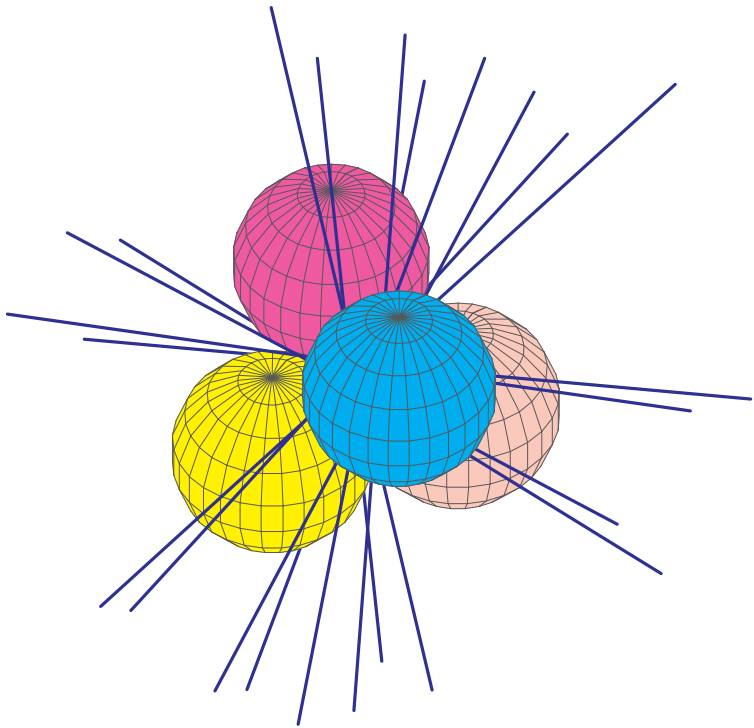}}
$$
\caption{Four spheres with 12 common tangents\label{F7:12lines}.\vspace{-10pt}}
\end{figure}

 How many common tangent lines are there to $2n{-}2$ spheres in $\R^n$?
For example, when $n=3$, how many common tangent lines are there to four
spheres in $\R^3$?
(The number $2n{-}2$ is the dimension of the space
of lines in $\R^n$ and is necessary for there to be finitely many
common tangents.)
Despite its simplicity, this question does not seem to have been asked
classically, but rather arose in computational
geometry.
Macdonald, Pach, and Theobald~\cite{MPT01} gave an elementary argument that four
spheres with the same radius in $\R^3$ can have at most 12 common tangents.
Then they considered the symmetric configuration where the spheres are centered at the
vertices of a regular tetrahedron.
If the spheres overlap pairwise, but no three have a common point, then there will be
exactly 12 common real tangents, as illustrated in Figure~\ref{F7:12lines}.

The general case was established soon after~\cite{STh02}.

\begin{thm}\label{thm:12lines}
 $2n-2$ general spheres in $\R^n$ $(n\geq 3)$ have $3\cdot 2^{n-1}$
 complex common tangent lines, and there are $2n-2$ such spheres with \DeCo{all}
 common tangent lines real.
\end{thm}

The same elementary arguments of Macdonald, Pach, and Theobald give a bound valid for all
$n$ and for spheres of any radius, and a generalization of the symmetric configuration of
Figure~\ref{F7:12lines}  gives a configuration of $2n-2$ spheres having $3\cdot 2^{n-1}$
common real tangents.

\begin{figure}[htb]
$$
  \includegraphics[height=2.4in]{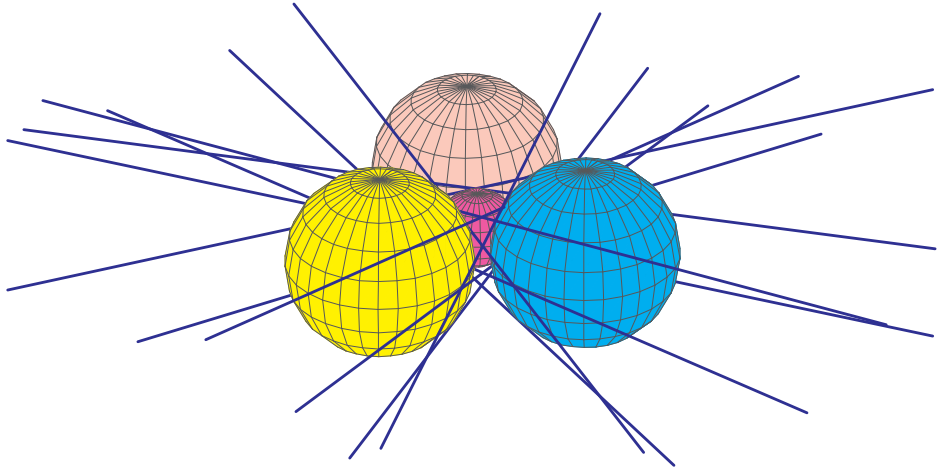}
$$
\caption{Four spheres with coplanar centers and 12 common tangents\label{F7:megyesi}.}
\end{figure}

\begin{figure}[htb]
$$
  \includegraphics[height=1.7in]{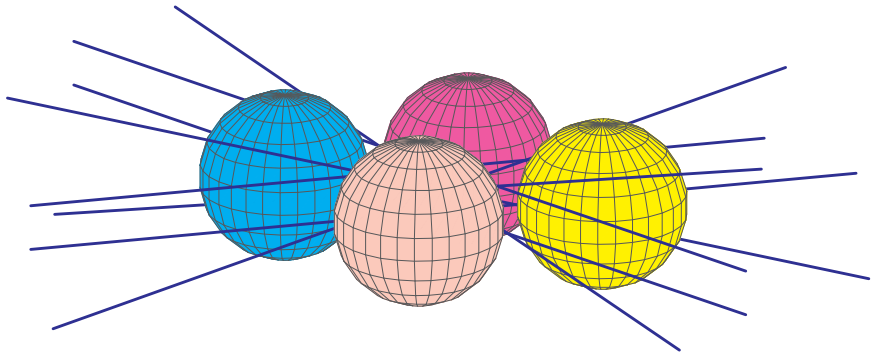}
$$
\caption{Four equal spheres with coplanar centers and 8 common tangents\label{F7:eight}.}
\end{figure}

Megyesi~\cite{Me01} showed that this result for $n=3$ remains true if the
spheres have coplanar centers (Figure~\ref{F7:megyesi}), 
but that there can
only be 8 common real tangents (out of 12 complex ones) if the spheres have the
same radii (Figure~\ref{F7:eight}).
 
The spheres in Figures~\ref{F7:12lines} and~\ref{F7:megyesi} are not disjoint, in fact
their union is connected.
Fulton asked if it were possible for 4 \Red{disjoint} spheres to have 12 common real
tangents.
It turns out that a perturbation of the configuration of Figure~\ref{F7:megyesi}
gives four pairwise disjoint spheres with 12 common tangents, as we show 
below.
$$
  \includegraphics[height=2.6in]{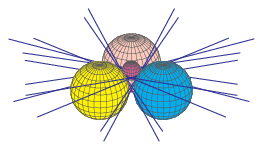}
$$
The three large spheres have radius $4/5$ and are centered at the vertices of an
equilateral triangle of side length $\sqrt{3}$, while the smaller sphere has
radius $1/4$ and is centered on the axis of symmetry of the triangle, but at a
distance of $35/100$ from the plane of the triangle.
It remains an open question whether it is possible for four disjoint
\Red{{\sl unit}} spheres to have 12 common tangents.

\section{Schubert calculus}\label{S7:SC}
The largest class of problems which have been studied from the perspective of having all
solutions real come
from the classical Schubert calculus of enumerative geometry, which involves
linear spaces meeting other linear spaces.
The simplest nontrivial example illustrates some of the vivid geometry behind
this class of problems.
Consider the following question:
\begin{center}
 How many line transversals are there to four given lines in space?
\end{center}

\begin{figure}[htb]
$$
  \includegraphics[height=2in]{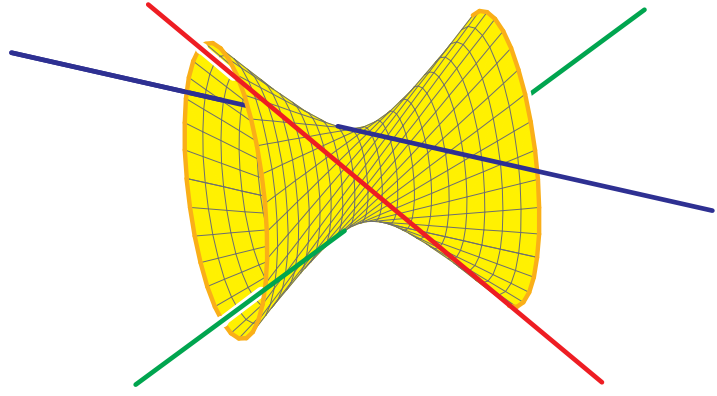}
$$
\caption{Hyperboloid containing three lines\label{F7:4lines}.}
\end{figure}

To answer this, first consider \DeCo{three} lines.
They lie on a unique hyperboloid.
(See Figure~\ref{F7:4lines}.)
This hyperboloid has two rulings by lines.
The three lines are in one ruling, and the other ruling (which is drawn on the
hyperboloid in Figure~\ref{F7:4lines}) 
consists of the lines which meet the three given lines.

The fourth line will meet the hyperboloid in two points (the hyperboloid is
defined by a quadratic polynomial).
Through each point of intersection there will be one line in the second family,
and that line will meet our four given lines.
In this way, we see that the answer to the question is \DeCo{2}.
Note that the fourth line may be drawn so that it meets the hyperboloid in
two real points, and both solution lines will be real when this
happens.\bigskip 

Let $\Gr(p,m{+}p)$ be the \DeCo{{\sl Grassmannian}} of $p$-dimensional linear subspaces 
(\DeCo{{\sl $p$-planes}}) in an $(m{+}p)$-dimensional vector space,
which is an algebraic manifold of dimension $mp$.
The Schubert calculus involves fairly general incidence conditions imposed on
$p$-planes $H$.
These general conditions are imposed by \DeCo{{\sl flags}}, which are sequences
of linear subspaces, one contained in the next.
More specifically, a flag \DeCo{{$\Fdot$}} is a sequence
\[
  \Fdot\ \colon\ 
   F_1\ \subset\ F_2\ \subset\ \dotsb\ \subset\ F_{m+p-1}\ \subset\ F_{m+p}\,,
\]
where $F_i$ is a linear subspace having dimension $i$.

These general conditions are called \DeCo{{\sl Schubert conditions}} and 
are indexed by sequences $\alpha\colon 1\leq\alpha_1<\alpha_2<\dotsb<\alpha_p<m{+}p$ of  
integers.
Write \DeCo{$\binom{[m+p]}{p}$} for the set of all such sequences.
The set of all $p$-planes satisfying the condition $\alpha$ imposed by the
flag $\Fdot$ is a \DeCo{{\sl Schubert variety}}, defined by
\[
  \DeCo{X_\alpha\Fdot}\ :=\ 
  \{ H\in\Gr(p,m{+}p)\ \mid\ \dim H\cap F_{\alpha_j}\geq j
      \quad\mbox{\rm for}\quad j=1,\dotsc,p\}\,.
\]
This is a subvariety of the Grassmannian of dimension
$\DeCo{|\alpha|}:=\alpha_1-1+\alpha_2-2+\dotsb+\alpha_p-p$.
Its codimension is $mp-|\alpha|$.

If $K$ is a linear subspace of codimension $p$ and $\Fdot$ is a flag
with $F_{m}=K$, then 
\[
   \DeCo{X_{\I}\Fdot}\ :=\ X_{(m,m+2,\dotsc,m+p)}\Fdot
   \ =\ \{H\ \mid\ H\cap K\neq\{0\}\}\,.
\]
The condition that $H$ meets $K$ is is called a \DeCo{{\sl simple Schubert
    condition}} and  
$X_{\I}\Fdot$ is a \DeCo{{\sl simple Schubert variety}}.
More generally, a \DeCo{{\sl special Schubert variety}} of codimension $a$ is
\[
  \DeCo{X_a\Fdot}\ :=\ X_{(m+1-a,m+2,\dotsc,m+p)}\Fdot
  \ =\ \{H\ \mid\ H\cap F_{m+1-a}\neq\{0\}\}\,.
\]

\begin{thm}\label{T7:real}
 Given Schubert conditions $\alpha^1,\alpha^2,\dotsc,\alpha^n$ in $\binom{[m+p]}{p}$ with 
 $\sum_i(mp-|\alpha^i|)=mp$, there exist real flags 
 $\Fdot^1,\Fdot^2,\dotsc,\Fdot^n$ such that the intersection
 \begin{equation}\label{Eq7:Schubert_intersection}
   \bigcap_{i=1}^n X_{\alpha^i}\Fdot^i
 \end{equation}
 is transverse with \DeCo{all} points real.
\end{thm}

The numerical condition that $\sum_i(mp-|\alpha^i|)=mp$ implies that the expected
dimension of the intersection is zero.
If the flags $\Fdot^i$ are in general position, then the Kleiman-Bertini
Theorem~\cite{KL74} implies that the intersection is in fact 
transverse and zero-dimensional, so there will be 
finitely many complex $p$-planes satisfying the incidence conditions $\alpha^i$ imposed
by general flags, and these $p$-planes are exactly the points in the
intersection. 

Theorem~\ref{T7:real} was proved in several stages.
First, when $p=2$~\cite{So97a}, and then for any $p$, but only for 
special Schubert conditions~\cite{So99}, and then finally for general Schubert
conditions by Vakil~\cite{Va06}.

\section{Quantum Schubert calculus}

Given points $s_1,s_2,\dotsc,s_{d(m+p)+mp}\in\P^1$ and 
$m$-planes $K_1,K_2,\dotsc,K_{d(m+p)+mp}$ in $\C^{m{+}p}$, there are finitely
many rational curves $\gamma\colon \P^1\to\Gr(p,m{+}p)$ of degree $d$ so that 
 \begin{equation}\label{E5:QSchCal}
   \gamma(s_i)\ \cap\ K_i\ \neq\ \{0\}\qquad 
    i=1,2,\dotsc,d(m{+}p)+mp\,.
 \end{equation}
These are \DeCo{{\sl simple quantum Schubert conditions}}.
More generally, one could (but we will not) impose the condition that the
$p$-plane $\gamma(s_i)$ lie in some predetermined Schubert variety.
The number of solutions to such problems are certain Gromov-Witten invariants of the
Grassmannian, and may be computed by the quantum Schubert
calculus~\cite{Be97,In91,ST97,Va92}. 

\begin{thm}[\cite{So00a}]\label{T7:quantum}
 There exist real points $s_1,s_2,\dotsc,s_{d(m+p)+mp}\in\R\P^1$ and real 
 $m$-planes $K_1,K_2,\dotsc,K_{d(m+p)+mp}$ in $\R^{m{+}p}$ so that \DeCo{every}
 rational curve  $\gamma\colon \P^1\to\Gr(p,m{+}P)$ of degree $d$
 satisfying~$\eqref{E5:QSchCal}$ is real.
\end{thm}

\section{Theorem of Mukhin, Tarasov, and Varchenko}

In May of 1995, Boris Shapiro communicated to the author a remarkable conjecture that
he and his brother Michael had made concerning reality in the Schubert calculus.
They conjectured that there would only be real points in a zero-dimensional intersection
of Schubert varieties given by flags osculating the rational normal curve.
Subsequent computation~\cite{RoSo98,So00b} gave strong evidence for the conjecture and
revealed that the intersection should be transverse.
Partial results were obtained~\cite{So99,EG02}, and the full conjecture was proven by
Mukhin, Tarasov, and Varchenko~\cite{MTV_Sh}.
They later gave a second proof~\cite{MTV_R}, which different from their original proof and
gave a proof of transversality.

This \DeCo{{\sl Shapiro Conjecture}} has been a motivating conjecture for the study of
reality in the Schubert calculus with several interesting (and as-yet-unproven)
generalizations that we will discuss in subsequent chapters.
Let $\gamma$ be the rational normal (or moment) curve in $\C^{m+p}$, which we
will take to be the image of the map
\[
  \gamma(t)\ =\ (1,\ t,\ t^2,\ \dotsc,\ t^{m+p-1})\ \in\ \C^{m+p}\,,
\]
defined for $t\in\C$.
Given a point $t\in\C$, the \DeCo{{\sl osculating flag}}
\DeCo{$\Fdot(t)$} is the flag of subspaces whose $i$-plane is the linear span of 
the first $i$ derivatives of $\gamma$, evaluated at $t$
\[
  F_i(t)\ :=\ \Span\{\gamma(t),\ \gamma'(t),\ 
    \gamma''(t),\ \dotsc,\ \gamma^{(i-1)}(t)\}\,.
\]
This flag makes sense for $t\in\P^1$.
Here is the strongest form of the Shapiro conjecture that has been proven~\cite{MTV_R}.

\begin{thm}\label{Th7:MTV}
 If $\alpha^1,\alpha^2,\dotsc,\alpha^n\in\binom{[m+p]}{p}$ are Schubert conditions satisfying
 $\sum_i(mp-|\alpha^i|)=mp$, then, for every choice of $n$ distinct
 points $s_1,s_2,\dotsc,s_n\in\R\P^1$, the intersection 
\[
   \bigcap_{i=1}^n X_{\alpha^i}\Fdot(s_i)
\]
 is transverse with \DeCo{all} points real.
\end{thm}

It is instructive to consider this for the problem of two
lines that we saw in Section~\ref{S7:SC}.
An osculating line will be a tangent line to the rational normal curve $\gamma$.
The Shapiro Conjecture asserts that given four lines that are
tangent to the rational normal curve at real points, there will be
two lines meeting all four, and the two lines will be real. 
As before, first consider three lines tangent to the rational normal curve.
They will lie in one ruling of a quadric surface, and the other ruling consists
of the lines meeting all four.

\begin{figure}[htb]
 \[
  \begin{picture}(372,172.2)(0,6)
   \put(0,0){\includegraphics[height=190pt]{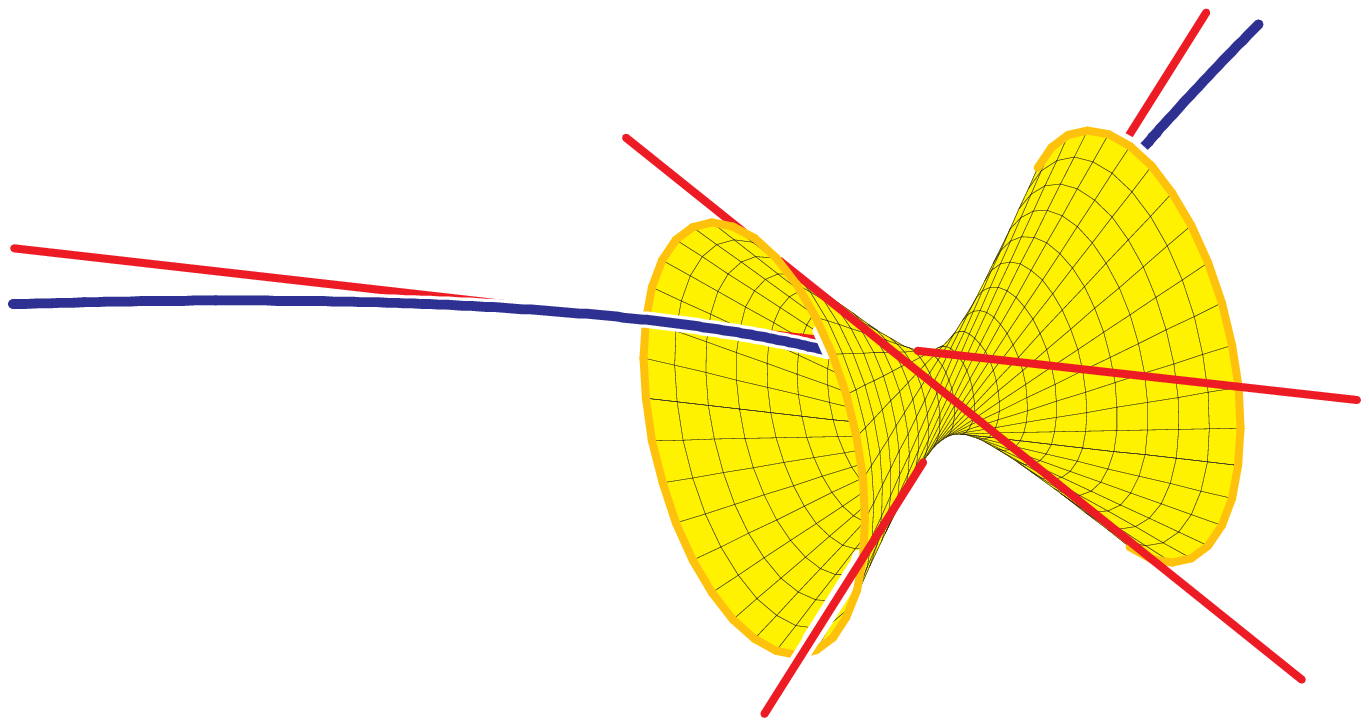}}
   \put(171,0){\Red{$\ell(-1)$}} \put(354,0){\Red{$\ell(0)$}}
   \put(0,138.2){\Red{$\ell(1)$}} \put(0,92){\DeCo{$\gamma$}}
   \put(324,138.2){\Brown{$Q$}}
  \end{picture}
 \]
 \caption{Quadric containing three lines tangent to the rational normal
 curve\label{F5:TanQuad}.} 
 \end{figure}

For example, if we let 
$\gamma(t)=(6t^2-1, \frac{7}{2}t^3+\frac{3}{2}t, -\frac{1}{2}t^3+\frac{3}{2}t)$, 
and consider tangent lines $\ell(-1)$, $\ell(0)$, and $\ell(1)$ to $\gamma$,
\[
   \ell(\pm 1)\ :\ (5,\pm 5,1) + u(\pm 1,1,0) \qquad
   \ell(0)\ :\ (-1, 0, 0) + u(0, 1,1)\qquad u\in\R\,,
\]
then the quadric is $x^2-y^2+z^2=1$.
We display this in Figure~\ref{F5:TanQuad}, together with the rational normal
curve. 
Since any three real points on any real rational normal curve may be carried to 
any three real points on any other real rational normal curve by an automorphism
of projective 3-space, we can assume
that we are in the situation of Figure~\ref{F5:TanQuad}.

The rational normal curve $\gamma$ (a cubic) meets the quadric tangentially at
the points where it is tangent to the lines.  
Since the total multiplicity of its intersection with the quadric is six, these
are its only points of contact.
In particular, $\gamma$ always lies on one side of the quadric, looping around
on the inside.
In this way we can see that a fourth tangent line to the rational normal curve must meet
the quadric in two 
real points, giving the two real transversals asserted by Shapiro's Conjecture. 

\begin{figure}[htb]
 \[
  \begin{picture}(332, 182)(-23,-5)
   \put(0,0){\includegraphics[height=160pt]{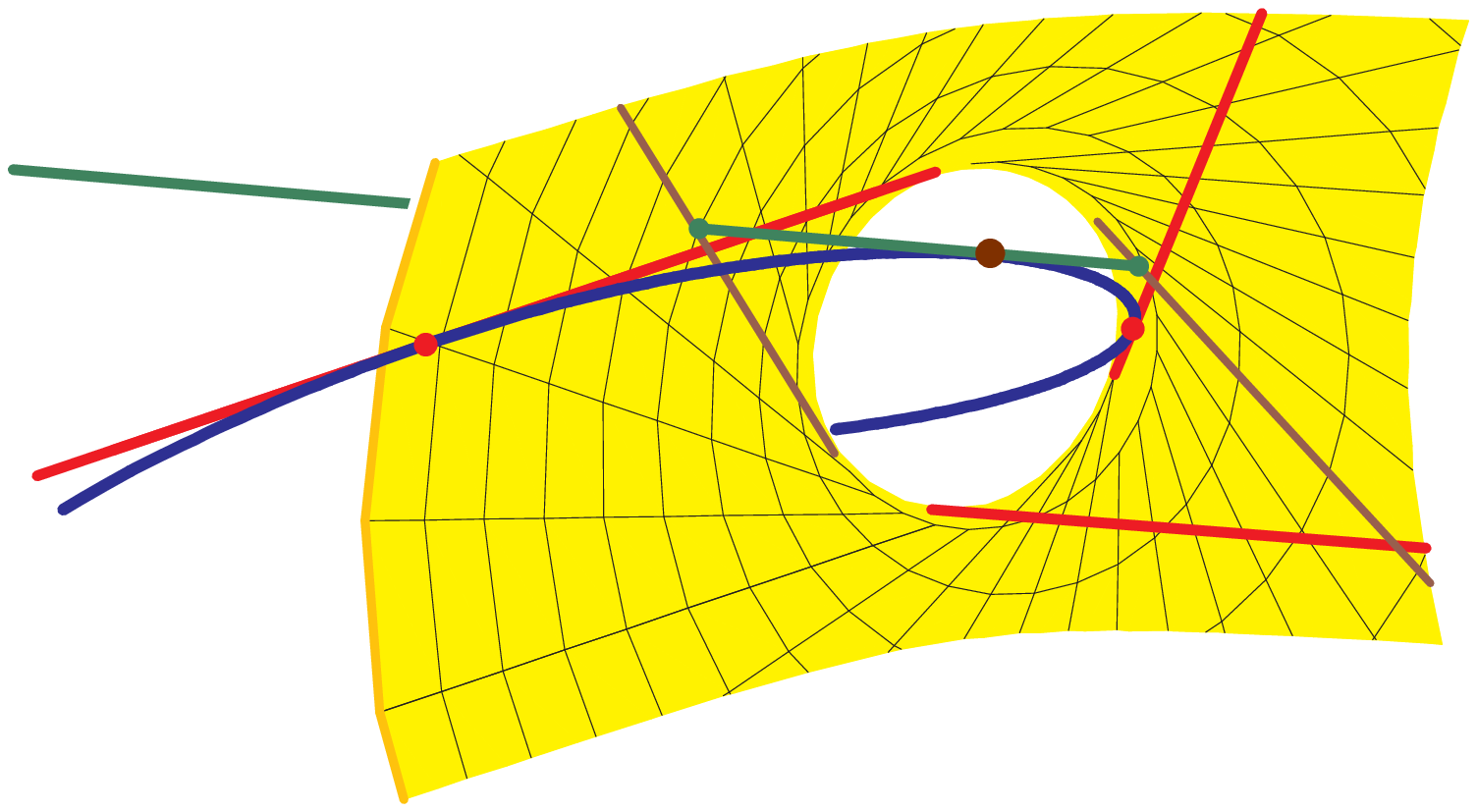} }
   \put(-30,63.5){\Red{$\ell(-1)$}} \put(286,48){\Red{$\ell(1)$}} 
   \put(244,166){\Red{$\ell(0)$}}
   \put(-2,47){\DeCo{$\gamma$}} \put(-8,135){\ForestGreen{$\ell(v)$}}
  \end{picture}
\] 
\caption{Configuration in the throat of quadric\label{F5:throat}.}
 \end{figure}

It is no loss to assume that the fourth tangent line is tangent at
some point \Brown{$\gamma(v)$} on the arc between points 
\Red{$\gamma(-1)$} and \Red{$\gamma(0)$}.
As illustrated in Figure~\ref{F5:throat}, the tangent line 
\ForestGreen{$\ell(v)$} does indeed
meet the hyperboloid in two real points.
%

%
%
%
%
\chapter{The Shapiro Conjecture for Grassmannians}\label{Ch:ScGr}


In Chapter~\ref{S:ERAG}, we considered the rational normal curve $\gamma$ in $\C^{m{+}p}$,
which we took to be the image of the map
\[
  \gamma(t)\ =\ (1,\ t,\ t^2,\ \dotsc,\ t^{m{+}p-1})\ \in\ \C^{m{+}p}\,,
\]
defined for $t\in\C$.
For a point $t\in\C$, the osculating flag
$\Fdot(t)$ is the flag of subspaces whose $i$-plane is the linear span of 
$\gamma(t)$ and the the first $i{-}1$ derivatives of $\gamma$, evaluated at $t$,
\[
  F_i(t)\ :=\ \Span\{\gamma(t),\ \gamma'(t),\ 
    \gamma''(t),\ \dotsc,\ \gamma^{(i-1)}(t)\}\,.
\]
We may also define $\Fdot(\infty)$ to be the limit as $s\to 0$ of $\Fdot(\frac{1}{s})$, to
get a family of flags $\Fdot(t)$ for $t\in\P^1$.

We work in the Grassmannian $\Gr(p,m{+}p)$ of $p$-planes in $\C^{m{+}p}$.
Recall from Section~\ref{S7:SC} that Schubert varieties are indexed by Schubert
conditions $\alpha\in\binom{[m{+}p]}{p}$, which are integer sequences  
$\alpha\colon 1\leq\alpha_1<\alpha_2<\dotsb<\alpha_p\leq m{+}p$.
Given a flag $\Fdot$, the corresponding Schubert variety is
\[
   X_\alpha\Fdot\ =\ 
  \{ H\in\Gr(p,m{+}p)\ \mid\ \dim H\cap F_{\alpha_j}\geq j
      \quad\mbox{\rm for}\quad j=1,\dotsc,p\}\,.
\]
This has dimension $\DeCo{|\alpha|}:=\sum_j (\alpha_j-j)$.

We investigate the Shapiro conjecture for Grassmannians (Theorem of Mukhin, Tarasov, and
Varchenko).
\medskip

\noindent{\bf  Theorem~\ref{Th7:MTV}}\ 
{\it 
 If $\alpha^1,\alpha^2,\dotsc,\alpha^n$ are Schubert conditions with 
 $\sum_i(mp-|\alpha^i|)=mp$, then, for every choice of $n$ distinct
 points $s_1,s_2,\dotsc,s_n\in\R\P^1$, the intersection 
\[
   \bigcap_{i=1}^n X_{\alpha^i}\Fdot(s_i)
\]
 is transverse with \DeCo{all} points real.
}

In particular, we show that a special case of this theorem is equivalent to the statement
(Theorem~\ref{T1:Shap_Conj}) of the Shapiro conjecture from the Introduction, and then
prove this special case in an asymptotic sense.

\section{The Wronski map}\label{S8:Wronski}

When all the Schubert conditions are simple (so that $\alpha=\bI=m,m+2,\dotsc,m{+}p$ and
$|\alpha|=mp-1$), the Shapiro conjecture has another formulation in terms of the Wronski
map. 
The \DeCo{{\sl Wronskian}} of a list $f_1(t),f_2(t),\dotsc,f_m(t)$ of polynomials of
degree $m{+}p{-}1$ is the determinant
\[
   \DeCo{\Wr(f_1,f_2,\dotsc,f_m)}\ := \ 
   \det\left( \left(\frac{d}{dt}\right)^{i-1} f_j(t)\right)
   _{i,j=1,\dotsc,m}\ ,
\]
which is a polynomial of degree $mp$, when the polynomials
$f_1(t),f_2(t),\dotsc,f_m(t)$ are generic among all polynomials of degree $m{+}p{-}1$.

Up to a scalar factor, this Wronskian depends only upon the linear span of
the polynomials $f_1(t),f_2(t),\dotsc,f_m(t)$.
Removing these ambiguities, we get the \DeCo{{\sl Wronski map}}
 \begin{equation}\label{E8:Wronski_map}
  \DeCo{\Wr}\ \colon\ \Gr(m,m{+}p)\ \longrightarrow\ \P^{mp}\ ,
 \end{equation}
where $\Gr(m,p)$ is the Grassmannian of $m$-planes in the space of
polynomials of degree $m{+}p{-}1$, and $\P^{mp}$ is the space of polynomials of
degree $mp$, modulo scalars.

\subsection{Some linear algebra}
Let us begin with the moment (rational normal) curve.
For $t\in\C$, set
\[
   \gamma(t)\ =\ (1,t,t^2,\dotsc,t^{m{+}p-1})\ \in\ \C^{m{+}p}\ .
\]
Let $\DeCo{\Gamma}=\Gamma(t)\colon\C^m\to\C^{m{+}p}$ be the map such that
 \begin{equation}\label{E8:Gamma}
   \Gamma(\e_i)\ =\ \gamma^{(i-1)}(t)\,,
 \end{equation}
the ($i{-}1$)-th derivative of $\gamma$.
(We will take $\e_1,\e_2,\dotsc$ to be the standard basis vectors of
the vector space in which we are working.)

A polynomial $f$ corresponds to a linear form (also written $f$):
\[
  f\ \colon\ \C^{m{+}p}\ \longrightarrow\ \C
  \qquad\mbox{so that}\qquad
  f\circ\gamma(t)\ =\ f(t)\ .
\]
The matrix in the definition of the Wronskian is the matrix of the composition
 \begin{equation}\label{E8:Wronski_Comp}
   \C^m\ \xrightarrow{\ \Gamma(t)\ }\ 
   \C^{m{+}p}\ \xrightarrow{\ \Psi\ }\ \C^m\,,
 \end{equation}
where the columns of $\Psi$ are the linear forms defining the polynomials 
$f_1,f_2,\dotsc,f_m$.
Now suppose that $H$ is a basis for the kernel of the map $\Psi$.
We consider $H$ to be a map
$\C^{p}\ \xrightarrow{\ H\ }\ \C^{m{+}p}$,
which we sum with $\Gamma(t)$ to get a map
\[
  \C^{m{+}p}\ =\ \C^m\oplus\C^{p}\ 
   \xrightarrow{\ [\Gamma(t):H]\ }\ \C^{m{+}p}\,.
\]
This map is invertible if and only if the composition~\eqref{E8:Wronski_Comp}
is invertible.
Thus, up to a constant, we have
 \begin{equation}\label{E8:Wronski_Concat}
  \Wr(f_1,f_2,\dotsc,f_m)\ =\ \det[\Gamma(t)\;:\;H]\,,
 \end{equation}
as both are polynomials of the same degree with the same roots.
(Strictly speaking, we need the Wronskian to have distinct roots for this
argument.
The general case follows via a limiting argument.)

We obtain a useful formula for the Wronskian when we expand the
determinant~\eqref{E8:Wronski_Concat} along the columns of $\Gamma(t)$
 \begin{equation}\label{E8:Laplace}
   \Wr(f_1,f_2,\dotsc,f_m)\ =\ 
   \det[\Gamma(t)\;:\;H]\ =\ \sum_{\alpha} 
    (-1)^{|\alpha|}p_{\alpha}(\Gamma(t))\cdot p_{\alpha^c}(H)\,.
 \end{equation}
Here, the sum is over all $\alpha\in\binom{[m{+}p]}{m}$, which are
choices $\{\alpha_1<\alpha_2<\dotsb<\alpha_m\}\subset[m{+}p]$ of 
$m$ distinct rows of the matrix $\Gamma(t)$.
Also, $\alpha^c=[m{+}p]-\alpha$ are the complimentary rows of $H$, and 
$p_\alpha(\Gamma(t))$ is the $\alpha$th maximal minor of $\Gamma(t)$, which is the
determinant of the submatrix of $\Gamma(t)$ formed by the rows in $\alpha$,
and similarly for $p_{\alpha^c}(H)$.
Observe that $\alpha^c\in\binom{[m{+}p]}{p}$.\smallskip

There is a similar expansion for the Wronskian using the
composition~\eqref{E8:Wronski_Comp}.
Take the top  exterior power ($\wedge^m$) of this composition, 
\[
   \C\ =\ \wedge^m\C^m\ \xrightarrow{\ \wedge^m\Gamma(t)\ }\ 
   \wedge^m\C^{m{+}p}\ \xrightarrow{\ \wedge^m\Psi\ }\ 
    \wedge^m\C^m\ =\ \C\,,
\]
where we have used the ordered basis of $\C^m$ so that 
$\wedge^m\C^m=\C\cdot\e_1\wedge\e_2\wedge\dotsc\wedge\e_m$, which identifies
$\wedge^m\C^m$ with $\C$.
Then $\wedge^m\Gamma(t)$ is simply a vector in $\wedge^m\C^{m{+}p}$ and 
$\wedge^m\Psi$ is a covector for $\wedge^m\C^{m{+}p}$.
If we use the basis 
$\e_\alpha:=\e_{\alpha_1}\wedge\e_{\alpha_2}\wedge\dotsc\wedge\e_{\alpha_m}$
for $\wedge^m\C^{m{+}p}$, where 
$\alpha\in\binom{[m{+}p]}{m}$, then
we see that the Wronskian has the form
\[
   \Wr(f_1,f_2,\dotsc,f_m)\ =\ 
   \sum_{\alpha\in\binom{[m{+}p]}{m}} 
   p_\alpha(\Gamma(t))\cdot p_\alpha(\Psi)\,.
\]
Here, $p_\alpha(\Gamma(t))$ and $p_\alpha(\Psi)$ are the $\alpha$th coordinates of
the corresponding vector/covector
(which are the $\alpha$th maximal minors of the corresponding matrices).
Equating these two expressions for the Wronskian gives the equality
(again up to a constant)
\[
    \sum_{\alpha\in\binom{[m{+}p]}{m}} 
    (-1)^{|\alpha|}p_{\alpha}(\Gamma(t))\cdot p_{\alpha^c}(H)
   \ =\ 
   \sum_{\alpha\in\binom{[m{+}p]}{m}} 
   p_\alpha(\Gamma(t))\cdot p_\alpha(\Psi)\,.
\]

 This argument does not use much about the matrix $\Gamma(t)$, besides that it
 depends upon an indeterminate $t$.
 Replacing $\Gamma(t)$ by a matrix of indeterminates proves an interesting (and well-known)
 matrix identity.

\begin{prop}\label{P6:matrix}
  Suppose that $\Psi$ and $H$ are matrices of format $n$ by $k$ and 
  $(n-k)$ by $n$ respectively such that $\Psi\circ H=0$ (the image of $H$ is the kernel
  of\/ $\Psi$).
  Then there is a constant $C$ so that
 \[
    p_{\alpha}(\Psi)\ =\ C\cdot (-1)^{|\alpha|} p_{\alpha^c}(H)\ ,
 \]
  for all $\alpha\in\binom{[n]}{k}$.
\end{prop}

\subsection{Connection to Schubert calculus}\label{S8:WrtoSC}
We explore some geometric consequences of the determinantal
formulas~\eqref{E8:Wronski_Concat} and~\eqref{E8:Laplace}.
Let $H$ be the column space of the matrix $H$, which is the kernel of the map
$\Psi$.
Then $H$ is a point in the Grassmannian $\Gr(p,m{+}p)$.
From the definition~\eqref{E8:Gamma} of $\Gamma$, we see that the column space
of the matrix $\Gamma(t)$ is the $m$-plane $F_m(t)$ osculating the rational
normal curve $\gamma$ at the point $\gamma(t)$. 
From~\eqref{E8:Wronski_Concat} and~\eqref{E8:Laplace}, we see that $s$ is a zero
of the Wronskian $\Phi(t):=\Wr(f_1(t),f_2(t),\dotsc,f_m(t))$ if and only if 
\[
    0\ =\ \det[\Gamma(s)\;:\;H]\,.
\]
This implies that there is a linear dependence among the columns of this matrix
and thus there is a nontrivial intersection between the subspaces $F_m(s)$ and $H$.

Suppose that a polynomial $\Phi(t)$ has distinct zeroes 
$s_1,s_2,\dotsc,s_{mp}$.
Then the columns of the matrix $H$ are linear forms cutting out
the linear span of polynomials  $f_1(t),f_2(t),\dotsc,f_m(t)$ of
degree $m{+}p{-}1$ having Wronskian $\Phi(t)$ 
 \begin{enumerate}
  \item[$\Leftrightarrow$]
     the $p$-plane $H$ meets the $m$-plane $F_m(s_i)$ nontrivially
     for each $i=1,2,\dotsc,mp$,
  \item[$\Leftrightarrow$]
     $H$ lies in the Schubert variety $X_\I\Fdot(s_i)$ for each
     $i=1,2,\dotsc,mp$. 
 \end{enumerate}
If the roots $s_1,\dotsc,s_{mp}$ are all real, then the Shapiro Conjecture
(Theorem~\ref{Th7:MTV}) asserts that all such $p$-planes $H$ are real, and there are the
expected number of them.
In this way, the second part of Theorem~\ref{T1:Shap_Conj} is a consequence of Theorem~\ref{Th7:MTV}.
\medskip

\noindent{\bf Second part of Theorem~\ref{T1:Shap_Conj} }
{\it 
 If the polynomial $\Phi(t)\in\P^{mp}$ has simple real roots 
 then there are $\#_{mp}$ real points in $\Wr^{-1}(F)$.
}\medskip

Recall from Chapter~\ref{ch:1} the formula for the degree of the Wronski map,
\[
  \tag{\ref{E1:WronDeg}}
   \#_{mp}\ =\ 
   \frac{1!2!\dotsb(m{-}1)!\cdot[mp]!}%
        {p!(p{+}1)!\dotsb(m{+}p{-}1)!}\ ,
\]
which is the number of inverse images of a regular value of the Wronski map.
The first part of Theorem~\ref{T1:Shap_Conj}, which asserts that all points are real in a
fiber of the Wronski map over a polynomial with only real roots, follows from the second
by a limiting argument that we give in Section~\ref{S8:reduction}.

\section{Asymptotic form of the Shapiro Conjecture is true}\label{S8:Asymp}

It is not too hard to show that the conclusion of the Shapiro conjecture when all
conditions $\alpha^i$ are simple, $\alpha^i=\bI$, holds for 
{\sl some} $s_1,\dotsc,s_{mp}\in\R$.
We give a (by now) standard argument for the following asymptotic form of the Shapiro
conjecture, which was given in~\cite{So99}, and independently in~\cite{EG01}.

\begin{thm}\label{T8:asymptotic}
 There exist real numbers $s_1,s_2,\dotsc,s_{mp}$ such that 
 \[
    X_\I\Fdot(s_1)\cap  X_\I\Fdot(s_2)\cap \dotsb \cap  X_\I\Fdot(s_{mp})
 \]
 is transverse with all points real.
\end{thm}

The proof is a version of Schubert's principle of degeneration to special position and the
same ideas can be used to establish similar results for other flag manifolds and related
varieties, such as Theorem~\ref{T7:quantum} on rational curves in Grassmannians.

Interchanging $\alpha$ with $\alpha^c$, the expansion~\eqref{E8:Laplace}
becomes (up to a sign)
\[
   \det[\Gamma(t)\;:\;H]\ =\ \sum_{\alpha\in\binom{[m{+}p]}{p}} 
    p_{\alpha}(H)\cdot(-1)^{|\alpha|}p_{\alpha^c}(\Gamma(t)) \,.
\]
We convert this into a very useful form by expanding the minor $p_{\alpha^c}(\Gamma(t))$.
Let $\alpha\in\binom{[m{+}p]}{m}$.
Observe that the determinant
\[
  \det\left[ \begin{matrix}
   t^{\alpha_1-1} & (\alpha_1-1)t^{\alpha_1-2}& \dotsc & 
      {\textstyle\frac{(\alpha_1-1)!}{(\alpha_1-m)!}}t^{\alpha_1-m}\\
   t^{\alpha_2-1} & (\alpha_2-1)t^{\alpha_2-2}& \dotsc &  \rule{0pt}{15pt}
      {\textstyle\frac{(\alpha_2-1)!}{(\alpha_2-m)!}}t^{\alpha_2-m}\\
    \vdots & \vdots& \ddots& \vdots\\
   t^{\alpha_m-1} & (\alpha_m-1)t^{\alpha_m-2}& \dotsc & 
      {\textstyle\frac{(\alpha_m-1)!}{(\alpha_m-m)!}}t^{\alpha_m-m}
   \end{matrix}\right]
\]
is equal to 
\[
   t^{|\alpha|} \cdot 
   \det\left[ \begin{matrix}
   1 & (\alpha_1-1)& \dotsc &{\textstyle\frac{(\alpha_1-1)!}{(\alpha_1-m)!}}\\
   1 & (\alpha_2-1)& \dotsc &{\textstyle\frac{(\alpha_2-1)!}{(\alpha_2-m)!}}
      \rule{0pt}{15pt}\\
    \vdots & \vdots& \ddots & \vdots\\
   1 & (\alpha_m-1)& \dotsc &{\textstyle\frac{(\alpha_m-1)!}{(\alpha_m-m)!}}
   \end{matrix}\right]
 \quad=\quad
  t^{|\alpha|} \cdot 
   \det\left[ \begin{matrix}
   1 & \alpha_1& \dotsc &\alpha_1^{m-1}\\
   1 & \alpha_2& \dotsc &\alpha_2^{m-1}  \rule{0pt}{15pt}\\
    \vdots & \vdots& \ddots & \vdots\\
   1 & \alpha_m& \dotsc &\alpha_m^{m-1}
   \end{matrix}\right]\ .
\]
(The second equality is via column operations.)
We recognize this last determinant as the Van der Monde, 
$\prod_{i<j}(\alpha_j-\alpha_i)$.
Write $V_\alpha$ for the product 
\[
   (-1)^{|\alpha|}\cdot\prod_{i<j}(\alpha_j-\alpha_i)\,.
\]
Since $|\alpha^c|=mp-|\alpha|$, we obtain the expansion for the Wronskian
(up to a global sign).
 \begin{equation}\label{E8:Wronski_Plucker}
   \det[\Gamma(t)\;:\; H]\ =\ 
    \sum_{\alpha\in\binom{[m{+}p]}{p}} 
    t^{mp-|\alpha|}\, V_{\alpha^c}\, p_\alpha(H)\,.
 \end{equation}

Observe that if we write $\DeCo{y_\alpha}:=V_{\alpha^c}p_\alpha(H)$,
then~\eqref{E8:Wronski_Plucker} becomes
\[
    \sum_{\alpha\in\binom{[m{+}p]}{p}} 
    t^{mp-|\alpha|}\, y_\alpha\,.
\]
In particular, the coefficient of $y_\alpha$ depends only upon the
rank, $|\alpha|$, of $\alpha\in\binom{[m{+}p]}{p}$.
Compare this to the Wronski polynomial~\eqref{E6:Wr_simple} of
Section~\ref{S6:PolyPoset}, and compare the original
polynomial~\eqref{E8:Wronski_Plucker}
to~\eqref{E6:Wronski_general}.  

%
%
\subsection{Schubert varieties}\label{S6:Schub_Vars}
We transpose all matrices, replacing
column vectors by row vectors.
Let $H\in\Gr(p,m{+}p)$ be represented as the row space of a $p$ by
$(m{+}p)$-matrix, and apply Gaussian elimination with pivoting from bottom to top and
right to left to $H$ to obtain a unique
representative matrix of the form
 \begin{equation}\label{Eq8:Schubert_cell}
  H\ =\ \Span\ \left(\begin{matrix}
   \Magenta{*\ \dotsb\ *}&\Blue{1}&0\ \dotsb\ 0&
   0&0\ \dotsb\ 0&0&0\ \dotsb\ 0\\
   \Magenta{*\ \dotsb\ *}&0&\Magenta{*\ \dotsb\ *}&
   \Blue{1}&0\ \dotsb\ 0&0&0\ \dotsb\ 0\\
   \Magenta{\vdots\ \qquad\vdots}&\vdots&\Magenta{\vdots\ \qquad\vdots}&
   0&\ \Blue{\ddots}\ &\vdots&\vdots\ \qquad\vdots\\
   \Magenta{*\ \dotsb\ *}&0&\Magenta{*\ \dotsb\ *}&
   \vdots&\Magenta{*\ \dotsb\ *}&\Blue{1}&0\ \dotsb\ 0\\
   \end{matrix}\right)\ .
 \end{equation}
Here, the entries $\Magenta{*}$ indicate an unspecified element of
our field ($\R$ or $\C$).

The set of columns containing the leading 1s (pivots) is a discrete invariant of 
the linear subspace $H$.
Let $\alpha\colon \alpha_1<\alpha_2<\dotsb<\alpha_p$ be the positions of the
pivots, that is, $\alpha_i$ is the column of the leading 1 in row $i$.
Observe that $p_\beta(H)=0$ unless $\beta_i\leq\alpha_i$ for every $i$.
This coordinatewise comparison defines the \DeCo{{\sl Bruhat order}} on the
indices $\alpha\in\binom{[m{+}p]}{p}$. 
The set of linear spaces whose row reduced echelon forms~\eqref{Eq8:Schubert_cell}
have pivots in the columns of $\alpha$ forms a topological cell of dimension
$|\alpha|$, called the \DeCo{{\sl Schubert cell}} and written $X^\circ_\alpha$.
The undetermined entries $\Magenta{*}$ in~\eqref{Eq8:Schubert_cell} show that
it is isomorphic to $\A^{|\alpha|}$, the affine space of dimension
$|\alpha|$. 

We ask: Which linear spaces are in the closure of the Schubert cell?
For the answer, let \DeCo{{$M_\alpha$}} be the set of matrices with full rank
$p$ where the entries in row $i$ are undetermined up to column $\alpha_i$, and
are 0 thereafter.
These matrices have the form
\[
  \left(\begin{matrix} 
   \Magenta{*\ \dotsb\ *}&0\ \dotsb\ 0&
    0\ \dotsb\ 0&0\ \dotsb\ 0\\
   \Magenta{*\ \dotsb\ *}&\Magenta{*\ \dotsb\ *}&
    0\ \dotsb\ 0&0\ \dotsb\ 0\\
   \Magenta{\vdots\ \qquad\vdots}&  \Magenta{\vdots\ \qquad\vdots}&
   \ \Magenta{\ddots}\ &0\ \dotsb\ 0 \\
   \Magenta{*\ \dotsb\ *}&\Magenta{*\ \dotsb\ *}&
   \Magenta{*\ \dotsb\ *}&0\ \dotsb\ 0 \\
  \end{matrix}\right)\ ,
\]
where the last undetermined entry $\Magenta{*}$ in row $i$ occurs in column
$\alpha_i$. 
This is a closed subset of the set of $p$ by $(m{+}p)$-matrices of full rank $p$.
The pivots $\beta$ of a matrix $M$ in $M_\alpha$ will occur weakly
to the left of the columns indexed by $\alpha$, so that
$\beta\leq\alpha$, and all possibilities $\beta$ can occur.

In particular, this shows that the set of $p$-planes $H$ parameterized by
matrices in $M_\alpha$ is the union of the Schubert cells indexed by $\beta$ for
$\beta\leq\alpha$ in the Bruhat order.
This is a closed subset of the Grassmannian, in fact it is one of the 
\DeCo{{\sl Schubert varieties}} defined in Chapter~\ref{S:ERAG}.
To see this, let $\e_1,\e_2,\dotsc,\e_{m{+}p}$ be basis vectors corresponding to the columns of our
matrices.
For each $i=1,\dotsc,m{+}p$ let $F_i$ be the linear span of the vectors
$\e_1,\e_2,\dotsc,\e_i$.
From the form of matrices in $M_\alpha$, we see that if $H$ is the row space of
a matrix in $M_\alpha$, then we have
\[
 \dim H\cap F_{\alpha_j}\ \geq\  j\qquad\mbox{for}\ j=1,\dotsc,p\,.
\]
These dimension inequalities define the Schubert variety \DeCo{$X_\alpha\Fdot$}.
Note that if $H\in X_\alpha\Fdot$, then $p_\beta(H)=0$ unless $\beta\leq\alpha$.
Write $X^\circ_\alpha\Fdot$ for the Schubert cell consisting of those $H$ of the
form~\eqref{Eq8:Schubert_cell}. 

The key lemma in our proof of Theorem~\ref{T8:asymptotic} is due
essentially to Schubert~\cite{Sch1886b}.

\begin{lemma} \label{L8:Schuberts_int} 
 For any $\alpha\in\binom{[m{+}p]}{p}$, 
 \[
   X_\alpha\Fdot\cap\{H\mid p_\alpha(H)=0\}\ =\ 
   \bigcup_{\beta\lessdot\alpha} X_\beta\Fdot\,.
 \]
\end{lemma}

Here $\beta\lessdot\alpha$ means that $\beta<\alpha$, but there is no index $\mu$
in the Bruhat order with $\beta<\mu<\alpha$.
This is easy to see set-theoretically, as for 
$H\in X_\alpha\Fdot$ we have $p_\beta(H)=0$ unless $\beta\leq \alpha$.

It is also easy to see that this is true on the generic point of each 
Schubert variety $X_\beta\Fdot$ for $\beta\lessdot\alpha$.
Fix some index $\beta$ with $\beta\lessdot\alpha$.
Then there is a unique index $k$ with $\beta_k=\alpha_k-1$, and for all other
indices $i$, $\beta_i=\alpha_i$.
Consider the subset of the matrices $M_\alpha$, where we require the entries
in row $i$ and column \DeCo{$\beta_i$} to be 1, and write $x_{k,\alpha_k}$ for
the entry in row $k$ and column $\alpha_k$ 
 \begin{equation}\label{E8:Mbeta}
  \left(\begin{matrix}
   \Magenta{*\ \dotsb\ *}&\DeCo{1}&0\ \dotsb\ 0&
   0&0&0\ \dotsb\ 0&0&0\ \dotsb\ 0\vspace{-3pt}\\
   \Magenta{\vdots\ \qquad\vdots}&0&\ \Blue{\ddots}\ &
   \vdots&\vdots&\vdots\ \qquad\vdots&\vdots&\vdots\ \qquad\vdots\\
   \Magenta{*\ \dotsb\ *}&0&\Magenta{*\ \dotsb\ *}&
   \Blue{1}&x_{k,\alpha_k}&0\ \dotsb\ 0&0&0\ \dotsb\ 0\\
   \Magenta{\vdots\ \qquad\vdots}&\vdots&\Magenta{\vdots\ \qquad\vdots}&
   0&0&\ \Blue{\ddots}\ &\vdots&\vdots\ \qquad\vdots\\
   \Magenta{*\ \dotsb\ *}&0&\Magenta{*\ \dotsb\ *}&
   \vdots&\vdots&\Magenta{*\ \dotsb\ *}&\Blue{1}&0\ \dotsb\ 0\\
   \end{matrix}\right)\ .
 \end{equation}
The row spans of these matrices form a dense subset of the Schubert variety
$X_\alpha\Fdot$, and therefore define a coordinate patch for $X_\alpha\Fdot$.
If we set $x_{k,\alpha_k}=0$, then we get all matrices of the
form~\eqref{Eq8:Schubert_cell}, but for the index $\beta$.

If $H$ is the row space of a matrix in this set~\eqref{E8:Mbeta}, then 
$p_\alpha(H)=x_{k,\alpha_k}$.
Thus, on this coordinate patch for $X_\alpha\Fdot$, the vanishing of the Pl\"ucker
coordinate $p_\alpha$ cuts out the Schubert variety $X_\beta\Fdot$,
scheme-theoretically.
Repeating this local argument for each $\beta\lessdot\alpha$,
proves~\eqref{L8:Schuberts_int}, at least at the generic point of each
component $X_\beta\Fdot$ (which is sufficient for our purposes).
More careful arguments show this is true even at the level of their homogeneous
ideals.

%
%
\subsection{Asymptotic form of Shapiro Conjecture}

We now have everything that we need to prove Theorem~\ref{T8:asymptotic}.
We will prove a stronger statement using induction on the Bruhat order (sometimes called 
\DeCo{{\sl Schubert induction}}).

\begin{lemma}\label{L8:schubert_induction}
   There exist real numbers $s_1,s_2,\dotsc,s_{mp}$ such that for all
   $\alpha\in\binom{[m{+}p]}{p}$,  
 \begin{equation}\label{Eq8:Inductive_step}
   X_{\alpha}\Fdot\ \cap\ \bigcap_{i=1}^{|\alpha|} X_\I\Fdot(s_i)
 \end{equation}
is transverse with all points of intersection real.
\end{lemma}

The statement of Theorem~\ref{T8:asymptotic} is the case 
$\alpha=m{+}1,m{+}2,\dotsc, m{+}p$, when $X_\alpha\Fdot$ is the Grassmannian.

\begin{rmk}\label{R8:For_7}
 It is not hard to see (it is equivalent to the Pl\"ucker
 formula~\cite{Plucker} for rational curves and was noted by Eisenbud and
 Harris~\cite[Theorem 2.3]{EH83}) that the intersection~\eqref{Eq8:Inductive_step} lies in
 the Schubert cell $X_\alpha^\circ\Fdot$  for the index $\alpha$.
 That is, every point $H$ in the intersection~\eqref{Eq8:Inductive_step} has
 row-reduced echelon form~\eqref{Eq8:Schubert_cell}, for the index $\alpha$.
 \QED 
\end{rmk}

Observe that when $\alpha=1,2,\dotsc,p$, then $|\alpha|=0$ and the Schubert
variety $X_\alpha\Fdot$ consists of the single point $\{F_p\}$.
Thus the base case of the induction to prove Lemma~\ref{L8:schubert_induction} is trivial,
as there is no intersection to contend with. 

Suppose that we have real numbers $s_1,\dotsc,s_j$ such that, for each $\alpha$ with
$|\alpha|=j$ the intersection~\eqref{Eq8:Inductive_step} is transverse with all points
real.
Let $\alpha\in\binom{[m{+}p]}{p}$ with $|\alpha|=j{+}1$.
Observe that by~\eqref{E8:Wronski_Plucker} the intersection 
$X_\alpha\Fdot\cap X_\I\Fdot(t)$ is defined by the single polynomial
equation
\[
  \sum_\beta t^{mp-|\beta|} V_{\beta^c} p_\beta(H)\ =\ 0
 \qquad\quad\mbox{for}\ H\in X_\alpha\Fdot\,.
\]
Since $p_\beta(H)=0$ unless $\beta\leq\alpha$, this becomes
\[
    \sum_{\beta\leq\alpha} t^{mp-|\beta|} V_{\beta^c} p_\beta(H)\ =\ 0\,.
\]
Dividing by the lowest power $t^{mp-|\alpha|}$ of $t$, this becomes
\[
   \sum_{\beta\leq\alpha} t^{|\alpha|-|\beta|} V_{\beta^c} p_\beta(H)
  \ =\  V_{\alpha^c} p_\alpha(H)\ +\ 
     t\cdot \sum_{\beta< \alpha} t^{|\alpha|-|\beta|-1} V_{\beta^c} p_\beta(H)\
  =\ 0\, .
\]
Since $V_{\alpha^c}\neq 0$, we see that in the limit as $t\to 0$, this equation
becomes $p_\alpha(H)=0$.
Using the lemma of Schubert~\eqref{L8:Schuberts_int}, we conclude that 
 \begin{equation}\label{Eq8:Limit}
   \lim_{t\to 0}\left( X_\alpha\Fdot\cap X_\I\Fdot(t)\right)\ =\ 
   X_\alpha\Fdot\cap\{H\mid p_\alpha(H)=0\}\ =\ 
   \bigcup_{\beta\lessdot\alpha} X_\beta\Fdot\,.
 \end{equation}

By our induction assumption on $j$, each intersection
\[
    X_\beta\Fdot  \ \cap\ \bigcap_{i=1}^j X_\I\Fdot(s_i)
\]
is transverse with all points real, and by Remark~\ref{R8:For_7} the intersection is
contained in the Schubert cell $X_\beta^\circ\Fdot$.
Since the Schubert cells are disjoint, we conclude that the intersection
\[
  \Bigl(\bigcup_{\alpha\lessdot\beta} X_\beta\Fdot\Bigr)
   \ \cap\ \bigcap_{i=1}^j X_\I\Fdot(s_i)
\]
is transverse with all points real.
By the computation of the limit~\eqref{Eq8:Limit}, and the observation that
transversality is preserved by small perturbations, we see that there is a number
$0<\varepsilon_\alpha$ such that if $t\leq \varepsilon_\alpha$ then
\[
   X_\alpha\Fdot \ \cap\ X_\I\Fdot(t)
   \ \cap\ \bigcap_{i=1}^j X_\I\Fdot(s_i)
\]
is transverse with all points real.

We complete the induction by setting $s_{j+1}$ to be the minimum of the
numbers $\varepsilon_\alpha$ where $|\alpha|=j+1$.  
This proves Theorem~\ref{T8:asymptotic}.\QED\smallskip

Similar asymptotic arguments are behind the proof of Theorem~\ref{T7:quantum},
which proved reality in the quantum Schubert calculus, as well as results for 
the classical flag manifolds and for the orthogonal Grassmannian~\cite{So00c}. 

\begin{figure}[htb]
\[
  \begin{picture}(210,350)(-3,0)
    \put(0,0){\includegraphics[height=350pt]{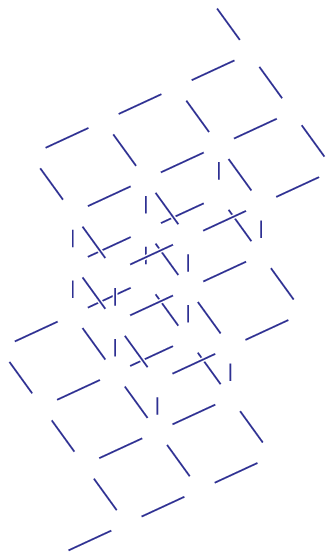}}

    \put(122,340){567} 
    \put(148,305){467} 
    \put(104,285){367} \put(172,270){457}
    \put(60,265){267}  \put(130,250){357} \put(198,235){456}
    \put(15,244){167}  \put(85,229){257}  \put(130,218){347} \put(154,214){356}
    \put(41,209){157}  \put(85,199){247}  \put(111,194){256} \put(154,183){346}
    \put(41,178){147}  \put(66,174){156}  \put(85,167){237}  \put(111,163){246} \put(179,148){345}
    \put(41,147){137}  \put(66,143){146}  \put(111,131){236} \put(135,127){245}
    \put(-3,126){127}  \put(66,111){136}  \put(91,108){145}  \put(135,96){235}
    \put(22,91){126}   \put(91,76){135}   \put(160,61){234}
    \put(47,56){125}   \put(116,41){134}
    \put(71,21){124}
    \put(28,1){123}
  \end{picture}
  \qquad
  \begin{picture}(210,350)(-3,0)
    \put(0,0){\includegraphics[height=350pt]{figures/8/G37.eps}}

    \put(122,340){462} 
    \put(147,305){462} 
    \put(104,285){252} \put(173,270){210}
    \put(62,265){84}  \put(129,250){168} \put(200,235){42}
    \put(17,244){14}  \put(87,229){70}  \put(131,218){56} \put(156,214){42}
    \put(43,210){14}  \put(87,198){35}  \put(113,194){21} \put(157,183){21}
    \put(46,177){9}  \put(71,174){5}  \put(87,167){10}  \put(111,163){16} \put(184,147){5}
    \put(46,147){4}  \put(71,143){5}  \put(115,132){6} \put(140,127){5}
    \put(1,126){1}   \put(71,111){3}  \put(96,108){2}  \put(140,96){3}
    \put(27,91){1}   \put(96,76){2}   \put(166,61){1}
    \put(52,56){1}   \put(122,41){1}
    \put(77,21){1}
    \put(33,1){1}
  \end{picture}
\]
\caption{Schubert's recursion for $\Gr(3,7)$.}
\label{F8:Schubert_Induction}
\end{figure}

\begin{rmk}\label{R8:construct}
 The proof of Theorem~\ref{T8:asymptotic} used induction to show that the 
 intersection~\eqref{Eq8:Inductive_step} is transverse with all points real.
 In fact, it gives an inductive method to construct all the points of intersection.
 The induction began with $\alpha=1,2,\dotsc,p$ so that $|\alpha|=0$ and the
 Schubert variety $X_\alpha\Fdot$ consists of the single point $\{F_p\}$.
 When $|\alpha|=j+1$, the limit
\[
  \Bigl(\lim_{t\to 0}X_\alpha\Fdot\cap X_\I\Fdot(t)\Bigr)
   \ \cap\ \bigcap_{i=1}^j X_\I\Fdot(s_i)\ =\ 
   \bigcup_{\beta\lessdot\alpha} X_\beta\Fdot
   \ \cap\ \bigcap_{i=1}^j X_\I\Fdot(s_i)\,
\]
 shows that each point in the intersection~\eqref{Eq8:Inductive_step}
 is connected to a point in 
 \begin{equation}\label{E8:start_solution}
   \bigcup_{\beta\lessdot\alpha} X_\beta\Fdot
   \ \cap\ \bigcap_{i=1}^j X_\I\Fdot(s_i)
    \ =\ 
   \bigcup_{\beta\lessdot\alpha} \Bigl( X_\beta\Fdot
     \ \bigcap_{i=1}^j X_\I\Fdot(s_i)\Bigr)
 \end{equation}
 along a path as $t$ ranges from $s_{j+1}$ to $0$, and the union on the right is disjoint.

 For the inductive construction, we may suppose that the 
 points in the set~\eqref{E8:start_solution} have been previously
 constructed as $\beta\lessdot\alpha$ implies that $|\beta|=j$.
 Starting at one of the points in~\eqref{E8:start_solution} and tracing the path from
 $t=0$ back to $t=s_{j+1}$ gives a point in the intersection~\eqref{Eq8:Inductive_step}
 for $\alpha$ with $|\alpha|=j{+}1$, and all such points in the
 intersection~\eqref{Eq8:Inductive_step} arise in this manner.
 Following paths along a general curve in $\C$ (as opposed to the line segment
 $[0,s_{j+1}]$) constructs points in the intersection~\eqref{Eq8:Inductive_step} where
 $s_{j+1}$ is any complex number.
 This is the idea behind the numerical \DeCo{{\sl Pieri homotopy algorithm}}, which 
 was proposed in~\cite{HSS98} and implemented in~\cite{HV00}.
 Its power was demonstrated in~\cite{Galois}, which used the Pieri homotopy algorithm to
 compute all solutions to a Schubert problem on $\Gr(3,9)$ with $17589$ solutions.

 If \DeCo{$\delta(\alpha)$} is the number of points in the
 intersection~\eqref{Eq8:Inductive_step}, then this limiting process also gives
 the recursion along the Bruhat order $\delta(\alpha)$,
 \begin{equation}\label{Eq8:delta-recursion}
  \begin{array}{rcl}
    \delta(1,2,\dotsc,p)&=& 1\,,\\
    \delta(\alpha)&=& \sum_{\beta\lessdot\alpha} \delta(\beta)\,.\rule{0pt}{14pt}
  \end{array}
 \end{equation}
 Schubert discovered this recursion~\cite{Sch1886c} and used it to
 compute the number $\delta(567)=462$ when $m=4$ and $p=3$.
 This is the number $\#_{4,3}$ given by the formula, 
\[
  \tag{\ref{E1:WronDeg}}
   \frac{1!2!\dotsb(m{-}1)!\cdot[mp]!}%
        {p!(p{+}1)!\dotsb(m{+}p{-}1)!}\ ,
\]
 which is also due to Schubert, as he solved his recursion to obtain a closed formula.
 This recursion shows that the number $\delta(\alpha)$ 
 may be interpreted combinatorially as the number of paths from the bottom of
 the Bruhat order to the element $\alpha$.
 In fact to each solution we constructed in~\eqref{Eq8:Inductive_step}, we may associate a
 path from $1,2,\dotsc,p$ to $\alpha$.

 Figure~\ref{F8:Schubert_Induction} shows the Bruhat order in this case when $m=4$ and $p=3$
 and the recursion for the numbers  $\delta(\alpha)$.
\QED
\end{rmk}

%
\subsection{Reduction to special case of Shapiro conjecture}\label{S8:reduction}

In Section~\ref{S8:WrtoSC}, we demonstrated that Theorem~\ref{T1:Shap_Conj} is equivalent
to Theorem~\ref{Th7:MTV} when all Schubert conditions are simple (each $\alpha^i=\bI$).
In fact this case of Theorem~\ref{Th7:MTV} implies a weak form of the general case, in
which we do not require transversality.
The main idea is to use the limit~\eqref{Eq8:Limit}, which we must first reinterpret.
The flag $\Fdot$ is the osculating flag $\Fdot(t)$ when $t=0$.
In fact, the limit~\eqref{Eq8:Limit} still holds if we replace $\Fdot$ by $\Fdot(s)$ and 
$0$ by $s$ for any point $s$ of $\P^1$.
That is,
 \begin{equation}\label{E8:newLimit}
   \lim_{t\to s} \bigl(X_\alpha\Fdot(s)\cap X_\I\Fdot(t)\bigr)\ =\ 
    \bigcup_{\beta\lessdot\alpha} X_\beta\Fdot(s)\,.
 \end{equation}
This is simply the limit~\eqref{Eq8:Limit} translated by the invertible matrix $\Phi(s)$
whose $i,j$-entry is 
\[
   \Phi(s)_{i,j}\ =\ \frac{1}{(i-1)!} \left(\frac{d}{dt}\right)^{\! i-1}
   t^{j-1}\Big|_{t=s}\,,
\]
as $\Phi(s).\Fdot(t)=\Fdot(s+t)$.

\begin{thm}\label{T8:real_limit}
 Suppose that Theorem~\eqref{Th7:MTV} holds for the Schubert problem in which all
 conditions $\alpha^i$ are simple $(\alpha^i=\bI)$.
 Then for any $\alpha^1,\dotsc,\alpha^n$ with $mp=\sum_i(mp-|\alpha^i|)$ and any 
 distinct $s_1,\dotsc,s_n\in\R\P^1$, the intersection
 \begin{equation}\label{E8:int}
   \bigcap_{i=1}^n X_{\alpha^i}\Fdot(s_i)
 \end{equation}
 has \DeCo{all} points real.
\end{thm}

We prove this by downward induction on the number $n$ of Schubert conditions
in~\eqref{E8:int}, using the limit~\eqref{E8:newLimit} and the simple idea that a limit of
a collection of real points is necessarily a collection of real points.

First, when $n=mp$ each $\alpha^i$ is simple and all points in the
intersection~\eqref{E8:int} are real as that is our hypothesis in Theorem~\ref{T8:real_limit}.
Suppose that $n<mp$.
Then we have $|\alpha^i|<mp-1$ for some $i$.
Suppose that $|\alpha^n|<mp-1$ and let $\beta\in\binom{[m{+}p]}{p}$ be a Schubert
condition with $\alpha^n\lessdot\beta$ so that $|\beta|=|\alpha|+1$.
Then
\[
   \lim_{t\to s_n} \Bigl[\,\bigcap_{i=1}^{n{-}1}X_{\alpha^i}\Fdot(s_i)\Bigr]\cap 
     \Bigl(X_{\beta}\Fdot(s_n)\cap X_{\I}\Fdot(t)\Bigr)\ =\ 
    \Bigl[\,\bigcap_{i=1}^{n{-}1}X_{\alpha^i}\Fdot(s_i)\Bigr]\cap 
    \Bigl(\bigcup_{\alpha\lessdot\beta}X_\alpha\Fdot(s_n)\Bigr)\,.
\]
The elementary inclusion $\subset$ of the limit in the set on the right is clear
from~\eqref{E8:newLimit}.
The equality of the two sides follows as the intersection on the right is
zero-dimensional, and therefore cannot contain any excess intersection.
 By induction, for general $t\in\R$, every point in the left-hand intersection is real,
 and so every point in the limit is real.
 Theorem~\ref{T8:real_limit} follows as $\alpha^n\lessdot\beta$ and so the
 intersection~\eqref{E8:int} is a subset of the right-hand side.

%
\section{Grassmann duality}\label{S8:duality}

In Section~\ref{S8:Asymp} we showed how the Wronski formulation of the Shapiro
Conjecture---$m$-dimensional spaces of polynomials of degree $m{+}p{-}1$ whose Wronskian
has distinct real roots---corresponds to an intersection of hypersurface Schubert
varieties in $\Gr(p,m{+}p)$ defined by flags that osculate the rational normal curve at
the roots of the Wronskian.
The first formulation concerns points in $\Gr(m,m{+}p)$ while the second concerns points
in $\Gr(p,m{+}p)$.
In our proof of this correspondence we considered an $m$-dimensional space of polynomials 
as a space of linear forms on $\C^{m+p}$, and associated this to the $p$-plane annihilated
by the linear forms.
This gives a natural bijection
 \begin{equation}\label{E8:GrassDuality}
   \Gr(m,(\C^{m{+}p})^*)\ \longrightarrow\ \Gr(p,\C^{m{+}p})\,.
 \end{equation}
Moreover, the annihilators of the subspaces in a flag $\Fdot$ in $\C^{m{+}p}$ form the
dual flag \DeCo{$\Fdot^*$} in $(\C^{m{+}p})^*$.
Under the identification of Grassmannians~\eqref{E8:GrassDuality}, the Schubert variety $X_\alpha\Fdot$ of
$\Gr(p,\C^{m{+}p})$ is identified with the Schubert variety
$X_{\alpha^*}\Fdot^*$ in $\Gr(m,(\C^{m{+}p})^*)$, where
\[
   \DeCo{\alpha^*}\ \colon\ 
    m{+}p{+}1-\alpha^c_m\ <\ 
    m{+}p{+}1-\alpha^c_{m-1}\ <\ \dotsb\ <\ 
    m{+}p{+}1-\alpha^c_2\ <\ m{+}p{+}1-\alpha^c_1\,,
\]
that is, to obtain $\alpha^*$, first form the complement $\alpha^c$ of $\alpha$ in
$[m{+}p]$, then subtract each component from $m{+}p{+}1$, and put the result in increasing order.
This is an exercise in combinatorial linear algebra, and the relation
$|\alpha|=|\alpha^*|$ is a combinatorial exercise.

Write \DeCo{$\C_{m{+}p{-}1}[t]$} for the space of polynomials of degree at most
$m{+}p{-}1$, which we identified as the dual space to $\C^{m{+}p}$. 
We describe the Schubert subvarieties of $\Gr(m,\C_{m{+}p{-}1}[t])$ that correspond to the Schubert
varieties $X_\alpha\Fdot(s)$.

Let $V\subset\C_{m{+}p{-}1}[t]$ be an $m$-dimensional space of polynomials.
For any $s\in\P^1$, $V$ has a distinguished basis $f_1,\dotsc,f_m$ whose orders of
vanishing at the point $s$ are strictly increasing,
\[
  \ord_s(f_1)\ <\ \ord_s(f_2)\ <\ \dotsb\ <\ \ord_s(f_m)\,.
\]
This follows by Gaussian elimination applied to any basis $f_1,\dotsc,f_m$ of $V$.
Suppose that $0\leq a_1$ is the minimal order of vanishing at $s$ of some $f_i$.
Reordering the basis, we may assume that $\ord_s(f_1)=a_1$.
Subtracting an appropriate multiple of $f_1$ from the subsequent elements 
gives a new basis, still written $f_1,\dots,f_n$, with $a_1<\ord_s(f_i)$ for $i>1$.
Suppose that $\ord_s(f_2)$ is minimal among the orders of vanishing at $s$ of $f_i$ for $i>1$ and
now subtract appropriate multiples of $f_2$ from the subsequent elements, and continue.
The resulting sequence $\DeCo{a}:=(a_1,\dotsc,a_m)=(\ord_s(f_1),\dotsc,\ord_s(f_m))$ is the 
\DeCo{{\sl ramification sequence}} of $V$ at $s$.

An elementary calculation shows that if $V$ has ramification sequence $a$ at a point
$s\in\P^1$, then the Wronskian of $V$ vanishes at $s$ to order
\[
   |a|\ =\ a_1{-}0 + a_2{-}1 +\ \dotsb\ + a_{m}{-}(m{-}1)\,.
\]

Define a flag $\Edot(s)\subset\C_{m{+}p{-}1}[t]$ where $E_i(s)$ is the space of all polynomials
that vanish to order at least $m{+}p{-}i$ at $s$.
With these definitions, we have the following lemma.

\begin{lemma}
 A space $V$ of polynomials in $\Gr(m,\C_{m{+}p{-}1}[t])$ has ramification sequence $a$ at $s$ if
 and only if
\[
   V\ \in\ X_{a^r}\Edot(s)\,,
\]
 where $\DeCo{a^r}\colon m{+}p{-}a_m<\dotsb<m{+}p{-}a_p\in \binom{[m{+}p]}{m}$.
\end{lemma}

A polynomial $f(t)\in E_i(s)$ if and only if $(t-s)^{m+p-i}$ divides $f$, if and only if
$f^{(j)}(s)=0$ for $j=0,1,\dotsc,m{+}p{-}1{-}i$.
If we view $f$ as a linear form on $\C^{m{+}p}$ so that $f(t)=f\circ\gamma(t)$, we see
that  $f^{(j)}(t)=f\circ \gamma^{(j)}(t)$, and therefore $f(t)\in E_i(s)$  if and only
if $f$ annihilates the osculating subspace $F_{m+p-i}(s)$ to $\gamma$ at $\gamma(s)$.
Thus $E_i(s)^\perp=F_{m+p-i}(s)$, and so $\Edot(s)$ is the dual flag to $\Fdot(s)$.
In particular, the Schubert variety $X_\alpha\Fdot(s)$ corresponds to
$X_{\alpha^*}\Edot(s)$ under Grassmann duality.

\begin{thm}\label{T8:GrassmannDuality}
  The identification of $\C_{m{+}p{-}1}[t]$ as the dual space to $\C^{m{+}p}$ induces an
  isomorphism of Grassmannians
\[
      \Gr(m,(\C^{m{+}p})^*)\ \longrightarrow\ \Gr(p,\C^{m{+}p})\,.
\]
  For any $s\in\P^1$, this restricts to an isomorphism of the Schubert varieties,
\[
   X_{\alpha^*}\Edot(s)\ \longrightarrow\ X_\alpha\Fdot(s)\,.
\]
\end{thm}

%
%
%
%
\chapter{The Shapiro Conjecture for Rational Functions}\label{Ch:EG}


We continue our study of the Shapiro Conjecture, which asserts that
if $f_1(t),f_2(t),\dotsc,f_m(t)$ are polynomials of degree $m{+}p{-}1$ whose
Wronskian is $\Phi(t)$ has degree $mp$ with distinct real roots,
then the linear span of  $f_1(t),f_2(t),\dotsc,f_m(t)$ is real.
We consider an apparently degenerate case, that of rational
functions, which is when $m=2$.

Eremenko and Gabrielov~\cite{EG02} originally gave a proof of the Shapiro
Conjecture when $m=2$ using essentially the uniformization theorem from complex
analysis. 
This proof was also complicated and challenging for {\sl this} reader.
They subsequently found a second, significantly more elementary proof~\cite{EG05}.
We will first discuss that proof, and then a generalization concerning rational functions that
are constant on prescribed sets, which leads to a generalization of the Shapiro Conjecture that
we will discuss in Chapter~\ref{Ch:Frontier}.

\section{The Shapiro Conjecture for rational functions}\label{S:SCRF}

The Shapiro conjecture for $m=2$ asserts that if $f_1(t)$ and $f_2(t)$ are univariate
polynomials whose Wronskian
 \begin{equation}\label{E9:Wronski}
   \Wr(f_1,f_2)\ =\ f'_1(t)f_2(t)\ -\ f_1(t)f'_2(t)
 \end{equation}
has only real roots, then the complex linear span $\langle f_1,f_2\rangle$ is real in that
there are real polynomials $g_1$ and $g_2$ with $\langle f_1,f_2\rangle=\langle g_1,g_2\rangle$.

This has a natural interpretation in terms of rational functions.
The quotient of two univariate polynomials $f_1$ and $f_2$ defines a rational function
$\rho\colon\P^1\to\P^1$ which on $\C\subset\P^1$ is
\[
   \rho\ \colon\ t\ \longmapsto\ f_1(t)/f_2(t)\ .
\]
Since
\[
   \rho'(t)\ =\ \frac{f'_1(t)f_2(t)\ -\ f_1(t)f'_2(t)}{f_2(t)^2}
   \ =\ \frac{\Wr(f_1,f_2)}{f_2(t)^2}\,,
\]
if $f_1$ and $f_2$ are relatively prime, the critical points of $\rho$ in $\P^1$ are
the roots of their Wronskian.

Two rational functions $\rho_1,\rho_2\colon\P^1\to\P^1$ are
\DeCo{{\sl equivalent}} if they differ by a fractional linear transformation on the target $\P^1$.
Since a fractional linear transformation on $f_1/f_2$ is a change of basis in the 
linear span $\langle f_1,f_2\rangle$, an equivalence class of rational functions is
simply a two-dimensional space of polynomials.
An equivalence class is \DeCo{{\sl real}} if the corresponding linear space is real.
We state the theorem of Eremenko and Gabrielov.

\begin{thm}[Shapiro's conjecture for rational functions~\cite{EG02,EG05}]\label{T9:EG}
  If a rational function $\rho\colon\P^1\to\P^1$ has only real critical points, then
  $\rho$ is equivalent to a real rational function.
\end{thm}

Theorem~\ref{T8:asymptotic} and the Grassmann duality of Theorem~\ref{T8:GrassmannDuality}
ensure the existence of a polynomial $\Phi_0(t)\in\R_{2p}[t]$ with only real roots such
that every space of polynomials with Wronskian $\Phi_0(t)$ is real, and there are exactly
$\#_{2,p}=\frac{1}{p+1}\binom{2p}{p}$ such spaces of polynomials.
The elementary proof of Theorem~\ref{T9:EG} analytically continues these
$\#_{2,p}$ real spaces of polynomials as the $2p$ distinct real roots of $\Phi_0(t)$ vary.
This continuation will produce fewer than $\#_{2,p}$ real spaces of polynomials for any
$\Phi(t)$ only if some of the spaces become complex during the continuation.
But this can happen only if two spaces of polynomials first become equal during the
continuation. 

The proof shows that such a collision cannot occur by associating discrete objects called nets,
to the real rational functions that are distinct from each of the $\#_{2,p}$ spaces of 
polynomials with Wronskian $\Phi_0(t)$, and which are preserved under a continuation that
varies the roots of $\Phi_0(t)$.
Thus no collisions are possible, which will imply Theorem~\ref{T9:EG}.

\subsection{Continuity and nets of rational functions}

A point $p$ in the Grassmannian $\DeCo{\Gr(2,\C_{p{+}1}[t])}\simeq\Gr(2,p{+}2)$ is a
two-dimensional space $V$ of univariate polynomials of degree at most $p{+}1$.
Each such point gives an equivalence class of rational functions
$\rho\colon\P^1\to\P^1$ of degree $p{+}1{-}\delta$, where $\delta$ is the maximum
degree of a common factor of the polynomials in $V$.
Working with this equivalence is awkward, so we will instead use the
\DeCo{{\sl Stiefel manifold}}, \DeCo{$\St_{p{+}1}$}, which is a $GL(2,\R)$-fiber bundle
over $\Gr(2,\C_{p{+}1}[t])$. 

The points of $\St_{p{+}1}$ are pairs of nonproportional univariate polynomials of
degree at most $p{+}1$.
Hence $\St_{p{+}1}$ is an open subset of $\R^{2p+4}$, with coordinates the
coefficients of the polynomials $f$ and $g$.
We give $\St_{p{+}1}$ the subspace topology.
The association $\St_{p{+}1}\ni(f,g)\longmapsto f/g$ defines a map (written
\DeCo{$\pi$}) from $\St_{p{+}1}$ to the
space of rational functions.
While this map is not continuous as a map of spaces, it does have the weak
continuity property given in Proposition~\ref{P9:Weak_cont} below.

Let $\DeCo{Z}\subset \St_{p{+}1}$ be the locus of pairs $(f,g)$ with either
 \begin{center}
   $\deg \gcf (f,g)>0$ \quad or \quad $\deg(f)<d$ and $\deg(g)<p{+}1$.
 \end{center}
That is, $f$ and $g$ either have a common root in $\C$ or else a common root at
$\infty$, and thus define a rational function $f/g$ of degree less than $p{+}1$.

\begin{prop}\label{P9:Weak_cont}
 Let $\{\varphi_j\mid j\in\N\}\subset \St_{p{+}1}\setminus Z$ be a sequence of points that
 converges to some point $\varphi=(f,g)\in Z$.
 Let $z_1,\dotsc,z_k$ be the common roots of $f$ and $g$ (including $\infty$ if 
 $\deg(f)$ and $\deg(g)$ are both less than $d$).
 Then the sequence of functions $\{\pi(\varphi_j)\mid j\in\N\}$ converges to
 $\pi(\varphi)$ uniformly on compact subsets of\/ $\P^1\setminus\{z_1,\dotsc,z_k\}$.
\end{prop}

 We give an elementary proof of this proposition.
 Let $K\subset\P^1$ be a compact subset disjoint from the common roots
 $\{z_1,\dotsc,z_k\}$ of $f$ and $g$.
 We may cover $\P^1$ by the standard affine charts $\C_0$ and $\C_\infty$ whose
 coordinates are $t$ and $1/t$, respectively.
 Then $K=K_0\cup K_\infty$, where $K_0\subset\C_0$ and $K_\infty\subset\C_\infty$ are
 compact subsets of the two affine charts.
 It suffices to show that the sequence of functions $\{\pi(\varphi_i)\}$ converges
 uniformly to $\pi(\varphi)$ on each set $K_0$ and $K_\infty$.

 Now $K_0$ is itself covered by compact sets $K_0^f$ and $K_0^g$, where $K_0^f$ contains
 no root of $f$ and $K_0^g$ contains no root of $g$.
 Removing finitely many members of the sequence 
 $\{\varphi_i=(f_i,g_i)\mid i\in\N\}$, we may assume that no $f_i$ has a root in $K_0^f$
 and no $g_i$ has a root in $K_0^g$.
 As $(f_i,g_i)$ converges to $(f,g)$ in $\St_{p{+}1}$, and no $g_i$ has a root
 in $K_0^g$, both sequences of functions
\[
   \{ f_i(t)\mid i\in\N\}\qquad\mbox{and}\qquad
   \{(g_i(t))^{-1}\mid i\in\N\}
\]
 are uniformly bounded in $K_0^g$.
 Therefore, the sequence of functions 
\[
    \frac{f_i(t)}{g_i(t)}\ \colon\ K_0^g\ \longrightarrow\ \C
\]
 is uniformly bounded and converges pointwise on the compact set $K_0^g$ to 
 $f(t)/g(t)$.
 Thus this convergence is uniform on $K_0^g$, and it remains
 uniform under $\C\xrightarrow{\,\sim\,}\C_0\subset\P^1$.
 The same arguments work for $K_0^f$, as well as $K_\infty$, which proves the proposition.
\hfil\QED

This proposition is half of the engine of this proof
of Eremenko and Gabrielov.
The other half is the asymptotic proof of Shapiro's conjecture, Theorem~\ref{T8:asymptotic}. 
\medskip

We now explain how to associate an embedded graph with distinguished vertices to each real
rational function.
Let \DeCo{$R_{p{+}1}$} be the set of nonconstant real rational functions of degree at
most $p{+}1$, all of whose critical points are real.
If $\rho\in R_{p{+}1}$, then $\rho^{-1}(\R\P^1)\subset\P^1$ defines an embedded (multi-)
graph $\Gamma$ with the following properties:
 \begin{enumerate}
  \item[(i)] $\Gamma$  is stable under complex conjugation and 
    $\R\P^1\subset \Gamma$.
 \end{enumerate}
Call any edge in $\Gamma\setminus\R\P^1$ an \DeCo{{\sl interior edge}}.
 \begin{enumerate}
  \item[(ii)] The vertices of $\Gamma$ lie on $\R\P^1$ and are the critical points
    of the rational function $\rho$.
    The valence of a vertex is even and it equals twice the order of ramification of $\rho$
    at the critical point, which we call the \DeCo{{\sl local degree}} of $\Gamma$ at the vertex.
 \end{enumerate}
The set-theoretic difference $\P^1\setminus\Gamma$ is a union of $2d$ cells, 
where $d$ is the degree of $\rho$.  
The closure of each cell is homeomorphic to a disc, and the boundary of each
cell maps homeomorphically to $\R\P^1$.
This is because the cells (and their closures) are the inverse images of one of 
 the two discs in $\P^1\setminus \R\P^1$ (or their closures), and there are no
 critical points in the interior of any cell.
 We deduce the following additional property of these multi-graphs.
 \begin{enumerate}
  \item[(iii)]  No interior edge of $\Gamma$ can begin and end at the same vertex.
 \end{enumerate}
 Indeed, if an interior edge $e$ begins and ends at the same vertex, then $r(e)=\R\P^1$ as $v$ is
 the only critical point on $e$.
 But then $e$ must be the boundary of any cell adjacent to $e$, which implies that 
 $\Gamma$ consists of only two cells and one edge $e$ and so $\rho$ has degree 1,
 and in fact $e$ was not an interior edge after all.

Here are three pictures of such embedded (multi-) graphs for quintic rational functions with
evenly spaced critical points.
We have drawn $\R\P^1$ as a circle with the upper half plane in its interior.
The point $\sqrt{-1}$ is at the center of the circle, 
$-\sqrt{-1}$ is the point at infinity, and complex
conjugation is inversion in the circle.
\[
   \includegraphics[height=80pt]{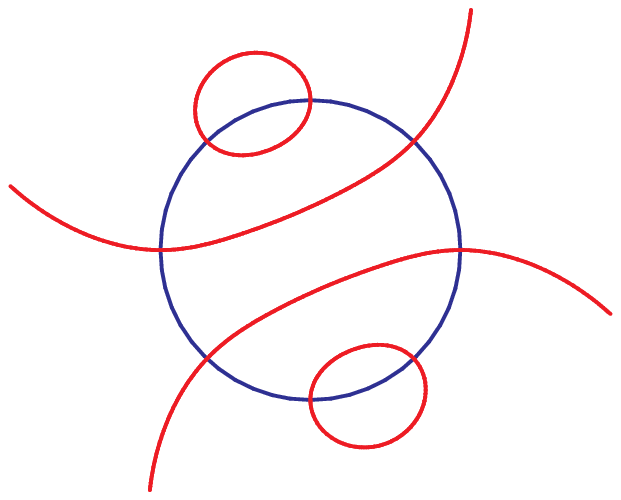}\qquad
   \includegraphics[height=80pt]{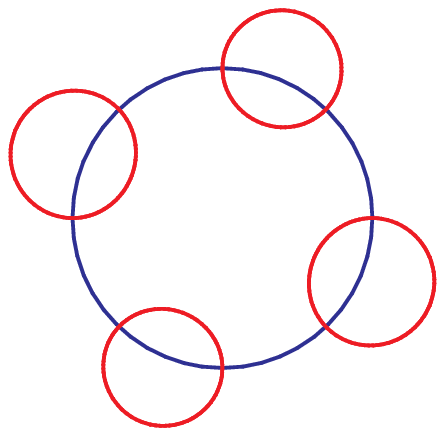}\qquad
   \includegraphics[height=80pt]{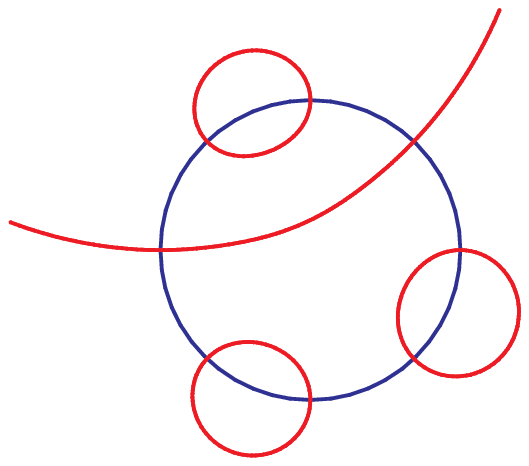}
\]

We seek to analytically continue rational functions whose Wronskians
lie in a curve of polynomials $\{\Phi_z(t)\mid z\in[0,1]\}$ where each
$\Phi_z(t)$ has degree $2p$ with distinct real roots 
$s_1(z), s_2(z), \dotsc,s_{2p}(z)$, and where each $s_i$ is a continuous
function of $z$.
The vertices of the graph $\rho^{-1}(\R\P^1)$ associated to a rational function $\rho$ with
Wronskian $\Phi_z(t)$ are labeled by these roots, or equivalently by the numbers
$1,2,\dotsc,2p$. 
Since the relative order of these roots $s_1(z), s_2(z), \dotsc,s_{2p}(z)$
does not change as $z$ varies (because each polynomial $\Phi_z(t)$ has distinct
roots), we may capture this information by labeling only one root, say $s_1(z)$
(which is a vertex of the corresponding graph), and assuming
that the roots are ordered in a manner consistent with a fixed orientation of
$\R\P^1$.
It is these labeled graphs that we wish to consider up to isotopy (deformation
in $\P^1$).

\begin{defn}
 A \DeCo{{\sl net}} is an (isotopy) equivalence
 class of such embedded multi-graphs in $\P^1$ satisfying (i), (ii), and
 (iii), with a  distinguished vertex.
\end{defn} 

 Here are the five nets with 6 vertices, each with local degree 2 at every vertex.
 These correspond to rational functions of degree four with simple ramification.
 \begin{equation}\label{Eq:quartic_nets}
  \raisebox{-30pt}{%
   \includegraphics[height=65pt]{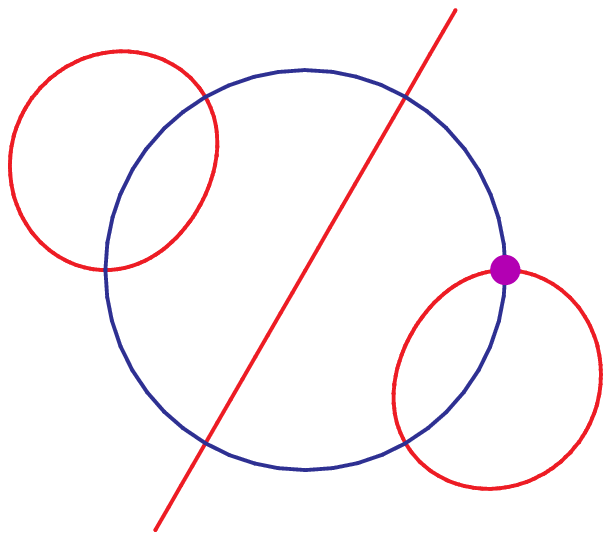}\quad
   \raisebox{-7.5pt}{\includegraphics[height=80pt]{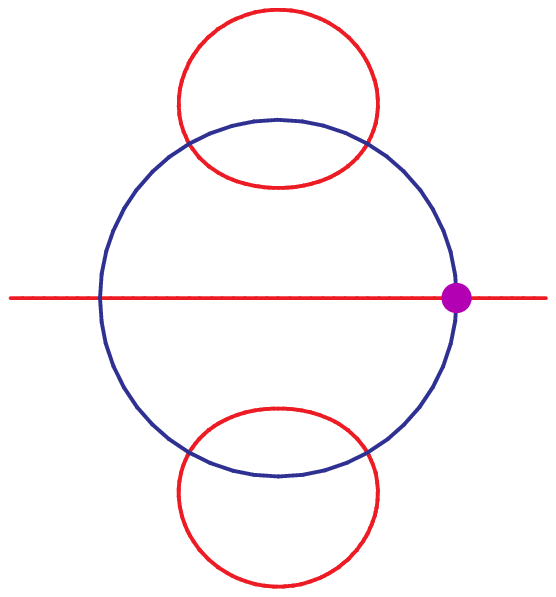}}\quad
   \includegraphics[height=65pt]{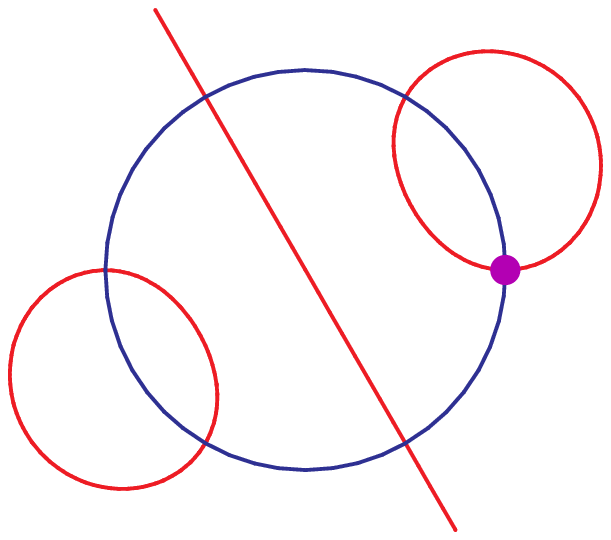}\quad
   \raisebox{4.2pt}{\includegraphics[height=65pt]{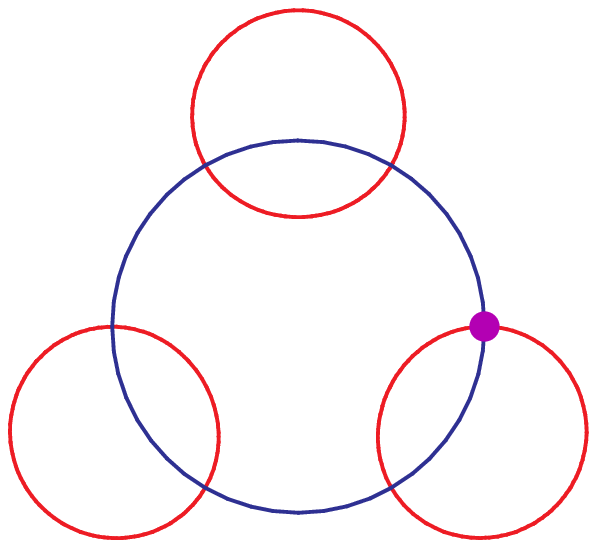}}\quad
   \raisebox{-9.8pt}{\includegraphics[height=65pt]{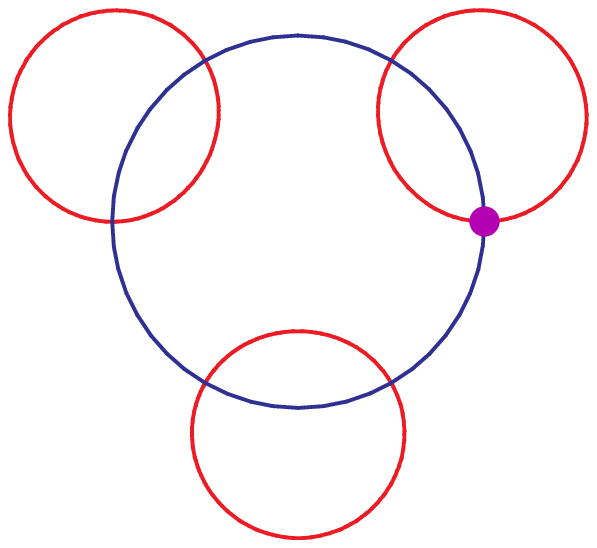}}}
 \end{equation}

The uniform convergence of Proposition~\ref{P9:Weak_cont} implies a certain continuity
of nets.
Two subsets $X,Y\subset\P^1$ have \DeCo{{\sl Hausdorff distance}} $\epsilon$ if every point of
$X$ lies within a distance $\epsilon$ of $Y$ and vice-versa.
This gives the \DeCo{{\sl Hausdorff metric}} on subsets of $\P^1$.

\begin{prop}\label{P9:Conv_nets}
 Let $\{\varphi_j\}\subset\St_{p{+}1}$ be a convergent sequence with limit $\varphi$.
 Then the sets $\{\pi(\varphi_j)^{-1}(\R\P^1)\}$ converge in the Hausdorff metric to
 the set $\{\pi(\varphi)^{-1}(\R\P^1)\}$.
\end{prop}

We deduce two corollaries from this proposition.

\begin{cor}\label{C9:analytic_continuation}
 Suppose that $\{\rho_z\mid z\in[0,1]\}$ is a continuous path in $R_{2p}$ where each
 $\rho_z$ has $2p$ critical points.
 Let $v_1(z)$ be a continuous function of $z$ which is equal to a critical point
 of $\rho_z$, for each $z$.
 Then the net of the pair
\[
   \bigl(\rho_z^{-1}(\R\P^1),\, v_1(z)\bigr)
\]
 does not depend upon $z$.
\end{cor}

\begin{cor}\label{C9:create_edge}
 Suppose that $\{\varphi_z\mid z\in[0,1]\}$ is a continuous path in the Stiefel
 manifold $\St_{2p}$.
 Suppose that for $0<z$ we have $\pi(\varphi_z)\in R_{2p}$ and $\pi(\varphi_z)$ has distinct
 critical points $v_1(z),v_2(z),\dotsc,v_n(z)$, 
 where $v_i(z)$ is a continuous function of $z$, 
 and also that, at $z=0$, we have $v_1(0)=v_2(0)$, but all other 
 critical points are distinct.
 Then the degree of $\pi(\varphi_z)$ is constant for $z\in(0,1]$, \Blue{and}
 $\deg(\pi(\varphi_0))<\deg(\pi(\varphi_1))$ if and only if the net of $\pi(\varphi_1)$
 has an interior edge between $v_1(1)$ and $v_2(1)$.
\end{cor}

Here are two nets for quartic rational functions~\eqref{Eq:quartic_nets}
as two of their vertices collide.
\[
  \begin{picture}(286,165)
   \put(0,88){\begin{picture}(281,77)(0,-2)
     \put(0,0){\includegraphics{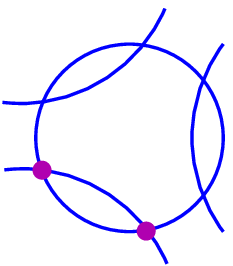}}
      \put(3,18){$v_2$}     \put(33,1){$v_1$}
      \put(75,35){$\longrightarrow$}
     \put(110,1){\includegraphics{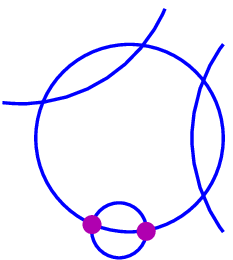}}
      \put(124,5){$v_2$}     \put(153,1){$v_1$}
      \put(188,35){$\longrightarrow$}
     \put(220,9){\includegraphics{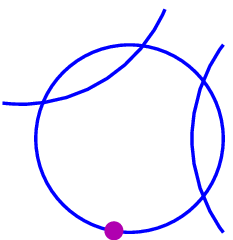}}
      \put(249,1){$v_1$}
   \end{picture}}
   \put(5,0){\begin{picture}(281,77)(0,-2)
     \put(0,0){\includegraphics{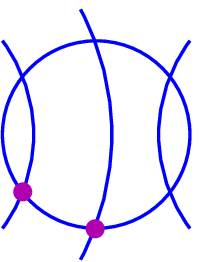}}
      \put(4,8){$v_2$}     \put(28,1){$v_1$}
      \put(70,35){$\longrightarrow$}
     \put(107,0){\includegraphics{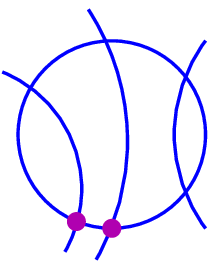}}
      \put(115,9){$v_2$}     \put(140,1){$v_1$}
      \put(182,35){$\longrightarrow$}
     \put(218,0){\includegraphics{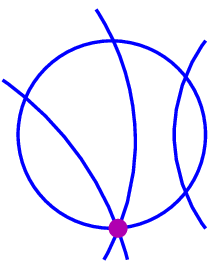}}
      \put(257,2){$v_1$}
   \end{picture}}
  \end{picture}
\]
 In the first row, there is an edge (in fact two) of $\Gamma\setminus\R\P^1$ between the
 vertices $v_2$ and $v_1$.
 This edge collapses in the limit as $v_2$ approaches $v_1$, which causes two regions of
 $\P^1\setminus\Gamma$ to collapse, so that the limiting net has two fewer regions in its
 complement, and thus corresponds to a rational function of degree $3$.
 There is no such edge in the nets of the second row, and the limiting net still has eight
 regions and thus its rational function still has degree four.\medskip 

\noindent{\it Proof of Corollary$~\ref{C9:create_edge}$.}
 The degree of a rational function $\rho\in R_{2p}$ is one-half the number of
 cells in the complement $\P^1\setminus \rho^{-1}(\R\P^1)$ of the net of $\rho$.
 Set $\rho_z:=\pi(\varphi_z)$.
 The only way for the number of cells in the complement of the net of $\rho_z$ to
 change at some $z_0\in[0,1]$ would be if some edge of $\rho_z^{-1}(\R\P^1)$
 collapsed as $z\to z_0$.  
 Since the vertices of $\rho_z^{-1}(\R\P^1)$ are the critical points 
 $v_1(z),\dotsc,v_n(z)$, which are distinct for $z\in(0,1]$,
 we see that the degree of $\rho_z$ is constant for $z\in(0,1]$.

 If the degree of $\rho_0$ is less than that of $\rho_1$, then some edges which bound a cell must
 disappear in the limit as $z\to 0$.
 But this can only happen if a cell is bounded by edges between $v_1(z)$ and $v_2(z)$,
 as they are the only critical points which collide in the limit as $z\to 0$.
 By (iii), such a cell must be bounded by more than one edge which implies that there was
 an interior edge between $v_1(1)$ and $v_2(1)$ outside of $\R\P^1$.
 This shows the necessity of such an interior edge between $v_1(1)$ and $v_2(1)$ for the degree
 to drop. 

 For sufficiency, note that if there is an interior edge between $v_1(1)$ and $v_2(1)$, then it
 must collapse in the limit as $z\to 0$ for otherwise 
 condition (iii) for nets would be violated.
\QED

%
\subsection{Schubert induction for rational functions}
In  Chapter~\ref{Ch:ScGr}, we used Schubert induction to construct a sequence of points 
$s_1,\dotsc,s_{mp}\in\R$ and sufficiently many real points in each
Schubert variety $X_\alpha\Fdot(0)$ which also lie in
$X_\I\Fdot(s_i)$ for $i=1,\dotsc,|\alpha|$.
Without re-running that proof, we will describe what that
construction gives for rational functions.

The construction of Theorem~\ref{T8:asymptotic} relevant for rational functions was
in the Grassmannian $\Gr(p,p{+}2)$.
Under the Grassmann duality of Theorem~\ref{T8:GrassmannDuality}, this becomes a 
construction in $\Gr(2,\C_{p+1}[t])$ and involves a Schubert variety 
$X_\alpha\Edot(t)$ with $\alpha\in\binom{[p{+}2]}{2}$.
In fact, the statements become identical after replacing $\Fdot(t)$ by $\Edot(t)$.
We will briefly recall these statements in this setting.

A point in the Schubert cell $X^\circ_\alpha\Edot(s)$ where $\alpha\colon\alpha_1<\alpha_2$
is a two-dimensional subspace $\langle f,g\rangle$ of polynomials of degree $p{+}1$ such
that
 \begin{equation}\label{eq:ord}
  \ord_s(f)\ =\ p{+}2{-}\alpha_2
   \qquad\mbox{and}\qquad
  \ord_s(g)\ =\ p{+}2{-}\alpha_1\,.
 \end{equation}
In particular,
\[
  (t-s)^{p+2-\alpha_2}|| f
   \qquad\mbox{and}\qquad
  (t-s)^{p+2-\alpha_1}|| g\,.
\]
(Here, $a^k||b$ means that $a^k$ divides $b$, but $a^{k+1}$ does not divide $b$.)

Given $\langle f,g\rangle\in X^\circ_\alpha\Edot(s)$, a consequence of~\eqref{eq:ord} is
that 
\[
   \Wr(f,g)\ =\ f'(t)g(t)-f(t)g'(t)
\]
vanishes to order $p{+}1{-}\alpha_2{+}p{+}2{-}\alpha_1=2p{-}|\alpha|$ at $s$.

By Lemma~\ref{L8:schubert_induction}, there exist real numbers $s_1,\dotsc,s_{2p}$
(in fact, we have $s_1>\dotsb>s_{2p}>0$) such that for all
$\alpha\in\binom{[p{+}2]}{2}$,  the intersection
 \begin{equation}\label{Eq9:Inductive_step}
   X_{\alpha}\Edot(0)\ \cap\ \bigcap_{i=1}^{|\alpha|} X_\I\Edot(s_i)
 \end{equation}
is transverse, and it consists of $\delta(\alpha)$ real points.
Any point $\langle f,g\rangle$ in the intersection~\eqref{Eq9:Inductive_step} will have
Wronskian
\[
  f'(t)g(t)-f(t)g'(t)\ =\ 
  \mbox{constant}\cdot t^{2p-|\alpha|}\cdot\prod_{i=1}^{|\alpha|}(t-s_i)
\]

In fact, as noted in Remark~\ref{R8:construct}, the proof of
Lemma~\ref{L8:schubert_induction} did much more.
Suppose that $|\alpha|>0$, and define
\[
  \DeCo{\beta^1}\ :=\ \alpha_1{-}1<\alpha_2
   \qquad\mbox{and}\qquad
  \DeCo{\beta^2}\ :=\ \alpha_1<\alpha_2{-}1\,,
\]
when possible.
($\beta^1$ is only defined if $1<\alpha_1$ and $\beta^2$ is only defined if
 $\alpha_1{+}1<\alpha_2$) 
Then the proof constructed $\delta(\alpha)=\delta(\beta^1)+\delta(\beta^2)$
families $\{\langle f_z,g_z\rangle\mid z\in[0,s_{|\alpha|}]\}$ of polynomials such that 

\begin{enumerate}
 \item For $z\neq 0$, $(f_z,g_z)\in X^\circ_\alpha$.
 \item ${\displaystyle  f_z'(t)g_z(t)-f_z(t)g_z'(t)\ =\ 
   \mbox{constant}\cdot t^{2p-|\alpha|}\cdot \Bigl(\ 
     \prod_{i=1}^{|\alpha|-1}(t-s_i)\Bigr)\cdot
      (t-z)}$.
 \item Exactly $\delta(\beta^i)$ of these families began in $X_{\beta^i}\Edot(0)$.
       That is,  for $\delta(\beta^i)$ of these families, we have 
       $(f_0,g_0)\in X_{\beta^i}\Edot(0)$.
\end{enumerate}

\subsection{Schubert induction for nets}

 The main idea in the proof is that the rational functions constructed in
 Lemma~\ref{L8:schubert_induction} each have different nets.

\begin{thm}\label{T9:Nets_Diff}
 The $\delta(\alpha)$ rational functions in $X^\circ_\alpha\Edot(0)$  constructed in
 Lemma~$\ref{L8:schubert_induction}$ each have  
 different nets. 
\end{thm}

 Suppose that $\langle f,g\rangle$ is a point in the intersection~\eqref{Eq9:Inductive_step}
 where $f$ and $g$ satisfy~\eqref{eq:ord} for $s=0$. 
 Then its Wronskian vanishes to order $2p{-}|\alpha|$ at $0$ and to order 1 at the points  
 $s_1,\dotsc,s_{|\alpha|}$.
 In particular, 0 is the only common zero of $f$ and $g$.
 Removing the common factor $t^{p+2-\alpha_2}$ from both $f$ and $g$ gives relatively prime
 polynomials of degree at most $\alpha_2{-}1$.
 Indeed, if $f$ and $g$ had a common root $s$, then a linear combination of them would vanish
 to order at least 2 and so their Wronskian would vanish to order at least $2$ at $s$.
 Then the rational function $\DeCo{r}:=f/g$ has degree $\alpha_2{-}1$ with Wronskian
\[
  \mbox{constant}\cdot t^{\alpha_2-\alpha_1-1}\cdot \prod_{i=1}^{|\alpha|}(t-s_i)\,.
\]

 The point  $\langle f,g\rangle$ corresponds to a unique path in the Bruhat order from
 $12$ to $\alpha$ in the Bruhat order.
 We claim that this path may be recovered from the net $\rho^{-1}(\R\P^1)$ of $\rho$.

 Indeed, consider the $i$th step in the construction, when the critical point $s_i$ was
 created. 
 By Corollary~\ref{C9:create_edge}, the interior edge from $s_i$ has other endpoint $0$ if the
 degree of the rational function increased at the $i$th step, and if its degree did not
 increase, then the other endpoint of that edge is at some critical point $s_k$ with $s_k>s_i$
 and so $k<i$.
 Subsequent steps in the construction will not affect an edge from $s_i$ to $s_k$ with $k<i$,
 but an edge between $0$ and $s_i$ may be moved to an edge between $s_i$ and $s_j$ where $j>i$.

 Thus, the degree of the rational function increased at step $i$ if and only if the other
 endpoint of an interior edge from $s_i$ is at $s_j$ with $j>i$.
 If $\beta\lessdot\beta'$ is the $i$th step in the chain corresponding to our rational function
 $\rho$, then 
\begin{enumerate}
 \item $\beta_2+1=\beta'_2$, so the degree of the rational function increased, if the interior
     edge from $s_i$ has endpoint $s_j$ with $j>i$ (so $s_j<s_i$), and
 \item $\beta_1+1=\beta'_1$, so the degree of the rational function did not increase, if the
     interior edge from $s_i$ has endpoint $s_k$ with $k<i$ (so $s_k>s_i$).
\end{enumerate}
 This completes the proof.
\QED\smallskip

\begin{figure}[!ht]
\[
\begin{picture}(420,585)(4,-17)
\thicklines

\put(335,-17){
  \begin{picture}(80,160)
   \put(-2,0){\includegraphics{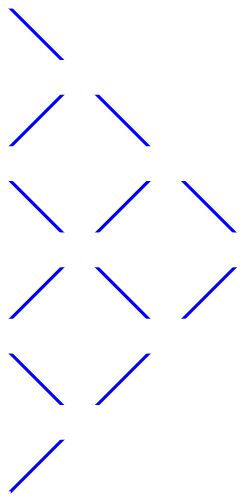}}
   \put( 0,150){45}
            \put(25,125){35}
   \put( 0,100){35}     \put(50,100){25}
            \put(25, 75){24}     \put(75, 75){15}
   \put( 0, 50){23}     \put(50,50){14}
            \put(25, 25){13}
   \put( 0,  0){12}
  \end{picture}
}

\put( 54,-13){\includegraphics{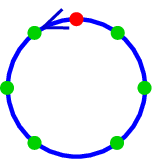}}%
\put(85, 29){$s_1$}
\put(99,  8){$s_2$}
\put(86,-16){$s_3$}
\put(55,-16){$s_4$}
\put(42,  8){$s_5$}
\put(53, 29){$s_6$}

\put( 94,24){\Blue{\vector(3,2){50}}}

\put(147,57){\includegraphics{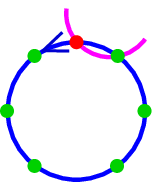}}%

\put(142,95){\Blue{\vector(-1,1){40}}}

\put(187,95){\Blue{\vector(3,2){60}}}

\put( 67.5,137){\includegraphics{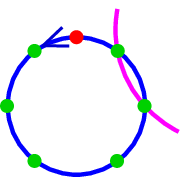}}%
\put(249,137){\includegraphics{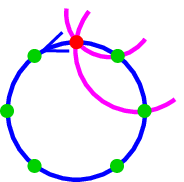}}%

\put(102,185){\Blue{\vector(1,2){19}}}
\put(252,185){\Blue{\vector(-1,1){42}}}
\put(298,185){\Blue{\vector(3,2){69}}}

\put(108,228.75){\includegraphics{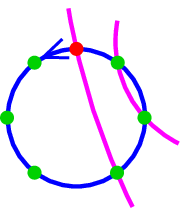}}%
\put(169,228){\includegraphics{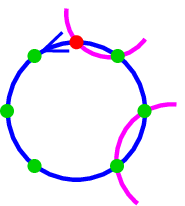}}%

\put(372,233.9){\includegraphics{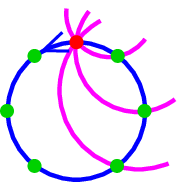}}%

\put(114,285){\Blue{\vector(-4,3){56}}}%
\put(173,285){\Blue{\vector(-4,3){56}}}%
\put(151,285){\Blue{\vector(4,3){56}}}%
\put(210,285){\Blue{\vector(4,3){56}}}%

\put(387,285){\Blue{\vector(-1,2){21}}}%

\put( 16,327){\includegraphics{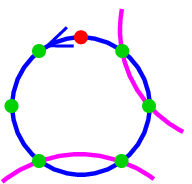}}%
\put( 77,327){\includegraphics{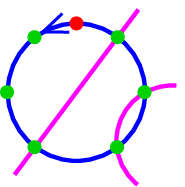}}%

\put(190,327){\includegraphics{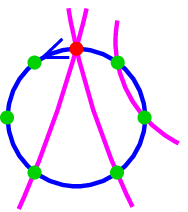}}%
\put(251,327){\includegraphics{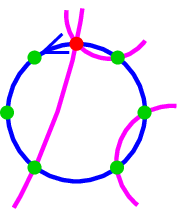}}%
\put(312,327){\includegraphics{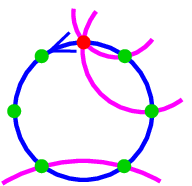}}%

\put( 39,384){\Blue{\vector(1,2){18}}}%
\put(100,384){\Blue{\vector(1,2){18}}}%

\put(193,384){\Blue{\vector(-1,2){18}}}%
\put(254,384){\Blue{\vector(-1,2){18}}}%
\put(315,384){\Blue{\vector(-1,2){18}}}%

\put( 30,423){\includegraphics{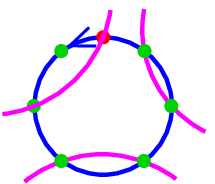}}%
\put( 91,419){\includegraphics{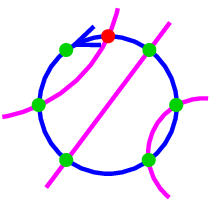}}%
\put(152,418){\includegraphics{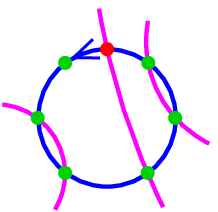}}%
\put(213,419){\includegraphics{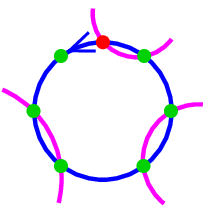}}%
\put(274,423){\includegraphics{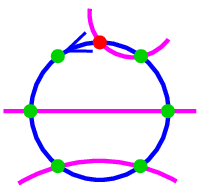}}%

\put( 51,472){\Blue{\vector(-1,2){18}}}%
\put(112,472){\Blue{\vector(-1,2){18}}}%
\put(173,472){\Blue{\vector(-1,2){18}}}%
\put(234,472){\Blue{\vector(-1,2){18}}}%
\put(295,472){\Blue{\vector(-1,2){18}}}%

\put(  4,510){\includegraphics{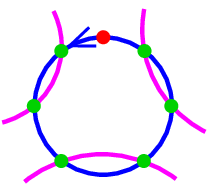}}%
\put( 65,507){\includegraphics{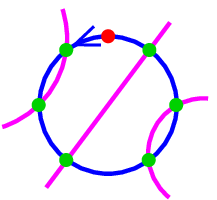}}%
\put(126,507){\includegraphics{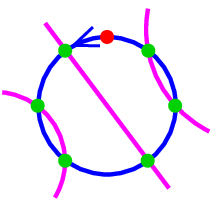}}%
\put(187,509){\includegraphics{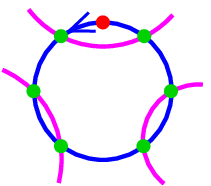}}%
\put(248,513){\includegraphics{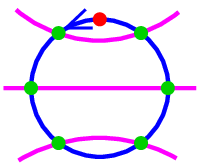}}%
\end{picture}
\]
\caption{Formation of nets during Schubert induction.}\label{F9:nets}
\end{figure}

Figure~\ref{F9:nets} illustrates the formation of the nets during the Schubert induction
for quartic rational functions, as well as the recursion for $\delta(\alpha)$.

We complete the proof of the Shapiro conjecture for rational curves.

\begin{thm}
 Let $\Phi(t)$ be a real polynomial of degree $2p$ whose roots are all real.
 Then there are exactly $\delta(p{+}1,p{+}2)$ real equivalence classes of rational
 functions with Wronskian $\Phi(t)$.
\end{thm}

\noindent{\it Proof.}
 Let $s_1>s_2>\dotsb>s_{2p}\in\R$ be the numbers such that 
 the intersection
\[
   \bigcap_{i=1}^{2p} X_\bI\Edot(s_i)
\]
 transverse with all points real.
 Each point in the intersection is an equivalence class of rational functions with
 Wronskian
\[
   \Phi_0(t)\ =\ \prod_{i=1}^{2p}(t-s_i)\,.
\]
 
 Let $\{\Phi_z\mid z\in[0,1]\}$ be a continuous family of  polynomials of degree $2p$ all
 with distinct real roots and with $\Phi_1(t)=\Phi(t)$. 
 We attempt to analytically continue each point in the fiber $\Wr^{-1}(\Phi_z)$ from
 $z=0$ to $z=1$.
 The only way this continuation could fail would be if it encountered a fiber 
 $\Wr^{-1}(\Phi_z)$ with a multiple  point.
 Then two of the rational functions in this fiber would have to coincide.
 In particular two would have the same net, (where we have labeled the
 nets by the root of the $\Phi_z(t)$ corresponding to $s_1$).
 This implies that two of the original rational functions in $\Wr^{-1}(\Phi_0)$ have
 the same net, by Corollary~\ref{C9:analytic_continuation}.
 But this contradicts Theorem~\ref{T9:Nets_Diff}.
\QED

%
%
\section{Rational functions with prescribed coincidences}
The results of Section~\ref{S:SCRF} can be used to prove a result about real rational functions
that satisfy a certain interpolation condition.
This is is due to Eremenko, Gabrielov, Shapiro, and Vainstein~\cite{EGSV}, and 
may be interpreted in the Grassmannian $\Gr(p,p{+}2)$ as an appealing
generalization of the Shapiro conjecture.
We will discuss this generalization, the Secant Conjecture, in Chapter~\ref{Ch:Frontier}. 

We first describe the interpolation problem.
Let $A_1,\dotsc,A_n$ be disjoint finite subsets of $\P^1$ where the set $A_i$ has $1+a_i$
elements with $1\leq a_i\leq p$ and $a_1+\dotsb+a_n=2p$.
Write $\ba$ for this sequence $(a_1,\dotsc,a_n)$ of numbers.
The interpolation problem is to determine the equivalence classes of rational functions $\rho$ of
degree $p{+}1$ that satisfy
\[
   \rho|_{A_i}\mbox{ is constant for }i=1,\dotsc,n\,.
\]
There are in fact finitely many such equivalence classes of rational functions when the sets
$A_i$ are general.
We will later prove this finiteness and show the number of equivalence classes is a Kostka
number $K_\ba$~\cite[p.25]{Fu97},\cite[I,6]{Mc95}.
Heuristically, there are finitely many equivalence classes because the
condition that a rational function is constant on a finite set of $1{+}a$ elements gives $a$
equations in the Stiefel coordinates for rational functions.

A collection of sets $A_i\subset\R\P^1$ for $i=1,\dotsc,n$ is \DeCo{{\sl separated}} if there
exist disjoint intervals $I_1,\dotsc,I_n$ of $\R\P^1$ with $A_i\subset I_i$ for
$i=1,\dotsc,n$.

\begin{thm}[\cite{EGSV}]\label{T9:EGSV}
 Let $\ba=(a_1,\dotsc,a_n)$ with $1\leq a_i\leq p$ and $a_1+\dotsb+a_b=2p$.
 For general separated subsets of $\R\P^1$, $A_1,\dotsc,A_n$ with 
 $|A_1|=1{+}a_i$, there are exactly $K_{\ba}$ real equivalence classes of rational
 functions $\rho$ such that 
 \begin{equation}\label{E9:constant}
   \rho|_{A_i}\mbox{ is constant for }i=1,\dotsc,n\,.
 \end{equation}
\end{thm}

Given separated subsets $A_1,\dotsc,A_n$ of $\R\P^1$, the proof first constructs such a
real rational function satisfying~\eqref{E9:constant} for every net with a certain
property, described below~\eqref{E9:hypothesis}.
Next, it relates this interpolation problem to a problem in the Schubert calculus that has
$K_\ba$ solutions, and finally shows that $K_\ba$ is the number of nets with the
property~\eqref{E9:hypothesis}. 
This will imply that we have constructed all the solutions for general sets $A_i$.

Theorem~\ref{T9:EGSV} generalizes Theorem~\ref{T9:EG}.
Suppose that we have subsets $\{A_z\mid z\in(0,1]\}$ of $\R\P^1$ depending continuously on
$z$, each of cardinality $a+1$, whose limit as $z\to 0$ consists of a single point, 
\[
  \lim_{z\to 0} A_z\ =\ \{s\}\,.
\]
Suppose further that we have a family $\{\rho_z\mid z\in[0,1]\}$ of rational functions that
depends continuously on $z$, and such that for $z>0$, $\rho_z$ is constant on $A$.
Then $\rho_0$ will have a critical point at $s$ of order at least $a$.

In this way, Theorem~\ref{T9:EGSV} implies Theorem~\ref{T9:EG} by simply considering the limit
as the points in each set $A_i$ collide.
In fact, this analysis will enable us to deduce a stronger form of Theorem~\ref{T9:EG}

\begin{thm}\label{T9:Pieri}
 Let $a_1,\dotsc,a_n$ be positive integers with $1\leq a_i\leq p$ and $a_1+\dotsb+a_n=2p$.
 Then every rational function of degree $p+1$ with $n$ real critical points of multiplicities 
 $a_1,\dotsc,a_n$ is real.
 There are exactly $K_\ba$ classes of such rational functions, and the corresponding Schubert
 varieties meet transversally.
\end{thm}

%
\subsection{Proof of Theorems~\ref{T9:EGSV} and~\ref{T9:Pieri}}

Let $\calR_{p+1}$ be the set of rational function of degree $p{+}1$ with exactly $2p$ critical
points of multiplicity 1.
 We use two consequences of the work in Section~\ref{S:SCRF}.
\begin{enumerate}
 \item If $\rho_1$ and $\rho_2$ are rational functions in $\calR_{p+1}$ with the same critical points
       and isotopic nets (where we use the same critical point for each net), then 
       $\rho_1$ is equivalent to $\rho_2$.

 \item For every net $\Gamma\subset\P^1$ with a given vertex set $V$ (and distinguished vertex
 $v_1\in V$), there is a unique equivalence class of rational functions in $\calR_{p+1}$ with
 critical set $V$ and net (with distinguished vertex $v_1$) isotopic to $\Gamma$.
\end{enumerate}

\begin{figure}[htb]
\[
  \begin{picture}(152,136)(-2,1)
   \put( 15, 13){\includegraphics{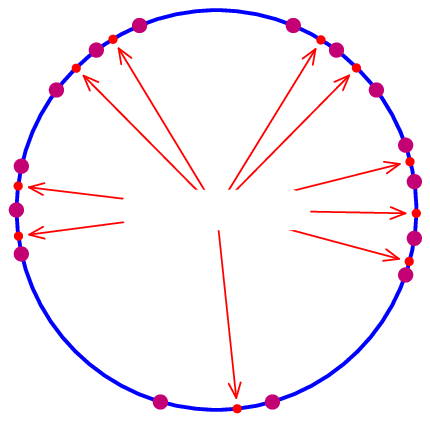}}
   \put( 49.5, 73){$s_1,\dotsc,s_{10}$}
   \put(137, 75){$A_1$}
   \put(108,126){$A_2$}
   \put( 21,125){$A_3$}
   \put( -2, 73){$A_4$}
   \put( 70,  2){$A_5$}
 \end{picture}
\]
\caption{Interlacing points.}
\label{F9:interlace}
\end{figure}

Actually, we only showed the second point for nets with local degree 2 at each critical
point.
The stronger statement follows from results in~\cite{EG02}.
To prove Theorem~\ref{T9:EGSV}, fix separated subsets $A_1,\dotsc,A_n$ of $\R\P^1$
satisfying the hypotheses.
Choose $2p$ additional points $s_1,\dotsc,s_{2p}$ where, for each $i$,  $a_i$ of the points
interlace the $a_i{+}1$ points of $A_i$.
Figure~\ref{F9:interlace} shows an example when $p=5$ and $\ba=(3,2,2,2,1)$.
Each point $s_j$ lies between two points of some set $A_i$.
Write $[x_j,y_j]$ for the interval that contains $s_j$ and note that $x_j,y_j\in A_i$, 
for some $i$.

Consider nets with the vertices $s_1,\dotsc,s_{2p}$ that have local degree $2$ at each
vertex and satisfy the additional hypothesis:
 \begin{equation}\label{E9:hypothesis}
  \mbox{There are no edges between points $s_j$ and $s_k$ that interlace the same set
  $A_i$}
 \end{equation}
There are five nets satisfying~\eqref{E9:hypothesis} for the points $A_i$ of
Figure~\ref{F9:interlace} (we only draw the edges in the upper half plane, which is the
interior of the circle). 
 \begin{equation}\label{Eq9:interpolate_nets}
   \raisebox{-35pt}{
   \includegraphics[height=70pt]{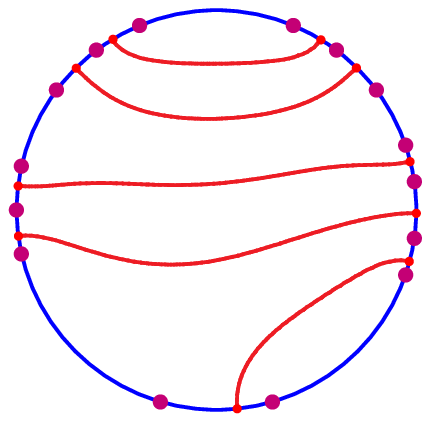}\quad
   \includegraphics[height=70pt]{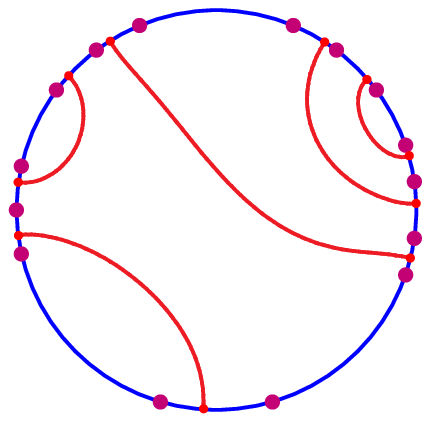}\quad
   \includegraphics[height=70pt]{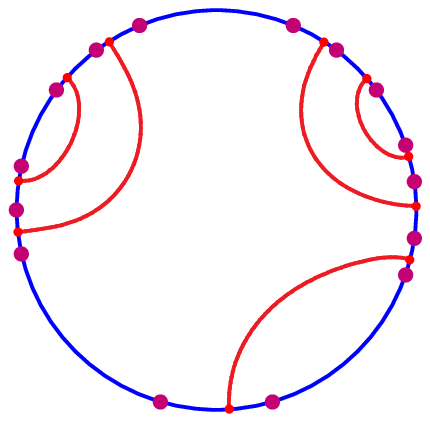}\quad
   \includegraphics[height=70pt]{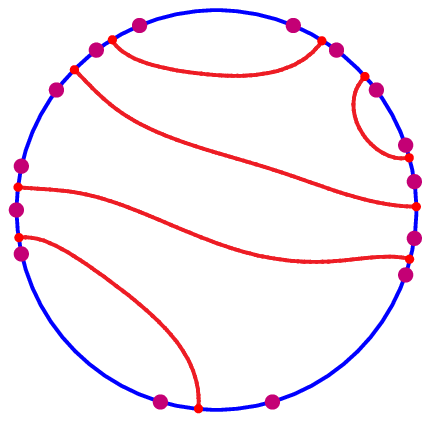}\quad
   \includegraphics[height=70pt]{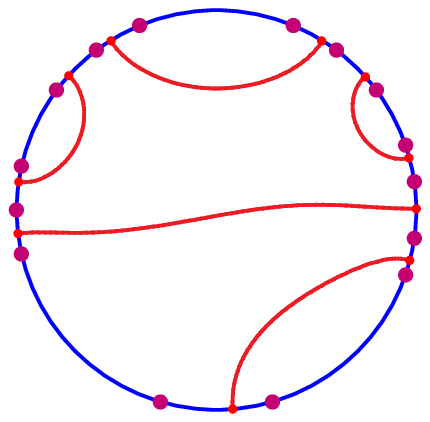}}
 \end{equation}

Suppose that we have critical points $s_1,\dotsc,s_{2p}$ interlacing the sets
$A_i$, where all the points $s_i$ are fixed and chosen arbitrarily within their intervals,
except for one, $s_j$, which is allowed to vary within its interval $[x_j,y_j]$.
Also fix an isotopy class $\Gamma$ of nets satisfying~\eqref{E9:hypothesis}.
For each $s\in[x_j,y_j]$ let $\rho_s\in\calR_{p+1}$ be a rational function with the critical
points $s_1,\dotsc,s_{2p}$ and net $\Gamma$.
We may suppose that $\{\rho_s\mid s\in[x_j,y_j]\}$ is a continuous family.

\begin{lemma}\label{L9:interpolate}
  There exists a point $s_j\in[x_j,y_j]$ such that $\rho_{s_j}(x_j)=\rho_{s_j}(y_j)$.
\end{lemma}

To see this, we may assume that $\rho_s$ is normalized so that $\rho_s(s)=0$ and $\rho_s$ maps the
interior edge of $\Gamma$ terminating in $s_j$ to the interval $[-\infty,0]$ in $\R\P^1$. 
Then the difference
\[
    \rho_s(x_j)\ -\ \rho_s(y_j)
\]
is positive when $s$ is near $y_j$ and negative when $s$ is near $x_j$, so it changes sign
on the interval $[x_j,y_j]$ and therefore takes the value zero at some $s_j\in[x_j,y_j]$.
We give a picture below.
The arrows point in the direction of increase of $\rho_s(t)$ for $t$ in the net.
\[
  \begin{picture}(120,76)(0,-13)
   \put(0,0){\includegraphics{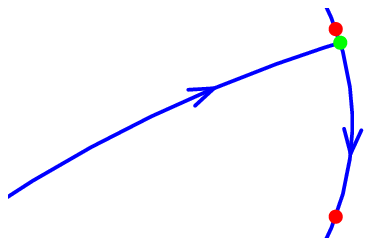}}
   \put(80,3){$x_j$} \put(80,60){$y_j$}
   \put(99,53){$s$}
   \put(0,-13){$\rho_s(x_j)-\rho_s(y_j)>0$}
  \end{picture}
   \qquad
  \begin{picture}(120,76)(0,-13)
   \put(0,0){\includegraphics{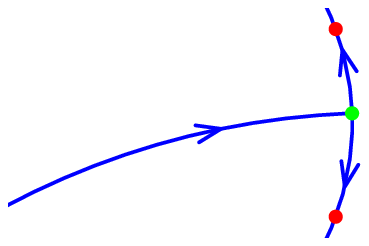}}
   \put(80,3){$x_j$} \put(80,60){$y_j$}
   \put(103,33){$s$}
  \end{picture}
   \qquad
  \begin{picture}(120,76)(0,-13)
   \put(0,0){\includegraphics{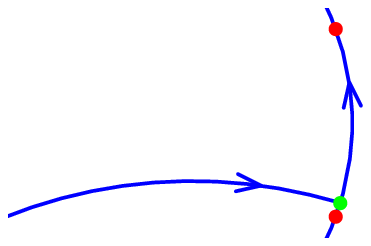}}
   \put(80,3){$x_j$} \put(80,60){$y_j$}
   \put(99,8){$s$}
   \put(0,-13){$\rho_s(x_j)-\rho_s(y_j)<0$}
  \end{picture}
\]

\begin{lemma}\label{L9:reality}
  If\/ $\Gamma$ is a net satisfying~\eqref{E9:hypothesis}, then there is a choice of critical
  points $s_1,\dotsc,s_{2p}$ interlacing the points of the sets $A_i$ such that every rational
  function $\rho$ of degree $p{+}1$ with the net $\Gamma$ and critical points $s_j$ satisfies
 \begin{equation}\label{E9:reprise}
   \rho|_{A_i}\mbox{ is constant for }i=1,\dotsc,n\,.
 \end{equation}
\end{lemma}

The set of possible critical points $\bs=(s_1,\dotsc,s_{2p})\in(\R\P^1)^{2p}$ interlacing the
sets $A_i$ forms the interior of a closed cube
\[
   Q\ :=\ [x_1,y_1]\times[x_2,y_2]\times\dotsb\times[x_{2p},y_{2p}]\,.
\]
By Lemma~\ref{L9:interpolate}, for every $j$, the function 
$\varphi_j(\bs):=\rho_{\bs}(x_j)-\rho_{\bs}(y_j)$ (defined as described in the proof of
Lemma~\ref{L9:interpolate}) is positive on the face $s_j=y_j$ and negative on the face
$s_j=x_j$.
This implies (via variant of Brower's fixed point theorem) that there is a point $\bs$ in the
interior of $Q$ where $\varphi_j(\bs)=0$ for all $j$, that is
$\rho_{\bs}(x_j)=\rho_{\bs}(y_j)$ for all $j$.
Since these intervals interlace the sets $A_i$, this implies~\eqref{E9:reprise}.
\QED

The next step is to show that the number of nets satisfying~\eqref{E9:hypothesis} for sets
$A_1,\dotsc,A_n$ where $A_i$ has $1{+}a_i$ members and $a_1+\dotsb+a_n=2p$ is the Kostka number
$K_{\ba}$. 
This Kostka number is the number of Young tableaux of shape $2\times p$ and content
$\ba$~\cite[p.25]{Fu}. 
These are arrays consisting of two rows of integers, each of length $p$ such that the integers
increase weakly across each row and strictly down each column.
For example, here are the five Young tableaux of shape $2\times p$ and content
$(3,2,2,2,1)$, showing that $K_{(3,2,2,2,1)}=5$.
 \begin{equation}\label{Eq9:tableaux}
  \raisebox{-11pt}{
  \begin{picture}(68,29)(-3.4,-3)
   \put(-3.4,-3){\includegraphics{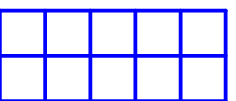}}
   \put(0,13){1} \put(13,13){1} \put(26,13){1} \put(39,13){2} \put(52,13){2}
   \put(0, 0){3} \put(13, 0){3} \put(26, 0){4} \put(39, 0){4} \put(52, 0){5}
  \end{picture}
 \quad
  \begin{picture}(63,27)(-3.4,-3)
   \put(-3.4,-3){\includegraphics{figures/9/YT.eps}}
   \put(0,13){1} \put(13,13){1} \put(26,13){1} \put(39,13){3} \put(52,13){4}
   \put(0, 0){2} \put(13, 0){2} \put(26, 0){3} \put(39, 0){4} \put(52, 0){5}
  \end{picture}
 \quad
  \begin{picture}(63,27)(-3.4,-3)
   \put(-3.4,-3){\includegraphics{figures/9/YT.eps}}
   \put(0,13){1} \put(13,13){1} \put(26,13){1} \put(39,13){3} \put(52,13){3}
   \put(0, 0){2} \put(13, 0){2} \put(26, 0){4} \put(39, 0){4} \put(52, 0){5}
  \end{picture}
 \quad
  \begin{picture}(63,27)(-3.4,-3)
   \put(-3.4,-3){\includegraphics{figures/9/YT.eps}}
   \put(0,13){1} \put(13,13){1} \put(26,13){1} \put(39,13){2} \put(52,13){4}
   \put(0, 0){2} \put(13, 0){3} \put(26, 0){3} \put(39, 0){4} \put(52, 0){5}
  \end{picture}
 \quad
  \begin{picture}(63,27)(-3.4,-3)
   \put(-3.4,-3){\includegraphics{figures/9/YT.eps}}
   \put(0,13){1} \put(13,13){1} \put(26,13){1} \put(39,13){2} \put(52,13){3}
   \put(0, 0){2} \put(13, 0){3} \put(26, 0){4} \put(39, 0){4} \put(52, 0){5}
  \end{picture}}
 \end{equation}

We only describe the map from nets to Young tableaux.
Given a net satisfying~\eqref{E9:hypothesis}, we will successively place 
integers into a left-justified two-rowed array while traversing $\R\P^1$.
This starts from the first (in the canonical ordering on $\R\P^1$) point in $A_1$ and 
begins with an empty array.
When a critical point $s_j$ is encountered, there will be an interior edge of the net
with endpoint $s_j$.
Place the integer $j$ in the second row if the other endpoint of that edge has already
been encountered, and in the first row if it has not been encountered.
For example, the tableaux in~\eqref{Eq9:tableaux} correspond, in order, to the nets
in~\eqref{Eq9:interpolate_nets}. 
(There, the order on $\R\P^1$ is counterclockwise on the circles.)

This bijection shows that we have constructed $K_{\ba}$ equivalence classes of rational
functions satisfying~\eqref{E9:constant}.
To complete the proof of Theorem~\ref{T9:EGSV}, we first show that the Kostka number $K_\ba$ is
the expected number of equivalence classes of rational functions
satisfying~\eqref{E9:constant}, and then that there is some choice of the sets $A_i$ for which
there are exactly $K_{\ba}$ equivalence classes of rational functions.
This last step will also prove Theorem~\ref{T9:Pieri}.

Recall that a polynomial $f$ of degree $p{+}1$ may be considered to be a linear map on 
$\C^{p+2}$ so that the composition with the rational normal curve
$\gamma(t)\colon\C\to\C^{p+2}$ gives the polynomial $f(t)$.
We used this to relate ramification to osculating flags as in Section~\ref{S8:duality}.
A two-dimensional space $\langle f,g\rangle$ of polynomials gives a map
\[
   \C\ \xrightarrow{\ \gamma(t)\ }\ 
   \C^{p+2}\ \xrightarrow{\ (f,g)\ }\ \C^2\,.
\]
The kernel $H$ of the map $\C^{p+2} \xrightarrow{\ ( f,g)\ }\C^2$ 
corresponds to $\langle f,g\rangle$  under Grassmann duality.

Suppose that the rational function $\rho=f/g$ is constant on some finite subset,
$A$.
This means that the line $(f(a),g(a))\subset\C^2$ is constant, for $a\in A$.
This implies that in $\C^{p+2}$ we have
\[
   \langle H, \gamma(a)\rangle\ =\  \langle H, \gamma(b)\rangle\,,
\]
for any $a,b\in A$.
In particular, $H$ has exceptional position with respect to the $|A|$-plane
$\DeCo{S(A)}:=\langle \gamma(a)\mid a\in A\rangle$ in that the two subspaces do not span
$\C^{p+2}$. 

Thus, the equivalence classes of rational functions $\rho$ of degree $p+1$ that
satisfy~\eqref{E9:reprise} correspond to the $p$-planes $H$ in $\C^{p+2}$ such that
 \begin{equation}\label{Eq9:span}
   \Span (H, S(A_i))\ \neq\ \C^{p+2}\qquad
   \mbox{for}\ i=1,\dotsc,n\,.
 \end{equation}
Those $H$ which satisfy~\eqref{Eq9:span} are the points of an intersection of Schubert
varieties.
Let $\DeCo{\sigma_a}:=(2,3,\dotsc,a{+}1,\; a{+}3,\dotsc,p{+}2)\in\binom{[p{+}2]}{p}$.
Then $X_{\sigma_a}\Fdot$ consists of those $H\in\Gr(p,p{+}2)$ such that 
\[
  \Span (H, F_{a+1})\ \neq\ \C^{p+2}\,.
\]
We will also write \DeCo{$X_{\sigma_a}F_{a{+}1}$} for this Schubert variety, which has
dimension $|\sigma_a|=2p-a$.
Thus the solutions to the interpolation problem~\eqref{E9:reprise} correspond to the
intersection of Schubert varieties
\[
   X_{\sigma_{a_1}}S(A_1)\;\cap\;
   X_{\sigma_{a_2}}S(A_2)\;\cap\; \dotsb\;\cap\;
   X_{\sigma_{a_n}}S(A_n)\,,
\]
which is expected to be zero dimensional.
These are special Schubert varieties, so the expected number of points in this
intersection may be computed using the Pieri formula, and it is the Kostka number
$K_{\ba}$~\cite[p.25]{Fu}. 

All that remains to show is that there is some choice of the sets $A_i$ for which there
are finitely many equivalence classes of rational functions satisfying the interpolation
condition~\eqref{E9:reprise}. 
We show that indirectly, by passing to the limit as each set $A_i$ collapses to a single
point, $s_i$.
If we consider the rational functions for a given net in this limit, we see that the
limiting rational function still has degree $p{+}1$, by Corollary~\ref{C9:create_edge} as no
interior edges were collapsed in the limit, by Condition~\eqref{E9:hypothesis}.
The limiting rational function has a critical point at each $s_i$ of multiplicity $a_i$, 
and is necessarily real.

It is easy to see that there are still $K_\ba$ nets that have critical points of
multiplicity $a_i$ at points $s_i$---the same bijection works.
For example, here are the nets of rational functions of degree five with critical points
having multiplicities $(3,2,2,2,1)$.
\[
  \includegraphics[height=75pt]{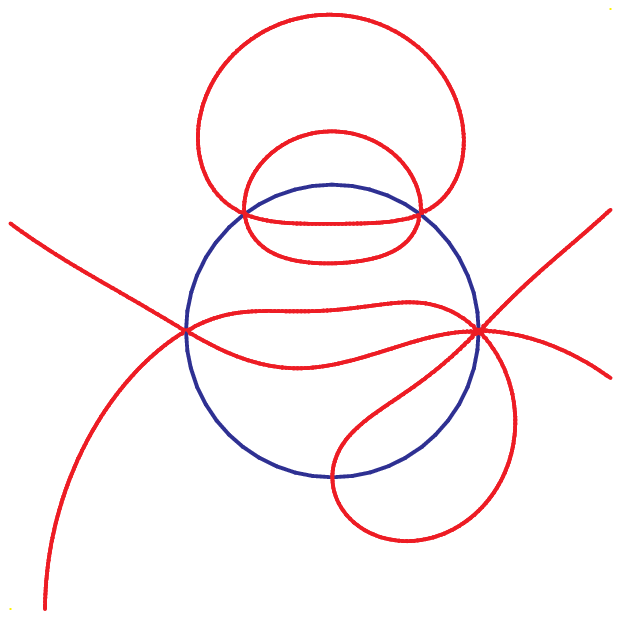}\quad
  \includegraphics[height=75pt]{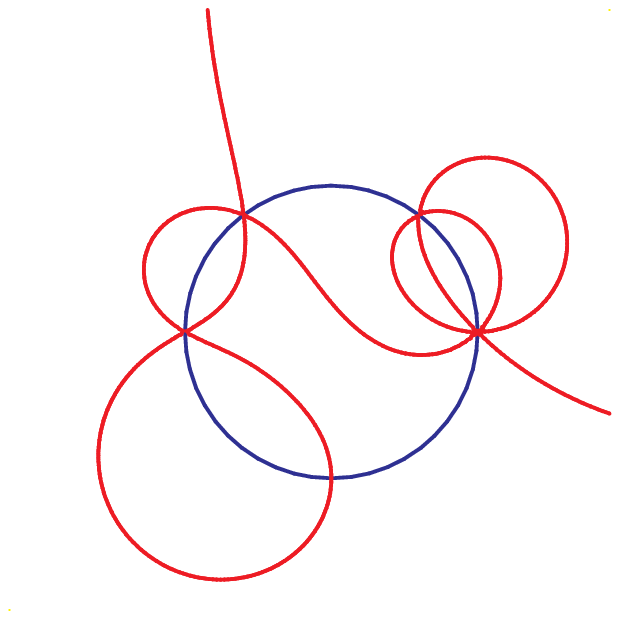}\quad
  \includegraphics[height=75pt]{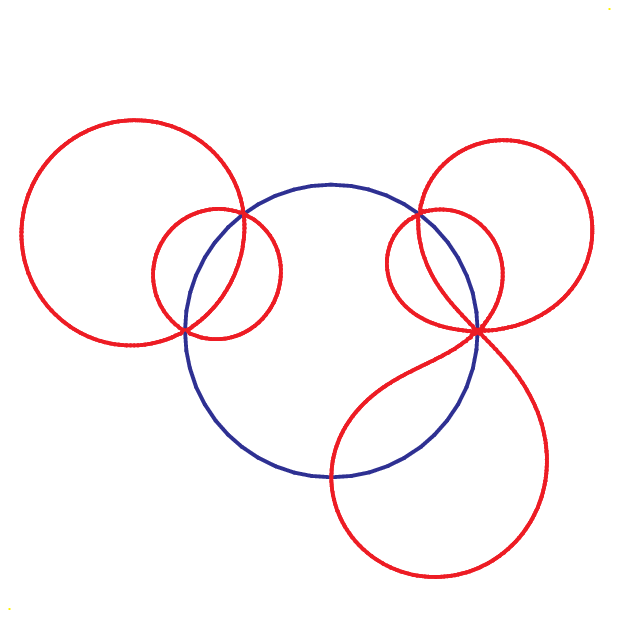}\quad
  \includegraphics[height=75pt]{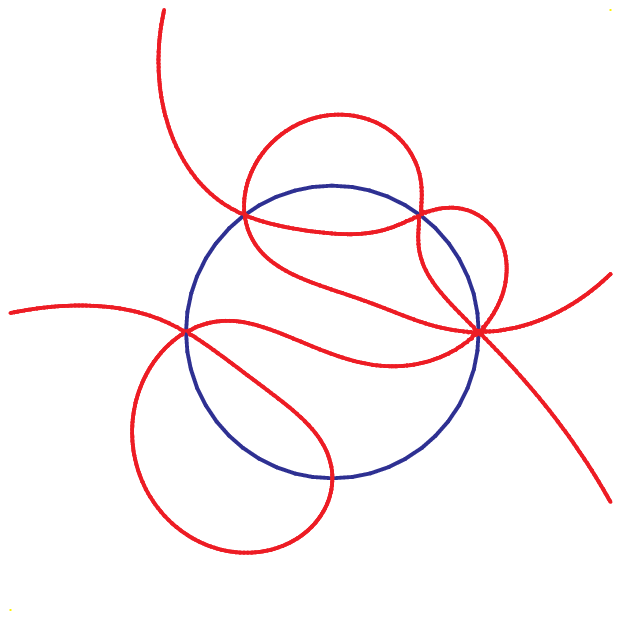}\quad
  \includegraphics[height=75pt]{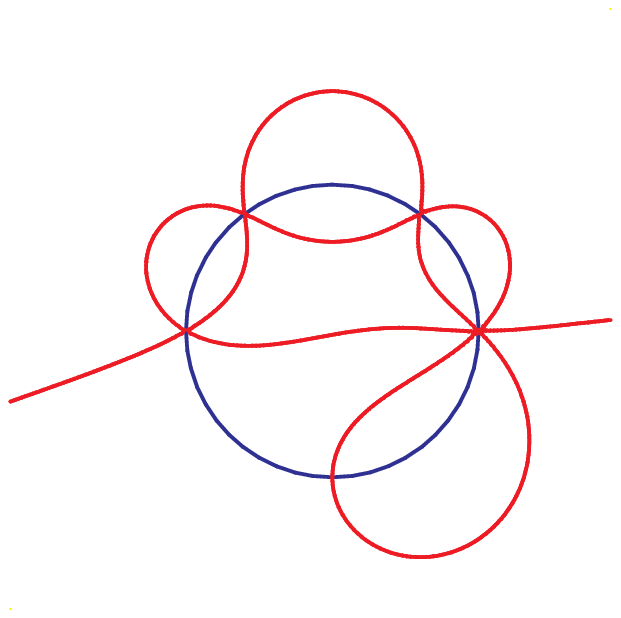}
\]
Moreover, the corresponding intersection of Schubert varieties is expected to have 
$K_\ba$ points.
Thus we have constructed the expected number of real rational functions with the desired
critical points. 
The proof of Theorem~\ref{T9:Pieri} now follows from the results of~\cite{EG02}, which gave a
bijective correspondence between nets and rational functions.

\newcommand{\bx}{{\bf x}}
\newcommand{\calS}{{\mathcal{S}}}
\newcommand{\dlog}{{\rm ln}'}
\newcommand{\frakh}{{\mathfrak{h}}}
\newcommand{\frakn}{{\mathfrak{n}}}
\newcommand{\sing}{{\rm sing}}
\newcommand{\glm}{{\mathfrak{gl}_m}}
\newcommand{\bM}{{\bf M}}
\newcommand{\bz}{{\bf z}}

\chapter{Proof of the Shapiro conjecture}\label{ch:MTV}

The Shapiro Conjecture was proven by Mukhin, Tarasov, and Varchenko in a preprint of
November 2005~\cite{MTV_Sh}.
Like the proofs of the special case of $m=2$ by Eremenko and Gabrielov~\cite{EG02,EG05}
(discussed in Chapter~\ref{Ch:EG}) the proof in the general case did not use much
algebraic geometry.
Instead it used results in mathematical physics, specifically integrable systems, with
some representation theory.
Unlike the rest of the material in this book, this chapter is not elementary.
It does contain a sketch of some of the main ideas in their proof, but it by no means
complete, and we recommend that the serious reader go to the original sources.
An expanded account of the Shapiro conjecture and its proof appeared in 
the Bulletin of the AMS~\cite{FRSC}.
In fact, what follows is a mildly revised version of Sections 2, 3, and 4 of~\cite{FRSC}.

%
\section{Spaces of polynomials with given Wronskian}\label{S:polys}

By Theorem~\ref{T8:real_limit}, the general case of the Shapiro conjecture follows from
the special case when all the Schubert conditions are equal to $\bI$\hspace{.5pt}, and this
case is equivalent to the Wronski formulation of Theorem~\ref{T1:Shap_Conj}.
A further reduction is possible, as the Wronski map
$\Wr\colon\Gr(m,m{+}p)\to\P^{mp}$ is a finite map in that it has finite fibers, 
a standard limiting argument (given, for example, in
Section 1.3 of~\cite{MTV_Sh} or Remark 3.4 of~\cite{So00c}) shows that it suffices to prove 
Theorem~\ref{T1:Shap_Conj} when the Wronskian has distinct real roots that are  
sufficiently general.
Since  $\#_{m,p}$ is the upper bound for the number of spaces of polynomials 
with a given Wronskian, it suffices to construct this number of distinct spaces of real
polynomials with a given Wronskian, when the Wronskian has distinct real roots that are 
sufficiently general.  
In fact, this is exactly what Mukhin, Tarasov, and Varchenko do~\cite{MTV_Sh}.\medskip 

\noindent{\bf Theorem~\ref{T1:Shap_Conj}\boldmath{$'$}.}
 {\it
  If $s_1,\dotsc,s_{mp}$ are generic real numbers, there are 
  $\#_{m,p}$ real spaces of polynomials in $\Gr(m,\C_{m{+}p{-}1}[t])$ with 
  Wronskian $\prod_{i=1}^{mp}(t-s_i)$.
}\medskip

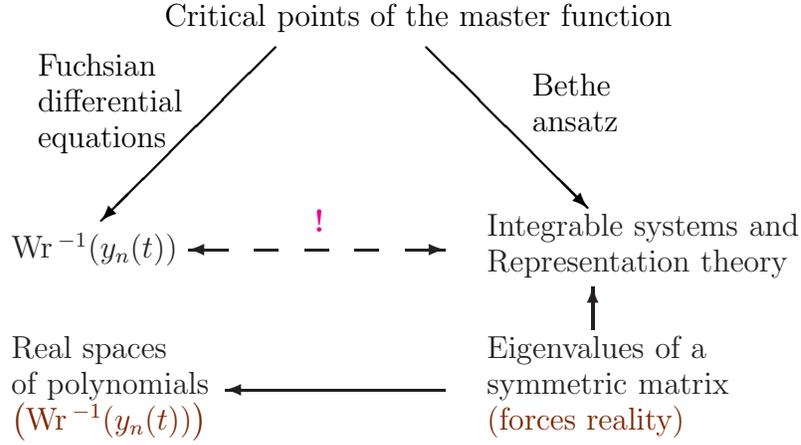
\begin{figure}[t]
\[
  \begin{picture}(286,160)(0,15)\thicklines   

    \put(58,168){Critical points of the master function}  

    \put(100,160){\vector(-1,-1){66}}\put(157,160){\vector(1,-1){61}}

    \put(10,122){\begin{minipage}[b]{70pt}Fuchsian\vspace{-1pt}\newline 
                          differential\vspace{-1pt}\newline equations\end{minipage}}

    \put(197,128){\begin{minipage}[b]{70pt}Bethe\vspace{-1pt}\newline  ansatz\end{minipage}}

    \put(114.5,90){{\bf \Magenta{!}}}

    \put(  0,80){$\Wr^{-1}(y_n(t))$}

    \put( 82,83){\vector(-1,0){15}}\put( 92,83){\line(1,0){9}}
    \put(111,83){\line(1,0){9}}    \put(130,83){\line(1,0){9}}
    \put(149,83){\vector(1,0){15}}

    \put(180,75){\begin{minipage}[b]{130pt} Integrable systems and\vspace{-1pt}\newline
                   Representation theory\end{minipage}}

    \put(220,53){\vector(0,1){16}}

    \put(180, 15){\begin{minipage}[b]{100pt} Eigenvalues of a\vspace{-1pt}\newline symmetric matrix\vspace{-1pt}\newline
                  \Brown{(forces reality)}\end{minipage}}

    \put(164,30){\vector(-1,0){83}}

    \put( 0, 15){\begin{minipage}[b]{75pt} Real spaces\vspace{-1pt}\newline of polynomials\vspace{-1pt}\newline  
             \Brown{$\bigl(\Wr^{-1}(y_n(t))\bigr)$}\end{minipage}}

  \end{picture}
\]
\caption{Schematic of proof of Shapiro conjecture.}
\label{F:schematic}
\end{figure}

The proof first constructs $\#_{m,p}$ distinct spaces of
polynomials with a given Wronskian having generic complex roots, 
which we describe in Section~\ref{S:construction}.
This uses a Fuchsian differential equation given by the critical points of a remarkable
symmetric function, called the master function.
The next step uses the Bethe ansatz in a certain representation $V$ of $\slm\C$: 
each critical point of the master function gives a Bethe eigenvector of the Gaudin Hamiltonians
which turns out to be a highest weight vector for an irreducible submodule of $V$. 
This is described in Section~\ref{S:BAGM}, where the eigenvalues of the 
Gaudin Hamiltonians on a Bethe vector are shown to be the coefficients of the Fuchsian
differential equation giving the corresponding spaces of polynomials.
This is the germ of the new, deep connection between representation theory and Schubert
calculus that led to the full statement of Theorem~\ref{Th7:MTV} (reality and transversality).
Finally, the Gaudin Hamiltonians are
real symmetric operators when the Wronskian has only real roots, so their
eigenvalues are real, and thus the Fuchsian differential equation
has real coefficients and the corresponding space of polynomials is also real.
Figure~\ref{F:schematic} presents a schematic of this extraordinary proof.

%
\subsection{Critical points of master functions}\label{S:construction}
The construction of $\#_{m,p}$ spaces of polynomials with a given Wronskian
begins with the critical points of a symmetric rational function that 
arose in the study of hypergeometric solutions to the Knizhnik-Zamolodchikov
equations~\cite{SV} and the Bethe ansatz method for the Gaudin model.

The master function depends upon parameters
$\DeCo{\bs}:=(s_1,\dotsc,s_{mp})$, which are 
the roots of our Wronskian $\Phi$, and an additional $p\binom{m}{2}$ 
variables 
\[
    \DeCo{\bx}\ :=\ (x_1^{(1)},\dotsc,x_{p}^{(1)},\,
              x_1^{(2)},\dotsc,x_{2p}^{(2)},\, \dotsc\,,\,
              x_1^{(m-1)},\dotsc,x_{(m-1)p}^{(m-1)})\,.
\]
Each set of variables $\bx^{(i)}:=(x_1^{(i)},\dotsc,x_{ip}^{(i)})$ will turn out to be the
roots of certain intermediate Wronskians.

Define the \DeCo{{\sl master function} $\Xi(\bx;\bs)$} by the (rather formidable) formula
 \begin{equation}\label{Eq:MasterFunction}
   \frac{\displaystyle \prod_{i=1}^n\ 
    \prod_{1\leq j<k\leq ip}(x_j^{(i)}-x_k^{(i)})^2
    \ \cdot\ \prod_{1\leq j<k<mp}(s_j-s_k)^2}{\displaystyle
   \prod_{i=1}^{n-1}\ \prod_{j=1}^{ip}\ \prod_{k=1}^{(i+1)p}
   (x_j^{(i)}-x_k^{(i+1)})\ \cdot\ 
   \prod_{j=1}^{(m-1)p}\ \prod_{k=1}^{mp}(x_j^{(m-1)}-s_k)}\ .
 \end{equation}
This is separately symmetric in each set of variables $\bx^{(i)}$.
The  Cartan matrix for $\slm$ appears in the exponents of the factors
$(x_*^{(i)}-x_*^{(j)})$ in~\eqref{Eq:MasterFunction}.
This hints at the relation of these master functions to Lie theory, which
we do not discuss.

The critical points of the master function are solutions to the system of equations
 \begin{equation}\label{Eq:Critical}
   \frac{1}{\Xi}\frac{\partial}{\partial x_j^{(i)}}\Xi(\bx;\bs)\ =\ 0
   \qquad\mbox{for}\quad i=1,\dotsc,m{-}1,\quad
                         j=1,\dotsc,ip\,.
 \end{equation}
When the parameters $\bs$ are generic, these \DeCo{{\sl Bethe ansatz equations}} turn out
to have finitely many solutions. 
The master function is invariant under the group 
\[
   \DeCo{\calS}\ :=\ 
    \calS_p\times\calS_{2p}\times\,\dotsb\,\times\calS_{(m-1)p}\,,
\]
where $\calS_N$ is the group of permutations of $\{1,\dotsc,N\}$, and 
the factor $\calS_{ip}$ permutes the variables in $\bx^{(i)}$.
Thus $\calS$ acts on the critical points.
The invariants of this action are polynomials whose roots are 
the coordinates of the critical points.

Given a critical point $\bx$, define monic polynomials $\bg_\bx:=(g_1,\dotsc,g_{m-1})$ 
where the components $\bx^{(i)}$ of $\bx$ are the roots of $g_i$,
 \begin{equation}\label{eq:crit_polys}
    \DeCo{g_i}\ :=\ \prod_{j=1}^{ip} (t-x_j^{(i)})
    \qquad\mbox{for}\quad i=1,\dotsc,m{-}1\,.
 \end{equation}
Also write $\DeCo{g_m}$ for the Wronskian, the monic polynomial with roots $\bs$.
The master function is greatly simplified by this notation. 
The discriminant \DeCo{$\mbox{Discr}(f)$} of a polynomial $f$ is the square of
the product of differences of its roots and the resultant \DeCo{$\mbox{Res}(f,h)$} is the
product of all differences of the roots of $f$ and $h$~\cite{CLO}.
Then the formula for the master function~\eqref{Eq:MasterFunction} becomes
 \begin{equation}\label{Eq:Resultant}
  \Xi(\bx;\bs)\ =\ \prod_{i=1}^m \mbox{Discr}(g_i) 
                    \Bigg/  \prod_{i=1}^{m-1} \mbox{Res}(g_i,g_{i+1}) \ .
 \end{equation}

The connection between the critical points of $\Xi(\bx;\bs)$ and spaces of
polynomials with Wronskian $\Phi$ is through a Fuchsian differential equation.
Given (an orbit of) a critical point $\bx$ represented by the list of polynomials
$\bg_\bx$, define the
 \DeCo{{\sl fundamental differential operator $D_\bx$ of the critical point $\bx$}} by
 \begin{equation}\label{Eq:FundDiffOp}
   \Bigl(\frac{d}{dt} - \dlog\Bigl(\frac{\Phi}{g_{m-1}}\Bigr)\Bigr)
   \,\dotsb\,
   \Bigl(\frac{d}{dt} - \dlog\Bigl(\frac{g_{2}}{g_1}\Bigr)\Bigr)
   \Bigl(\frac{d}{dt} - \dlog(g_1)\Bigr)\,,
 \end{equation}
where $\dlog(f):=\frac{d}{dt}\ln f$.
The kernel \DeCo{$V_\bx$} of $D_\bx$ is the 
\DeCo{{\sl fundamental space of the critical point $\bx$}}.

\begin{ex}\label{Ex:solutions}
 Since 
\[
   \Bigr( \frac{d}{dt}\ -\ \dlog(f)\Bigr)\; f\ =\ 
   \Bigr( \frac{d}{dt}\ -\ \frac{f'}{f}\Bigr)\; f\ =\ 
   f' - \frac{f'}{f}f\ =\ 0\,,
\]
 we see that $g_1$ is a solution of $D_\bx$.
 It is instructive to look at $D_\bx$ and $V_\bx$ when $m=2$.
 Suppose that $g$ a solution to $D_\bx$ that is linearly independent from $g_1$.
 Then
\[
  0\ =\ 
   \Bigl(\frac{d}{dt} - \dlog\Bigl(\frac{\Phi}{g_{1}}\Bigr)\Bigr)
   \Bigl(\frac{d}{dt} - \dlog(g_1)\Bigr)\;g
  \ =\  \Bigl(\frac{d}{dt} - \dlog\Bigl(\frac{\Phi}{g_{1}}\Bigr)\Bigr)
     \bigl( g' - \frac{g_1'}{g_1} g \bigr)\,.
\]
 This implies that
\[
   \frac{\Phi}{g_1}\ =\ g' - \frac{g_1'}{g_1} g\,,
\]
 so $\Phi=\Wr(g,g_1)$, and the kernel of $D_\bx$ is a two-dimensional space of functions with
 Wronskian $\Phi$. 
\end{ex}

What we just saw is always the case.
The following result is due to Scherbak and Varchenko~\cite{ScVa} for $m=2$ and to Mukhin
and Varchenko~\cite[\S 5]{MV04} for all $m$.

\begin{thm}\label{Th:fundSpace}
 Suppose that $V_\bx$ is the fundamental space of a critical point $\bx$ of
 the master function $\Xi$ with generic parameters $\bs$ which are the roots of $\Phi$. 
 \begin{enumerate}
  \item Then $V_\bx\in\Gr(m,\C_{m{+}p{-}1}[t])$ has Wronskian $\Phi$.
  \item The critical point $\bx$ is recovered from $V_\bx$ in some cases as follows.
        Suppose that $f_1,\dotsc,f_m$ are monic polynomials in $V_\bx$ with 
        $\deg f_i=p{-}1{+}i$, each $f_i$ is square-free, and that the pairs $f_i$ and
        $f_{i+1}$ are relatively prime.
        Then, up to scalar multiples, the polynomials $g_1,\dotsc,g_{m-1}$ in the sequence
        $\bg_\bx$ are
\[
    f_1\,,\ \Wr(f_1,f_2)\,,\ \Wr(f_1,f_2,f_3)\,,\ \dotsc\,,\ 
    \Wr(f_1,\dotsc,f_m)\,.
\]
 \end{enumerate} 
\end{thm}

Statement~(2) includes a general result about factoring a linear differential operator into
differential operators of degree 1.
Linearly independent  $C^\infty$ functions $f_1,\dotsc,f_m$ 
span the kernel of the differential operator of degree $m$,
\[
   \det\left(\begin{matrix}
          f_1 & f_2 & \dotsb& f_m &  1 \\
          f'_1& f'_2& \dotsb& f'_m&\frac{d}{dt}\\
        \vdots&\vdots&\ddots&\vdots&\vdots\\
     f_1^{(m)}& f_2^{(m)}& \dotsb& f_m^{(m)}&\frac{d^m}{dt^m}
    \end{matrix}\right)\ .
\]
If we set $g_i:=\Wr(f_1,\dotsc,f_i)$, then~\eqref{Eq:FundDiffOp} is a factorization
over $\C(t)$ of this determinant into differential operators of degree 1.
This follows from some interesting identities among Wronskians shown in the Appendix
of~\cite{MV04}.

Theorem~\ref{Th:fundSpace} is deeper than this curious fact.
When the polynomials $g_1,\dotsc,g_m$ and $\Phi$ are square-free,
consecutive pairs are relatively prime, and $\bs$ is generic, it implies that the kernel $V$ of
an operator of the form~\eqref{Eq:FundDiffOp} is a space of polynomials with Wronskian $\Phi$
having roots $\bs$ if and only if the polynomials $g_1,\dotsc,g_m$ come from the
critical points of the master function~\eqref{Eq:MasterFunction} corresponding to $\Phi$.

This gives an injection from $\calS$-orbits of critical
points of the master function $\Xi$ with parameters $\bs$ to spaces of polynomials 
whose Wronskian has roots $\bs$.
Mukhin and Varchenko showed that this is a bijection when $\bs$ is generic. 

\begin{thm}[Theorem 6.1 in~\cite{MV05}]\label{Th:distinct}
  For generic complex numbers $\bs$, the master function $\Xi$ has nondegenerate critical
  points that form $\#_{m,p}$ distinct orbits. 
\end{thm}

The structure (but not of course the details) of their proof is remarkably similar to the
structure of the proof of Theorem~\ref{T8:asymptotic}; they allow the
parameters to collide one-by-one, and study how the orbits of critical points behave.
Ultimately, they obtain the same recursion as in~\eqref{Eq8:delta-recursion}, which mimics the
Pieri formula for the branching rule for tensor products of representations of $\slm$ with its
last fundamental representation $V_{\omega_{m-1}}$.
This same structure is also found in the main argument in~\cite{EG01}.
In fact, this is the same recursion in $\alpha$ that Schubert established for intersection
numbers $\delta(\alpha)$, and then solved to obtain the
formula~\eqref{E1:WronDeg} in~\cite{Sch1886c}. 

%
\section{The Bethe  ansatz for the Gaudin model}\label{S:BAGM}

The Bethe ansatz is a general (conjectural) method to find pure states, called 
\DeCo{{\sl Bethe vectors}}, of quantum integrable systems. 
The (periodic) Gaudin model is an integrable system consisting of a family of commuting
operators called the Gaudin Hamiltonians that act on a representation $V$ of $\slm\C$.
In this Bethe ansatz, a vector-valued rational function 
is constructed so that for certain values of the parameters it yields a complete
set of Bethe vectors.
As the Gaudin Hamiltonians commute with the action of $\slm\C$, the Bethe vectors turn out to
be highest weight vectors generating irreducible submodules of $V$, and so this also gives a
method for decomposing some representations $V$ of $\slm\C$ into irreducible submodules. 
The development, justification, and refinements of this Bethe ansatz are the subject of a large
body of  work, a small part of which we mention.

%
\subsection{Representations of  $\mathfrak{sl}_m\C$}
The Lie algebra $\slm\C$ (or simply $\slm$) is the space of $m\times m$-matrices
with trace zero.
It has a decomposition
\[
   \slm\ =\ \frakn_-\oplus\frakh\oplus\frakn_+\,,
\]
where $\frakn_+$ $(\frakn_-)$ are the strictly upper (lower) triangular matrices,
and $\frakh$ consists of the diagonal matrices with zero trace.
The universal enveloping algebra $U\slm$ of $\slm$ is the associative algebra generated by
$\slm$ subject to the relations $uv-vu=[u,v]$ for $u,v\in\slm$, where $[u,v]$ is the Lie
bracket in $\slm$.

We consider only finite-dimensional representations of $\slm$ (equivalently, of $U\slm$).
For a more complete treatment, see~\cite{FuHa}.
Any representation $V$ of $\slm$ decomposes into joint
eigen\-spaces of $\frakh$, called \DeCo{{\sl weight spaces}},
\[
   V\ =\ \bigoplus_{\mu\in\frakh^*} V[\mu]\,,
\]
where, for $v\in V[\mu]$ and $h\in\frakh$, we have $h.v=\mu(h)v$.
The possible weights $\mu$ of representations lie in the integral 
\DeCo{{\sl weight lattice}}.
This has a distinguished basis of \DeCo{{\sl fundamental weights}}
$\DeCo{\omega_1},\dotsc,\DeCo{\omega_{m-1}}$ that generate the cone of 
\DeCo{{\sl dominant weights}}.

An irreducible representation $V$ has a unique one-dimensional weight space that is annihilated
by the nilpotent subalgebra $\frakn_+$ of $\slm$.
The associated weight $\mu$ is dominant, and it is called the \DeCo{{\sl highest weight}} of $V$.
Any nonzero vector with this weight is a highest weight vector of $V$, and it generates $V$.
Furthermore, any two irreducible modules with the same highest weight are isomorphic.
Write \DeCo{$V_\mu$} for the \DeCo{{\sl highest weight module}} with highest weight $\mu$.
Lastly, there is one highest weight module for each dominant weight.

More generally, if $V$ is any representation of $\slm$ and $\mu$ is a weight, 
then the \DeCo{{\sl singular vectors}} in $V$ of weight $\mu$, 
written \DeCo{$\sing(V[\mu])$}, are the vectors in $V[\mu]$ annihilated by $\frakn_+$. 
If $v\in \sing(V[\mu])$ is nonzero, then the submodule $U\slm.v$ it generates is
isomorphic to the highest weight module $V_\mu$.
Thus $V$ decomposes as a direct sum of submodules generated by the singular
vectors, 
 \begin{equation}\label{Eq:sing_decomp}
   V\ =\ \bigoplus_\mu  U\slm.\sing (V[\mu])\,,
 \end{equation}
so that the multiplicity of the highest weight module $V_\mu$ in $V$ is simply the
dimension of its space of singular vectors of weight $\mu$.

When $V$ is a tensor product of highest weight modules, the Littlewood-Richard\-son
rule~\cite{Fu} gives formulas for the dimensions of the spaces of singular vectors.
Since this is the same rule for the number of points in an
intersection~\eqref{Eq7:Schubert_intersection} of Schubert varieties from a
Schubert problem,  these geometric intersection numbers are equal to the dimensions of 
spaces of singular vectors.
In particular, if $V_{\omega_1}\simeq\C^m$ is the defining representation of $\slm$
and $\DeCo{V_{\omega_{m-1}}}=\bigwedge^{m-1} V_{\omega_1}=V_{\omega_1}^*$ (these are the first and last
fundamental representations of $\slm$), then 
 \begin{equation}\label{Eq:SC=RT}
   \dim \sing( V_{\omega_n}^{\otimes mp}[0])\ =\ 
    \#_{m,p}\,.
 \end{equation}
It is important to note that this equality of numbers is purely formal, in that the same
formula governs both numbers.
A direct connection remains to be found.

%
\subsection{The (periodic) Gaudin model}
The Bethe ansatz is a conjectural method to obtain a complete set of eigenvectors for the
integrable system on $V:=V_{\omega_{m-1}}^{\otimes n}$ given by the Gaudin Hamiltonians (defined
below). 
Since these Gaudin Hamiltonians commute with $\slm$, the Bethe ansatz has the additional
benefit of giving an explicit basis for $\sing(V[\mu])$, thus explicitly giving the 
decomposition~\eqref{Eq:sing_decomp}.

The Gaudin Hamiltonians act on $V_{\omega_{m-1}}^{\otimes n}$ and depend upon $n$
distinct complex numbers $s_1,\dotsc,s_n$ and a complex variable $t$.
Let $\glm$ be the Lie algebra of $m\times m$ complex matrices.
For each $i,j=1,\dotsc,m$, let $E_{i,j}\in\glm$ be the matrix whose only nonzero entry is a
1 in row $i$ and column $j$.
For each pair $(i,j)$ consider the differential operator $X_{i,j}(t)$ acting on 
$V_{\omega_{m-1}}^{\otimes n}$-valued functions of $t$,
\[
   \DeCo{X_{i,j}(t)}\ :=\ \delta_{i,j}\frac{d}{dt}\ -\ 
    \sum_{k=1}^n \frac{E_{j,i}^{(k)}}{t-s_k}\ ,
\]
where $E_{j,i}^{(k)}$ acts on tensors in $V_{\omega_{m-1}}^{\otimes n}$  by $E_{j,i}$ in the
$k$th factor and by the identity in other factors.
Define a differential operator acting on $V_{\omega_{m-1}}^{\otimes n}$-valued functions of
$t$, 
\[
   \DeCo{\bM}\ := \sum_{\sigma\in\calS} \mbox{sgn}(\sigma)\; 
    X_{1,\sigma(1)}(t)\; 
    X_{2,\sigma(2)}(t)\; \dotsb\; X_{m,\sigma(m)}(t)\ ,
\]
where $\calS$ is the group of permutations of $\{1,\dotsc,m\}$ and
$\mbox{sgn}(\sigma)=\pm$ is the sign of a permutation $\sigma\in\calS$.
Write $\bM$ in standard form
\[
   \bM\ =\ \frac{d^m}{dt^m}\ +\ M_1(t) \frac{d^{m-1}}{dt^{m-1}}\ +\ 
    \dotsb\ +\ M_{m}(t)\,.
\]
These coefficients $M_1(t),\dotsc,M_{m}(t)$ are called the (higher) 
\DeCo{{\sl Gaudin Hamiltonians}}.
They are linear operators that depend rationally on $t$ and act on  
$V_{\omega_{m-1}}^{\otimes n}$.
We collect together some of their properties.

\begin{thm}
 Suppose that $s_1,\dotsc,s_n$ are distinct complex numbers.
 Then
 \begin{enumerate}
   \item[{$1.$}] The Gaudin Hamiltonians commute, that is, $[M_i(u),M_j(v)]=0$ for all
      $i,j=1,\dotsc,m$ and $u,v\in\C$.
   \item[{$2.$}]  The Gaudin Hamiltonians commute with the action of\/ $\slm$ on
      $V_{\omega_{m-1}}^{\otimes n}$.
 \end{enumerate}
\end{thm}

Proofs are given in~\cite{KuS}, as well as
Propositions 7.2 and 8.3 in~\cite{MTV_06}, and are based on results of Talalaev~\cite{Ta}.
A consequence of the second assertion is that the Gaudin Hamiltonians preserve the weight
space decomposition of the singular vectors of $V_{\omega_{m-1}}^{\otimes n}$.
Since they commute, the singular vectors of $V_{\omega_{m-1}}^{\otimes n}$
have a basis of common eigenvectors of the Gaudin Hamiltonians.
The Bethe ansatz is a method to write down joint eigenvectors and their eigenvalues. 

%
%
\subsection{The Bethe ansatz for the Gaudin model}\label{S:BAGMwv}

This begins with a rational function that takes
values in a weight space $V_{\omega_{m-1}}^{\otimes n}[\mu]$, 
\[
   v\ \colon\ \C^l\times\C^n\ \longmapsto\  V_{\omega_{m-1}}^{\otimes n}[\mu]\ .
\]
This \DeCo{{\sl universal weight function}} was introduced in~\cite{SV} to solve the
Knizhnik-Zamo\-lod\-chi\-kov equations with values in $V_{\omega_{m-1}}^{\otimes n}[\mu]$.
When $(\bx,\bs)$ is a critical point of a master function, the vector 
$v(\bx,\bs)$ is both singular and an eigenvector of the Gaudin Hamiltonians.
(This master function is a generalization of the one defined
by~\eqref{Eq:MasterFunction}.) 
The Bethe ansatz conjecture for the periodic Gaudin model asserts that the 
vectors $v(\bx,\bs)$ form a basis for the space of singular vectors.

Fix a highest weight vector $\DeCo{v_m}\in V_{\omega_{m-1}}[\omega_{m-1}]$.
Then $v_m^{\otimes n}$ generates $V_{\omega_{m-1}}^{\otimes n}$ as a
$U\slm^{\otimes n}$-module.  
In particular, any vector in $V_{\omega_{m-1}}^{\otimes n}$ is a linear combination of vectors that
are obtained from $v_m^{\otimes n}$ by applying a sequence of operators $E^{(k)}_{i+1,i}$, for
$1\leq k\leq n$ and $1\leq i\leq m{-}1$.
The universal weight function is a linear combination of such vectors of weight
$\mu$.

When 
$n=mp$, $l=p\binom{m}{2}$, and $\mu=0$,
the universal weight function is a map
\[
  v\ \colon\ \C^{p\binom{m}{2}}\times\C^{mp} \ 
    \longrightarrow\  V_{\omega_{m-1}}^{\otimes mp}[0]\,.
\]
To describe it, note that a vector $E_{a+1,a}E_{b+1,b}\dotsb E_{c+1,c}.v_{n+1}$ is nonzero only if
\[
  (a,b,\dotsc,c)\ =\ (a, a{+}1,\dotsc,m{-}2,m{-}1)\,.
\]
Write \DeCo{$v_a$} for this vector.
The vectors $v_1,\dotsc,v_m$ form a basis of $V_{\omega_{m-1}}$.
Thus only some sequences of operators $E^{(k)}_{i+1,i}$ applied to $v_m^{\otimes mp}$ 
give a nonzero vector.
These sequences are completely determined once we know the weight of the result.
The operator $E^{(k)}_{i+1,i}$ lowers the weight of a weight vector by the root $\alpha_i$.
Since
 \begin{equation}\label{eq:weight_sum}
    m\omega_{m-1}\ =\ \alpha_1+2\alpha_2+\dotsb+(m{-}1)\alpha_{m-1}\,,
 \end{equation}
there are $ip$ occurrences of $E^{(k)}_{i+1,i}$, which is the number of
variables in $\bx^{(i)}$.

Let \DeCo{$\calB$} be the set of all sequences $(b_1,b_2,\dotsc,b_{mp})$,
where $1\leq b_k\leq m$ for each $k$, and we have
\[
   \#\{k\mid b_k\leq i\}\ =\ ip\,.
\]
Given a sequence $B$ in $\calB$, define
 \begin{eqnarray*}
   \DeCo{v_B} &:=& v_{b_1}\otimes v_{b_2}\otimes\dotsb\otimes v_{b_{mp}}\\
               &=& \bigotimes_{k=1}^{mp} 
           \bigl( E^{(k)}_{b_k+1,b_k} \dotsb E^{(k)}_{m-1,m-2}\cdot
           E^{(k)}_{m,m-1}\bigr).v_m\,,
 \end{eqnarray*}
where the operator $ E^{(k)}_{b_k+1,b_k}\dotsb E^{(k)}_{m-1,m-2}\cdot E^{(k)}_{m,m-1}$ is
the identity if $b_k=m$.
Then $v_B$ is a vector of weight $0$, by~\eqref{eq:weight_sum}.
The universal weight function is a linear combination of these vectors $v_B$,
 \begin{equation}\label{Eq10:ratfn}
   v(\bx;\bs)\ =\ \sum_{B\in\calB} w_B(\bx;\bs)\cdot v_B\,,
 \end{equation}
where the function $w_B(\bx,\bs)$ is separately symmetric in each set of variables
$\bx^{(i)}$. 

To describe $w_B(\bx;\bs)$, suppose that 
\[
   \bz\ =\ (\bz^{(1)},\bz^{(2)},\dotsc,\bz^{(mp)})
\]
is a partition of the variables $\bx$ into $mp$ sets of variables where the $k$th
set $\bz^{(k)}$ of variables has exactly one variable from each set $\bx^{(i)}$ with $b_k\leq i$
(and is empty when $b_k=m$).
That is, if $b_k\leq m{-}1$, then
 \begin{equation}\label{eq:yk}
   \bz^{(k)}\ =\ (x^{(b_k)}_{c_{b_k}}, x^{(b_k+1)}_{c_{b_k+1}},\dotsc, x^{(m-1)}_{c_{m-1}})\,,
 \end{equation}
for some indices $c_{b_k},\dotsc,c_{m-1}$. 
If $b_k=m$, set $w_k(\bz):=1$, and otherwise
\[
   w_k(\bz;\bs)\ :=\ \frac{1}{x^{(b_k)}_{c_{b_k}}-x^{(b_k+1)}_{c_{b_k+1}}}\dotsb
               \frac{1}{x^{(m-2)}_{c_{m-2}}-x^{(m-1)}_{c_{m-1}}}\cdot
               \frac{1}{x^{(m-1)}_{c_{m-1}}-s_k}\,,
\]
in the notation~\eqref{eq:yk}.
Then we set 
\[
   w(\bz;\bs)\ :=\ \prod_{k=1}^{mp} w_k(\bz;\bs)\,.
\]
Finally, $w_B(\bx;\bs)$ is the sum of the rational functions $w(\bz;\bs)$ over all 
such partitions $\bz$ of the variables $\bx$.
(Equivalently, the symmetrization of any single $w(\bz;\bs)$.)

While $v(\bx,\bs)$~\eqref{Eq10:ratfn} is a rational function of $\bx$ and hence not globally
defined,  if the coordinates of $\bs$ are distinct and $\bx$ is a critical point of the master 
function~\eqref{Eq:MasterFunction}, then the vector 
$v(\bx,\bs)\in V^{\otimes mp}_{\omega_{m-1}}[0]$ is well-defined, nonzero and it is in fact
a singular vector (Lemma~2.1 of~\cite{MV05}).
Such a vector $v(\bx,\bs)$ when $\bx$ is a critical point of the master function is called a 
\DeCo{{\sl Bethe vector}}.
Mukhin and Varchenko also prove the following, which is the second part of Theorem~6.1
in~\cite{MV05}.

\begin{thm}\label{Th:BV_basis}
 When $\bs\in\C^{mp}$ is general, the Bethe vectors form a basis of  the
 space $\sing\bigl( V^{\otimes mp}_{\omega_{m-1}}[0]\bigr)$.
\end{thm}

These Bethe vectors are the joint
eigenvectors of the Gaudin Hamiltonians.

\begin{thm}[Theorem 9.2 in \cite{MTV_06}]\label{Th:BV_EV}
 For any critical point $\bx$ of the master function $\Xi$~$\eqref{Eq:MasterFunction}$, the
 Bethe vector $v(\bx,\bs)$ is a joint eigenvector of the Gaudin Hamiltonians
 $M_1(t),\dotsc,M_{n+1}(t)$.
 Its eigenvalues $\mu_1(t),\dotsc,\mu_{n+1}(t)$ are given by the formula
 \begin{multline}\label{Eq:ev}
   \quad \frac{d^m}{dt^m}\ +\ \mu_1(t)\frac{d^{m-1}}{dt^{m-1}}\ +\ 
      \dotsb\ +\ \mu_{m-1}(t)\frac{d}{dt}\ +\ \mu_{m}(t)\ =\ \\ \rule{0pt}{17pt}
     \Bigl(\frac{d}{dt}+\dlog(g_1)\Bigr)
     \Bigl(\frac{d}{dt}+\dlog\Bigl(\frac{g_2}{g_1}\Bigr)\Bigr)
       \ \dotsb\ 
     \Bigl(\frac{d}{dt}+\dlog\Bigl(\frac{g_{m-1}}{g_{m-2}}\Bigr)\Bigr)
     \Bigl(\frac{d}{dt}+\dlog\Bigl(\frac{\Phi}{g_{m-1}}\Bigr)\Bigr)\,,\quad
 \end{multline}
 where $g_1(t),\dotsc,g_{m-1}(t)$ are the polynomials~\eqref{eq:crit_polys} associated to
 the critical point $\bx$ and $\Phi(t)$ is the polynomial with roots $\bs$.
\end{thm}

Observe that~\eqref{Eq:ev} is similar to the formula~\eqref{Eq:FundDiffOp} for the
differential operator $D_\bx$ of the critical point $\bx$.
This similarity is made more precise if we replace the Gaudin Hamiltonians by a different
set of operators.
Consider the differential operator formally conjugate to $(-1)^mM$,
 \begin{eqnarray*}
   \DeCo{K}&=& \frac{d^m}{dt^m}\ -\ \frac{d^{m-1}}{dt^{m-1}}M_1(t)\ +\ 
    \dotsb\ +\ (-1)^{m-1}\frac{d}{dt} M_{m-1}(t)\ +\ (-1)^m M_m(t)\\
   &=& \frac{d^m}{dt^m}\ +\ \DeCo{K_1(t)}\frac{d^{m-1}}{dt^{m-1}}\ +\ 
    \dotsb\ +\  \DeCo{K_{m-1}(t)}\frac{d}{dt}\ +\ \DeCo{K_{m}(t)}\ .
 \end{eqnarray*}
These coefficients $K_i(t)$ are operators on $V_{\omega_{m-1}}^{\otimes mp}$
that depend rationally on $t$, and are also called the Gaudin Hamiltonians.
Here are the first three, 
 \begin{eqnarray*}
   K_1(t)&=&-M_1(t)\,,\qquad\qquad\qquad
   K_2(t)\ =\ M_2(t)\ -\ n M_1'(t)\,,\\
   K_3(t)&=& -M_3(t)\ +\ (n{-}1)M_2''(t)\ -\ \binom{n}{2}M_1'''(t)\,,
 \end{eqnarray*}
and in general $K_i(t)$ is a differential polynomial in $M_1(t),\dotsc,M_i(t)$.

Like the $M_i(t)$, these operators commute with each other
and with $\slm$, and the Bethe vector $v(\bx,\bs)$ is a joint eigenvector of these new Gaudin 
Hamiltonians $K_i(t)$. 
The corresponding eigenvalues $\lambda_1(t),\dotsc,\lambda_{m}(t)$ are given
by the formula
 \begin{multline}\label{Eq:Kev}
   \quad \frac{d^m}{dt^m}\ +\ \lambda_1(t)\frac{d^{m-1}}{dt^{m-1}}\ +\ 
      \dotsb\ +\ \lambda_{m-1}(t)\frac{d}{dt}\ +\ \lambda_m(t)\ =\ \\
   \Bigl(\frac{d}{dt} - \dlog\Bigl(\frac{\Phi}{g_{m-1}}\Bigr)\Bigr)   \rule{0pt}{18pt}
   \Bigl(\frac{d}{dt} - \dlog\Bigl(\frac{g_{m-1}}{g_{m-2}}\Bigr)\Bigr)
   \,\dotsb\,
   \Bigl(\frac{d}{dt} - \dlog\Bigl(\frac{g_{2}}{g_1}\Bigr)\Bigr)
   \Bigl(\frac{d}{dt} - \dlog(g_1)\Bigr)\,,
 \end{multline}
which is $(!)$ the fundamental differential operator $D_\bx$ of the critical point $\bx$.  

\begin{cor}\label{Co:GH_simple}
 Suppose that $\bs\in\C^{mp}$ is generic.
 \begin{enumerate}
  \item[{$1.$}]   The Bethe vectors form an eigenbasis of\/ 
       $\sing(V^{\otimes mp}_{\omega_{m-1}}[0])$ for the Gaudin Hamiltonians 
       $K_1(t),\dotsc,K_m(t)$.
  \item[{$2.$}]  The Gaudin Hamiltonians $K_1(t),\dotsc,K_m(t)$ have simple spectrum in that
       their eigenvalues separate the basis of eigenvectors.
\end{enumerate}
\end{cor}

Statement (1) follows from Theorems~\ref{Th:BV_basis} and~\ref{Th:BV_EV}.
For Statement (2), suppose that two Bethe vectors $v(\bx,\bs)$ and $v(\bx',\bs)$ have the
same eigenvalues.
By~\eqref{Eq:Kev}, the corresponding fundamental differential operators would be equal,
$D_{\bx}=D_{\bx'}$.
But this implies that the fundamental spaces coincide, $V_\bx=V_{\bx'}$.
By Theorem~\ref{Th:fundSpace} the fundamental space determines the orbit of critical
points, so the critical points $\bx$ and $\bx'$ lie in the same orbit, which implies that 
 $v(\bx,\bs)=v(\bx',\bs)$.

%
%
\section{Shapovalov form and the proof of the Shapiro conjecture}\label{S:Shapovalov}
The last step in the proof of Theorem~\ref{T1:Shap_Conj} is to show that if
$\bs\in\R^{mp}$ is generic and $\bx$ is a critical point of the master
function~\eqref{Eq:MasterFunction}, then the fundamental space $V_\bx$ of the critical
point $\bx$ has a basis of real polynomials.
The reason for this reality is that the eigenvectors and
eigenvalues of a symmetric matrix are real.

We begin with the Shapovalov form.
The map $\tau\colon E_{i,j}\mapsto E_{j,i}$ induces an antiautomorphism on $\slm$.
Given a highest weight module $V_\mu$, and a highest weight vector $v\in V_\mu[\mu]$
the \DeCo{{\sl Shapovalov form} $\langle \cdot,\cdot\rangle$} on $V_\mu$ is defined
recursively by 
\[
  \langle v,v\rangle\ =\ 1\qquad\mbox{and}\qquad
  \langle g.u,v\rangle\ =\ \langle u,\tau(g).v\rangle\,,
\]
for $g\in\slm$ and $u,v\in V$.
In general, this Shapovalov form is nondegenerate on $V_\mu$ and positive definite on the
real part of $V_\mu$.

For example, the Shapovalov form on $V_{\omega_{m-1}}$ is the standard Euclidean inner
product, $\langle v_i,v_j\rangle=\delta_{i,j}$, in the basis $v_1,\dotsc,v_{m}$ of
Section~\ref{S:BAGMwv}.
This induces the symmetric (tensor)
Shapovalov form on the tensor product $V_{\omega_{m-1}}^{\otimes mp}$,
which is  positive definite on the real part of 
$V_{\omega_{m-1}}^{\otimes mp}$.

\begin{thm}[Proposition 9.1 in~\cite{MTV_06}]
  The Gaudin Hamiltonians are symmetric with respect to the tensor Shapovalov form,
\[
   \langle K_i(t).u,\, v\rangle \ =\ \langle u,\, K_i(t).v\rangle\,,
\]
  for all $i=1,\dotsc,m$, $t\in\C$, and $u,v\in V_{\omega_{m-1}}^{\otimes mp}$.
\end{thm}

We give the most important consequence of this result for our story.

\begin{cor}\label{C:symmetric}
  When the parameters $\bs$ and variable $t$ are real, the Gaudin Hamiltonians 
  $K_1(t),\dotsc, K_m(t)$ are real linear operators with real spectrum.
\end{cor}

\noindent{\it Proof.}
 The Gaudin Hamiltonians $M_1(t),\dotsc,M_m(t)$ are real linear operators which
 act on the real part of  $V_{\omega_{m-1}}^{\otimes mp}$, by their definition.
 The same is then also true of the Gaudin Hamiltonians $K_1(t),\dotsc,K_m(t)$.
 But these are symmetric with respect to the Shapovalov form
 and thus have real spectrum.
\hfill\QED

\noindent{\it Proof of Theorem~$\ref{T1:Shap_Conj}$.}
 Suppose that $\bs\in\R^{mp}$ is general.
 By Corollary~\ref{C:symmetric}, the Gaudin Hamiltonians for $t\in\R$ acting on 
 $\sing(V_{\omega_{m-1}}^{mp}[0])$ are symmetric operators on a Euclidean space, and
 so have real eigenvalues.
 The Bethe vectors $v(\bx,\bs)$ for critical points $\bx$ of the master function with
 parameters $\bs$ form an eigenbasis for the Gaudin Hamiltonians.
 As $\bs$ is general, the eigenvalues are distinct by Corollary~\ref{Co:GH_simple} (2),
 and so the Bethe vectors must be real.

 Given a critical point $\bx$, the eigenvalues $\lambda_1(t),\dotsc,\lambda_m(t)$ of
 the Bethe vectors are then real rational functions, and so the fundamental differential
 operator $D_\bx$ has real coefficients.
 But then the fundamental space $V_\bx$ of polynomials is real.

 Thus each of the $\#_{m,p}$ spaces of polynomials $V_\bx$
 whose Wronskian has roots $\bs$ that were constructed in Section~\ref{S:polys} is in fact
 real.
 This proves Theorem~\ref{T1:Shap_Conj}.
\hfill\QED


%
%
%
\chapter{Beyond the Shapiro Conjecture}\label{Ch:Frontier}


Here, we will touch on further topics related to the Shapiro
Conjecture, including\vspace{-5pt}
\begin{enumerate}
 \item[I] Transversality and Discriminants.\vspace{-5pt}
 \item[II] Maximally Inflected Curves.\vspace{-5pt}
 \item[III] The Shapiro Conjecture for flag manifolds (Monotone Conjecture).\vspace{-5pt}
 \item[IV] The Secant Conjecture and the Monotone Secant Conjecture\vspace{-5pt}
 \item[V]  The Shapiro Conjecture for Lagrangian and Orthogonal Grassmannians.
\end{enumerate}

\Red{{\sf This chapter currently only in a rough draft, and will require a complete
    rewrite.}} 

%
%
%

\section{Transversality}
In Chapter~\ref{Ch:EG}, we presented a proof of the Shapiro Conjecture for rational
functions by Eremenko and Gabrielov~\cite{EG05}.
Its main point was that there is no obstruction to analytically continuing the
rational functions that were constructed in Chapter~\ref{Ch:ScGr} to give
rational functions with any given Wronskian having distinct real zeroes.
The key to this was the association of a net to each rational function with only
real critical points, which showed that the analytic continuation was possible.

A consequence of this proof is the statement that when $\min(m,p)=2$, the
Wronski map 
\[
  \Wr\ \colon\ \Gr(p,\C^{m+p})\ \longrightarrow\ \C\P^{mp}
\]
is unramified over the locus of hyperbolic polynomials with distinct (real) roots.
This is in fact true for all Grassmannians, as Mukhin, Tarasov, and Varchenko
showed~\cite{MTV_R}.

\begin{thm}\label{T11:trans}
 The Wronski map is unramified over the locus of hyperbolic polynomials
 with distinct roots, for any $m$ and $p$.
\end{thm}

 More generally, given called Schubert data, 
 $\alpha_1,\alpha_2,\dotsc,\alpha_n$ we can consider
 intersections of the form 
 \begin{equation}\label{E8:Shapiro}
   X_{\alpha_1}\Fdot(s_1)\ \cap\ 
   X_{\alpha_2}\Fdot(s_2)\ \cap\ \dotsb \ \cap\ 
   X_{\alpha_n}\Fdot(s_n)\ ,
 \end{equation}
where $s_1,s_2,\dotsc,s_n\in\P^1$.
Recall that a collection of subvarieties meet \DeCo{{\sl transversally}} if at
all points of their intersection, their tangent spaces meet properly (have the
expected dimension of intersection).
We rephrase Theorem~\ref{T11:trans}.\medskip

\noindent{\bf Transversality Theorem.}
 {\em If the points $s_1,s_2,\dotsc,s_n$ are real and distinct, then the
 intersection~$\eqref{E8:Shapiro}$ is transverse.
}\medskip

This Transversality Theorem has a strengthening, which is stated in
terms of real algebra.

\begin{defn}
  The \DeCo{{\sl discriminant}} is the locus in $\P^{mp}$ of critical values of
  the Wronski map (points over which it is ramified).  
  It is an algebraic hypersurface and defined by a single polynomial, also
  called the discriminant.
\end{defn}

 The Transversality Theorem asserts that the discriminant does not meet the
 set of polynomials with distinct roots.

More generally, if we consider a given family of
intersections of the form~\eqref{E8:Shapiro}, then the discriminant is the set
of points $(s_1,s_2,\dots,s_n)$ for which this intersection is not transverse.
Again, this discriminant is a polynomial $\Delta(s_1,s_2,\dots,s_n)$ in the
parameters $(s_1,s_2,\dots,s_n)\in(\P^1)^n$, and the Transversality Conjecture 
asserts that $\Delta$ does not vanish when the parameters are real and
distinct, that is, when $s_i\neq s_j$, for all $i,j$.
We conjecture something much stronger.

\begin{conj}[Discriminant Conjecture]\label{C8:discr}
 The discriminant polynomial $\Delta$ is a sum of squares.
 Each term in the sum is a monomial in the differences $(s_i-s_j)$.
\end{conj}

There is some evidence for this conjecture.
For the problem of two lines in 3-space meeting four lines $\ell(s)$, $\ell(t)$,
$\ell(u)$, and $\ell(v)$, each tangent to the rational normal curve, the
discriminant is a constant multiple of 
\[
    (s-t)^2(u-v)^2 \ +\ (s-u)^2(t-v)^2\ +\ (s-v)^2(t-u)^2\,.
\]

For $\Gr(2,5)$, the discriminant of the Schubert intersection
\[
  X_\I\Fdot(0)\ \cap\ X_\I\Fdot(\Blue{s})\ \cap\ 
  X_\I\Fdot(\Blue{t})\ \cap\ X_\I\Fdot(\Blue{u})\ \cap\ 
  X_\I\Fdot(\Blue{v})\ \cap\ X_\I\Fdot(\infty)\,,
\]
has degree 20 in the variables $\Blue{s},\Blue{t},\Blue{u}$, and $\Blue{v}$, and
it has 711 different terms. 
Ad hoc methods~\cite{So00b} showed that it was a sum of squares.\fauxfootnote{Give more details?}
This was quite surprising.
Hilbert~\cite{Hi1888} showed that not every polynomial of even degree greater
than 2 in four variables that is nonnegative can be written as a sum of squares.

We feel that it is a useful project to further investigate these discriminant
polynomials. 

\section{Maximally inflected curves}

This section represents joint work with Kharlamov~\cite{KhS03}.
A list $f_1(t),f_2(t),\dotsc,f_k(t)$ of polynomials with degree $d$
defines a map $\varphi\colon\P^1\to\P^{k-1}$ as follows,
\[
   \varphi\ \colon\ \P^1\ni t\ \longrightarrow\ 
          [f_1(t):f_2(t):\dotsb:f_k(t)]\in\P^{k-1}\,.
\]
The image of the curve is \DeCo{{\sl convex}} at a point $\varphi(t)$ if and
only if the first $k{-}1$ derivatives of $\varphi(t)$ are linearly independent. 
The failure to be convex is measured exactly by the vanishing of
the Wronskian of the polynomials $f_1(t),f_2(t),\dotsc,f_k(t)$.
We say that $\varphi$ is \DeCo{{\sl ramified}} at a point $t$ if 
this Wronskian vanishes at $t$.
Another geometric term is that the curve $\varphi$ has an 
\DeCo{{\sl inflection point}} (or flex) at $t$.
This corresponds exactly to the usual notion of an inflection point for plane
curves. 
A rational curve $\varphi$ of degree $d$ in $\P^{k-1}$ has $k(d{+}1{-}k)$
flexes, counted with multiplicity.

The connection between the Schubert calculus and rational curves in
projective space (linear series on ${\mathbb P}^1$) originated in
work of Castelnuovo~\cite{Ca1889} on $g$-nodal rational curves.
This led to the use of Schubert calculus in Brill-Noether
theory (see Chapter~5 of~\cite{HaMo} for an elaboration).
In turn, the theory of limit linear series of Eisenbud and
Harris~\cite{EH83,EH87}
provides essential tools to show reality of the special Schubert
calculus~\cite{So99}.
(That result on the special Schubert calculus was a generalization of
Theorem~\ref{T8:asymptotic} of Chapter~\ref{Ch:ScGr}.)
\medskip

\noindent{\bf Mukhin, Tarasov, Varchenko Theorem for rational curves.}
{\em 
  If a rational curve in $\P^{k-1}$ has all of its flexes real, then it
  \Blue{must} be real.}\medskip

\begin{defn}
  A real rational curve with all of its flexes real is said to be 
  \DeCo{{\sl maximally inflected}}.
\end{defn}

The Shapiro conjecture asserts that there are lots (maximally many, in fact) of
these maximally inflected curves.

\begin{figure}[htb]
 $$
    \includegraphics[height=70pt]{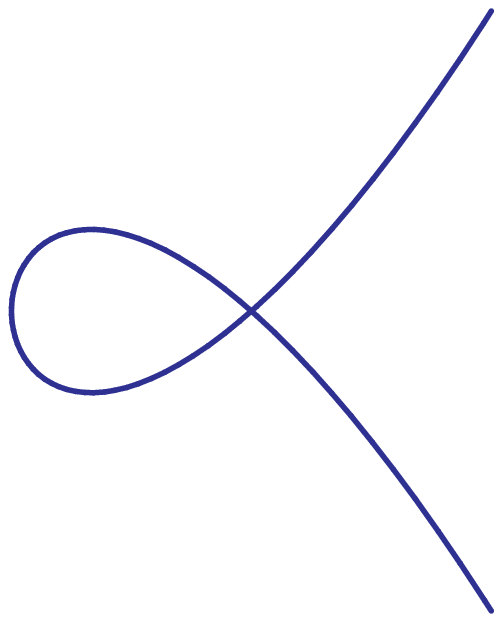}\qquad\qquad\qquad %
    \includegraphics[height=70pt]{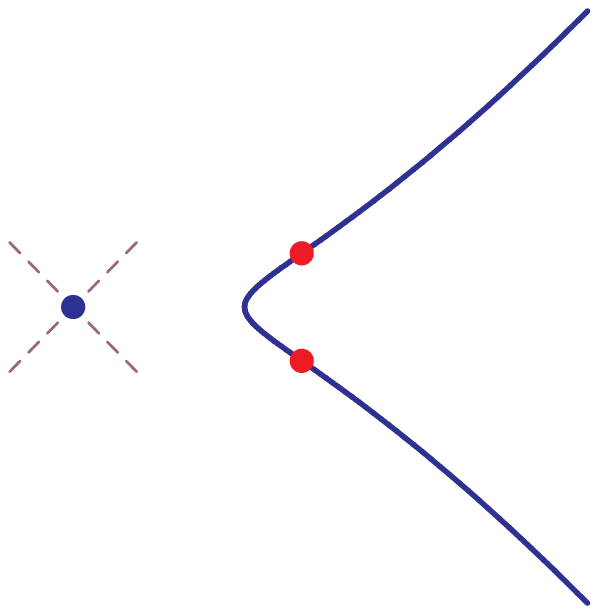}\qquad\qquad\qquad %
    \includegraphics[height=70pt]{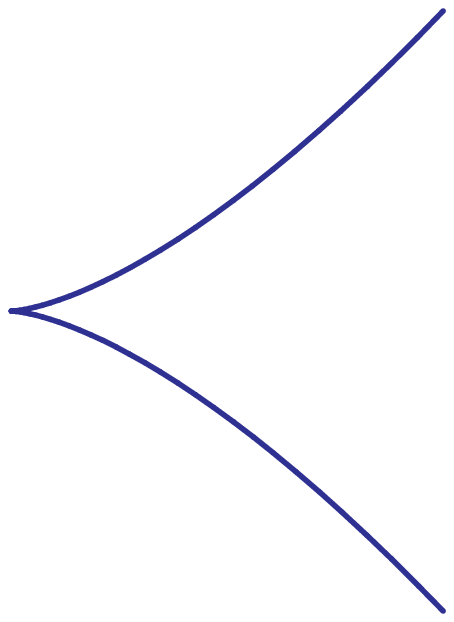}
 $$
 \caption{Real rational cubics.}\label{F8:cubics}
\end{figure}

Up to projective transformation and reparameterization,
there are only three real rational plane cubic curves.
They are represented by the
equations 
\[
   y^2\ =\ x^3+x^2\,,\qquad 
   y^2\ =\ x^3-x^2\,, \qquad\mbox{and}\qquad y^2\ =\ x^3\,,
\]
and they have the shapes shown in Figure~\ref{F8:cubics}.
All three have a real flex at infinity and are singular at the origin.
The first has a real node and no other real flexes, the second has
a solitary point and two real flexes at
$(\frac{4}{3},\pm\frac{4}{3\sqrt{3}})$ (we indicate these with dots
and the complex conjugate tangents at the
solitary point with dashed lines), and the
third has a real cusp.
The last two are maximally inflected, while the
first is not.

\begin{figure}[htb]
\[
   \includegraphics[height=70pt]{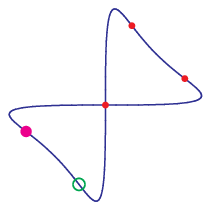}\qquad
   \includegraphics[height=70pt]{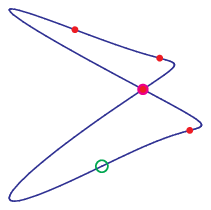}\qquad
   \includegraphics[height=70pt]{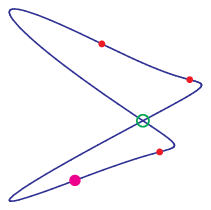}\qquad
   \includegraphics[height=70pt]{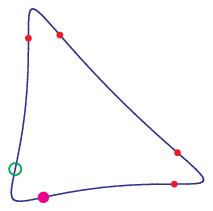}\qquad
   \includegraphics[height=70pt]{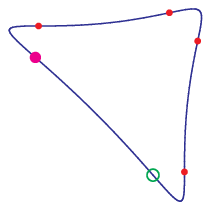}
\]
   \caption{The 5 curves with flexes at
            $\{-3,-1,0,1,3,\infty\}$.\label{F8:6flex}}
 \end{figure}

The Schubert calculus gives 5 rational quartics with 6 given points
of inflection and Figure~\ref{F8:6flex} shows 5 maximally inflected curves with
flexes at $\{-3,-1,0,1,3,\infty\}$.
(Each nodal curve has 2 flexes at its node, which is a consequence of the
 symmetry in the choice of flexes.)
We indicate the differences in the parameterizations of these curves,
labeling the flex at $-3$ by the larger dot and the flex at $-1$ by the
circle.
The solitary points are not drawn.
The first three curves have two solitary points, while the last 2 have three
solitary points.

Here are a few quintics.
The flexes are indicated (the symmetric curves have one additional flex at infinity), but
we do not draw solitary points. 
Also, the open circles represent two flexes which have merged into a \DeCo{{\sl
    planar point}}.
\[
  \includegraphics[height=90pt]{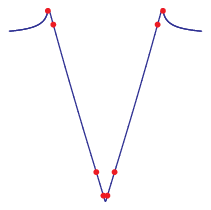}\qquad
  \includegraphics[height=90pt]{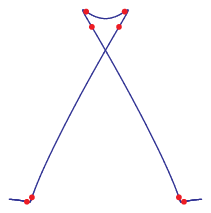}\qquad
  \includegraphics[height=90pt]{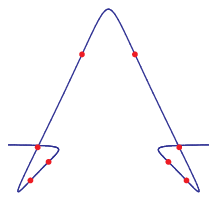}\qquad
  \includegraphics[height=90pt]{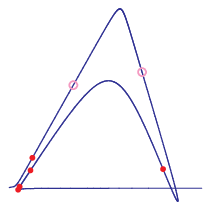}
\]
\[
  \includegraphics[height=90pt]{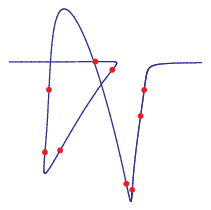}\qquad
  \includegraphics[height=90pt]{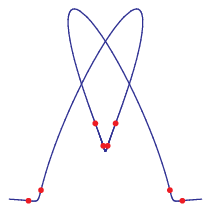}\qquad
  \includegraphics[height=90pt]{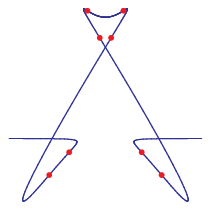}\qquad
  \includegraphics[height=90pt]{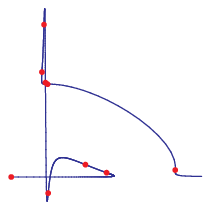}
\]
Finally, here are four singular quintics.
The solitary point is drawn on the second. 
\[
   \includegraphics[height=90pt]{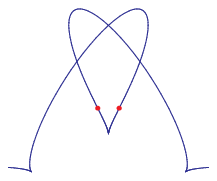}\qquad
   \includegraphics[height=90pt]{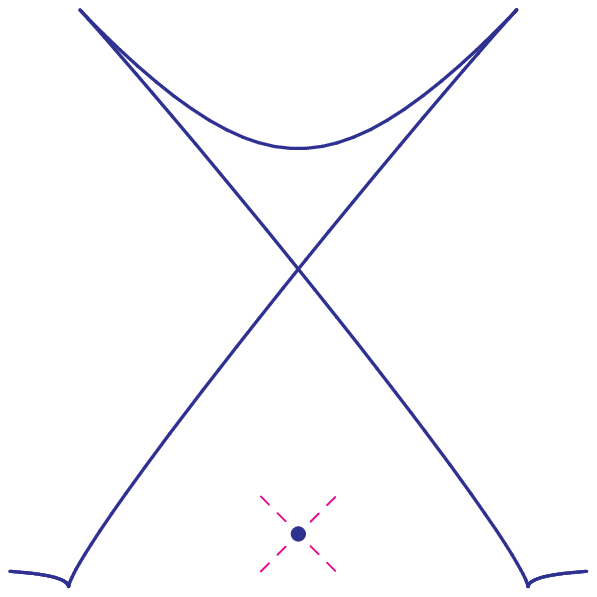}\qquad
   \includegraphics[height=90pt]{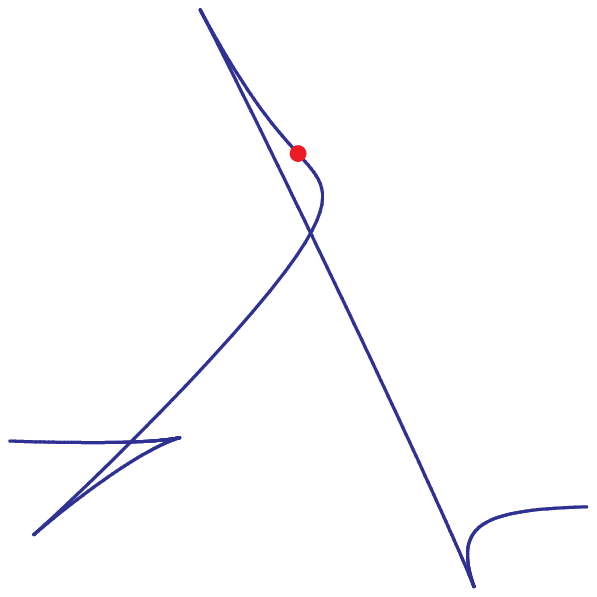}\qquad
   \includegraphics[height=90pt]{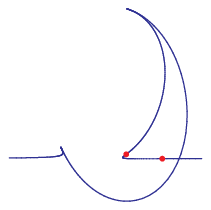}
\]

None of these maximally
inflected curves had many nodes.
Indeed, rational quartics typically have 3 double points and quintics have
6, yet we saw quartics with at most 1 node and quintics with at most 3 nodes.
More generally, 
consider a maximally inflected curve with only flexes and cusps, and whose other
singularities are ordinary double points.
Let $\iota$ be its number of flexes and $\kappa$ be its number of cusps.
Then $\iota+2\kappa=3(d-2)$.
By the genus formula, it has $\binom{d-1}{2}-\kappa$ double
points.
The following theorem is an easy consequence of the Klein~\cite{Klein} and
Pl\"ucker~\cite{Plucker} formulas.

\begin{thm}[Topological Restrictions~\cite{KhS03}]\label{T8:restrictions}
 Suppose that a maximally inflected curve has only flexes, cusps, and ordinary
 double points, and let $\iota$ and $\kappa$ be its number of flexes and cusps,
 respectively. 
 Then it has at least $d-2-\kappa$ solitary points and at most $\binom{d-2}{2}$
 nodes.
\end{thm}

Thus maximally inflected cubics have at most $\binom{3-2}{2}=0$ nodes, quartics
have at most $\binom{4-2}{2}=1$ node, and quintics have at most
$\binom{5-2}{2}=3$ nodes.

The existence of curves satisfying the hypotheses of
Theorem~\ref{T8:restrictions} is not guaranteed, even if we know the Shapiro
conjecture. 
For example, the construction in Theorem~4 of Chapter~\ref{Ch:ScGr} does not guarantee that
the maximally inflected curve has only ordinary double points. 
There are, however two constructions which guarantee curves 
having only double points.
The first uses Shustin's patchworking of singular curves~\cite{Sh85} to obtain
degenerate Harnack curves with $\binom{d-1}{2}$ solitary points (the figure on
the left illustrates this patchworking).
The second perturbs $d-2$ lines tangent to a conic to obtain maximally inflected
curves of degree $d$ with the minimal number of solitary points (and up to $d-2$
cusps). 
\[
  \raisebox{-25pt}{\includegraphics[height=150pt]{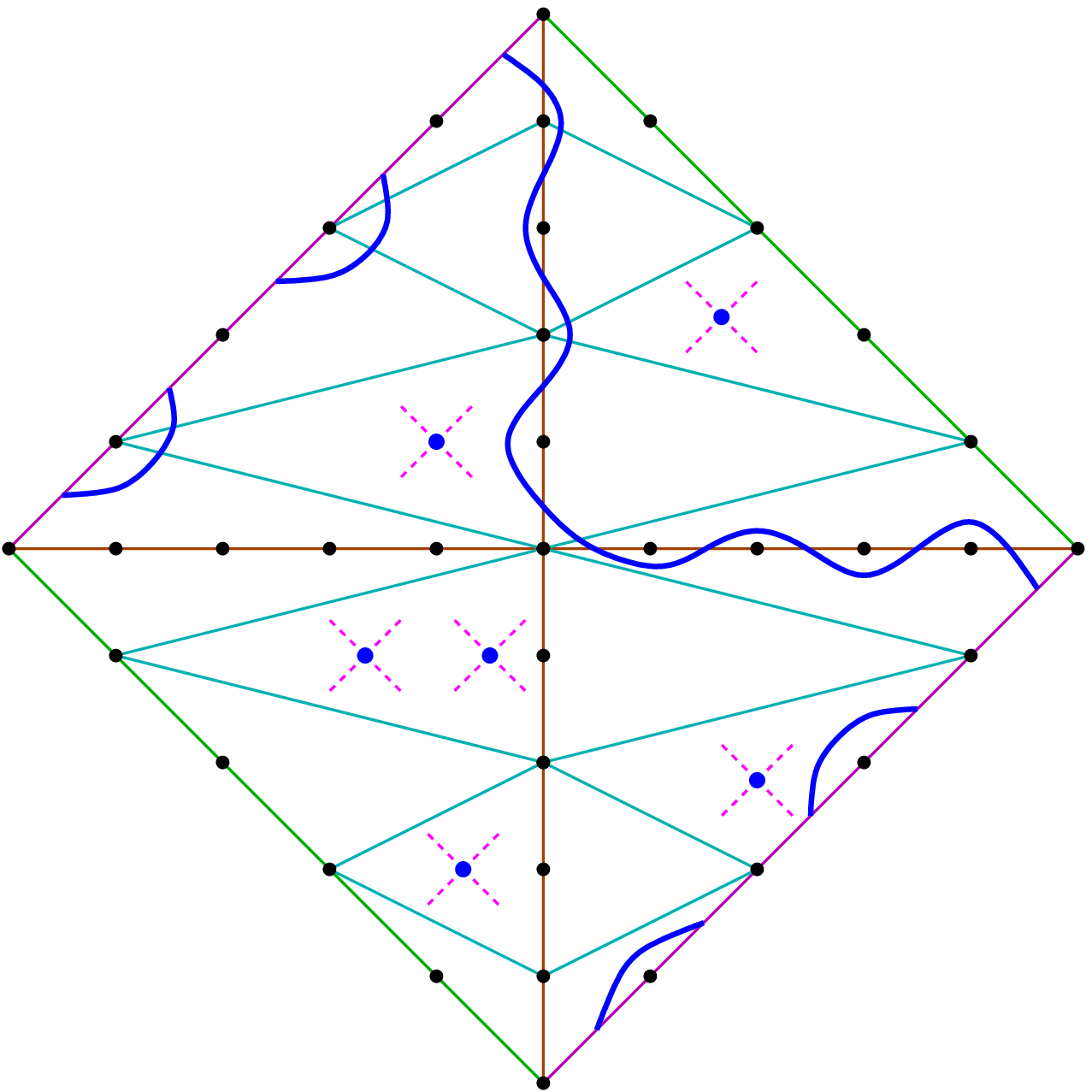}}\qquad
  \includegraphics[height=100pt]{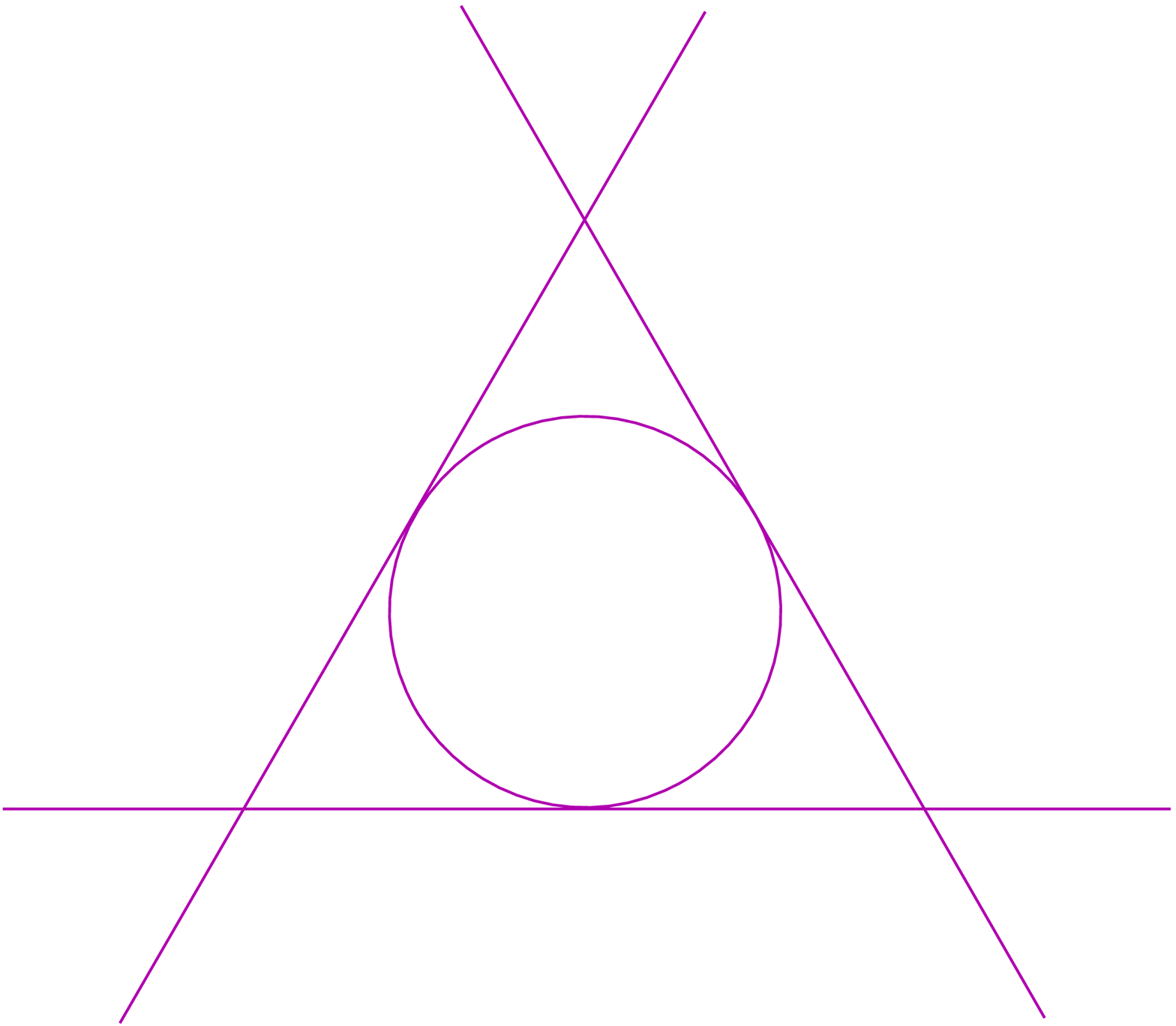}
   \raisebox{45pt}{\Large$\Longrightarrow$}
  \includegraphics[height=100pt]{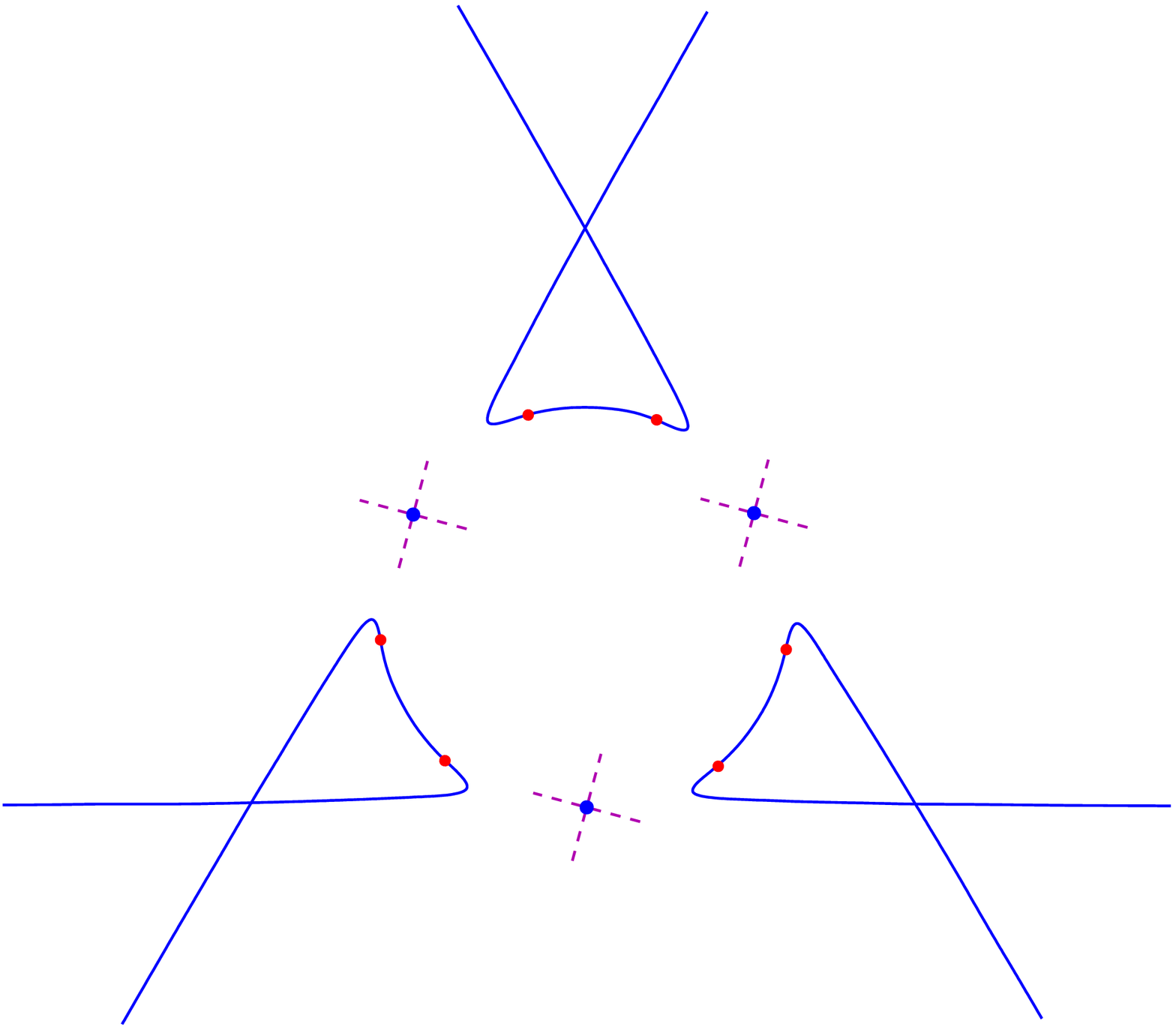}
\]

The topological classification of maximally inflected plane quintics is
open.\fauxfootnote{Explain a couple of open questions}
Also, maximally inflected curves in higher-dimensional space have not been
investigated.
For example, which knot types can occur for maximally inflected space curves?
We speculate that there should be some restrictions on the number of real
quadrisecants of a maximally inflected space curve.

\section{The Shapiro Conjecture for flag manifolds}

The Shapiro conjecture involved linear spaces satisfying incidence conditions
imposed by flags osculating the
rational normal curve at real points.
One natural variant is to consider flags satisfying incidence conditions imposed
by osculating flags.

\begin{ex}
 Consider partial flags of subspaces $m\subset H$, where $m$ is a line which
 lies on a plane $H$ in $\P^3$.
 Here are two typical incidence conditions on these flags.
\begin{itemize}
 \item $m$ meets a fixed line $\ell(t)$ tangent to the rational normal curve
       $\gamma$ at the point $\gamma(t)$. 
 \item $H$ contains a fixed point $\gamma(t)$ of the rational normal curve.
\end{itemize}

We consider the problem of  partial flags  $m\subset H$ where $m$ meets three
tangent lines $\Red{\ell(-1)}$,  $\Red{\ell(0)}$, and $\Red{\ell(1)}$ to 
$\gamma$, and $H$
contains two points $\gamma(v)$ and $\gamma(w)$.
First observe that $H$ contains the line $\lambda(v,w)$ spanned by these two
points.
But then  $m$ must also meet this secant line $\lambda(v,w)$, as $m\subset H$. 
In this way, we see that this problem reduces to our old friend, the problem of
four lines.
Now, however, we seek lines $m$ which meet three tangent lines  $\ell(-1)$,
$\ell(0)$, and $\ell(1)$ to the rational normal curve, and one secant line
$\lambda(v,w)$. 
Given such a line $m$, we may recover $H$ as the linear span of $m$ and
$\lambda(v,w)$.  

As in Section 4 of Chapter~\ref{S:ERAG}, we begin with the quadric 
\Brown{$Q$} containing the three lines 
\Red{$\ell(-1)$},  \Red{$\ell(0)$}, and \Red{$\ell(1)$} tangent to  the rational
normal curve \Blue{$\gamma$}. 
 \begin{figure}[htb]
 \[
  \begin{picture}(372,172.2)(0,6)
   \put(0,0){\includegraphics[height=190pt]{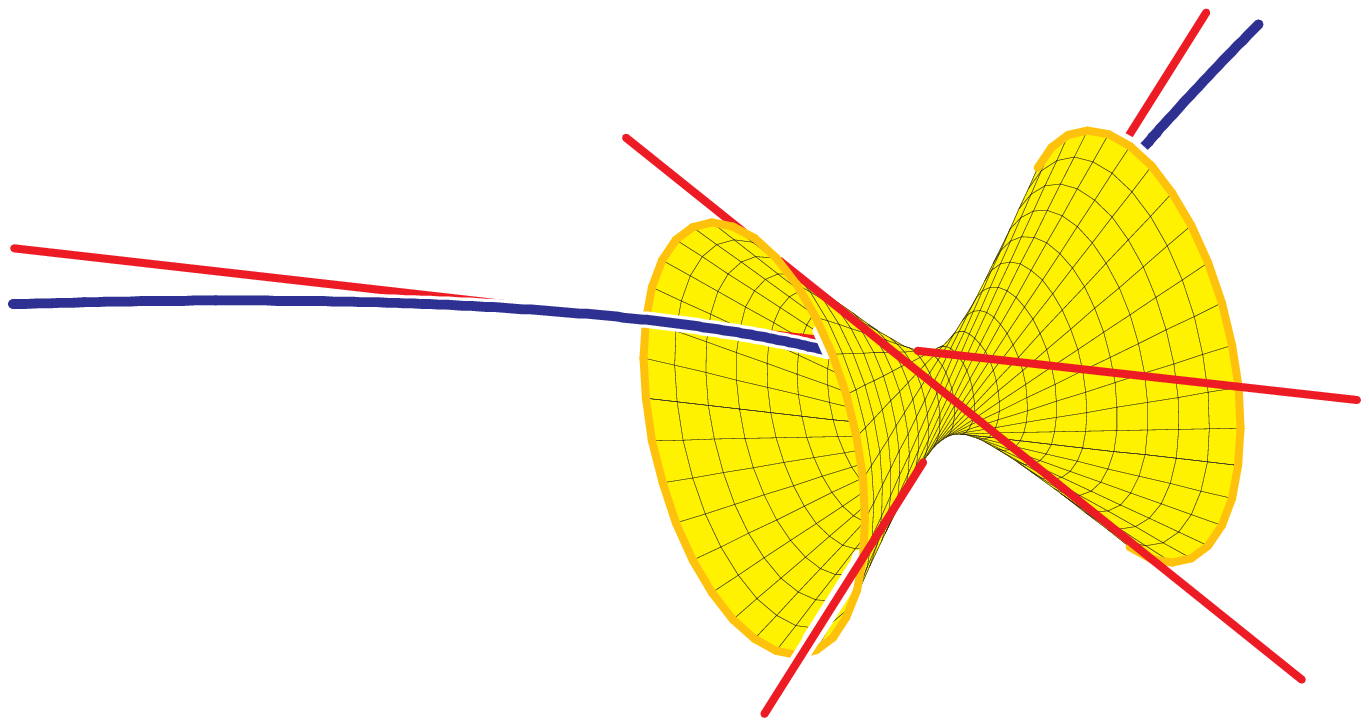} }
   \put(170,0){\Red{$\ell(-1)$}} \put(354,0){\Red{$\ell(0)$}}
   \put(0,136){\Red{$\ell(1)$}} \put(0,94){\Blue{$\gamma$}}
   \put(320,135){\Brown{$Q$}}
  \end{picture}
 \]
 \caption{Quadric containing three lines tangent to the rational normal
 curve\label{F8:3TanQuad}.} 
 \end{figure}
The lines meeting $\Red{\ell(-1)}$, $\Red{\ell(0)}$, $\Red{\ell(1)}$, 
and the secant line $\lambda(v,w)$ correspond to the points where $\lambda(v,w)$
meets the quadric \Brown{$Q$}. 
In Figure~\ref{F8:NC}, we display a secant line $\lambda(v,w)$ which meets the
hyperboloid in two real points, and therefore these choices for $v$ and $w$
give two real solutions to our geometric problem.
 \begin{figure}[htb]
 \[
  \begin{picture}(360,170)(-22,5)
   \put(0,0){\includegraphics[height=175pt]{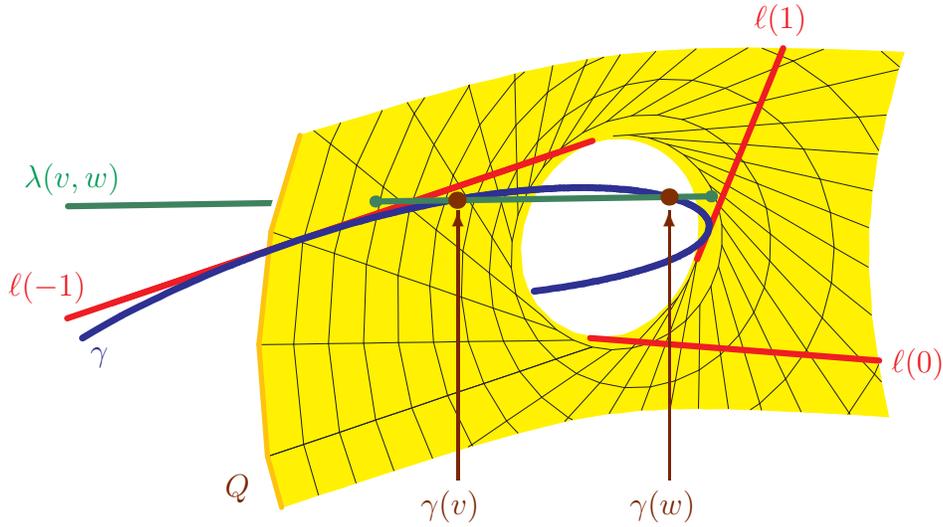} }
   \put(-22,81){\Red{$\ell(-1)$}} \put(312,52){\Red{$\ell(0)$}} 
   \put(260,182){\Red{$\ell(1)$}}
   \put(9,56){\Blue{$\gamma$}} \put(-16,122){\ForestGreen{$\lambda(v,w)$}}
   \put(60,5){\Brown{$Q$}}
 \thicklines
    \put(134,-2){\Brown{$\gamma(v)$}}  \put(148,11){\Brown{\vector(0,1){102}}}
    \put(213,-2){\Brown{$\gamma(w)$}}  \put(228,11){\Brown{\vector(0,1){102}}} 
  \end{picture}
\] \caption{A secant line meeting \Brown{$Q$}\label{F8:NC}.}
 \end{figure} 
There is also a secant line which meets the hyperboloid $\Brown{Q}$ in no real
points, and hence in two complex conjugate points. 
For this secant line, neither flag solving our problem is real.
 \begin{figure}[ht]
 \[
  \begin{picture}(440,220)(10,20)
   \put(0,0){\includegraphics[height=250pt]{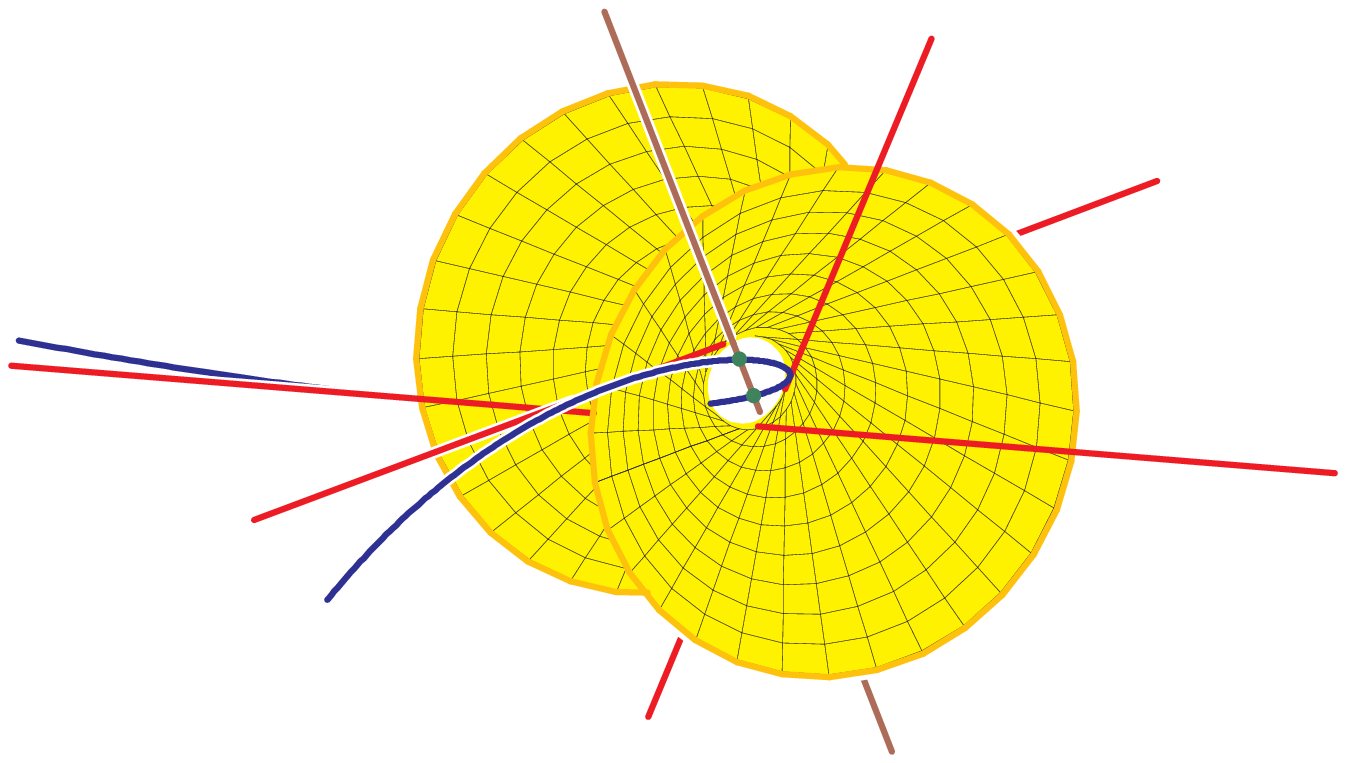} }
   \put(54,69){\Red{$\ell(-1)$}} \put(428,80){\Red{$\ell(0)$}} 
   \put(315,235){\Red{$\ell(1)$}}
   \put(100,40){\Blue{$\gamma$}}\put(163,238){\Brown{$\lambda(v,w)$}}
   \thicklines
   \put(265,243){\ForestGreen{$\gamma(v)$}}  
       \put(274,238){\ForestGreen{\vector(-1,-4){25.2}}}
   \put(171,11){\ForestGreen{$\gamma(w)$}}  
       \put(186,22){\ForestGreen{\vector(2,3){63.5}}} 
    \put(420,-20){\includegraphics[height=19pt]{figures/qed.eps}}
  \end{picture}
 \]
 \caption{A secant line not meeting \Brown{$Q$}\label{F8:CR}.}
 \end{figure}
We show this configuration in Figure~\ref{F8:CR}.
\end{ex}

\begin{figure}\label{F:mon}
\[
  \includegraphics[height=120pt]{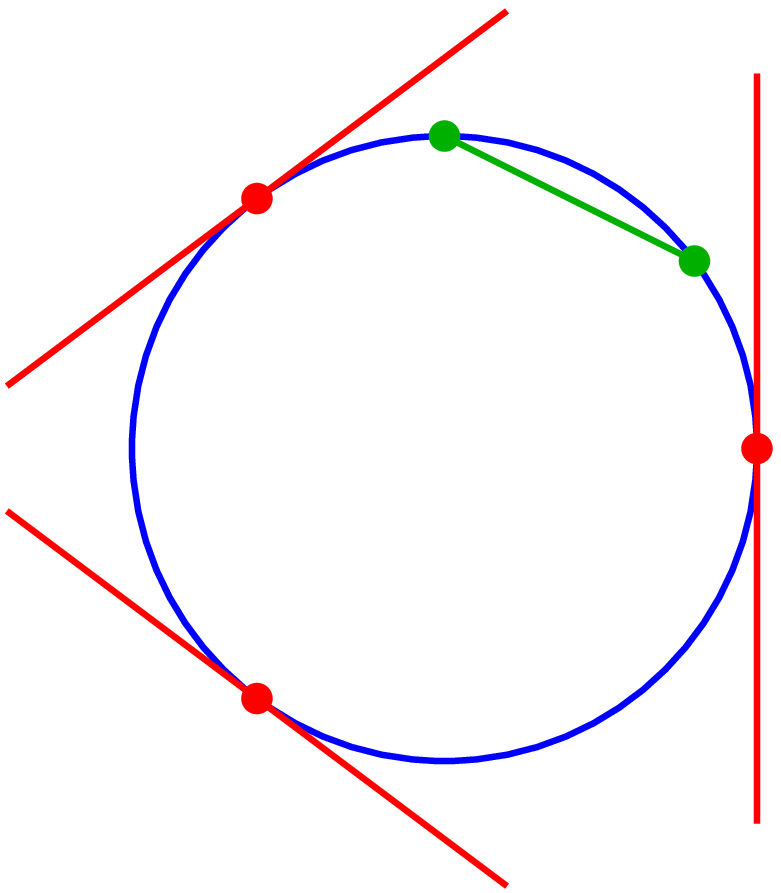}
   \qquad\qquad
   \includegraphics[height=120pt]{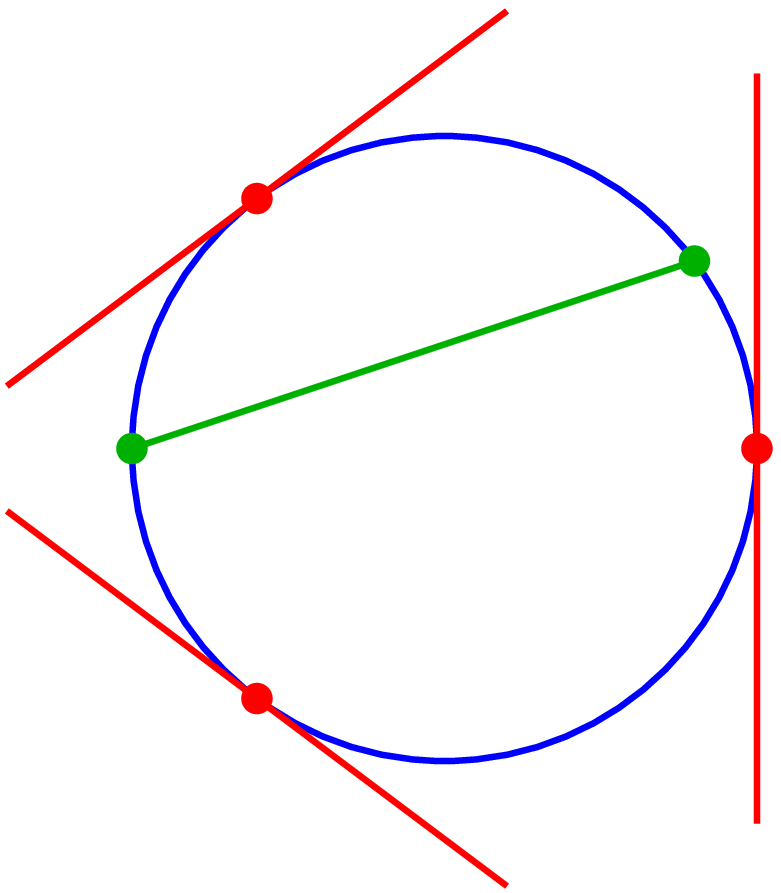}
\]
 \caption{\hspace{1.5cm}Always real.\hspace{3cm}Not always real.\hspace{3cm}}
\end{figure}

The original Shapiro conjecture concerned partial flags of a given type
meeting flags osculating the rational normal curve (at real points), and it
asserted that all such partial flags would be real.
This example shows that the Shapiro conjecture fails for flags, but not too
badly.
It is useful to consider this failure schematically.
In Figure~\ref{F:mon} we represent the rational normal curve $\gamma$ as a
circle, and indicate the relative positions of three tangent lines and the
secant line. 
Observe that in the diagram on the left, one may travel along the circle, first
encountering the three points where the lines are tangent and then the two points
where the line is secant.
If we record the dimension of the piece of the flag $m\subset H$ which
is affected by the tangent line or point on the rational normal curve as we
travel along the circle, we get
the sequence $(1,1,1,2,2)$, which is weakly increasing.
The reading the diagram on the right either gives the sequence $(1,1,2,1,2)$,
or its reverse $(2,1,2,1,1)$, or $(1,2,1,2,1)$, none of which are monotone.

\subsection{The Monotone Conjecture}

A \DeCo{{\sl partial flag $\Edot$}} of 
\DeCo{{\sl type $a=(a_1,a_2,\dotsc,a_k)$}} is a sequence of linear subspaces 
\[
  \Edot\ \colon\ E_{a_1}\ \subset\ E_{a_2}\ 
   \subset\ \dotsb\ \subset\ E_{a_k}\ \subset\ \C^n\,,
\]
where $a_i=\dim E_i$.
Its possible positions with respect to a fixed (complete) flag are encoded by
certain permutations, $\alpha$. 
The original Shapiro conjecture concerned partial flags of a given type
meeting flags osculating the rational normal curve (at real points), and it
asserted that all such partial flags would be real.
The example that we just gave shows that it was wrong\fauxfootnote{In fact, it
  fails for the first nontrivial problem on a flag manifold.}.
However, it can be repaired.

A position $\alpha$ for flags $\Edot$ of type $a$ which may be expressed in
terms of only one piece $E_{a_i}$ of the flag $\Edot$ is called 
\DeCo{{\sl Grassmannian}} of \DeCo{{\sl index}} $a_i$.
When this occurs, write $\iota(\alpha)=a_i$.
In the counterexample to the Shapiro conjecture, both conditions are
Grassmannian---one has index 2 (for the requirement that $m$ meet a tangent
line) and the other has index 3 (for the requirement that $H$ meet a point on
the rational normal curve.)  
These differ from the dimension of the linear spaces $m$ and $H$, as we
were working in projective 3-space, rather than a 4-dimensional vector space.

\begin{conj}[Monotone Conjecture]
 Let $\alpha_1,\alpha_2,\dotsc,\alpha_m$ be Grassmannian conditions on flags of
 type $a$ whose indices satisfy
 $\iota(\alpha_1)\leq\iota(\alpha_2)\leq\dotsb\leq\iota(\alpha_m)$.
 Then if $s_1<s_2<\dotsb<s_m$ are real numbers, the intersection
\[
  \bigcap_{i=1}^m X_{\alpha_i}\Fdot(s_i)
\]
 is transverse with \Blue{all} points real.
\end{conj}

As with the Shapiro conjecture for Grassmannians, it is known that transversality
will imply that the points of intersection are real.
There is also a similar discriminant conjecture.
The \DeCo{{\sl preprime}} generated by polynomials $g_1,g_2,\dotsc,g_m$ 
is the collection of polynomials $f$ of the form
\[
   f\ =\ \sigma_0\ +\ \sigma_1 g_1\ +\ \sigma_2 g_2\ +\ \dotsb\ +\ 
   \sigma_m g_m\,,
\]
where the polynomials $\sigma$ are sums of squares.
Any polynomial $f$ in this preprime is positive on the set
\[
   K\ :=\  \{x\mid g_i(x)>0\}
\]
and a representation of $f$ as an element of this preprime is a certificate of
its positivity on $K$.

\begin{conj}[Discriminant Conjecture]
  The discriminant of the polynomial system associated to an intersection of
  Schubert varieties defined by Grassmannian conditions as in the monotone
  conjecture lies in the preprime generated by the differences
  $(s_i-s_j)$ for $i>j$.
\end{conj}

 The point of this conjecture is that not only is the discriminant positive on
 the set $\{(s_1,s_2,\dotsc,s_m)\mid s_1<s_2<\dotsb<s_m\}$ of monotone
 parameters, but that it has a special form from which this positivity is
  transparent. 
 It is known that not all polynomials which are positive on a set of this form
 can lie in the preprime generated by the differences 
 $s_i-s_j$ for $i>j$.

 There is a continually unfolding story of the Shapiro conjecture and
 its generalizations, both to secant flags (by Eremenko, Gabrielov, Shapiro, and 
 Vainshtein~\cite{EGSV}) and for other flag manifolds, for other groups.
 Very few computations have been done (but see~\cite{So00c} for the beginnings
 for the Lagrangian Grassmannian and orthogonal Grassmannian.)

 We end these notes with the remark that, just as they were being completed, 
 Mukhin, Tarasov, and Varchenko~\cite{MTV_Sh} gave a proof the the Shapiro
conjecture for 
Grassmannians, using the Bethe Ansatz and the Gaudin model (a connection between
differential equations with polynomial solutions, representation theory of the
special linear group, and intersections of Schubert varieties given by flags
osculating the rational normal curve).
The discriminant conjecture remain open, however.


\Magenta{{\sf Topics to work on for next rewrite:}}

\begin{enumerate}

\item
Discuss some of the evidence for the Monotone Conjecture

\item
Explain how the result of Eremenko, Gabrielov, Shapiro, and Vainshtein (Section 9.2) 
establishes the monotone conjecture for certain flag manifolds.

\item
Pose the Secant Conjecture and discuss its the evidence for it.
Include gaps, and relate to Chapter 6.

\item
Pose the Monotone secant conjecture and give a few words about current computations.

\item Explain (briefly) the Shapiro conjecture for Lagrangian and Orthogonal
  Grassmannians, including Purbhoo's proof for Orthogonal Grassmannians.

\end{enumerate}

\bibliographystyle{amsplain}
\bibliography{bib}

\end{document}